\documentclass{amsart}
\usepackage{tikz}
\usepackage{xcolor}
\usepackage{amssymb,latexsym,amsmath,extarrows}
\usepackage{graphicx,mathrsfs,comment}
\usepackage[margin=1in, centering]{geometry}
\usepackage{hyperref,url}

\makeatletter
\newcommand{\tocpartheading}[2]{\l@part{Part #1.\quad #2}{}}
\newcommand{\tocsectionindent}[1]{%
  \def\l@section{\@tocline{1}{0pt}{#1}{}{}}}
\makeatother

\usepackage{esint}
\usepackage{amstext}
\usepackage{bbm}

\numberwithin{equation}{section}

\renewcommand{\theequation}{\arabic{section}.\arabic{equation}}

\newtheorem{theorem}{Theorem}[section]
\newtheorem{lemma}[theorem]{Lemma}

\newtheorem{proposition}[theorem]{Proposition}
\newtheorem{remark}[theorem]{Remark}

\newtheorem{definition}[theorem]{Definition}
\newtheorem{corollary}[theorem]{Corollary}
\newtheorem{question}[theorem]{Question}

\newtheorem{convention}[theorem]{Convention}

\newcommand{\e}{\epsilon}

\newcommand{\R}{\mathbb{R}}

\newcommand{\B}{\mathbb{B}}

\newcommand{\supp}{\textup{supp}}
\definecolor{replycolor}{RGB}{140,55,75}

\begin{document}

\title[]{Amplitude dependent wave envelope type estimates for moment curves and applications}

\author{Sangwon Jeong} \address{Sangwon Jeong\\Department of Mathematical Sciences, Seoul National University, Republic of Korea  }\email{hawl1209@snu.ac.kr}

\author{Changkeun Oh}\address{ Changkeun Oh\\ Department of Mathematical Sciences and RIM, Seoul National University, Republic of Korea} \email{changkeun.math@gmail.com}

\maketitle

\setcounter{tocdepth}{1}

\begin{abstract}
We prove amplitude dependent wave envelope type estimates for moment curves  in arbitrary dimensions. The proof uses the sharp square function estimates of Maldague \cite{Mal} and Guth-Maldague \cite{GMSF}. As an application, we establish sharp small cap decoupling inequalities for the moment curve in every dimension \(n\ge4\). This proves the conjecture of Demeter-Guth-Wang \cite{DGW} on exponential sums in all dimensions.
\end{abstract}

\tableofcontents



\section{Introduction}\label{sec:results}

A central problem in Fourier restriction theory is to estimate
$\|f\|_{L^p}$ when the Fourier support of $f$ lies in a small
neighborhood of a manifold $M$. Wave packet decomposition is a standard tool for estimating \(\|f\|_{L^p}\).
We briefly recall the decomposition. First, we partition the Fourier support of \(f\) into rectangular boxes.
To each frequency box we associate a dual box in physical space with reciprocal side lengths, as suggested by the uncertainty principle.
We then write the corresponding part of \(f\) as a sum of functions, each essentially supported on a translate of the dual box.
These functions are called wave packets.
This leads to the following informal question.

\begin{question}\label{1118}
Let $M\subset\R^n$ be a manifold, and let $f$ be a Schwartz
function with Fourier support in a small neighborhood of $M$.
How sharply can $\|f\|_{L^p}$ be bounded in terms of the geometry
and amplitudes of its wave packets?
\end{question}

Guth-Wang-Zhang \cite{GWZ} first provided a satisfactory answer for the parabola in $\R^2$ and the cone
in $\R^3$ by introducing the notion of wave envelopes and proving an $L^4$ wave envelope estimate. Their $L^4$ wave envelope estimate led to a resolution
of the local smoothing conjecture for the wave equation in $2+1$
dimensions. Guth-Maldague \cite{GMAD} subsequently considered the following
informal question, which also takes the amplitude of $f$ into account.

\begin{question}
Let $M\subset\R^n$ be a manifold, and let $f$ be a Schwartz
function with Fourier support in a small neighborhood of $M$.
For $\alpha>0$, how sharply can the measure
\[
 \bigl|\{x\in\R^n:\alpha\le |f(x)|<2\alpha\}\bigr|
\]
be bounded in terms of the geometry and amplitudes of its
wave packets?
\end{question}

Guth-Maldague \cite{GMAD} gave
satisfactory answers for the parabola in $\R^2$ and the cone in $\R^3$
through amplitude dependent wave envelope estimates. Their estimates provide a unified proof of several results in harmonic
analysis that had previously been established by different methods based
on wave packet decomposition. They yield Bochner--Riesz estimates,
square function estimates, decoupling inequalities, small cap decoupling
inequalities, and restriction estimates for both the parabola in $\R^2$ and the cone in $\R^3$.
Guth-Maldague subsequently used their small cap decoupling for the cone in \(\mathbb R^3\) to prove sharp
small cap decoupling for the moment curve in \(\mathbb R^3\) \cite{GMSC}.
The same inequality for the cone was also used by Li-Yang \cite{LY23},
and more recently by Cen-Gao-Zhang-Zhang \cite{CGZZ26}, to obtain
improved bounds in the Gauss circle problem.
\medskip

In this paper, we answer Question~1.2 for the moment curve in
\(\mathbb R^n\) by proving an amplitude dependent wave envelope type estimate. A key ingredient in our proof is the sharp square function estimate for
the moment curve, proved by Maldague \cite{Mal} for \(n=3\) and by
Guth-Maldague \cite{GMSF} for \(n\ge4\). Our amplitude dependent wave envelope type estimate implies
the following known results in harmonic analysis for the moment curve
in \(\mathbb R^n\):
\begin{enumerate}
 \item square function estimates \cite{Mal,GMSF};
 \item the wave envelope estimate of Gan-Maldague-Oh \cite{GMO};
 \item Bochner--Riesz type estimates \cite{GMO};
 \item restriction estimates \cite{Drury1985}.
\end{enumerate}
Moreover, our amplitude dependent wave envelope type estimate implies the following new result:
\begin{enumerate}
 \setcounter{enumi}{4}
 \item sharp small cap decoupling inequalities for the moment curve in every dimension \(n\ge4\),
\end{enumerate}
thereby giving an affirmative answer to the conjecture of Demeter-Guth-Wang \cite{DGW}.
The applications are proved in Part~III. 

\subsection{Small cap decoupling inequalities}\label{subsec:mainresults}
Our main result is the following exponential sum estimate.

\begin{theorem}\label{cor:sums}
Let $n\ge4$, let $N\ge1$ be an integer, and let $2\le p<\infty$.
For $0\le\sigma \leq n$ and an interval $H\subset\R$ of length
$N^{-\sigma}$,
\begin{equation}
 \int_{[0,1]^{n-1}\times H}
 \Bigl|\sum_{k=1}^N a_k e^{2\pi i(kx_1+k^2x_2+\cdots+k^nx_n)}\Bigr|^p dx
 \lesssim_\e N^\e
 \bigl(N^{p/2-\sigma}+N^{p-n(n+1)/2}\bigr)
 \label{eq:sumsmain}
\end{equation}
for every complex sequence $(a_k)_{k=1}^N$ satisfying $|a_k|\le1$.
\end{theorem}

Theorem~\ref{cor:sums} is sharp up to $N^{\epsilon}$ losses: the two terms on the right hand side are obtained by square root cancellation and constructive interference.

The case $\sigma=0$ is the Vinogradov mean value estimate, proved for
$n=3$ by Wooley \cite{Wooley2016Cubic}, and in general by
Bourgain-Demeter-Guth \cite{BDG} and Wooley \cite{Wooley2019}.
The same bound at $\sigma=n$ follows from the
$(n-1)$ dimensional estimate, applied for each fixed $x_n$.
For $0<\sigma\le n-1$, taking $a_k=1$ and $H=[0,N^{-\sigma}]$
resolves the conjecture of Demeter-Guth-Wang
\cite[Conjecture 2.5]{DGW} in all dimensions. For $n-1<\sigma\le n$, the estimate with $a_k=1$ is predicted by Wooley's extended main conjecture in Vinogradov's mean value theorem \cite[Conjecture 8.2]{Wooley2023Subconvexity}.

Demeter-Guth-Wang \cite{DGW} proved the conjecture for $n=3$ in the range $0\le\sigma\le3/2$. Jung \cite[Theorem~1.2]{Jung2023} proved sharp \(L^{10}\) small cap decoupling for the moment curve in \(\mathbb R^3\) at cap width \(R^{-1/2}\), which recovers the tenth moment bound in \eqref{eq:sumsmain} for cubic exponential sums when \(\sigma=1\).
Guth-Maldague \cite{GMSC} subsequently proved
\eqref{eq:sumsmain} for $n=3$ and $0\le\sigma\le2$,
resolving the full cubic case of the conjecture.
Maldague-Oh \cite[Theorem 1.1]{MaldagueOh2026} considered the remaining
range $2<\sigma<3$. For $0\le\sigma<3$, $2\le p<\infty$,
$|a_k|\le1$, and an interval $H$ of length $N^{-\sigma}$, their result gives
\[
 \int_{[0,1]^2\times H}
 \Bigl|\sum_{k=1}^N a_k e^{2\pi i(kx_1+k^2x_2+k^3x_3)}\Bigr|^p dx
 \lesssim_\e N^\e
 \bigl(N^{p/2-\sigma}+N^{1+\sigma(p-7)/3}+N^{p-6}\bigr).
\]
The middle term is dominated by the other two for $\sigma\le5/2$,
but is necessary when $5/2<\sigma<3$. Note that Theorem~\ref{cor:sums} is stated only for $n\ge4$, where no analogous middle term is needed.
In dimension four, the range $3\le\sigma<4$ of \eqref{eq:sumsmain} follows from \cite[Corollary 1.4]{Glidewell}.
\\

Theorem~\ref{cor:sums} follows from a small cap decoupling estimate
which we now describe. Let $\gamma_n(t)=(t,t^2,\ldots,t^n)$ for
$0\le t\le1$. For $n\ge2$, $R\ge1$, and $\beta\ge1/n$, consider
the anisotropic neighborhood
\[
 \mathcal M^n(R^\beta,R)=
 \left\{\gamma_n(t)+\sum_{j=2}^n\frac{a_j}{j!}\gamma_n^{(j)}(t):
 \begin{array}{l}
 t\in[0,1],\quad |a_j|\le R^{-j\beta}\ (2\le j\le n-1),\\
 |a_n|\le R^{-1}
 \end{array}\right\}.
\]
For $n=1$, set $\mathcal M^1(R)=[0,1]$ and use the canonical 
partition.
For $0<\delta\le1$, let $\mathcal I(\delta)$ be the dyadic partition of $[0,1]$
with interval lengths in $(\delta/2,\delta]$.
For $I\in\mathcal I(R^{-\beta})$, define the small cap
\[
 \gamma_I:=\left\{\gamma_n(t)+\sum_{j=2}^n
 \frac{a_j}{j!}\gamma_n^{(j)}(t):
 \begin{array}{l}
 t\in I,\quad |a_j|\le R^{-j\beta}\ (2\le j<n),\\
 |a_n|\le R^{-1}
 \end{array}\right\}.
\]
Since the first coordinate is $t$, these caps form a disjoint
partition of $\mathcal M^n(R^\beta,R)$.
Each $\gamma_I$ is contained in a fixed dilate of an
$R^{-\beta}\times\cdots\times R^{-(n-1)\beta}\times R^{-1}$
box.
When $\beta=1/n$, this is the canonical partition of
$\mathcal M^n(R):=\mathcal M^n(R^{1/n},R)$, whose caps have dimensions
$R^{-1/n}\times R^{-2/n}\times\cdots\times R^{-1}$.
For $n=3$, this is the small cap  considered by Guth-Maldague \cite{GMSC}.  For each $I\in\mathcal I(R^{-\beta})$, define $f_{\gamma_I}$  by $\widehat{f_{\gamma_I}}=\mathbf 1_{\gamma_I}\widehat f$.

\begin{theorem}[Small cap decoupling]\label{thm:SC}
Let $n\ge4$, $\beta\ge1/n$, and $2\le p<\infty$.
For every $\e>0$,
\begin{equation}
 \|f\|_p^p\lesssim_\e R^\e
 \bigl(R^{\beta(p/2-1)}+
 R^{\beta(p-n(n-1)/2-1)-1}\bigr)
 \sum_{I\in\mathcal I(R^{-\beta})}\|f_{\gamma_I}\|_p^p
 \label{eq:SCmain}
\end{equation}
for every Schwartz function $f$ with
$\supp\widehat f\subset\mathcal M^n(R^\beta,R)$.
\end{theorem}
Theorem~\ref{thm:SC} implies Theorem~\ref{cor:sums} for $0\le\sigma<n$ by the argument of Guth-Maldague \cite[Section~1.1]{GMSC}.
Both powers of $R$ are sharp. The first comes from
square root cancellation and the second from constructive
interference of wave packets centered at the origin. The two powers agree at the critical exponent $p=n(n-1)+\frac2\beta$.
At $\beta=1/n$, Theorem~\ref{thm:SC} recovers the canonical
$\ell^p(L^p)$ decoupling estimate with critical exponent $n(n+1)$.
As $\beta$ increases, the critical exponent decreases to $n(n-1)$,
the canonical critical exponent for the moment curve projected
onto the first $n-1$ coordinates.
\medskip

\subsection{Amplitude dependent wave envelope type estimates}
The key ingredient in the proof of Theorem~\ref{thm:SC} is an
amplitude dependent wave envelope type estimate for the moment curve in
$\R^n$. Since its formulation is involved, we first
introduce the ordinary wave envelope type estimate.
As mentioned above, Guth-Wang-Zhang \cite{GWZ} introduced the wave
envelope estimate for the cone in $\R^3$, which
is equivalent to the $L^4$ square function estimate. For the moment curve, the sharp square function estimates
were proved by Maldague \cite{Mal} in $\R^3$ and by
Guth-Maldague \cite{GMSF} in higher dimensions. We formulate a wave
envelope type estimate, which is equivalent to their square function estimate. This ``wave envelope type'' estimate is different from wave envelope estimates in the literature. We explain more after introducing Proposition \ref{thm:intro-ordinary}. 

 For $0<s\le1$, let $\Theta(s)$ denote the canonical partition of
$\mathcal M^n(s^{-n})$. For each cap $\tau=\gamma_I$ with
$I\in\mathcal I(s)$, write $|\tau|=s$ and let $t_\tau$ be the left endpoint
of $I$. We write $\theta$ for a cap in $\Theta(R^{-1/n})$.
For cones, use
\[
 \phi_d(t)=(1,t,t^2/2!,\ldots,t^d/d!)
\]
and the neighborhood of the cone over the moment curve
\[
 \Gamma^{d+1}(R)=
 \left\{\sum_{i=0}^d a_i\phi_d^{(i)}(t):
 \begin{array}{l}
 t\in[0,1],\quad 1\le a_0\le2,\\
 |a_i|\le R^{-i/d}\ (1\le i\le d)
 \end{array}\right\}.
\]
For $0<s\le1$ and $\xi=(\xi_0,\ldots,\xi_d)\in\Gamma^{d+1}(s^{-d})$, let $u(\xi)$ be the point of $[0,1]$ nearest to $\xi_1/\xi_0$. Define
\[
 \tau_I=\{\xi\in\Gamma^{d+1}(s^{-d}):u(\xi)\in I\},\qquad
 \Xi(s)=\{\tau_I:I\in\mathcal I(s)\}.
\]
These caps partition $\Gamma^{d+1}(s^{-d})$. Each $\tau_I$ is contained in a fixed dilate of a $1\times s\times\cdots\times s^d$ box. For each cap $\tau=\tau_I\in\Xi(s)$, write $|\tau|=s$ and let $t_\tau$ be the left endpoint of $I$.
In the wave envelope type estimates below, $\theta$ ranges over $\Theta(R^{-1/n})$ for the curve and $\Xi(R^{-1/d})$ for the cone. In either case, $f_\theta$ denotes the sharp Fourier projection defined by $\widehat{f_\theta}=\mathbf 1_\theta\widehat f$, so that $f=\sum_\theta f_\theta$.\footnote{The same estimates hold if the pieces
$f_\theta$ are defined using a smooth partition of unity subordinate to
fixed enlargements of the caps.
This follows by grouping the caps into $O(1)$ families
and applying the sharp projection estimate to each family
at a comparable scale.} 

For $R^{-1/n}\le s\le1$ and $\tau\in\Theta(s)$, define the centered box
\begin{equation}
 V_{\tau,R}:=\bigl\{x\in\R^n:
 |x\cdot\gamma_n^{(j)}(t_\tau)|\le R^{2/n}s^{2-j}
 \text{ for }1\le j\le n\bigr\}.
 \label{eq:curvebox}
\end{equation}
For $R^{-1/d}\le s\le1$ and $\tau\in\Xi(s)$, set $\phi_d^{(0)}=\phi_d$ and define
\begin{equation}
 U_{\tau,R}:=\bigl\{x\in\R^{d+1}:
 |x\cdot\phi_d^{(j)}(t_\tau)|\le R^{2/d}s^{2-j}
 \text{ for }0\le j\le d\bigr\}.
 \label{eq:conebox}
\end{equation}
We call these the averaging boxes for the estimates below. For curves with $n\ge3$ and cones with $d\ge3$, the averaging boxes at $s=1$ need not contain a full wave packet box and hence need not be wave envelopes in the usual sense.

For a spatial box $B\subset\R^n$, $W_B$ denotes a nonnegative smooth weight
adapted to $B$, normalized by $\int W_B=|B|$.
Write $V\parallel V_{\tau,R}$ and $U\parallel U_{\tau,R}$
for the translates in fixed tilings of the ambient space.
We now state the ordinary wave envelope type estimates.

\begin{proposition}[Wave envelope type estimate for curves]
\label{thm:intro-ordinary}
Let $n\ge1$ and $p_n=\frac{n(n+1)}{2}+1$. For any $\e>0$,
\begin{equation}
 \int_{\R^n}|f|^{p_n}
 \lesssim_\e R^\e
 \sum_{\substack{R^{-1/n}\le s\le1\\ s\ \mathrm{dyadic}}}\sum_{\tau\in\Theta(s)}
 \sum_{V\parallel V_{\tau,R}}
 |V|\biggl[\frac1{|V|}
 \int\bigl(\sum_{\theta\subset\tau}|f_\theta|^2\bigr)^{p_n/4}
 W_V\biggr]^2
 \label{eq:intro-ordinary}
\end{equation}
for every Schwartz function $f$ with $\supp\widehat f\subset\mathcal M^n(R)$.
\end{proposition}

\begin{proposition}[Wave envelope type estimate for cones over moment curves]
\label{thm:intro-ordinary-cone}
Let $d\ge1$ and $p_d=\frac{d(d+1)}{2}+1$. For any $\e>0$,
\begin{equation}
 \int_{\R^{d+1}}|f|^{p_d}
 \lesssim_\e R^\e
 \sum_{\substack{R^{-1/d}\le s\le1\\ s\ \mathrm{dyadic}}}\sum_{\tau\in\Xi(s)}
 \sum_{U\parallel U_{\tau,R}}
 |U|\biggl[\frac1{|U|}
 \int\bigl(\sum_{\theta\subset\tau}|f_\theta|^2\bigr)^{p_d/4}
 W_U\biggr]^2
 \label{eq:intro-ordinary-cone}
\end{equation}
for every Schwartz function $f$ with $\supp\widehat f\subset\Gamma^{d+1}(R)$.
\end{proposition}

To answer Question \ref{1118} for the moment curve, we examined all the examples known to us of functions \(f\) with Fourier support near the curve, and this suggested the formulation of Proposition \ref{thm:intro-ordinary}. Since the estimate is different from a wave envelope estimate in the sense that $V_{\tau,R}$ is comparable to a ball of radius $R^{2/n}$ when $s=1$, it was not so clear if this estimate is useful to the study of small cap decouplings. Interestingly, the amplitude dependent form of this estimate (Theorem~\ref{thm:main}) turned out to be precisely what we needed to prove sharp small cap decoupling for the moment curve.

Proposition~\ref{thm:intro-ordinary} and the sharp square function estimate for the moment curve \cite{Mal,GMSF} are equivalent up to \(R^\epsilon\) losses.
To obtain the square function estimate from \eqref{eq:intro-ordinary}, note that H\"older's inequality bounds the contribution at each scale by a constant times \(\int(\sum_\theta |f_\theta|^2)^{p_n/2}\). For the reverse implication, it suffices to prove
\begin{equation}
 \int\bigl(\sum_\theta|f_\theta|^2\bigr)^{p_n/2}
 \lesssim_\e R^\e
 \sum_{\substack{R^{-1/n}\le s\le1\\ s\ \mathrm{dyadic}}}\ \sum_{\tau\in\Theta(s)}
 \sum_{V\parallel V_{\tau,R}}
 |V|\biggl[\frac1{|V|}
 \int\bigl(\sum_{\theta\subset\tau}|f_\theta|^2\bigr)^{p_n/4}
 W_V\biggr]^2.
 \label{eq:intro-equivalence}
\end{equation}
The cases $n=1,2$ are obvious.
For \(n\ge3\), apply Lemma~\ref{lem:curvefine} to \(h_\theta=|f_\theta|^2\), whose Fourier transform is supported in a fixed dilate of \(\theta-\theta\). Enlarging the sum to all \(V\parallel V_{\tau,R}\) and letting \(\alpha\downarrow0\) gives \eqref{eq:intro-equivalence}.
\\

\label{app:ADstart}
For the application to the small cap decoupling inequality, we need the amplitude dependent
form of the estimate. For each dyadic $s\in[R^{-1/n},1]$, let
\[
 \#\tau=\#\bigl\{\tau'\in\Theta(s):
 f_\theta\not\equiv0\text{ for some }\theta\subset\tau'\bigr\}.
\]
For the cone, use the same definition with $\Xi(s)$ in place of $\Theta(s)$ for each dyadic $s\in[R^{-1/d},1]$. 
For caps define
\begin{equation}
 \mathcal G_\tau^{\mathcal M}(\alpha)=
 \biggl\{V\parallel V_{\tau,R}:
 \frac1{|V|}\int
 \bigl(\sum_{\theta\subset\tau}|f_\theta|^2\bigr)^{p_n/4}W_V
 \ge\Bigl(\frac{\alpha}{C_\e R^\e\#\tau}\Bigr)^{p_n/2}
 \biggr\},
 \label{eq:ADfunctional}
\end{equation}
and for caps in $\R^{d+1}$ define
\[
 \mathcal G_\tau^{\mathcal C}(\alpha)=
 \biggl\{U\parallel U_{\tau,R}:
 \frac1{|U|}\int
 \bigl(\sum_{\theta\subset\tau}|f_\theta|^2\bigr)^{p_d/4}W_U
 \ge\Bigl(\frac{\alpha}{C_\e R^\e\#\tau}\Bigr)^{p_d/2}
 \biggr\}.
\]

We now state the amplitude dependent wave envelope type estimates for the moment curve and the cone over the moment curve.

\begin{theorem}
\label{thm:main}
Let $n\ge1$ and $p_n=\frac{n(n+1)}{2}+1$. For every $\e>0$ and $\alpha>0$,
\begin{equation}
 \int_{\{|f|>\alpha\}}|f|^{p_n}
 \le C_\e R^\e
 \sum_{\substack{R^{-1/n}\le s\le1\\ s\ \mathrm{dyadic}}}\sum_{\tau\in\Theta(s)}
 \sum_{V\in\mathcal G_\tau^{\mathcal M}(\alpha)}
 |V|\biggl[\frac1{|V|}
 \int\bigl(\sum_{\theta\subset\tau}|f_\theta|^2\bigr)^{p_n/4}
 W_V\biggr]^2
 \label{eq:AD}
\end{equation}
for every Schwartz function $f$ with $\supp\widehat f\subset\mathcal M^n(R)$.
\end{theorem}

\begin{theorem}
\label{thm:main-cone}
Let $d\ge1$ and $p_d=\frac{d(d+1)}{2}+1$. For every $\e>0$ and $\alpha>0$,
\begin{equation}
 \int_{\{|f|>\alpha\}}|f|^{p_d}
 \le C_\e R^\e
 \sum_{\substack{R^{-1/d}\le s\le1\\ s\ \mathrm{dyadic}}}\sum_{\tau\in\Xi(s)}
 \sum_{U\in\mathcal G_\tau^{\mathcal C}(\alpha)}
 |U|\biggl[\frac1{|U|}
 \int\bigl(\sum_{\theta\subset\tau}|f_\theta|^2\bigr)^{p_d/4}
 W_U\biggr]^2
 \label{eq:AD-cone}
\end{equation}
for every Schwartz function $f$ with $\supp\widehat f\subset\Gamma^{d+1}(R)$.
\end{theorem}

The cases \(n=1\) and \(d=1\) follow from the locally constant property and
\(L^2\) orthogonality at the finest scale. The case \(n=2\)
for curves and the case \(d=2\) for cones follow from
\cite[Theorems~4 and~2]{GMAD}, respectively. Actually, these
theorems estimate \(\alpha^4|\{|f|>\alpha\}|\); the argument
in Appendix~\ref{subsec:curveweakstrong} replaces this left hand
side by \(\int_{|f|>\alpha}|f|^4\).
Since \(p_2/4=1\), the convolutions can be absorbed into the
weighted averages over \(V\) and \(U\) by \eqref{eq:weightedYoung}.
Theorems~\ref{thm:main} and~\ref{thm:main-cone} are new for \(n\geq 3\) and \(d\geq 3\), respectively.

\subsection{Proof strategy and organization}

We prove the amplitude dependent wave envelope type estimates for moment curves and cones over moment curves in Theorems~\ref{thm:main} and~\ref{thm:main-cone} by induction on the ambient dimension. Building on the wave envelope estimate of Guth-Wang-Zhang \cite[Theorem~1.3]{GWZ} for the cone in \(\mathbb R^3\), Guth-Maldague proved an amplitude dependent wave envelope estimate for this cone \cite[Theorem~2]{GMAD}. Together with their amplitude dependent estimate for the parabola in \(\mathbb R^2\) \cite[Theorem~4]{GMAD}, this gives the base cases of the induction, as explained after Theorem~\ref{thm:main-cone}.

We divide the induction step into two parts. In Part I (Sections~\ref{sec:curve}--\ref{sec:curvesecondrecurrence}), we use the estimate for the cone over the moment curve in \(\mathbb R^n\) to prove the estimate for the moment curve in the same ambient dimension. In Part II (Sections~\ref{sec:highercone}--\ref{sec:hcbase}), we use the estimates for moment curves in \(\mathbb R^k\), \(2\le k\le d\), to prove the estimate for the cone over the moment curve in \(\mathbb R^{d+1}\). Starting with the cone in \(\mathbb R^3\), we alternate between these two steps to obtain
\[
 \text{cone in }\mathbb R^3
 \ \Longrightarrow\ \text{curve in }\mathbb R^3
 \ \Longrightarrow\ \text{cone in }\mathbb R^4
 \ \Longrightarrow\ \text{curve in }\mathbb R^4
 \ \Longrightarrow\ \cdots .
\]
Sections~\ref{sec:curve} and~\ref{sec:highercone} outline the proofs of the two steps, respectively. Readers familiar with wave envelope estimates can follow the main ideas of the proof from these two sections alone.
\medskip

In Part I, we split the superlevel set \(\{x\in\mathbb R^n:|f(x)|>\alpha\}\) into a narrow part, which we estimate by rescaling a single cap, and a broad part, where several separated caps contribute. To estimate the broad part, we use a high low decomposition of  \(\sum_\tau\bigl|\sum_{\theta\subset\tau}f_\theta\bigr|^2*w_\tau\) and their pruned versions, where  \(w_\tau\) is an \(L^1\) normalized weight adapted to the dual box of \(\tau\). The low frequency part is controlled by the sum over the averaging boxes at the coarsest angular scale \(s=1\), while the high frequency pieces, after grouping and rescaling, have Fourier support in caps of the cone over the moment curve, as in the cubic case \cite[Proposition~21]{GMSC}. To estimate these high frequency pieces, we prune wave packets as in \cite[Section~2.1]{GMSC}, thereby obtaining improved \(L^\infty\) bounds for the pruned functions \(f_\tau^j\). Together with the hypothesis of Proposition~\ref{prop:conetocurve} on the cone over the moment curve, this bounds the broad contribution at the given amplitude by a sum over averaging boxes for the pruned functions. To turn this into a wave envelope type estimate for the original function, we adapt Maldague's unpruning algorithm \cite[proof of Proposition~4.6]{Mal}, which bounds the sum for the pruned functions by the corresponding sum for the original functions \(f_\theta\). We group caps of width \(\sigma\) inside caps of width \(s\). At each step, we either obtain a bound in terms of the original functions \(f_\theta\), pass to a finer scale, or advance the pruning index. Since there are only finitely many scales and pruning stages, the iteration terminates.

Part II uses the estimates for moment curves obtained in this way to prove the estimate for the cone over the moment curve in the next dimension. The three dimensional precedent is the proof of Guth-Maldague \cite[Section~8]{GMAD}, who apply the estimate for the parabola to slices of a cone. In higher dimensions we again use a high low decomposition, but the Fourier support of the rescaled high frequency part need not lie in the neighborhood of the cone over the moment curve. Indeed, the radial coefficient can be small while the coefficient of the tangent vector to the generating curve has full size at its scale. Following Guth-Maldague \cite{GMSF}, we include Taylor cones in the induction and use their square function estimates for these sets. More precisely, the high frequency part at Taylor order \(m<d-2\) lies, after rescaling, in a Taylor cone of order \(m+1\). We begin at \(m=d-2\), where \(L^2\) orthogonality controls the high frequency part, and descend to \(m=0\), the cone over the moment curve. At each order \(m\), the estimate for the moment curve in \(\mathbb R^{d-m}\) gives a bound in terms of the sums \(\sum_{\theta\subset\tau}f_\theta\) over caps \(\tau\) at a coarser angular scale, so we still need to bound these averages in terms of the original functions \(f_\theta\). In their proof for the cone in \(\mathbb R^3\), Guth-Maldague \cite[Section~7]{GMAD} prune the fine pieces at successive scales, keeping their contributions on envelopes that meet the amplitude condition. In place of this pruning, we exploit the multiscale structure of the wave envelope type estimate itself to argue by induction on the product \(\sigma s\). The idea comes from the unpruning in Part I, where we reduce either \(\sigma\) or \(s\) to bound the averages of the pruned pieces by those of the original functions \(f_\theta\).

For the cone in \(\mathbb R^3\), the wave envelope associated with \(\tau\in\Xi(s)\) in the work of Guth-Wang-Zhang \cite{GWZ} and Guth-Maldague \cite{GMAD} is comparable to the convex hull of the centered wave packet boxes for all \(\theta\in\Xi(R^{-1/2})\) with \(\theta\subset\tau\). In higher dimensions, our estimates for moment curves and for cones over moment curves require smaller averaging boxes at coarser angular scales. The difficulty in both parts of the proof is that the estimate for the high frequency part, when brought back to the original coordinates, can give averages over boxes larger than these averaging boxes. To recover the required scales, we start with the low frequency cutoff that gives the average over the coarsest box. For a Taylor cone of order \(m\), this cutoff restricts the first \(m+2\) Taylor coordinates, so we take dyadic differences in those coordinates to isolate the high frequency pieces. We track the spatial widths of the resulting kernels under rescaling and use the locally constant property at those widths. When an averaging box is still too large, we subdivide it at the required scale and select the smaller averaging boxes on which the amplitude condition holds.

In Section~\ref{chap:smallcap}, we prove sharp small cap decoupling for the moment curve in \(\mathbb R^n\) for \(n\ge 4\). Motivated by Guth-Maldague's work on the cone in \(\mathbb R^3\) \cite[\S1.1]{GMAD}, we use the amplitude dependent wave envelope type estimate to bound the level set \(\{|f|>\alpha\}\). Roughly speaking, when estimating the integral over each box \(V\) on the right hand side of \eqref{eq:AD}, we need only consider points \(x\in V\) for which \(\sum_{\theta\subset\tau}|f_\theta(x)|^2>c_\epsilon R^{-\epsilon}(\alpha s)^2\), where \(s\) is the width of \(\tau\). When this threshold is sufficiently high, we first apply decoupling for the cone over the moment curve in \(\mathbb R^{n-1}\) and then decoupling for the moment curve in the same dimension. Together these estimates give a stronger bound for the square function than decoupling for the curve alone, and this gain controls the contribution of these scales to the small cap decoupling inequality.

We give two further applications of the amplitude dependent wave envelope type estimate for moment curves. In Section~13 we recover the wave envelope estimate of Gan-Maldague-Oh \cite[Theorem~2.1]{GMO} for moment curves in dimensions \(n\ge 3\), and in Section~\ref{sec:restriction} we derive  restriction estimates for the moment curve. Appendix~\ref{sec:curveauxiliary} gives the interpolation for the cone over the moment curve and Taylor cones. Appendix~\ref{subsec:curveweakstrong} gives the passage from weak type to strong type for the moment curve.

\subsection{Notations}

Throughout the paper, dimensions and exponents, including $\beta$
when it occurs, are fixed. The notation $A\lesssim B$ means $A\le CB$ for
an absolute constant $C$; $A\sim B$ means both inequalities hold.
A subscript records an additional dependence, and $\e$ denotes
an arbitrarily small positive loss exponent.
We write $A\lessapprox B$ if $A\lesssim_\e R^\e B$ for every $\e>0$. We also use \(A\ll B\) to mean that \(A\) is sufficiently small relative to \(B\), with the required smallness depending on the context.

For a spatial box $B\subset\R^n$, let $\phi_B:\R^n\to\R^n$
be an affine map sending $B$ onto $[-1/2,1/2]^n$.
Define $W_B(x)=w(\phi_B(x))$ and $w_B=|B|^{-1}W_B$,
where $w\ge0$, $\int w=1$, and
$w(z)\sim(1+|z|)^{-M}$ with $M$ sufficiently large.
Thus $\int W_B=|B|$ and $\int w_B=1$.
When necessary, we write $w_B=w_{B,E}$ to indicate
the decay exponent $E$.
 For a spatial box \(V\subset\mathbb{R}^n\) and a function \(F:\mathbb{R}^{n}\rightarrow\mathbb{C}\), write
\[
 \fint_V F
 :=\frac1{|V|}\int_{\R^n} F(x)W_V(x)\,dx,
 \qquad
 \|F\|_{L^p_\#(V)}
 :=\biggl(\fint_V |F|^p\biggr)^{1/p}.
\]

\subsection{Acknowledgements}
Changkeun Oh was
supported by the POSCO Science Fellowship of POSCO TJ Park Foundation, and the National Research
Foundation of Korea (NRF) grant funded by the Korea government (MSIT) RS-2024-00341891. The second author learned the high-low method through conversations with Larry Guth and Dominique Maldague. He is grateful for his past collaborations with them.

\subsection{Use of AI}

The authors used ChatGPT-6 for discussions during the development of the
arguments and for assistance in preparing the manuscript. The authors take
full responsibility for the mathematical claims and the final text.

\addtocontents{toc}{\protect\tocpartheading{I}{From cones over moment curves to moment curves}\protect\tocsectionindent{3pc}}
\section{\texorpdfstring{\except{toc}{From cones over moment curves to moment curves: }}{From cones over moment curves to moment curves: }Overview of the proof}
\label{chap:cone-to-curve}\label{sec:curve}

We prove Theorem~\ref{thm:main}  for the curve in $\R^n$ assuming Theorem~\ref{thm:main-cone} for the cone over the moment curve in $\R^n$. In this section, we give overview of the proof.

\begin{proposition}\label{prop:conetocurve}
Let \(n\ge3\). Assume the amplitude dependent wave envelope type estimate
of Theorem~\ref{thm:main-cone} for the cone over the moment curve in \(\R^n\),
at exponent \(p_{n-1}\). Then   the amplitude dependent wave envelope type estimate
of Theorem~\ref{thm:main} for the moment curve in $\R^n$  holds  at exponent \(p_n\).
\end{proposition}

The proof of Proposition \ref{prop:conetocurve} is involved and consisted of several steps. Our proof relies on the induction on the dimension $n$, and the base case is  an amplitude dependent estimate for the cone in \(\R^3\), proved in \cite[Theorem~2]{GMAD}.
\medskip

Our goal is to prove the following 
\begin{equation}
 \int_{\{|f|>\alpha\}}|f|^{p_n}
 \le C_\e R^\e
 \sum_{R^{-1/n}\le s\le1}\sum_{\tau\in\Theta(s)}
 \sum_{V\in\mathcal G_\tau^{\mathcal M}(\alpha)}
 |V|\biggl[\frac1{|V|}
 \int\bigl(\sum_{\theta\subset\tau}|f_\theta|^2\bigr)^{p_n/4}
 W_V\biggr]^2.
\end{equation}
for  every Schwartz
function $f$ with $\supp\widehat f\subset\mathcal M^n(R)$.
Recall that for $\tau\in\Theta(s)$,  $V_{\tau,R}$ is the box
centered at the origin whose side lengths in the dual derivative
coordinates at $t_\tau$ are
\begin{equation}
 R^{2/n}s\times R^{2/n}\times R^{2/n}s^{-1}
 \times\cdots\times R^{2/n}s^{2-n}.
\end{equation}
Recall also that
\begin{equation}
 \mathcal G_\tau^{\mathcal M}(\alpha)=
 \biggl\{V\parallel V_{\tau,R}:
 \frac1{|V|}\int
 \bigl(\sum_{\theta\subset\tau}|f_\theta|^2\bigr)^{p_n/4}W_V
 \ge\Bigl(\frac{\alpha}{C_\e R^\e\#\tau}\Bigr)^{p_n/2}
 \biggr\}.
\end{equation}

{For a cap \(\tau\), write
\(w_\tau=w_{\tau^*}=|\tau^*|^{-1}W_{\tau^*}\),
where \(\tau^*=\tau-\tau\) is its dual box centred at the origin.} To prove our goal, we first  estimates the function $f$ in terms of weighted averages of
\(\sum_\theta|f_\theta|^2*w_\theta\) instead of $\sum_\theta|f_\theta|^2$. Since $|f_{\theta}|$ is ``essentially constant'' on a dual set $\theta$, naively speaking, $|f_{\theta}|^2 *w_{\theta}$ behaves similarly to $w_{\theta}$. To this end, we introduce the following definitions. 

\begin{definition}\label{def2.2} Let  $1 \leq r \leq R$ and $R^{-1/n} \leq s \leq 1$.
For a Schwartz function \(f:\R^n\to\mathbb C\) with
\(\supp\widehat f\subset \mathcal M^n(r)\),
call \(\tau\in\Theta(s)\) active if
\(f_\theta\not\equiv0\) for some \(\theta\subset\tau\), and write
\[
 N(s)=\#\{\tau\in\Theta(s):\tau\text{ is active}\},
 \qquad
 N_\rho(s)=\#\{\tau\in\Theta(s):\tau\subset\rho,
                    \ \tau\text{ is active}\}.
\]
 For cones over moment curves, we use the same
notation with \(\Xi(s)\) in place of \(\Theta(s)\). 
The counts are always taken from the finest scale family under
consideration, and sums over parents omit inactive caps. Define the set of good averaging boxes corresponding to $\tau$ with $|\tau|=s$ by
\[
 \mathcal G^{\mathcal M,w}_\tau(\alpha;f)
 :=\Biggl\{V\parallel V_{\tau,r}:
 \biggl\|\bigl(\sum_{\theta\subset\tau}
 |f_\theta|^2*w_\theta\bigr)^{1/2}\biggr\|_{L^{p_n/2}_\#(V)}
 \ge\frac{\alpha}{N(s)}\Biggr\}.
\]
The superscript $\mathcal M$ denotes the moment curve case, as in \eqref{eq:ADfunctional}. 
\end{definition}

\begin{definition}
Define \(S_1(R)\ge1\) to be the least constant for which
\begin{equation}
 \begin{split}
 \int_{\{|f|>\alpha\}}|f|^{p_n}
 &\le S_1(R)^{p_n}
 \sum_{r^{-1/n}\le s\le1}\sum_{\tau\in\Theta(s)}
 \sum_{V\in\mathcal G^{\mathcal M,w}_\tau(\frac{\alpha}{S_1(R)};f)}
 |V|
 \biggl\|\bigl(\sum_{\theta\subset\tau}
 |f_\theta|^2*w_\theta\bigr)^{1/2}
 \biggr\|_{L^{p_n/2}_\#(V)}^{p_n}
 \end{split}
 \label{eq:curveSone}
\end{equation}
holds for every \(1\le r\le R\), every such Schwartz function with \(\supp\widehat f\subset \mathcal M^n(r)\),
and every \(\alpha>0\).

Define \(S_2(R)\ge1\) analogously, replacing
\(\mathcal G^{\mathcal M,w}_\tau\) with \(\mathcal G^{\mathcal M}_\tau\) and
\(|f_\theta|^2*w_\theta\) with \(|f_\theta|^2\)
in \eqref{eq:curveSone}.
\end{definition}

Our proof has a lot of similarities to those of \cite{GMAD,Mal,GMSF}. By the nature of the problems, the definitions and statements of \cite{GMAD} are involved. Unfortunately, ours are more involved than theirs as we are combining and developing the tools in \cite{Mal,GMSC}. To simplify the notations and follow the key ideas, we use the following convention.

\begin{convention}\label{1102.conv24}
When the scale and family are clear, we abbreviate
\(\mathcal G^{\mathcal M,w}_\tau(\alpha;f)\) to
\(\mathcal G_\tau(\alpha)\).
More generally, in a sum over averaging boxes of the form
\[
 \sum_{\substack{B\parallel B_0\\B\in\mathcal G_\tau(\alpha)}}
 |B|\,\|F_\tau\|_{L^{p_n/2}_\#(B)}^{p_n},
\]
the choice of $B$ is always understood to use the same function \(F_\tau\)
as the summand:
\[
 \mathcal G_\tau(\alpha)
 =\Bigl\{B\parallel B_0:
 \|F_\tau\|_{L^{p_n/2}_\#(B)}
 \ge\frac{\alpha}{N(s)}\Bigr\},\qquad |\tau|=s.
\]
The count \(N(s)\) is always computed from the underlying
finest scale family specified above.
\end{convention}

Our goal is to prove $S_2(R) \lessapprox 1$. However,
as in
\cite[Sections~1.2 and 6.3]{Mal}, we first prove the estimate with
these weights and then remove them by applying it at a smaller
scale. The following proposition can be proved by following the argument of \cite[Proposition~6.8]{Mal}.

\begin{proposition}\label{prop:curveStworecurrence}
Under the hypothesis of Proposition~\ref{prop:conetocurve} on the cone over the moment curve,
\begin{equation}
 S_2(R)\lessapprox S_1(R^{1/n})S_2(CR^{1-1/n}).
 \label{eq:curveD2}
\end{equation}
\end{proposition}

By the proposition, the inequality $S_2(R) \lessapprox 1$ immediately follows from $S_1(R) \lessapprox 1$. So it suffices to prove $S_1(R) \lessapprox 1$. The following is the main step of the proof of Proposition~\ref{prop:conetocurve}.

\begin{proposition}
\label{lem:curveSoneRecurrence}
Under the hypothesis of Proposition~\ref{prop:conetocurve} on the cone over the moment curve, the following holds true.
Let \(\eta>0\) be sufficiently small, with \(N=\eta^{-1}\in\mathbb N\), let \(1\le N_0<N\), and let \(0<\kappa\ll\eta/(1000n^2)\). Then
\begin{equation}
\begin{split}
 S_1(R)\lesssim_{\eta}
 \max\big\{&R^{C\eta N_0},\,
 R^{C\eta}S_1(R^\eta)^{N-N_0+C},\,R^{C\eta^{10}}
S_1(R^{1-n\kappa})\big\}.
 \label{eq:curveD1}
\end{split}
\end{equation}
\end{proposition}

The above proposition easily implies the following corollary.

\begin{corollary}\label{prop:curveSonebound}
Under the hypothesis of Proposition~\ref{prop:conetocurve} on the cone over the moment curve, $S_1(R) \lessapprox 1$.
\end{corollary}

Proposition~\ref{lem:curveSoneRecurrence} and
Proposition~\ref{prop:curveStworecurrence} are proved in
Sections~\ref{sec:curvefirstrecurrence} and~\ref{sec:curvesecondrecurrence},
respectively. As mentioned above, the proof of Proposition \ref{prop:curveStworecurrence} is similar to that for Proposition~6.8 of Maldague~\cite{Mal}. So the main task is to prove Proposition \ref{lem:curveSoneRecurrence}.

\begin{proof}[Proof of Corollary~\ref{prop:curveSonebound}]
Set $d_1=\inf\bigl\{d\ge0:
\sup_{R\ge1}\frac{S_1(R)}{R^d}<\infty\bigr\}$,
and suppose that \(d_1>0\). Then for every \(\delta>0\), $S_1(u)\lesssim_\delta u^{d_1+\delta},~ u\ge1.$
Put
\(N_0=\eta^{-1/2}\) and \(\kappa=\eta^2\).
Substituting into \eqref{eq:curveD1}, we obtain
\[
S_1(R)\lesssim_{\eta,\delta}
\max\left\{
\begin{aligned}
&R^{C\sqrt\eta},\\
&R^{(d_1+\delta)(1-\sqrt\eta+C\eta)+C\eta},\\
&R^{(d_1+\delta)(1-n\eta^2)+C\eta^{10}}
\end{aligned}
\right\}.
\]
Taking \(\eta\) sufficiently small in terms of \(d_1\),
and then \(\delta>0\) sufficiently small in terms of
\(d_1\) and \(\eta\), makes all three exponents strictly
less than \(d_1\), contradicting its definition.
Hence \(d_1=0\), and
\(S_1(R)\lessapprox1\).
\end{proof}

\subsection{Roadmap of Proposition \ref{lem:curveSoneRecurrence}} {Throughout Sections~\ref{sec:curve}--\ref{sec:curvesecondrecurrence},
we assume the hypothesis of Proposition~\ref{prop:conetocurve} on the cone over the moment curve.
Here we outline the proof of Proposition~\ref{lem:curveSoneRecurrence}
and explain how it follows from the lemmas and propositions stated below.
We give the detailed proof in
Sections~\ref{sec:curvefineestimates}--\ref{sec:curvefirstrecurrence}.}
Fix a Schwartz function \(f:\R^n\to\mathbb C\) with
\(\supp\widehat f\subset \mathcal M^n(R)\), and let \(\alpha>0\).
By 
standard reductions, it suffices to bound
\[
 \alpha^{p_n}\bigl|\{|f|>\alpha\}\bigr|
\]
in place of the integral on the left hand side of
\eqref{eq:curveSone}, at the cost of a logarithmic factor.  We give a proof of the reduction in details in {Appendix~\ref{subsec:curveweakstrong}} for the completeness of the paper. 
We now describe the estimates needed for this bound. We split the contribution into two cases:
\[
 \alpha^{p_n}
 \Big|\Big\{|f|>\alpha,\quad
 \sum_\theta|f_\theta|^2*w_\theta>\alpha^2
 \Big\}\Big|
 \quad\text{or}\quad
 \alpha^{p_n}
 \Bigl|\Big\{|f|>\alpha,\quad
 \sum_\theta|f_\theta|^2*w_\theta\le\alpha^2
 \Big\}\Big|.
\]
{We first treat the first case, where
\(\alpha^{p_n}\le(\sum_\theta |f_\theta|^2*w_\theta)^{p_n/2}\).
The high frequency pieces of this sum have Fourier supports in
caps after rescaling. Since the hypothesis on the cone is at
exponent \(p_{n-1}>p_n/2\), we can lower the exponent through
interpolation and apply the estimate for the cone over the moment curve in a
high low decomposition, without introducing pruning.}

\medskip
\noindent\textit{Case 1: Bounding \(~ \alpha^{p_n}
 \big|\big\{|f|>\alpha,\quad
 \sum_\theta|f_\theta|^2*w_\theta>\alpha^2
 \big\}\bigr|\).}

By  Lemma~\ref{lemmalargesqaure}, this case already satisfies the required bound.
We  prove Lemma~\ref{lemmalargesqaure} in
Section~\ref{sec:curvefineestimates}.

\begin{lemma}\label{lem:curvefine}
{Let \(\{h_\theta\}_{\theta\in\Theta(R^{-1/n})}\)
be a family of nonnegative functions in \(L^1(\mathbb R^n)\)
satisfying \(\supp\widehat h_\theta\subset C(\theta-\theta)\).}
For every \(\e>0\) and \(\alpha>0\),
\begin{equation}
\begin{split}
 &\int_{\sum_\theta h_\theta>\alpha^2}
       \bigl(\sum_\theta h_\theta\bigr)^{p_n/2}
 \lesssim_\e R^\e
 \sum_{R^{-1/n}\le s\le1}\sum_{\tau\in\Theta(s)}
 \sum_{V\in\mathcal G_\tau(\frac{\alpha}{R^\e})}
 |V|\,
 \biggl\|\bigl(\sum_{\theta\subset\tau}h_\theta\bigr)^{1/2}
 \biggr\|_{L^{p_n/2}_\#(V)}^{p_n}.
 \label{eq:curvefine}
\end{split}
\end{equation}
{Here \(N(s)\) in the selection condition is computed from the family \(\{h_\theta\}\).}

\label{lemmalargesqaure}
\end{lemma}
{
We use this lemma repeatedly in Sections~4--7.}

Let us assume the lemma and finish the proof for the case. Taking \(h_\theta=|f_\theta|^2*w_\theta\)
in Lemma~\ref{lem:curvefine} gives
\begin{equation}
\begin{aligned}
 &\alpha^{p_n}
 \Big|\Big\{|f|>\alpha,\quad
 \sum_\theta|f_\theta|^2*w_\theta>\alpha^2\Big\}\Big|
 ~\le~
 \int_{\sum_\theta|f_\theta|^2*w_\theta>\alpha^2}
 \bigl(\sum_\theta|f_\theta|^2*w_\theta\bigr)^{p_n/2}
 \\
 &\quad\lesssim_\e R^\e
 \sum_{R^{-1/n}\le s\le1}\sum_{\rho\in\Theta(s)}
 \sum_{V\in\mathcal G_\rho(\frac{\alpha}{R^\e})}
 |V|\,
 \biggl\|\bigl(\sum_{\theta\subset\rho}|f_\theta|^2*w_\theta
 \bigr)^{1/2}\biggr\|_{L^{p_n/2}_\#(V)}^{p_n}.
\end{aligned}
\label{eq:curvelargesquareguide}
\end{equation}
This finishes the proof.
\\

\noindent\textit{Case 2: Bounding \(~\alpha^{p_n}
 \big|\big\{|f|>\alpha,\quad
 \sum_\theta|f_\theta|^2*w_\theta\le\alpha^2
 \big\}\big|\).}

On \(\{|f|>\alpha\}\), by Cauchy--Schwarz inequality and the locally constant property, we obtain
\[
 \alpha^2<|f(x)|^2
 \le N(R^{-1/n})\sum_\theta|f_\theta(x)|^2
 \lesssim N(R^{-1/n})\sum_\theta|f_\theta|^2*w_\theta(x).
\]
Hence, in case 2, we have
\[
 \frac{\alpha^2}{N(R^{-1/n})}
 \lesssim\sum_\theta|f_\theta|^2*w_\theta(x) \le\alpha^2.
\]
By applying the dyadic pigeonholing to the above inequality, it suffices to bound \(\alpha^{p_n}|E_{\alpha,\beta}|\), where
\begin{equation}
 E_{\alpha,\beta}
 =\Big\{x:|f(x)|>\alpha,\quad
 \beta\le
 \bigl(\sum_\theta|f_\theta|^2*w_\theta\bigr)(x)
 \le2\beta\Big\}, \quad \frac{\alpha^2}{N(R^{-1/n})}\lesssim\beta\le\alpha^2.
 \label{beta reduction}
\end{equation}
To bound this contribution, we prune the wave packets of \(f\), with
square root cancellation guiding the choice of threshold. The pruned functions
satisfy the improved \(L^\infty\) bound
\eqref{eq:curveprunebounds}. We then decompose \(E_{\alpha,\beta}\)
into the sets \(\Omega_j\) and a remaining low set \(L\), according to
the first scale at which the square function exceeds its threshold
as we pass from fine to coarse scales. We describe this construction
below.

\subsection{Pruning wave packets}\label{sec:curvepruning}

We first introduce cutoff functions for the wave packet decomposition.
Choose \(\varphi\in C_c^\infty(\mathbb R^n)\) with
\(\|\varphi\|_2=1\) and support in a sufficiently small
neighborhood of the origin, and put
\[
 \psi_0(x)=\int_{[-1/2,1/2]^n}|\check\varphi(x-y)|^2\,dy.
\]
For a cap \(\tau\in\Theta(s)\), let \(A_\tau\) be an invertible
linear map such that
\(\tau^*=A_\tau([-1/2,1/2]^n)\). For \(T\parallel\tau^*\), set
\[
 \psi_{\tau,T}(x)=\psi_0(A_\tau^{-1}(x-c_T)).
\]
Then this function satisfies the following.
\[
 \sum_T\psi_{\tau,T}=1,\qquad
 \sum_T\psi_{\tau,T}^{1/2}\lesssim1,\qquad
 \int\psi_{\tau,T}^{1/2}\lesssim|T|,\qquad
 \supp\widehat\psi_{\tau,T}\subset(\tau-\tau).
\]
By the locally constant property, we obtain
\begin{equation}
 \|\psi_{\tau,T}^{1/2}g_\tau\|_\infty^r
 \lesssim_r\fint_T|g_\tau|^r,
 \qquad r\ge1,\quad\supp\widehat g_\tau\subset C\tau.
 \label{eq:packetBernstein}
\end{equation}

Choose \(\eta>0\) sufficiently small, with
\(N=\eta^{-1}\in\mathbb N\). We have a sequence of numbers.
\begin{equation}
    R^{0} \leq R^{\eta} \leq R^{2\eta} \leq R^{3\eta} \leq \cdots \leq R^{N\eta}=R.
\end{equation}
 We use the notations
\begin{itemize}
    \item $R_j=R^{j\eta}$ for $0 \leq j \leq N$

    \item  \(K=R^\kappa\), where \(0<\kappa\ll\eta/(1000n^2)\), and choose a sufficiently large
\(A=A(\eta)\).

\item  $\tau_j\in\Theta(R_j^{-1/n})$.
\end{itemize}
Set \(f^{N+1}_{\tau_N}=f^{N+1}_\theta=f_\theta\).
For \(j=N,\ldots,1\), define good planks
\begin{equation}
\begin{gathered}
 \mathbb T^g_{\tau_j}
 =
 \Bigl\{T\parallel\tau_j^*:
 \|\psi_{\tau_j,T}^{1/2}f^{j+1}_{\tau_j}\|_\infty
 \le K^{100n^2}A^{N-j+1}\frac{\beta}{\alpha}\Bigr\}.
\end{gathered}
\label{eq:curveprune}
\end{equation}
For each cap \(\tau\), write \(f_\tau=\sum_{\theta\subset\tau}f_\theta\).
Define the pruned functions by
\begin{equation}
     f^j_{\tau_j}
 =
 f^{j+1}_{\tau_j}
 \sum_{T\in\mathbb T^g_{\tau_j}}\psi_{\tau_j,T},
 \qquad
 f^j_{\tau_{j-1}}
 =
 \sum_{\tau_j\subset\tau_{j-1}}f^j_{\tau_j}.
\end{equation}
For \(|\tau|\ge R_j^{-1/n}\), write
\(f^j_\tau=\sum_{\tau_j\subset\tau}f^j_{\tau_j}\).

\begin{lemma}{\cite[Lemma~4.1]{GMSF}, \cite[Lemma~2.1]{Mal}}\label{lem:curveprunebounds}
For \(R\) sufficiently large in terms of \(\eta\), the pruned
functions satisfy
\begin{equation}
\begin{gathered}
 \supp\widehat{f^j_{\tau_j}}\subset3\tau_j,\qquad
 |f^N_\theta|\le|f_\theta|,\qquad
 |f^j_{\tau_j}|\le|f^{j+1}_{\tau_j}|,\\
 \|f^j_{\tau_j}\|_\infty
 \lesssim K^{100n^2}A^{N-j+1}\frac{\beta}{\alpha}.
\end{gathered}
\label{eq:curveprunebounds}
\end{equation}
\end{lemma}

Choose \(1\le N_0<N\), as in
Proposition~\ref{lem:curveSoneRecurrence}, and define for 
\[
\begin{aligned}
&\begin{aligned}[c]
\Omega_j=\biggl\{x\in E_{\alpha,\beta}:\;&
\bigl(\sum_{\tau_j}|f^{j+1}_{\tau_j}|^2*w_{\tau_j}\bigr)(x)
\ge A^{N-j}\beta,
\\
&\bigl(\sum_{\tau_l}|f^{l+1}_{\tau_l}|^2*w_{\tau_l}\bigr)(x)
< A^{N-l}\beta
\quad\text{for all }j<l<N\biggr\},
\end{aligned}
\qquad N_0\le j<N,
\\[6pt]
&L=E_{\alpha,\beta}\setminus
\bigcup_{j=N_0}^{N-1}\Omega_j.
\end{aligned}
\]
 Let \(\eta\in C_c^\infty(\R^n)\) be a radial bump function with
\(\eta=1\) on \(B(0,1)\) and \(\supp\eta\subset B(0,2)\), and
\[
 \eta_{\leq r}(\xi)=\eta(\xi/r),\qquad
 \eta_{>r}=1-\eta_{\leq r},\qquad
 \eta_r=\eta_{\leq r}-\eta_{\leq r/2}.
\]
The following estimates, for \(N_0\le j<N\), give high frequency
dominance on \(\Omega_j\), control the square function on \(L\),
and bound the pruning error.

\begin{lemma}{\cite[Lemma~4.4, Corollary~4.5,
and Lemma~4.6]{GMSF}}\label{11.01.lem210}
\begin{gather}
\begin{gathered}
\Bigl|
  \bigl(\sum_{\tau_j}|f^{j+1}_{\tau_j}|^2*w_{\tau_j}\bigr)
  *\check{\eta}_{\leq R_{j+1}^{-1/n}}
\Bigr|
\lesssim
\sum_{\tau_{j+1}}|f^{j+2}_{\tau_{j+1}}|^2*w_{\tau_{j+1}},
\\[3pt]
\sum_{\tau_j}|f^{j+1}_{\tau_j}|^2*w_{\tau_j}
\lesssim
\Bigl|
  \bigl(\sum_{\tau_j}|f^{j+1}_{\tau_j}|^2*w_{\tau_j}\bigr)
  *\check{\eta}_{>R_{j+1}^{-1/n}}
\Bigr|
\quad\text{on }\Omega_j.
\end{gathered}
\label{eq:curvelowcontrol}
\\
\sum_{\tau_{N_0}}
  |f^{N_0+1}_{\tau_{N_0}}|^2*w_{\tau_{N_0}}
\lesssim A^{N-N_0}\beta
\quad\text{on }L.
\label{eq:curvehighsets}
\\
\begin{aligned}
|f_\tau-f^{j+1}_\tau|
&\le \frac{\alpha}{A^{1/2}K^{100n^2}}
\quad\text{on }\Omega_j,
\\[3pt]
|f_\tau-f^{N_0}_\tau|
&\le \frac{\alpha}{A^{1/2}K^{100n^2}}
\quad\text{on }L,
\end{aligned}
\qquad |\tau|\ge R_1^{-1/n}.
\label{eq:curvepruneerror}
\end{gather}
\end{lemma}
\smallskip

Since the proof is already given in \cite[Lemma~4.4, Corollary~4.5, and Lemma~4.6]{GMSF}, we only give a sketch of the proof.

\subsubsection*{Sketch of the proof of Lemma \ref{11.01.lem210}}  \hfill

{We briefly outline the argument in \cite{GMSF}.
Expand each \(f_{\tau_j}^{j+1}\) into its
\(R_{j+1}^{-1/n}\) children.}
For \(\tau_{j+1}\subset\tau_j\), we have
\(\tau_j^*\subset C\tau_{j+1}^*\), and the shortest side of
\(\tau_{j+1}^*\) is comparable to \(R_{j+1}^{1/n}\).
The rapid decay of \(\check\eta_{\leq R_{j+1}^{-1/n}}\) outside
the ball of this radius therefore gives
\[
 w_{\tau_{j+1}}*w_{\tau_j}*
 \bigl|\check\eta_{\leq R_{j+1}^{-1/n}}\bigr|
 \lesssim w_{\tau_{j+1}}.
\]
{Together with the locally constant property and
\(\lvert f_{\tau_{j+1}}^{j+1}\rvert
\le\lvert f_{\tau_{j+1}}^{j+2}\rvert\), this proves}
the first inequality in \eqref{eq:curvelowcontrol}.
The definitions of \(\Omega_j\) and \(L\){,
together with \(\sum_\theta|f_\theta|^2*w_\theta\le2\beta\)
on \(E_{\alpha,\beta}\), then}
give the second inequality
{in \eqref{eq:curvelowcontrol}}
and \eqref{eq:curvehighsets} for sufficiently large \(A\).

For \eqref{eq:curvepruneerror}, \eqref{eq:packetBernstein}
{with \(r=2\) and the pruning threshold give,
on \(\Omega_j\) for \(j<l\le N\), and on \(L\) for
\(N_0\le l\le N\),}
\[
 {
 \sum_{\tau_l\subset\tau}
 \bigl|f_{\tau_l}^{l+1}-f_{\tau_l}^{l}\bigr|
 \lesssim
 \frac{\alpha}{K^{100n^2}A^{N-l+1}\beta}
 \sum_{\tau_l\subset\tau}|f_{\tau_l}^{l+1}|^2*w_{\tau_l}
 \lesssim \frac{\alpha}{AK^{100n^2}}.
 }
\]
{Summing over \(l=j+1,\ldots,N\) on \(\Omega_j\),
and over \(l=N_0,\ldots,N\) on \(L\), bounds both errors in
\eqref{eq:curvepruneerror} by
\(\lesssim N\alpha/(AK^{100n^2})\).}
Taking \(A\) sufficiently large in terms of \(N\) proves
{\eqref{eq:curvepruneerror}}.

\subsection{Roadmap of Proposition \ref{lem:curveSoneRecurrence} continued}

We return to Case 2 of the roadmap, where we have
\(\sum_\theta|f_\theta|^2*w_\theta\le\alpha^2\). This subsection is intended to give brief ideas of the proof. So I'll omit some details here. 

We need to bound the set $\alpha^{p_n}|E_{\alpha,\beta}|$. By the standard broad narrow analysis, we may assume that all the points are ``broad''. The narrow case can be dealt by using rescaling as always in  literature.
By the  error estimate \eqref{eq:curvepruneerror}, we may replace $f_{\tau}$ by the pruned functions $f_{\tau}^{j+1}$.
On \(\Omega_j\), by computations, the square function is dominated by its
high frequency part. After dyadic decomposition on the frequency side and
making suitable affine changes of variables, the Fourier supports
of these pieces lie in a neighborhood of the cone over the moment curve in \(\R^n\).
We can therefore apply the hypothesis of Proposition~\ref{prop:conetocurve}
on the cone over the moment curve. 
More precisely, after the standard broad narrow decomposition,
\eqref{eq:curvepruneerror} allows us to apply this argument to
the broad set defined below.

For \(K^{-1}\ge R_1^{-1/n}\), call a tuple of caps in
\(\Theta(K^{-1})\) \(K^{-1}\) separated if the corresponding intervals in \([0,1]\) have pairwise distance at least \(K^{-1}\).
Fix a \(K^{-1}\) separated tuple
\(\tau_1,\ldots,\tau_n\in\Theta(K^{-1})\), and put
\[
\operatorname{Br}_{\alpha,\beta}^K
=
\Bigl\{
x\in E_{\alpha,\beta}:
\max_{\tau\in\Theta(K^{-1})}|f_\tau(x)|
\le {\alpha},
\quad
|f_{\tau_1}(x)\cdots f_{\tau_n}(x)|^{1/n}
\ge \frac{\alpha}{K^{100n}}
\Bigr\}.
\]
The hypothesis on the cone in the same ambient dimension gives
the following proposition below; see
\eqref{eq:curvebroadreturned} and \eqref{eq:curveLfinal}.
\begin{proposition}\label{prop:curvebroadpieces}
Let \(1\le N_0<N\), and let \(f:\mathbb R^n\to\mathbb C\)
be a Schwartz function with
\(\operatorname{supp}\widehat f\subset\mathcal M^n(R)\).
For each \(N_0\le k<N\), there is a dyadic \(\sigma\) with
\(R_{k+1}^{-1/n}<\sigma\le R_k^{-1/n}\) such that
\begin{equation}
\begin{aligned}
\alpha^{p_n}
|\operatorname{Br}_{\alpha,\beta}^{K}\cap\Omega_k|
&\lesssim_\eta
R^{C\eta}
\sum_{\sigma\le s\le1}
\sum_{\rho\in\Theta(s)}
\sum_{\substack{
 B\parallel(s/\sigma)V_{\rho,\sigma^{-n}}\\
 B\in\mathcal G_\rho(\frac{\alpha}{R^{C\eta}})
}}
|B|
\biggl\|
 \Bigl(
 \sum_{\substack{\tau\subset\rho\\|\tau|=\sigma}}
 |f_\tau^{k+1}|^2*w_\tau
 \Bigr)^{1/2}
\biggr\|_{L^{p_n/2}_\#(B)}^{p_n}.
\\
\alpha^{p_n}|L|
&\lesssim_\eta R^{C\eta N_0}
\sum_{R^{-1/n}\le s\le1}
\sum_{\rho\in\Theta(s)}
\sum_{V\in\mathcal G_\rho(\frac{\alpha}{R^{C\eta N_0}})}
|V|
\biggl\|
\bigl(
\sum_{\theta\subset\rho}|f_\theta|^2*w_\theta
\bigr)^{1/2}
\biggr\|_{L_\#^{p_n/2}(V)}^{p_n}.
\end{aligned}
\label{eq:curvebroadreturned or}
\end{equation}
\end{proposition}
Here the averaging boxes \((s/\sigma)V_{\rho,\sigma^{-n}}\), \(|\rho|=s\),
are obtained by pulling the averaging boxes \(U_{\rho,\sigma^{-(n-1)}}\) back to the original coordinates. 
The second inequality is always of the desired form. On the other hand, the right hand side of the first inequality is still far from what we are hoping for.
Proposition~\ref{prop:curveunpruning} below
\emph{unprunes} the right hand side of
\eqref{eq:curvebroadreturned or}, expressing it in terms of
the original pieces \(f_\theta\).

\begin{proposition}\label{prop:curveunpruning}
Let \(1\le j\le N\) and
\(R_j^{-1/n}\le\sigma\le R_{j-1}^{-1/n}\).
For \(s\in[\sigma,1]\) and \(\alpha>0\),
\begin{equation}
\begin{split}
 &\sum_{\rho\in\Theta(s)}
 \sum_{\substack{
 B\parallel(s/\sigma)V_{\rho,\sigma^{-n}}\\
 B\in\mathcal G_\rho(\alpha)
 }}
 |B|\,
 \biggl\|
 \Bigl(
 \sum_{\substack{\tau\subset\rho\\|\tau|=\sigma}}
 |f^j_\tau|^2*w_\tau
 \Bigr)^{1/2}
 \biggr\|_{L^{p_n/2}_\#(B)}^{p_n}
 \lesssim_\eta
 \bigl(R^{C\eta}
 S_1(R^\eta)^{N-j+C}\bigr)^{p_n}
 \\
 &\qquad{}\times
 \sum_{R^{-1/n}\le s'\le s}
 \sum_{\rho'\in\Theta(s')}
 \sum_{V\in\mathcal G_{\rho'}\bigl(
 \frac{\alpha}{
 R^{C\eta}S_1(R^\eta)^{N-j+C}}
 \bigr)}
 |V|\biggl\|
 \bigl(
 \sum_{\theta\subset\rho'}|f_\theta|^2*w_\theta
 \bigr)^{1/2}
 \biggr\|_{L^{p_n/2}_\#(V)}^{p_n}.
 \label{eq:curveunprunecost}
\end{split}
\end{equation}
{The counts \(N(s)\) and \(N(s')\) in the
selection conditions are computed from the fine scale  \(\{f_\theta\}_{\theta\in \Theta(R^{-1/n})}\).} 
\end{proposition}

Combining Propositions~\ref{prop:curvebroadpieces} and
\ref{prop:curveunpruning} gives 

\begin{equation}\label{eq:curvebroadbound}
\begin{aligned}
&\alpha^{p_n}|\operatorname{Br}_{\alpha,\beta}^{K}|
\lesssim_\eta R^{C\eta N_0}
 \sum_{R^{-1/n}\le s\le1}
 \sum_{\rho\in\Theta(s)}
 \sum_{V\in\mathcal G_\rho(\frac{\alpha}{R^{C\eta N_0}})}
 |V|
 \biggl\|
 \bigl(\sum_{\theta\subset\rho}|f_\theta|^2*w_\theta\bigr)^{1/2}
 \biggr\|_{L^{p_n/2}_\#(V)}^{p_n}
\\
&\quad{}+
 \bigl(
 R^{C\eta}
 S_1(R^\eta)^{N-N_0+C}
 \bigr)^{p_n}
 \sum_{R^{-1/n}\le s\le1}
 \sum_{\rho\in\Theta(s)}
 \sum_{V\in\mathcal G_\rho\bigl(
 \frac{\alpha}{
 R^{C\eta}S_1(R^\eta)^{N-N_0+C}}
 \bigr)}
 |V| \biggl\|
\bigl(\sum_{\theta\subset\rho}|f_\theta|^2*w_\theta\bigr)^{1/2}
 \biggr\|_{L^{p_n/2}_\#(V)}^{p_n}.
\end{aligned}
\end{equation}

Recall Proposition~\ref{lem:curveSoneRecurrence}:
\begin{equation}
\begin{split}
 S_1(R)\lesssim_{\eta}
 \max\big\{&R^{C\eta N_0},\,
 R^{C\eta}S_1(R^\eta)^{N-N_0+C},\,R^{C\eta^{10}}
S_1(R^{1-n\kappa})\big\}.
\end{split}
\end{equation}
The two sums in
\eqref{eq:curvebroadbound} give the first two terms in
\[
 \max\bigl\{
 R^{C\eta N_0},\,
 R^{C\eta}S_1(R^\eta)^{N-N_0+C},\,R^{C\eta^{10}}
S_1(R^{1-n\kappa})
 \bigr\}
\]
in \eqref{eq:curveD1}, where \(K=R^{\kappa}\).
The last term comes from the narrow contribution by induction, via rescaling. This is an outline of the proof of Proposition~\ref{lem:curveSoneRecurrence}.
\\

We finish this subsection by explaining what we prove in Section~\ref{sec:curverefinement}--~\ref{sec:curvesecondrecurrence}. Section~\ref{sec:curverefinement} develops the geometric lemmas
needed for Proposition~\ref{prop:curveunpruning}.
Section~\ref{sec:curveunpruning} then iterates these estimates
to prove Proposition~\ref{prop:curveunpruning}.
In Section~\ref{sec:curvefirstrecurrence},
we first prove Proposition~\ref{prop:curvebroadpieces}
and then use \eqref{eq:curvebroadbound} in a broad narrow
analysis to prove Proposition~\ref{lem:curveSoneRecurrence}.
 This finishes the proof of Proposition~\ref{lem:curveSoneRecurrence}.
Section~\ref{sec:curvesecondrecurrence} proves
Proposition~\ref{prop:curveStworecurrence}, following the strategy
of \cite[Section~6.3]{Mal}.

\subsection{\texorpdfstring{{Preliminaries}}{Preliminaries}}

{We record the rescaling maps and auxiliary
lemmas used in the subsequent sections.}

\subsubsection{Rescaling}

{We begin with} the affine rescaling for the moment curve; see
\cite[proof of Lemma~4.3]{Mal}. For \(0<s_0\le1\) and
\(0\le t_0\le1-s_0\), set
\[
 L_{t_0,s_0}
 =\Bigl(s_0\gamma_n'(t_0),
        \frac{s_0^2}{2!}\gamma_n''(t_0),\ldots,
        \frac{s_0^n}{n!}\gamma_n^{(n)}(t_0)\Bigr).
\]
Taylor's formula gives
\[
 L_{t_0,s_0}^{-1}
 \bigl(\gamma_n(t_0+s_0u)-\gamma_n(t_0)\bigr)
 =\gamma_n(u),\qquad 0\le u\le1.
\]
Thus the affine frequency map
\(\xi\mapsto L_{t_0,s_0}^{-1}(\xi-\gamma_n(t_0))\)
sends a cap of parameter length \(\sigma\le s_0\), contained
in the parent cap over \([t_0,t_0+s_0]\), to a cap
of length \(\sigma/s_0\). In particular, the fine scale
\(R^{-1/n}\) becomes \((s_0^nR)^{-1/n}\).
Difference sets \(\theta-\theta\) are transformed by the linear
map \(L_{t_0,s_0}^{-1}\). The corresponding change of spatial
variables is \(x=L_{t_0,s_0}^{-T}y\); returning to the original
coordinates multiplies the \(i\)th dual derivative length by
\(s_0^{-i}\).

For the cone generated by \(\gamma_n'\), the linear frequency
map is \(\xi\mapsto s_0L_{t_0,s_0}^{-1}\xi\), since
\[
 s_0L_{t_0,s_0}^{-1}\gamma_n'(t_0+s_0u)=\gamma_n'(u).
\]
This preserves the radial parameter and sends aperture \(\sigma\)
to \(\sigma/s_0\). The spatial change is
\(x=s_0L_{t_0,s_0}^{-T}y\), so the corresponding factors for the
dual derivative lengths are \(s_0^{1-i}\), \(1\le i\le n\).

\subsubsection{\texorpdfstring{{Auxiliary lemmas}}{Auxiliary lemmas}}
\begin{lemma}
For a measurable \(F:\R^n\to\mathbb C\), \(t>0\), and a tiling
\(B\parallel B_0\), we have
\begin{equation}
 \sum_{\substack{B\parallel B_0\\\fint_B |F|>t}}
 |B|\biggl(\fint_B |F|\biggr)^2
 \lesssim\int_{|F|>t/2}|F|^2.
 \label{eq:selaverage}
\end{equation}
\end{lemma}
\begin{proof}
If $\fint_B |F|>t$, then
\[
 \int_{|F|\le t/2}|F|W_B
 \le\frac t2|B|<\frac12\int |F|W_B.
\]
Since $\int W_B=|B|$, Cauchy--Schwarz gives
\[
 |B|\biggl(\fint_B |F|\biggr)^2
 \lesssim\frac1{|B|}\biggl(\int_{|F|>t/2}|F|W_B\biggr)^2
 \lesssim\int_{|F|>t/2}|F|^2W_B.
\]
Summing over $B$ proves \eqref{eq:selaverage}.
\end{proof}

The first lemma allows us to replace sums over larger averaging boxes
by sums over \emph{smaller} averaging boxes while preserving the
amplitude condition.

\begin{lemma}[Changing averaging boxes]\label{lem:envelopechange}
Let \(F\ge0\) and \(t>0\). Fix a tiling by translates
\(B\parallel B_0\) of a parallelepiped \(B_0\).
If \(U_0\subset CB_0\) has the same axes as \(B_0\), then
\begin{equation}
 \sum_{\substack{B\parallel B_0:\fint_B F>t}}
 |B|\biggl(\fint_B F\biggr)^2
 \lesssim
 \sum_{\substack{U\parallel U_0:\fint_U F>ct}}
 |U|\biggl(\fint_U F\biggr)^2.
 \label{eq:subdivide}
\end{equation}
Here \(c>0\) is fixed independently of \(F\), \(B_0\), and \(t\).
\end{lemma}
\begin{proof}
Set
\[
 G(x)=\sum_{U\parallel U_0}\biggl(\fint_U F\biggr)\mathbf1_U(x).
\]
Since \(U_0\subset CB_0\) has the same axes as \(B_0\), on each
\(U\) we have \(W_B\sim W_B(c_U)\) and \(W_U\gtrsim1\). Hence
\[
 \fint_B F
 =\frac1{|B|}\sum_{U\parallel U_0}\int_U F W_B
 \lesssim\sum_{U\parallel U_0}
 \frac{W_B(c_U)|U|}{|B|}\fint_U F
 \sim\fint_B G.
\]
Applying \eqref{eq:selaverage} to \(G\), for a sufficiently small
fixed \(c>0\), gives
\[
\begin{aligned}
 \sum_{\substack{B\parallel B_0:\fint_B F>t}}
 |B|\biggl(\fint_B F\biggr)^2
 &\lesssim\sum_{\substack{B\parallel B_0:\fint_B G>2ct}}
 |B|\biggl(\fint_B G\biggr)^2\\
 &\lesssim\int_{G>ct}G^2
 =\sum_{\substack{U\parallel U_0:\fint_U F>ct}}
 |U|\biggl(\fint_U F\biggr)^2.
\end{aligned}
\]
This proves \eqref{eq:subdivide}.
\end{proof}
\begin{remark} Let us briefly explain our notations. We use $V,U,B$ to denote averaging boxes.
{We write \(V\parallel V_{\tau,r}\) for the averaging boxes in the tiling for the moment curve
and \(U\parallel U_{\tau,r}\) for those in the tiling for the cone over the moment curve.
To estimate the high frequency part, we will apply the estimate for the cone over the moment curve
after rescaling. The boxes obtained from \(U\) on returning
to the original coordinates will be denoted by \(B\).} 
\end{remark}

The next lemma allows us to remove convolution kernels.
\begin{lemma}\label{lem:selection}
Let \(F:\R^n\to\mathbb C\) be measurable, \(t>0\), \(r\ge1\),
and let \(K\) be an integrable kernel with \(\|K\|_1\lesssim1\).
Fix a tiling \(B\parallel B_0\). Then
\begin{equation}
 \int_{|K*F|>t}|K*F|^r
 \lesssim_r\int_{|F|>ct}|F|^r.
 \label{eq:tailYoung}
\end{equation}
\begin{equation}
 \sum_{\substack{B\parallel B_0\\\fint_B|K*F|>t}}
 |B|\biggl(\fint_B|K*F|\biggr)^2
 \lesssim
 \sum_{\substack{B\parallel B_0\\\fint_B |F|>ct}}
 |B|\biggl(\fint_B |F|\biggr)^2.
 \label{eq:commonreturn}
\end{equation}
If \(|K|*W_B\lesssim W_B\), then
\begin{equation}\label{eq:weightedYoung}
 \fint_B|K*F|^r
 \lesssim_r\fint_B|F|^r.
\end{equation}
\end{lemma}
\begin{proof}
Choosing
$c>0$ sufficiently small, on $\{|K*F|>t\}$,
\[
 |K*F|\lesssim|K|*(|F|\mathbf1_{|F|>ct}).
\]
Young's inequality and $\|K\|_1\lesssim1$ give
\eqref{eq:tailYoung}.

For $y\in B$, the definition of the weights gives
\[
 W_{B_0}(x-y+c_{B_0})\sim W_B(x),\qquad
 \int |F(x)|w_{B_0}(x-y+c_{B_0})\,dx
 \sim\fint_B |F|.
\]
The comparison is uniform in $x$, $y$, and $B$. It also holds
with $|F|$ replaced by $|K|*|F|$.
Since $|K*F|\le |K|*|F|$, the left side of
\eqref{eq:commonreturn} is bounded by
\begin{align*}
 &C\int_{\int(|K|*|F|)(x)w_{B_0}(x-y+c_{B_0})\,dx>ct}
 \biggl[
 \int(|K|*|F|)(x)w_{B_0}(x-y+c_{B_0})\,dx
 \biggr]^2dy\\
 &\lesssim
 \int_{\int |F(x)|w_{B_0}(x-y+c_{B_0})\,dx>ct}
 \biggl[
 \int |F(x)|w_{B_0}(x-y+c_{B_0})\,dx
 \biggr]^2dy\\
 &\lesssim
 \sum_{\fint_B |F|>ct}
 |B|\biggl(\fint_B |F|\biggr)^2.
\end{align*}
The middle inequality is \eqref{eq:tailYoung} with $r=2$,
since convolution with $|K|$ commutes with the displayed
averaging integral.
The last inequality follows by integrating over the tiles $B$.
The fixed constant $c$ may decrease in these inequalities.

Finally, H\"older's inequality gives
\[
 |K*F|^r
 \le\|K\|_1^{r-1}\bigl(|K|*|F|^r\bigr)
 \lesssim_r |K|*|F|^r.
\]
By Fubini's theorem, using $W_B(2c_B-x)\sim W_B(x)$
in the last inequality,
\[
 \fint_B|K*F|^r
 \lesssim_r
 \frac1{|B|}\int \bigl(|K|*|F|^r\bigr)W_B
 =
 \frac1{|B|}\int |F(x)|^r
       \biggl(\int |K(y)|W_B(x+y)\,dy\biggr)dx
 \lesssim_r
 \fint_B|F|^r.
\]
This finishes the proof.
\end{proof}

\begin{lemma}\label{lem:selectedaverages}
Let \(F:\R^n\to\mathbb C\) be measurable, \(p\ge2\), and let
\(Q\parallel Q_0\) tile \(\R^n\) by translates of a fixed parallelepiped.
Suppose that, uniformly in \(Q\),
\[
 |F(x)|\lesssim\|F\|_{L^{p/2}_\#(Q)}
 \qquad\text{for a.e. }x\in Q.
\]
Then, for every \(a>0\),
\begin{equation}\label{eq:selectednormintegral}
 \int_{|F|>Ca}|F|^p
 \lesssim
 \sum_{\substack{Q\parallel Q_0\\\|F\|_{L^{p/2}_\#(Q)}\ge a}}
 |Q|\|F\|_{L^{p/2}_\#(Q)}^p
 \lesssim \int_{|F|>ca}|F|^p,
\end{equation}
for some constants \(C>1\) and \(0<c<1\).
\end{lemma}
\begin{proof}
The second inequality follows from \eqref{eq:selaverage} applied to \(|F|^{p/2}\).
For the first, integrate the pointwise bound over
\(Q\cap\{|F|>Ca\}\) and sum over \(Q\); this set is null unless
\(\|F\|_{L^{p/2}_\#(Q)}\ge a\).
\end{proof}

Finally, {in proving
Propositions~\ref{prop:curveStworecurrence}
and~\ref{lem:curveSoneRecurrence}}, we may assume the counting
inequalities in the following uniformization lemma. Recall Definition \ref{def2.2}.

\begin{lemma}[Uniformization]\label{lem:regular}
A subcollection of \(\{f_{\theta}\}_{\theta\in\Theta(R^{-1/n})}\) can be partitioned into
\(O_\e(R^\e)\) classes such that, in each class,
\begin{equation}
 N(s)N_\tau(\sigma)\lesssim_\e N(\sigma),
 \qquad \tau\in\Theta(s),\quad R^{-1/n}\le\sigma\le s\le1.
 \label{eq:regularcounts}
\end{equation}
The same inequality holds hold inside any dyadic parent \(\rho\), with
\(N_\rho(s)\) and \(N_\rho(\sigma)\) in place of
\(N(s)\) and \(N(\sigma)\), for
\(\tau\subset\rho\) and \(\sigma\le s\le|\rho|\).
\end{lemma}

\begin{proof}
Let $\mathcal A=
\{\theta\in\Theta(R^{-1/n}):f_\theta\not\equiv0\}$
and choose an integer $b$ such that $\log_2(b+1)/b<\epsilon$.
Starting from scale $1$, take every $b$th dyadic scale and include
the finest scale, obtaining the block scales
\[
1=s_0>s_1>\cdots>s_L=R^{-1/n},
\qquad
L=\Bigl\lceil\frac{\log_2 R}{nb}\Bigr\rceil.
\]
Each cap at scale $s_j$ has at most $2^b$ children at scale $s_{j+1}$. First consider the scales $s_{L-1},s_L$. For each $\tau\in\Theta(s_{L-1})$, write
\[
\#\{\theta\in\mathcal A:\theta\subset\tau\}
=
\sum_{k_1\in E_\tau}2^{k_1},
\qquad
E_\tau\subset\{0,\ldots,b\}.
\]
Partition these caps into sets $C_{\tau,k_1}$ of size $2^{k_1}$,
and put
\[
\mathcal A_{k_1}
=
\bigcup_{\tau:\,k_1\in E_\tau}C_{\tau,k_1}.
\]
The nonempty sets $\mathcal A_{k_1}$ partition $\mathcal A$.
In each $\mathcal A_{k_1}$, every active $s_{L-1}$ cap contains
exactly $2^{k_1}$ caps at scale $s_L=R^{-1/n}$.

Proceed upwards one block at a time.
After $r$ steps, fix a class $\mathcal A_{k_1,\ldots,k_r}$.
For each parent $\tau\in\Theta(s_{L-r-1})$, count its children
$\nu\in\Theta(s_{L-r})$ containing a member of this class.
Using the binary expansion of this count, partition these children
into groups $C_{\tau,k_{r+1}}$ of size $2^{k_{r+1}}$, with
$k_{r+1}\in E_\tau\subset\{0,\ldots,b\}$, and define
\[
\mathcal A_{k_1,\ldots,k_r,k_{r+1}}
=
\bigcup_{\tau:\,k_{r+1}\in E_\tau}
\ \bigcup_{\nu\in C_{\tau,k_{r+1}}}
\{\theta\in\mathcal A_{k_1,\ldots,k_r}:\theta\subset\nu\}.
\]
In each new class, every active $s_{L-r-1}$ cap has exactly
$2^{k_{r+1}}$ active $s_{L-r}$ children.
Each selected child keeps all its fine descendants in
$\mathcal A_{k_1,\ldots,k_r}$, so the branching numbers already
fixed at finer block scales are preserved.

Each step splits a class into at most $b+1$ classes.
Hence the number of final classes is at most
\[
(b+1)^L
=
(b+1)^{\lceil\log_2 R/(nb)\rceil}
\lesssim_\epsilon R^\epsilon.
\]
In each final class, all active parents at a block scale have
the same number of active children at the next block scale.
Multiplying these branching numbers gives
$N(s)N_\tau(\sigma)=N(\sigma)$ at block scales.

Rounding $s,\sigma$ to adjacent block scales changes the counts
by a factor depending only on $b$, which proves \eqref{eq:regularcounts}.
Applying the same argument inside $\rho$ gives the corresponding counting inequalities.
\end{proof}

In proving the bounds for \(S_1\) and \(S_2\), we assume from now on
that the original fine family \(\{f_\theta\}\) satisfies the counting
inequalities of Lemma~\ref{lem:regular}. To justify this reduction, write
\(f=\sum_{i=1}^J f^{(i)}\), where each \(f^{(i)}\) is the sum over
one of these classes. We have
\[
 \int_{\{|f|>\alpha\}} |f|^{p_n}
 \le J^{p_n}\sum_{i=1}^J
 \int_{\{|f^{(i)}|>\alpha/J\}}
 |f^{(i)}|^{p_n}.
\]
and for every class \(\mathcal A\),
\[
 N_{\mathcal A}(s)\le N(s),\qquad
 \sum_{\substack{\theta\subset\tau:\theta\in\mathcal A}}
       |f_\theta|^2*w_\theta
 \le\sum_{\theta\subset\tau}|f_\theta|^2*w_\theta.
\]
The same inequality holds without the convolutions. Hence the boxes
entering the sums for a class at amplitude \(\alpha/J\) also enter those
for the full family after multiplying the selection constant by \(J\). Taking the uniformization parameter
sufficiently small absorbs these factors into \(R^\e\).

\section{\texorpdfstring{\except{toc}{From cones over moment curves to moment curves: }}{From cones over moment curves to moment curves: }Proof of Lemma~\ref{lem:curvefine}}
\label{sec:curvepreliminaries}\label{sec:curvefineestimates}

In this section, we prove Lemma~\ref{lem:curvefine}. Here is what we need to prove: for any \(\alpha>0\)
\begin{equation}
\begin{split}
 &\int_{\sum_\theta h_\theta>\alpha^2}
       \bigl(\sum_\theta h_\theta\bigr)^{p_n/2}
 \lesssim_\e R^\e
 \sum_{R^{-1/n}\le s\le1}\sum_{\tau\in\Theta(s)}
 \sum_{V\in\mathcal G_\tau(\frac{\alpha}{R^\e})}
 |V|\,
 \biggl\|\bigl(\sum_{\theta\subset\tau}h_\theta\bigr)^{1/2}
 \biggr\|_{L^{p_n/2}_\#(V)}^{p_n}
\end{split}
\end{equation}
{for any nonnegative function \(\{h_\theta\}_{\theta\in\Theta(R^{-1/n})}\)
satisfying \(\supp\widehat h_\theta\subset C(\theta-\theta)\).}
Recall that we are assuming Theorem~\ref{thm:main-cone} for the cone over the moment curve in
\(\R^n\) at exponent $p_{n-1}$.  

Since the
high low analysis applies it to quadratic expressions in the
Fourier projections \(f_\theta\), we need the following estimate at exponent
\(p_n/2\), as in \cite[Proposition~4.4]{Mal}.
A fixed affine map takes \(\phi_{n-1}\) to \(\gamma_n'\).
This change of normalization affects only implicit constants
and introduces no additional difficulty in adapting Maldague's
argument under our hypothesis on the cone. After the corresponding change of spatial variables,
\(U_{\rho,\sigma^{-(n-1)}}\), for \(\rho\in\Xi(s)\), agrees up to a fixed dilation with
\[
 \bigl\{x\in\R^n: |x\cdot\gamma_n^{(i)}(t_\rho)|\le
 \sigma^{-2}s^{3-i}\text{ for }1\le i\le n\bigr\}.
\]
For the cone families on \(\Gamma^n\), the notation for good averaging boxes
\(\mathcal G_\rho\) follows the convention \ref{1102.conv24}.

The following estimate is very simlar to \eqref{eq:AD-cone}. The main difference is that we raise $p_n/2$th power than $p_{n-1}$ (note that we are using $d=n-1$ as  hypothesis). Since $2 \leq p_{n}/2 \leq p_{n-1}$, by interpolating the exponents $p_{n-1}$ and $2$, one may expect that the estimate is true for $p_n/2$ as well.

\begin{lemma}\label{lem:curvconetail}
Assume  Theorem~\ref{thm:main-cone} for the cone over the moment curve in
\(\R^n\). Let \(\alpha>0\). We have
\[
 \int_{|g|> \alpha^2}|g|^{p_n/2}
 \lesssim_\e\sigma^{-\e}
 \sum_{\sigma\le s\le1}\sum_{\rho\in\Xi(s)}
 \sum_{\substack{U\parallel U_{\rho,\sigma^{-(n-1)}}\\
 U\in\mathcal G_\rho(\frac{\alpha}{C_\e\sigma^{-\e}})}}
 |U|\biggl\|\bigl(\sum_{\tau\subset\rho}|g_\tau|\bigr)^{1/2}
 \biggr\|_{L^{p_n/2}_\#(U)}^{p_n}
\]
for a Schwartz function \(g\) with
\(\supp\widehat g\subset\Gamma^n(\sigma^{-(n-1)})\).
Here \(N(s)\) counts the active parents of
\(\{g_\tau\}_{\tau\in\Xi(\sigma)}\).
\end{lemma}

Since the proof of the lemma is rather standard,
 we postpone the proof to Appendix~\ref{subsec:momentconeinterpolation}. Let us begin the proof of Lemma~\ref{lem:curvefine}.

\begin{proof}[Proof of Lemma~\ref{lem:curvefine}]
By standard approximation, we may assume that the
\(h_\theta\) are Schwartz functions.
Using the ball cutoffs from Section~\ref{sec:curvepruning},
dyadic pigeonholing, we have one of two cases: 
\begin{equation}\label{1104.32}
     \int_{\sum_\theta h_\theta>\alpha^2}
 \bigl(\sum_\theta h_\theta\bigr)^{p_n/2}
 \lessapprox
 \int_{\lvert(\sum_\theta h_\theta)*\check\eta_\lambda\rvert
          \gtrapprox \alpha^2}
 \Bigl|\bigl(\sum_\theta h_\theta\bigr)*\check\eta_\lambda\Bigr|^{p_n/2}
\end{equation}
for some $R^{-\frac2n} \ll \lambda < R^{-\frac1n}$,
or
\begin{equation}\label{1104.33}
\int_{\sum_\theta h_\theta>\alpha^2}
 \bigl(\sum_\theta h_\theta\bigr)^{p_n/2}
 \lessapprox    \int_{\lvert(\sum_\theta h_\theta)*\check\eta_{\leq R^{-2/n}}\rvert
          \gtrapprox \alpha^2}
 \Bigl|\bigl(\sum_\theta h_\theta\bigr)
          *\check\eta_{\leq R^{-2/n}}\Bigr|^{p_n/2}.
\end{equation}
The first case and the second case are called the high frequency case and the low frequency case, respectively.
 The high frequency case
is treated by Lemma~\ref{lem:curvconetail} after rescaling; the low frequency case is treated by averaging over balls of radius \(R^{2/n}\) as usual.

\subsubsection*{Case 1: High frequency case}

Every \(\xi\in\supp(\widehat h_\theta\eta_\lambda)\) satisfies \(|\xi|\sim\lambda\) and has the expansion
\[
 \xi=\sum_{i=1}^n a_i\gamma_n^{(i)}(t_\theta),\qquad
 |a_1|\lesssim\lambda,\quad |a_i|\lesssim R^{-i/n}\quad(2\le i\le n).
\]
Under \(\lambda^{-1}\) dilation, the dimensions of
\(\supp(\widehat h_\theta\eta_\lambda)\) change as follows:
\[
 \Big(\lambda\times R^{-2/n}\times\cdots\times R^{-1} \Big)
 \ \xrightarrow{\ \lambda^{-1}\text{ dilation}\ }
 \biggl(1\times\frac{R^{-2/n}}\lambda\times\cdots\times\frac{R^{-1}}\lambda \biggr)
 \ \subset
 \biggl(1\times\frac{R^{-2/n}}\lambda\times\cdots\times
 \Bigl(\frac{R^{-2/n}}\lambda\Bigr)^{n-1} \biggr).
\]
Put \(\sigma=R^{-2/n}/\lambda\) and group the caps $\theta\in\Theta(R^{-1/n})$ according to their parents \(\tau\in\Theta(\sigma)\). Since \(|t_\theta-t_\tau|\le\sigma\),
Taylor expansion from \(t_\theta\) to \(t_\tau\), as in the proof of
Lemma~\ref{lem:smallcapdensity}, gives
\[
 \lambda^{-1}\supp
 \Bigl[\bigl(\sum_{\theta\subset\tau}h_\theta\bigr)
            *\check\eta_\lambda\Bigr]^{\wedge}
 \subset
 \Bigl\{\sum_{i=1}^n b_i\gamma_n^{(i)}(t_\tau):
 |b_1|\sim1,\ |b_i|\lesssim\sigma^{i-1}\ (i\ge2)\Bigr\}.
\]
Thus the Fourier supports of the dilated functions lie in
caps of aperture \(\sigma\), whose active count
at scale \(s\) is at most \(N(s)\).
For \(g(x)=[(\sum_\theta h_\theta)*\check\eta_\lambda](x/\lambda)\),
apply Lemma~\ref{lem:curvconetail} to the family
\(\bigl\{[(\sum_{\theta\subset\tau}h_\theta)*\check\eta_\lambda](x/\lambda)\bigr\}_{\tau\in\Theta(\sigma)}\)
with \(\sigma=R^{-2/n}/\lambda\) and \(\alpha\mapsto\frac{\alpha}{R^\e}\).
On returning to the original coordinates, the averaging boxes \(U\parallel U_{\rho,\sigma^{-(n-1)}}\)
become
\[
 B\parallel\lambda^{-1}U_{\rho,\sigma^{-(n-1)}}
 =\lambda sR^{2/n}V_{\rho,R},
 \qquad \sigma\le s\le1.
\]
Their dimensions are
\(\lambda R^{4/n}s^2\times\lambda R^{4/n}s\times\cdots
\times\lambda R^{4/n}s^{3-n}\).
The same Jacobian occurs in the integral in \eqref{1104.32}
and in \(|U|\).
Thus the right hand side of
\eqref{1104.32} is bounded by
\begin{equation}
 R^\e
 \sum_{\sigma\le s\le1}\sum_{\rho\in\Theta(s)}
 \sum_{\substack{B\parallel\lambda sR^{2/n}V_{\rho,R}\\
 B\in\mathcal G_\rho(\frac{\alpha}{R^\e})}}
 |B|\biggl\|
 \Bigl(\sum_{\tau\subset\rho}
 \Bigl|\bigl(\sum_{\theta\subset\tau}h_\theta\bigr)
                     *\check\eta_\lambda\Bigr|\Bigr)^{1/2}
 \biggr\|_{L^{p_n/2}_\#(B)}^{p_n}.
 \label{eq:curvefineconeoutput}
\end{equation}
Applying H\"older and \eqref{eq:commonreturn} in
Lemma~\ref{lem:selection} to each sum over \(B\) in
\eqref{eq:curvefineconeoutput} bounds \eqref{eq:curvefineconeoutput} by
\begin{equation}
 R^\e
 \sum_{\sigma\le s\le1}\sum_{\rho\in\Theta(s)}
 \sum_{\substack{B\parallel\lambda sR^{2/n}V_{\rho,R}\\
 B\in\mathcal G_\rho(\frac{\alpha}{R^\e})}}
 |B|\biggl\|\bigl(\sum_{\theta\subset\rho}h_\theta\bigr)^{1/2}
 \biggr\|_{L^{p_n/2}_\#(B)}^{p_n}.
 \label{eq:curvefinecleansum}
\end{equation}
Now apply Lemma~\ref{lem:envelopechange} to each sum over \(B\)
in \eqref{eq:curvefinecleansum}, with  \(F=(\sum_{\theta\subset\rho}h_\theta)^{p_n/4}\),
\(B_0=\lambda sR^{2/n}V_{\rho,R}\), and \(U_0=V_{\rho,R}\).
Since \(\lambda sR^{2/n}\ge1\), this bounds
\eqref{eq:curvefinecleansum} by
\begin{equation}
 R^\e
 \sum_{R^{-1/n}\le s\le1}\sum_{\rho\in\Theta(s)}
 \sum_{\substack{V\parallel V_{\rho,R}\\
 V\in\mathcal G_\rho(\frac{\alpha}{R^\e})}}
 |V|\biggl\|\bigl(\sum_{\theta\subset\rho}h_\theta\bigr)^{1/2}
 \biggr\|_{L^{p_n/2}_\#(V)}^{p_n}.
 \label{eq:curvefinefinalsum}
\end{equation}

\subsubsection*{Case 2: Low frequency case}

Tile by
\(Q\parallel V_{\rho,R}\), \(|\rho|=1\).
These boxes are comparable to balls of radius \(R^{2/n}\).
The decay of \(\check\eta_{\leq R^{-2/n}}\) and H\"older give
\[
 \sup_{x\in Q}
 \Bigl|\bigl(\sum_\theta h_\theta\bigr)*
                  \check\eta_{\leq R^{-2/n}}(x)\Bigr|^{p_n/4}
 \lesssim\fint_Q\bigl(\sum_\theta h_\theta\bigr)^{p_n/4}.
\]
Consequently, the right hand side of \eqref{1104.33} is bounded by the \(s=1\)
sum in \eqref{eq:curvefinefinalsum}. This proves \eqref{eq:curvefine}.
\end{proof}

\section{\texorpdfstring{\except{toc}{From cones over moment curves to moment curves: }}{From cones over moment curves to moment curves: }Auxiliary estimates for unpruning}\label{sec:curverefinement}

In this section and the next, we prove the unpruning estimate in Proposition~\ref{prop:curveunpruning}.
Here we establish two auxiliary estimates for the pruned square functions, which we combine in Section~\ref{sec:curveunpruning} to complete the proof.

\subsection{Wave envelope type estimates for the pruned functions}

The following lemma bounds the sum over \(Q\) in
\eqref{eq:curveproducerbox} by sums over \(V\) and \(B\), involving the
original functions \(f_\theta\) and finer pruned groups, respectively.

\begin{lemma}\label{lem:curveproducer}
Let \(1\le j\le N\) and
\(R_j^{-1/n}\le\sigma\le R_{j-1}^{-1/n}\).
For every \(s\in[\sigma,1]\) and \(\alpha>0\),
there exists a dyadic \(\sigma'\in[R^{-1/n},\sigma]\) such that
\begin{equation}
\begin{aligned}
 &\sum_{\rho\in\Theta(s)}
 \sum_{\substack{Q\parallel\rho^*\\
 Q\in\mathcal G_\rho(\alpha)}}
 |Q|\biggl\|
 \Bigl(\sum_{\substack{\tau\subset\rho\\|\tau|=\sigma}}
 |f_\tau^j|^2*w_\tau\Bigr)^{1/2}
 \biggr\|_{L^{p_n/2}_\#(Q)}^{p_n}
 \\
 &\lesssim_{\eta,\e} R^\e
 \sum_{R^{-1/n}\le s'\le s}\sum_{\rho'\in\Theta(s')}
 \sum_{\substack{V\parallel V_{\rho',R}\\
 V\in\mathcal G_{\rho'}(\frac{\alpha}{R^\e})}}
 |V|\biggl\|
 \bigl(\sum_{\theta\subset\rho'}|f_\theta|^2*w_\theta\bigr)^{1/2}
 \biggr\|_{L^{p_n/2}_\#(V)}^{p_n}
 \\
 &\quad+R^\e
 \sum_{\sigma'\le s'\le s}\sum_{\rho'\in\Theta(s')}
 \sum_{\substack{B\parallel(s'/\sigma')V_{\rho',(\sigma')^{-n}}\\
 B\in\mathcal G_{\rho'}(\frac{\alpha}{R^\e})}}
 |B|\biggl\|
 \Bigl(\sum_{\substack{\tau'\subset\rho'\\|\tau'|=\sigma'}}
 |f_{\tau'}^{j'}|^2*w_{\tau'}\Bigr)^{1/2}
 \biggr\|_{L^{p_n/2}_\#(B)}^{p_n}.
\end{aligned}
\label{eq:curveproducerbox}
\end{equation}
Here \(j\le j'\le N\) satisfies
\(R_{j'}^{-1/n}\le\sigma'\le R_{j'-1}^{-1/n}\).

For the original functions \(f_\tau\), the same inequality
holds for \(R^{-1/n}\le\sigma\le s\le1\), with \(j,j'\) omitted.
\end{lemma}

By Lemma~\ref{lem:selectedaverages}, for some fixed \(0<c<1\),
\[
 \text{LHS of }\eqref{eq:curveproducerbox}
 \lesssim
 \sum_{\rho\in\Theta(s)}
 \int_{\sum_{\substack{\tau\subset\rho\\|\tau|=\sigma}}
 |f_\tau^j|^2*w_\tau>(c\alpha/N(s))^2}
 \Bigl(\sum_{\substack{\tau\subset\rho\\|\tau|=\sigma}}
 |f_\tau^j|^2*w_\tau\Bigr)^{p_n/2}.
\]
We use a high low decomposition to the square function. The high frequency part is estimated
by Lemma~\ref{lem:curvconetail} after rescaling; in the low frequency case
we increase the pruning index and repeat the decomposition, applying
Lemma~\ref{lem:curvefine} at \(j=N\).

\begin{proof}
For each \(\rho\in\Theta(s)\), split its superlevel set in the
preceding display. The high frequency region consists of the points where
\[
 \biggl|\Bigl(\sum_{\substack{\tau\subset\rho\\|\tau|=\sigma}}
 |f_\tau^j|^2*w_\tau\Bigr)*\check\eta_{>R_j^{-1/n}}\biggr|(x)
 \ge\frac12
 \Bigl(\sum_{\substack{\tau\subset\rho\\|\tau|=\sigma}}
 |f_\tau^j|^2*w_\tau\Bigr)(x).
\]
On the  low frequency region, we have
\[
 \Bigl(\sum_{\substack{\tau\subset\rho\\|\tau|=\sigma}}
 |f_\tau^j|^2*w_\tau\Bigr)(x)
 \lesssim \biggl|\Bigl(\sum_{\substack{\tau\subset\rho\\|\tau|=\sigma}}
 |f_\tau^j|^2*w_\tau\Bigr)*\check\eta_{\leq R_j^{-1/n}}\biggr|(x).
\]
We distinguish which region contributes at least half of the sum of
integrals. In the high frequency case, decompose this part into dyadic
pieces using \(\eta_{\sigma'}\),
\(R_j^{-1/n}<\sigma'\le\sigma\), and pigeonhole \(\sigma'\) after
summing over \(\rho\). Let us first consider the case that the high frequency part dominates.

\subsubsection{Case 1: High frequency part}

For each \(\rho\in\Theta(s)\), expand every
\(f_\tau^j\) into its \(\sigma'\) pieces and write
\[
 F=\Bigl(\sum_{\substack{\tau\subset\rho\\|\tau|=\sigma}}
 |f_\tau^j|^2*w_\tau\Bigr)*\check\eta_{\sigma'}
 =\sum_{\substack{\tau'\subset\rho\\|\tau'|=\sigma'}}h_{\tau'}.
\]
Here, for \(\tau'\subset\tau\) with \(|\tau|=\sigma\),
\[
 h_{\tau'}=
 \sum_{\substack{\widetilde\tau'\subset\tau,\ |\widetilde\tau'|=\sigma'\\
 |t_{\tau'}-t_{\widetilde\tau'}|\lesssim\sigma'}}
 (f_{\tau'}^j\overline{f_{\widetilde\tau'}^j})
 *w_\tau*\check\eta_{\sigma'}.
\]
The ball cutoff selects only these neighbouring pairs.
Dyadic pigeonholing after summing over \(\rho\) bounds the
left hand side of \eqref{eq:curveproducerbox} by
\(\lessapprox\sum_{\rho\in\Theta(s)}
\int_{|F|\gtrapprox(\alpha/N(s))^2}|F|^{p_n/2}\)
for some \(R_j^{-1/n}<\sigma'\le\sigma\).

Fix \(\rho\in\Theta(s)\) in the preceding sum and let \(t_\rho\) be the left endpoint of the interval indexing \(\rho\).
As in the proof of Lemma~\ref{lem:curvefine}, every \(\xi\in\supp\widehat h_{\tau'}\)
satisfies \(|\xi|\sim\sigma'\) and can be written as
\[
 \xi=\sum_{i=1}^n a_i\gamma_n^{(i)}(t_{\tau'}),\qquad |a_i|\lesssim(\sigma')^i\quad(1\le i\le n).
\]
Since \(F\) is the sum over \(\tau'\subset\rho\), we follow
dilation by \((\sigma')^{-1}\) with the rescaling
\(sL_{t_\rho,s}^{-1}\) that sends \([t_\rho,t_\rho+s]\) to \([0,1]\). The widths of each piece then change as follows:
\[
 \bigl(\sigma'\times(\sigma')^2\times\cdots\times(\sigma')^n\bigr)
 \xrightarrow{\;(\sigma')^{-1}\text{ dilation}\;}
 \bigl(1\times\sigma'\times\cdots\times(\sigma')^{n-1}\bigr)
 \xrightarrow{\;sL_{t_\rho,s}^{-1}\;}
 \biggl(1\times\frac{\sigma'}s\times\cdots\times
       \Bigl(\frac{\sigma'}s\Bigr)^{n-1}\biggr).
\]
Thus the rescaled pieces have Fourier support in caps
of aperture \(\sigma'/s\).
For \(g(x)=F((s/\sigma')L_{t_\rho,s}^{-T}x)\), apply Lemma~\ref{lem:curvconetail}
to the family \(\bigl\{h_{\tau'}((s/\sigma')L_{t_\rho,s}^{-T}x)\bigr\}_{\tau'\subset\rho,\,|\tau'|=\sigma'}\)
at amplitude \(\alpha/(R^\epsilon N(s))\).

In the sum given by Lemma~\ref{lem:curvconetail}, let \(\rho''\)
be a parent of width \(s'/s\) and let
\(U\parallel U_{\rho'',(\sigma'/s)^{-(n-1)}}\).
Returning to the original coordinates sends \(\rho''\) to a cap
\(\rho'\subset\rho\) of width \(s'\) and \(U\) to
\(B=(s/\sigma')L_{t_\rho,s}^{-T}U\).
The \(i\)th side length of \(B\) is
\[
 \frac{s^{1-i}}{\sigma'}
 \left(\frac{\sigma'}s\right)^{-2}
 \left(\frac{s'}s\right)^{3-i}
 =(\sigma')^{-3}(s')^{3-i},\qquad 1\le i\le n.
\]
Thus \(B\parallel(s'/\sigma')V_{\rho',(\sigma')^{-n}}\).

The active parent count at width \(s'/s\) is at most
\(N_\rho(s')\), since each nonzero \(h_{\tau'}\) comes from
a cap \(\tau'\) containing a cap \(\theta\) for which \(f_\theta\not\equiv0\).
Using \(N(s)N_\rho(s')\lesssim_\e N(s')\) from
Lemma~\ref{lem:regular} in the selection condition and summing
over \(\rho\), we obtain
\begin{equation}\label{eq:curveproducerconeoutput}
\begin{aligned}
 \sum_{\rho\in\Theta(s)}
 \int_{|F|\gtrapprox(\alpha/N(s))^2}|F|^{p_n/2}
 &\lesssim_\e R^\e
 \sum_{\sigma'\le s'\le s}\sum_{\rho'\in\Theta(s')}
 \sum_{\substack{B\parallel(s'/\sigma')V_{\rho',(\sigma')^{-n}}\\
 B\in\mathcal G_{\rho'}(\frac{\alpha}{R^\e})}}
 |B|\biggl\|\Bigl(
 \sum_{\substack{\tau'\subset\rho'\\|\tau'|=\sigma'}}
 |h_{\tau'}|\Bigr)^{1/2}\biggr\|_{L^{p_n/2}_\#(B)}^{p_n}.
\end{aligned}
\end{equation}

Since \(\tau'\subset\tau\), we have
\( |f_{\tau'}^j|^2*w_\tau
\lesssim |f_{\tau'}^j|^2*w_{\tau'}\).
Bounding neighbouring products by separate squares and absorbing
neighbouring parents by bounded overlap, we may replace
\(\sum_{\tau'\subset\rho',\,|\tau'|=\sigma'}|h_{\tau'}|\)
on the right of \eqref{eq:curveproducerconeoutput} by
\(
 |\check\eta_{\sigma'}|*
 \sum_{\substack{\tau'\subset\rho'\\|\tau'|=\sigma'}}
 |f_{\tau'}^j|^2*w_{\tau'}.
\)
Apply H\"older and \eqref{eq:commonreturn} to each sum over \(B\) on the right of \eqref{eq:curveproducerconeoutput}, with
\(F=(\sum_{\tau'\subset\rho',\,|\tau'|=\sigma'}
|f_{\tau'}^j|^2*w_{\tau'})^{p_n/4}\) and
\(K=|\check\eta_{\sigma'}|\).
This bounds the right hand side of
\eqref{eq:curveproducerconeoutput} by the second sum on the right
of \eqref{eq:curveproducerbox}, with \(j'=j\).

\subsubsection{Case 2: The low frequency part dominates}

Using \eqref{eq:curvelowcontrol} for the \(\sigma\) pieces
inside \(\rho\) gives
\[
 \biggl|\Bigl(\sum_{\substack{\tau\subset\rho\\|\tau|=\sigma}}
 |f_\tau^j|^2*w_\tau\Bigr)*\check\eta_{\leq R_j^{-1/n}}\biggr|
 \lesssim
 \sum_{\tau_j\subset\rho}|f_{\tau_j}^j|^2*w_{\tau_j}.
\]
For \(j<N\), applying \(|f_{\tau_j}^j|\le|f_{\tau_j}^{j+1}|\)
yields
\begin{equation}\label{eq:curveproducerlowoutput}
 \text{LHS of \eqref{eq:curveproducerbox}}
 \lesssim_{\eta,\e} R^\e\sum_{\rho\in\Theta(s)}
 \sum_{\substack{Q\parallel\rho^*\\
 Q\in\mathcal G_\rho(\frac{\alpha}{R^\e})}}
 |Q|\biggl\|\bigl(
 \sum_{\tau_j\subset\rho}|f_{\tau_j}^{j+1}|^2*w_{\tau_j}
 \bigr)^{1/2}\biggr\|_{L^{p_n/2}_\#(Q)}^{p_n}.
\end{equation}
The sum over \(Q\) on the right of
\eqref{eq:curveproducerlowoutput} is the left hand side of
\eqref{eq:curveproducerbox} with \((j,\sigma,\alpha)\) replaced
by \((j+1,R_j^{-1/n},\alpha/R^\e)\), while \(s\) is unchanged.
Repeat the decomposition into Cases 1 and 2 on this sum over \(Q\)
with these parameters. If the high frequency part in Case 1
dominates at index \(j'\), repeat the argument leading to
\eqref{eq:curveproducerconeoutput} at that index and remove
\(|\check\eta_{\sigma'}|\) as above. This gives
the second sum on the right of \eqref{eq:curveproducerbox}, with
\(R_{j'}^{-1/n}<\sigma'\le R_{j'-1}^{-1/n}\).

If the low frequency part dominates at every step, we repeat Case~2
until the pruning index reaches \(N\). At \(j=N\), we use
\(|f_\theta^N|\le|f_\theta|\) to obtain
\begin{equation}\label{eq:curveproducerterminal}
 \text{LHS of \eqref{eq:curveproducerbox}}
 \lesssim_{\eta,\e} R^\e\sum_{\rho\in\Theta(s)}
 \sum_{\substack{Q\parallel\rho^*\\
 Q\in\mathcal G_\rho(\frac{\alpha}{R^\e})}}
 |Q|\biggl\|\bigl(
 \sum_{\theta\subset\rho}|f_\theta|^2*w_\theta
 \bigr)^{1/2}\biggr\|_{L^{p_n/2}_\#(Q)}^{p_n}.
\end{equation}
It remains to bound the right hand side of
\eqref{eq:curveproducerterminal} by the first sum on the right
of \eqref{eq:curveproducerbox}. For each \(\rho\in\Theta(s)\),
\eqref{eq:selaverage} bounds the sum over \(Q\) in
\eqref{eq:curveproducerterminal} by
\[
 \int_{\sum_{\theta\subset\rho}|f_\theta|^2*w_\theta
       \,\gtrsim\,(\alpha/(R^\e N(s)))^2}
 \Bigl(\sum_{\theta\subset\rho}|f_\theta|^2*w_\theta\Bigr)^{p_n/2}.
\]
The Taylor rescaling taking \(\rho\) to the unit cap sends each \(\theta\subset\rho\) to a cap of \(\mathcal M^n(s^nR)\) at scale \(R^{-1/n}/s\).
We therefore apply Lemma~\ref{lem:curvefine} to
\(\{(|f_\theta|^2*w_\theta)(L_{t_\rho,s}^{-T}\,\cdot)\}_{\theta\subset\rho}\),
with parameter \(s^nR\) and amplitude \(\alpha/(R^\e N(s))\).

Returning to the original coordinates sends a parent of width
\(s'/s\) to a cap \(\rho'\subset\rho\) of width \(s'\), where
\(R^{-1/n}\le s'\le s\). Its averaging box returns to
\(V\parallel V_{\rho',R}\), since its \(i\)th side length is
\[
 s^{-i}(s^nR)^{2/n}(s'/s)^{2-i}
 =R^{2/n}(s')^{2-i},\qquad 1\le i\le n.
\]
Using \(N(s)N_\rho(s')\lesssim_\e N(s')\) from
Lemma~\ref{lem:regular} in the selection condition and summing
over \(\rho\) gives the first sum on the right of
\eqref{eq:curveproducerbox}.

There are at most \(N-j+1\le\eta^{-1}\) steps. Pigeonholing
in \(j'\) and \(\sigma'\), with a sufficiently small loss at
each step, gives \eqref{eq:curveproducerbox}. For the original
functions, repeat the arguments of Cases~1 and~2 with \(f_\tau\)
in place of \(f_\tau^j\), taking the low frequency cutoff at
\(R^{-1/n}\).
\end{proof}

\subsection{Refining either \texorpdfstring{\(\sigma\)}{sigma} or \texorpdfstring{\(s\)}{s}}
Since the first sum on the right of \eqref{eq:curveproducerbox}
already involves the original functions \(f_\theta\), it remains
to estimate the second sum, which involves the pruned pieces.
Lemma~\ref{lem:curvehigh} bounds this sum over \(B\) by two
\(Q\) sums of the form on the left of \eqref{eq:curveproducerbox},
one with \(\sigma\) reduced to \(\sigma'\) and the other with \(s\) reduced to \(s'\).
In the proof of Proposition~\ref{prop:curveunpruning}, we can
therefore apply Lemma~\ref{lem:curveproducer} again at the smaller
scales. Since \(\sigma,s\ge R^{-1/n}\), neither width can decrease
indefinitely. In Section~\ref{sec:curveunpruning} we use this fact
to advance the pruning index \(j\) to \(N\) and bound the sum over
\(B\) in terms of the original functions \(f_\theta\).

\begin{lemma}\label{lem:curvehigh}
Let \(1\le j\le N\) and
\(R_j^{-1/n}\le\sigma\le R_{j-1}^{-1/n}\).
For every \(s\in[\sigma,1]\) and \(\alpha>0\),
\begin{equation}
\begin{aligned}
 &\sum_{\rho\in\Theta(s)}
 \sum_{\substack{B\parallel(s/\sigma)V_{\rho,\sigma^{-n}}\\
 B\in\mathcal G_\rho(\alpha)}}
 |B|\biggl\|
 \Bigl(\sum_{\substack{\tau\subset\rho\\|\tau|=\sigma}}
 |f_\tau^j|^2*w_\tau\Bigr)^{1/2}
 \biggr\|_{L^{p_n/2}_\#(B)}^{p_n}
 \\
 &\lesssim_{\e} R^\e
 \sum_{\rho\in\Theta(s)}
 \sum_{\substack{Q\parallel\rho^*\\
 Q\in\mathcal G_\rho(\frac{\alpha}{R^\e})}}
 |Q|\biggl\|
 \Bigl(\sum_{\substack{\tau'\subset\rho\\|\tau'|=\sigma'}}
 |f_{\tau'}^j|^2*w_{\tau'}\Bigr)^{1/2}
 \biggr\|_{L^{p_n/2}_\#(Q)}^{p_n}
 \\
 &\quad+R^\e
 \sum_{\rho'\in\Theta(s')}
 \sum_{\substack{Q\parallel(\rho')^*\\
 Q\in\mathcal G_{\rho'}(\frac{\alpha}{R^\e})}}
 |Q|\biggl\|
 \Bigl(\sum_{\substack{\tau\subset\rho'\\|\tau|=\sigma}}
 |f_\tau^j|^2*w_\tau\Bigr)^{1/2}
 \biggr\|_{L^{p_n/2}_\#(Q)}^{p_n},
\end{aligned}
\label{eq:curvehighreturn}
\end{equation}
where \(\sigma'=\max\{R_j^{-1/n},\sigma^2/s\}\) and
\(s'=s(\sigma/s)^{1/(n-1)}\in[\sigma,s]\).
For the original grouped functions \(f_\tau\), the same estimate holds
for \(R^{-1/n}\le\sigma\le s\le1\), with \(j\) omitted and
\(\sigma'=\max\{R^{-1/n},\sigma^2/s\}\).
\end{lemma}

The high frequency part of the proof uses the following square function estimate
for the cone over the moment curve.

\begin{proposition}[{\cite[Theorem~1 and Proposition~1.1]{GMSF}}]
\label{prop:cone-square-function}
For a Schwartz function \(g\) with
\(\supp\widehat g\subset\Gamma^n(\sigma^{-(n-1)})\)
and \(2\le p\le p_{n-1}\), we have
\[
 \int_{\R^n}|g|^p
 \lesssim_\e\sigma^{-\e}
 \int_{\R^n}
 \bigl(\sum_{\tau\in\Xi(\sigma)}|g_\tau|^2\bigr)^{p/2}.
\]
\end{proposition}

\begin{proof}[Proof of Lemma~\ref{lem:curvehigh}]
For each \(\rho\in\Theta(s)\), decompose
\(\sum_{\tau\subset\rho,\,|\tau|=\sigma}|f_\tau^j|^2*w_\tau\)
using the cutoffs \(\eta_{\leq\sigma'}\) and \(\eta_\lambda\),
where \(\sigma'<\lambda\lesssim\sigma\).
We split the sum over \(B\) on the left of
\eqref{eq:curvehighreturn} according to whether the
low frequency part controls the average over \(B\).

\subsubsection{Case 1: Low frequency part}

Consider the boxes \(B\in\mathcal G_\rho(\alpha)\),
\(B\parallel(s/\sigma)V_{\rho,\sigma^{-n}}\), for which
\[
 \fint_B\Bigl(
 \sum_{\substack{\tau\subset\rho\\|\tau|=\sigma}}
 |f_\tau^j|^2*w_\tau
 \Bigr)^{p_n/4}
 \lesssim
 \fint_B\biggl|\Bigl(
 \sum_{\substack{\tau\subset\rho\\|\tau|=\sigma}}
 |f_\tau^j|^2*w_\tau
 \Bigr)*\check\eta_{\leq\sigma'}\biggr|^{p_n/4}.
\]
When we expand each \(f_\tau^j\) into its \(\sigma'\) pieces,
the cutoff \(\eta_{\leq\sigma'}\) selects only products of
neighbouring pieces. As in Case~2 of the proof of
Lemma~\ref{lem:curveproducer}, this gives
\[
 \biggl|\Bigl(
 \sum_{\substack{\tau\subset\rho\\|\tau|=\sigma}}
 |f_\tau^j|^2*w_\tau
 \Bigr)*\check\eta_{\leq\sigma'}\biggr|
 \lesssim
 \sum_{\substack{\tau'\subset\rho\\|\tau'|=\sigma'}}
 |f_{\tau'}^j|^2*w_{\tau'}.
\]
Since \(B_0=(s/\sigma)V_{\rho,\sigma^{-n}}=(s/\sigma)^3\rho^*\)
and \(s\ge\sigma\), we apply Lemma~\ref{lem:envelopechange}
to the sum over these \(B\), with
\(F=(\sum_{\tau'\subset\rho,\,|\tau'|=\sigma'}
|f_{\tau'}^j|^2*w_{\tau'})^{p_n/4}\) and \(U_0=\rho^*\).
Summing over \(\rho\) gives the first sum over \(Q\) in \eqref{eq:curvehighreturn}.

\subsubsection{Case 2: High frequency part}

Fix \(\rho\in\Theta(s)\) and let \(t_\rho\) be the left endpoint of the interval indexing \(\rho\).
For each box \(B\in\mathcal G_\rho(\alpha)\) not covered by
Case~1, dyadic pigeonholing gives a width \(\lambda\),
\(\sigma'<\lambda\lesssim\sigma\), for which
\begin{equation}\label{ineqLem5.5}
 \fint_B\Bigl(
 \sum_{\substack{\tau\subset\rho\\|\tau|=\sigma}}
 |f_\tau^j|^2*w_\tau
 \Bigr)^{p_n/4}
 \lessapprox
 \fint_B\biggl|\Bigl(
 \sum_{\substack{\tau\subset\rho\\|\tau|=\sigma}}
 |f_\tau^j|^2*w_\tau
 \Bigr)*\check\eta_\lambda\biggr|^{p_n/4}.
\end{equation}
Fix \(\lambda\) and consider the boxes \(B\) satisfying
\eqref{ineqLem5.5}. To estimate its right hand side, we group
the \(\sigma\) caps in \(\rho\) into parents \(\rho'\) of width
\(s'=s(\sigma/s)^{1/(n-1)}\) and write
\[
 \Bigl(
 \sum_{\substack{\tau\subset\rho\\|\tau|=\sigma}}
 |f_\tau^j|^2*w_\tau
 \Bigr)*\check\eta_\lambda
 =
 \sum_{\substack{\rho'\subset\rho\\|\rho'|=s'}}h_{\rho'},
 \qquad
 h_{\rho'}=
 \Bigl(
 \sum_{\substack{\tau\subset\rho'\\|\tau|=\sigma}}
 |f_\tau^j|^2*w_\tau
 \Bigr)*\check\eta_\lambda.
\]

When we expand each \(f_\tau^j\) into its \(R_j^{-1/n}\) pieces,
the cutoff \(\eta_\lambda\) selects only pairs indexed by intervals at distance \(O(\lambda)\) from each other.
Since \(R_j^{-1/n}\le\lambda\), the Taylor expansion in the proof of Lemma~\ref{lem:curveproducer}
gives, for every \(\xi\) in the Fourier support of \(|f_\tau^j|^2*w_\tau*\check\eta_\lambda\),
\[
 |\xi|\sim\lambda,\qquad \xi=\sum_{i=1}^n a_i\gamma_n^{(i)}(t_\tau),\qquad
 |a_i|\lesssim\lambda\sigma^{i-1}\quad(1\le i\le n).
\]
For \(\tau\subset\rho'\), we have \(\sigma\le s'\) and \(|t_\tau-t_{\rho'}|\le s'\).
Expanding at \(t_{\rho'}\) as in Case~1 of the proof of Lemma~\ref{lem:curvefine}, we obtain
\[
 \xi\in\supp\widehat h_{\rho'}\quad\Longrightarrow\quad
 \xi=\sum_{i=1}^n b_i\gamma_n^{(i)}(t_{\rho'}),\qquad
 |b_i|\lesssim\lambda(s')^{i-1}\quad(1\le i\le n).
\]

As in Lemma~\ref{lem:curveproducer}, we follow dilation by
\(\lambda^{-1}\) with the rescaling \(sL_{t_\rho,s}^{-1}\)
that sends \([t_\rho,t_\rho+s]\) to \([0,1]\).
The widths of each piece change as follows:
\[
 \bigl(\lambda\times\lambda s'\times\cdots\times
       \lambda(s')^{n-1}\bigr)
 \xrightarrow{\;\lambda^{-1}\text{ dilation}\;}
 \bigl(1\times s'\times\cdots\times(s')^{n-1}\bigr)
 \xrightarrow{\;sL_{t_\rho,s}^{-1}\;}
 \biggl(1\times\frac{s'}s\times\cdots\times
       \Bigl(\frac{s'}s\Bigr)^{n-1}\biggr).
\]
Thus the rescaled \(h_{\rho'}\) have Fourier support in
caps of aperture \(s'/s\).

For each \(\rho'\subset\rho\) of width \(s'\), let
\(a_{\rho'}=|\check\eta_\lambda|*
\sum_{\tau\subset\rho',\,|\tau|=\sigma}|f_\tau^j|^2*w_\tau\),
so that \(|h_{\rho'}|\le a_{\rho'}\).
To obtain the second sum on the right of \eqref{eq:curvehighreturn},
we estimate the high frequency average on the right of
\eqref{ineqLem5.5}. We claim that, for each \(B\) satisfying
\eqref{ineqLem5.5},
\begin{equation}\label{eq:curvehighboxtail}
 \biggl[
 \fint_B\biggl|
 \Bigl(\sum_{\substack{\tau\subset\rho\\|\tau|=\sigma}}
 |f_\tau^j|^2*w_\tau\Bigr)*\check\eta_\lambda
 \biggr|^{p_n/4}
 \biggr]^2
 \lesssim_\epsilon \frac{R^\epsilon}{|B|}
 \sum_{\substack{\rho'\subset\rho\\|\rho'|=s'}}
 \int_{a_{\rho'}>[\alpha/(R^\epsilon N(s))]^2}
 a_{\rho'}^{p_n/2}\,W_B .
\end{equation}
Once we prove \eqref{eq:curvehighboxtail}, we combine it with
\eqref{ineqLem5.5}, multiply by \(|B|\), and sum over the boxes for which this \(\lambda\)
was chosen. Since \(\sum_B W_B\lesssim1\), we apply
\eqref{eq:tailYoung} with
\(F=\sum_{\tau\subset\rho',\,|\tau|=\sigma}|f_\tau^j|^2*w_\tau\),
\(K=|\check\eta_\lambda|\), and exponent \(p_n/2\) to remove the
convolution from \(a_{\rho'}\). The locally constant property on
\(Q\parallel(\rho')^*\) and \(N(s')\ge N(s)\) then give the second
sum on the right of \eqref{eq:curvehighreturn} after summing over
\(\rho\) and \(\lambda\).

It remains to prove \eqref{eq:curvehighboxtail}.
When \(n\ge4\), we use Proposition~\ref{prop:cone-square-function}
at \(p=p_n/4\), since \(2\le p_n/4\le p_{n-1}\).
When \(n=3\), Proposition~\ref{prop:cone-square-function} requires
\(p\ge2\), so we first use H\"older on \(B\) to pass from
\(p_3/4=7/4\) to \(p=7/3\).

First suppose \(n\ge4\). Choose \(\chi_B\sim W_B^{4/p_n}\) with
\(\operatorname{supp}\widehat{\chi_B}\subset B^*\).
By the definition of \(h_{\rho'}\), the left hand side of
\eqref{eq:curvehighboxtail} satisfies
\begin{equation}\label{eq:curvehighlocalization}
\begin{aligned}
 \biggl[
 \fint_B\biggl|
 \Bigl(\sum_{\substack{\tau\subset\rho\\|\tau|=\sigma}}
 |f_\tau^j|^2*w_\tau\Bigr)*\check\eta_\lambda
 \biggr|^{p_n/4}
 \biggr]^2
 &\sim
 \biggl[
 \frac{1}{|B|}\int
 \Bigl|\chi_B
 \sum_{\substack{\rho'\subset\rho\\|\rho'|=s'}}h_{\rho'}
 \Bigr|^{p_n/4}
 \biggr]^2 .
\end{aligned}
\end{equation}

To apply Proposition~\ref{prop:cone-square-function} to the integral
on the right of \eqref{eq:curvehighlocalization}, we check the Fourier
support of each \(\chi_Bh_{\rho'}\) after rescaling.
Multiplication by \(\chi_B\) enlarges the Fourier support of
\(h_{\rho'}\) by at most the widths of \(B^*\), namely
\(\sigma^3s^{i-3}\) in direction \(i\).
Under the frequency map \(\xi\mapsto(s/\lambda)L_{t_\rho,s}^{-1}\xi\),
the added width in each direction is at most
\[
 \frac{\sigma^3s^{i-3}}{\lambda s^{i-1}}
 =\frac{\sigma^3}{\lambda s^2}
 \le\frac{\sigma}{s}
 =\left(\frac{s'}s\right)^{n-1},
\]
where we used \(\lambda>\sigma'\ge\sigma^2/s\).
Since \((s'/s)^{n-1}\) is the smallest width of a cap of
aperture \(s'/s\), the Fourier support of
\((\chi_Bh_{\rho'})\bigl((s/\lambda)L_{t_\rho,s}^{-T}\,\cdot\bigr)\)
remains in the corresponding caps for \(\Gamma^n((s/s')^{n-1})\),
up to a fixed enlargement.

We may therefore apply Proposition~\ref{prop:cone-square-function}
at \(p=p_n/4\) to the sum in the integral on the right of
\eqref{eq:curvehighlocalization} after rescaling, that is, to
\(x\mapsto\bigl[\chi_B\sum_{\rho'\subset\rho,\,|\rho'|=s'}h_{\rho'}\bigr]
\bigl((s/\lambda)L_{t_\rho,s}^{-T}x\bigr)\), with aperture \(s'/s\).
On returning to the original coordinates, the square function integrand is
\(\displaystyle |\chi_B|^{p_n/4}
\bigl(\sum_{\rho'\subset\rho,\,|\rho'|=s'}|h_{\rho'}|^2\bigr)^{p_n/8}\).
Using \(|\chi_B|^{p_n/4}\sim W_B\) in \eqref{eq:curvehighlocalization},
we obtain
\begin{equation}\label{inequ2Lem5.5}
\begin{aligned}
 \biggl[
 \fint_B\biggl|
 \Bigl(\sum_{\substack{\tau\subset\rho\\|\tau|=\sigma}}
 |f_\tau^j|^2*w_\tau\Bigr)*\check\eta_\lambda
 \biggr|^{p_n/4}
 \biggr]^2
 &\lesssim_\epsilon R^\epsilon \biggl[
 \fint_B
 \Bigl(\sum_{\substack{\rho'\subset\rho\\|\rho'|=s'}}
 |h_{\rho'}|^2\Bigr)^{p_n/8}
 \biggr]^2 .
\end{aligned}
\end{equation}

On the right of \eqref{inequ2Lem5.5}, use \(|h_{\rho'}|\le a_{\rho'}\)
and split \(\sum_{\rho'}a_{\rho'}^2\) at the threshold
\([\alpha/(R^\epsilon N(s))]^2\) in \eqref{eq:curvehighboxtail}.
The sum below the threshold satisfies
\[
 \sum_{a_{\rho'}\le[\alpha/(R^\epsilon N(s))]^2}a_{\rho'}^2
 \le
 \left(\frac{\alpha}{R^\epsilon N(s)}\right)^2
 \sum_{\rho'}a_{\rho'}.
\]
Raising this inequality to the power \(p_n/8\), averaging over \(B\),
and applying H\"older gives
\[
\begin{aligned}
 \fint_B\Bigl(
 \sum_{a_{\rho'}\le[\alpha/(R^\epsilon N(s))]^2}a_{\rho'}^2
 \Bigr)^{p_n/8}
 &\lesssim
 \left(\frac{\alpha}{R^\epsilon N(s)}\right)^{p_n/4}
 \biggl[\fint_B\Bigl(\sum_{\rho'}a_{\rho'}\Bigr)^{p_n/4}\biggr]^{1/2}\\
 &=
 \left(\frac{\alpha}{R^\epsilon N(s)}\right)^{p_n/4}
 \biggl[\fint_B\Bigl(
 |\check\eta_\lambda|*
 \sum_{\substack{\tau\subset\rho\\|\tau|=\sigma}}
 |f_\tau^j|^2*w_\tau\Bigr)^{p_n/4}\biggr]^{1/2}\\
 &\lesssim
 \left(\frac{\alpha}{R^\epsilon N(s)}\right)^{p_n/4}
 \biggl[\fint_B\Bigl(
 \sum_{\substack{\tau\subset\rho\\|\tau|=\sigma}}
 |f_\tau^j|^2*w_\tau\Bigr)^{p_n/4}\biggr]^{1/2}.
\end{aligned}
\]
In the last inequality, we used \(\lambda^{-1}\le\sigma^{-3}s^2\),
the shortest side of \(B\), so that
\(|\check\eta_\lambda|*W_B\lesssim W_B\), and applied
\eqref{eq:weightedYoung} with
\(F=\sum_{\tau\subset\rho,\,|\tau|=\sigma}|f_\tau^j|^2*w_\tau\),
\(K=|\check\eta_\lambda|\), and \(r=p_n/4\).

Since \(B\in\mathcal G_\rho(\alpha)\), we have
\(\fint_B(\sum_{\tau\subset\rho,\,|\tau|=\sigma}|f_\tau^j|^2*w_\tau)^{p_n/4}
\ge(\alpha/N(s))^{p_n/2}\).
Squaring the preceding bound, we obtain
\[
\begin{aligned}
 \biggl[\fint_B\Bigl(
 \sum_{a_{\rho'}\le[\alpha/(R^\epsilon N(s))]^2}a_{\rho'}^2
 \Bigr)^{p_n/8}\biggr]^2
 &\lesssim R^{-\epsilon p_n/2}
 \biggl[\fint_B\Bigl(
 \sum_{\substack{\tau\subset\rho\\|\tau|=\sigma}}
 |f_\tau^j|^2*w_\tau\Bigr)^{p_n/4}\biggr]^2\\
 &\lesssim_\epsilon R^{-(p_n/2-3)\epsilon}
 \biggl[\fint_B\Bigl(\sum_{\rho'}|h_{\rho'}|^2\Bigr)^{p_n/8}\biggr]^2 ,
\end{aligned}
\]
where we used the square of \eqref{ineqLem5.5} and then
\eqref{inequ2Lem5.5} in the last inequality.
Since \(p_n/2-3>0\), for sufficiently large \(R\) we can absorb
the terms below the threshold and replace \(\sum_{\rho'}|h_{\rho'}|^2\)
on the right of \eqref{inequ2Lem5.5} by
\(\sum_{a_{\rho'}>[\alpha/(R^\epsilon N(s))]^2}a_{\rho'}^2\).

For the remaining terms, we use
\[
 \sum_{a_{\rho'}>[\alpha/(R^\epsilon N(s))]^2}a_{\rho'}^2
 \le
 \Bigl(\sum_{\rho'}a_{\rho'}\Bigr)
 \Bigl(\sum_{a_{\rho'}>[\alpha/(R^\epsilon N(s))]^2}
 a_{\rho'}^{p_n/2}\Bigr)^{2/p_n}.
\]
Raising this inequality to the power \(p_n/8\), applying H\"older
in the \(B\) average, and squaring gives
\[
\begin{aligned}
 &\biggl[\fint_B\Bigl(
 \sum_{a_{\rho'}>[\alpha/(R^\epsilon N(s))]^2}a_{\rho'}^2
 \Bigr)^{p_n/8}\biggr]^2\\
 &\quad\lesssim
 \fint_B\Bigl(\sum_{\rho'}a_{\rho'}\Bigr)^{p_n/4}
 \biggl[
 \fint_B\sum_{a_{\rho'}>[\alpha/(R^\epsilon N(s))]^2}
 a_{\rho'}^{p_n/2}
 \biggr]^{1/2}.
\end{aligned}
\]
As above, \eqref{eq:weightedYoung} bounds
\(\fint_B(\sum_{\rho'}a_{\rho'})^{p_n/4}\) by
\(\fint_B(\sum_{\tau\subset\rho,\,|\tau|=\sigma}
|f_\tau^j|^2*w_\tau)^{p_n/4}\).
By \eqref{ineqLem5.5}, this last average is
\(\lessapprox\fint_B|\sum_{\rho'}h_{\rho'}|^{p_n/4}\).
Substituting the resulting bound for the remaining terms on the right
of \eqref{inequ2Lem5.5}, we obtain
\[
 \biggl[\fint_B\Bigl|\sum_{\rho'}h_{\rho'}\Bigr|^{p_n/4}\biggr]^2
 \lessapprox
 \fint_B\Bigl|\sum_{\rho'}h_{\rho'}\Bigr|^{p_n/4}
 \biggl[
 \fint_B\sum_{a_{\rho'}>[\alpha/(R^\epsilon N(s))]^2}
 a_{\rho'}^{p_n/2}
 \biggr]^{1/2}.
\]
Dividing by \(\fint_B|\sum_{\rho'}h_{\rho'}|^{p_n/4}\)
and squaring gives \eqref{eq:curvehighboxtail} for \(n\ge4\).

Now suppose \(n=3\). Choose \(\chi_B\sim W_B^{3/7}\) with
\(\operatorname{supp}\widehat{\chi_B}\subset B^*\).
H\"older on \(B\), followed by Proposition~\ref{prop:cone-square-function}
at \(p=7/3\) with the same rescaling as in \eqref{inequ2Lem5.5}, gives
\[
 \biggl[\fint_B\Bigl|\sum_{\rho'}h_{\rho'}\Bigr|^{7/4}\biggr]^{4/3}
 \lesssim \frac{1}{|B|}\int
 \Bigl|\chi_B\sum_{\rho'}h_{\rho'}\Bigr|^{7/3}
 \lessapprox \fint_B\Bigl(\sum_{\rho'}|h_{\rho'}|^2\Bigr)^{7/6}.
\]
As in the \(n\ge4\) case, we use \eqref{eq:weightedYoung}, the selection
condition on \(B\), and \eqref{ineqLem5.5} to absorb the terms with
\(a_{\rho'}\le[\alpha/(R^\epsilon N(s))]^2\) into the square function
average on the right. For the terms above this threshold, the same
pointwise bound followed by H\"older with exponents \(3/2\) and \(3\)
gives
\[
 \biggl[\fint_B\Bigl|\sum_{\rho'}h_{\rho'}\Bigr|^{7/4}\biggr]^{4/3}
 \lessapprox
 \biggl[\fint_B\Bigl|\sum_{\rho'}h_{\rho'}\Bigr|^{7/4}\biggr]^{2/3}
 \biggl[\fint_B
 \sum_{a_{\rho'}>[\alpha/(R^\epsilon N(s))]^2}
 a_{\rho'}^{7/2}\biggr]^{1/3}.
\]
Dividing by the \(2/3\) power on the right and cubing gives
\eqref{eq:curvehighboxtail} for \(n=3\).
Together with Case~1, this proves \eqref{eq:curvehighreturn}.
For the original functions, the same argument applies with \(R^{-1/n}\)
in place of \(R_j^{-1/n}\) and \(f_\tau\) in place of \(f_\tau^j\)
throughout. This completes the proof of the lemma.
\end{proof}

\section{\texorpdfstring{\except{toc}{From cones over moment curves to moment curves: }}{From cones over moment curves to moment curves: }Proof of Proposition~\ref{prop:curveunpruning}}\label{sec:curveunpruning}
 
 In this section, we prove Proposition~\ref{prop:curveunpruning}.
We follow the algorithm in \cite[proof of Proposition~4.6]{Mal}
to bound the left hand side of \eqref{eq:curveunprunecost}, namely
\begin{equation}\label{eq:curveunpruneiteration}
 \sum_{\rho\in\Theta(s)}
 \sum_{\substack{B\parallel(s/\sigma)V_{\rho,\sigma^{-n}}\\
 B\in\mathcal G_\rho(\alpha)}}
 |B|\biggl\|
 \Bigl(\sum_{\substack{\tau\subset\rho\\|\tau|=\sigma}}
 |f_\tau^j|^2*w_\tau\Bigr)^{1/2}
 \biggr\|_{L^{p_n/2}_\#(B)}^{p_n}.
\end{equation}
At each step, \(j\) is the pruning index, \(\alpha\) is the amplitude, and caps of width \(\sigma\) are grouped inside caps of width \(s\).

Let us breifly explain the idea.
To obtain the right hand side of \eqref{eq:curveunprunecost},
we seek a bound involving the original functions \(f_\theta\),
of the form
\begin{equation}\label{eq:curveunpruneterminal}
 \sum_{R^{-1/n}\le u\le s}
 \sum_{\rho'\in\Theta(u)}
 \sum_{\substack{V\parallel V_{\rho',R}\\
 V\in\mathcal G_{\rho'}(\alpha)}}
 |V|\biggl\|
 \Bigl(\sum_{\theta\subset\rho'}
 |f_\theta|^2*w_\theta\Bigr)^{1/2}
 \biggr\|_{L^{p_n/2}_\#(V)}^{p_n}.
\end{equation}
Using Lemmas~\ref{lem:curveproducer} and \ref{lem:curvehigh}
with \(0<\epsilon\ll\eta\), together with
\eqref{eq:curveSone} after Taylor rescaling, we either obtain
a sum of the form \eqref{eq:curveunpruneterminal} and stop,
or return to \eqref{eq:curveunpruneiteration} with
\[
 (j,\sigma,s,\alpha)\longmapsto(j',\sigma',s',\alpha'),
 \qquad j'\ge j,\quad \sigma'\le\sigma,\quad s'\le s,
\]
where \(j'>j\), \(\sigma'\ll\sigma\), or \(s'\ll s\).
Since \(\sigma,s\ge R^{-1/n}\), the decreases in \(\sigma\)
and \(s\) cannot continue indefinitely. Continuing the iteration
therefore advances the pruning index to \(j=N\).
The final step then uses \(|f_\theta^N|\le|f_\theta|\)
to obtain the averages of the original functions \(f_\theta\)
in \eqref{eq:curveunpruneterminal}.

The factor outside the sum and the new amplitude \(\alpha'\)
will be specified in each case. We will bound the accumulated
factor by
\(\bigl(R^{C\eta}S_1(R^\eta)^{N-j+C}\bigr)^{p_n}\)
and keep the amplitude in \eqref{eq:curveunpruneterminal}
at least
\(\alpha/\bigl(R^{C\eta}S_1(R^\eta)^{N-j+C}\bigr)\),
with the initial values of \(j\) and \(\alpha\).
Let \(0<\epsilon\ll\eta\) be sufficiently small.
We will show that the iteration has \(O(\eta^{-3})\) steps,
so all the \(R^\epsilon\) losses arising during the iteration
can be absorbed into \(R^{C\eta}\), both outside the sum
and in the amplitude denominator.
\\

Let us now give a detailed proof.
First consider the case that \(\sigma=R^{-1/n}\). In this case, we have \(j=N\).
For each \(\rho\), \eqref{eq:selaverage}, with
\(F=(\sum_{\theta\subset\rho}|f_\theta^N|^2*w_\theta)^{p_n/4}\)
and \(t=(\alpha/N(s))^{p_n/2}\), together with
\(|f_\theta^N|\le|f_\theta|\), bounds its sum over \(B\)
in \eqref{eq:curveunpruneiteration} by
\begin{equation}\label{eq:curveunprunefinetail}
 \int_{\sum_{\theta\subset\rho}|f_\theta|^2*w_\theta
 \gtrsim(\alpha/N(s))^2}
 \bigl(\sum_{\theta\subset\rho}|f_\theta|^2*w_\theta\bigr)^{p_n/2}.
\end{equation}
The Taylor rescaling taking \(\rho\) to the unit cap sends each \(\theta\subset\rho\) to a cap of \(\mathcal M^n(s^nR)\) at scale \(R^{-1/n}/s\).
Applying Lemma~\ref{lem:curvefine} to
\(\{(|f_\theta|^2*w_\theta)(L_{t_\rho,s}^{-T}\,\cdot)\}_{\theta\subset\rho}\),
with parameter \(s^nR\) and amplitude \(\alpha/N(s)\),
then returning to the original coordinates and summing over
\(\rho\), bounds \eqref{eq:curveunpruneiteration} by
\[
 \lesssim_{\epsilon,\eta} R^\epsilon
 \sum_{R^{-1/n}\le u\le s}
 \sum_{\rho'\in\Theta(u)}
 \sum_{\substack{V\parallel V_{\rho',R}\\
 V\in\mathcal G_{\rho'}(\alpha/R^\epsilon)}}
 |V|\biggl\|
 \Bigl(\sum_{\theta\subset\rho'}|f_\theta|^2*w_\theta\Bigr)^{1/2}
 \biggr\|_{L_\#^{p_n/2}(V)}^{p_n},
\]
where we used
\(\frac{\alpha}{(s^nR)^\epsilon N(s)N_\rho(u)}
\gtrsim_\eta\frac{\alpha}{R^\epsilon N(u)}\)
by Lemma~\ref{lem:regular}.
This is \(R^\epsilon\) times \eqref{eq:curveunpruneterminal},
with amplitude \(\alpha/R^\epsilon\), and the iteration stops.

We may therefore assume \(\sigma>R^{-1/n}\).
We estimate \eqref{eq:curveunpruneiteration} separately for
\(s\le R^{\eta^3}\sigma\) and \(s>R^{\eta^3}\sigma\).
In the first range, each cap of width \(s\) contains at most \(R^{\eta^3}\) caps of width \(\sigma\), and we use \eqref{eq:curveSone} after
Taylor rescaling. In the second range, we use
Lemma~\ref{lem:curvehigh} to reduce either \(\sigma\) or \(s\).

\subsection{Case 1: \texorpdfstring{$s\le R^{\eta^3}\sigma$}{s <= R^(eta^3) sigma}}
\leavevmode\par\nobreak
\noindent
If \(\sigma=R_j^{-1/n}\), first use
\(|f_{\tau_j}^j|\le|f_{\tau_j}^{j+1}|\) in
\eqref{eq:curveunpruneiteration} and replace \(j\) by \(j+1\).
We may therefore assume \(\sigma>R_j^{-1/n}\).

As in the derivation of \eqref{eq:curveunprunefinetail},
\eqref{eq:selaverage} bounds the sum over \(B\) in
\eqref{eq:curveunpruneiteration}, for each \(\rho\), by the
integral on the left of \eqref{eq:curvescalarize} below.
Since each \(\rho\) contains at most \(R^{\eta^3}\) caps
of width \(\sigma\), we use H\"older to move the power
\(p_n/2\) inside the sum over \(\tau\), and then apply
\eqref{eq:tailYoung} with \(F=|f_\tau^j|^2\) and \(K=w_\tau\)
to remove the convolution, obtaining 
\begin{equation}\label{eq:curvescalarize}
 \int_{\sum_{\tau\subset\rho,\,|\tau|=\sigma}
 |f_\tau^j|^2*w_\tau\gtrsim(\alpha/N(s))^2}
 \Bigl(\sum_{\substack{\tau\subset\rho\\|\tau|=\sigma}}
 |f_\tau^j|^2*w_\tau\Bigr)^{p_n/2}
 \lesssim R^{\eta^3(p_n/2-1)}
 \sum_{\substack{\tau\subset\rho\\|\tau|=\sigma}}
 \int_{|f_\tau^j|\gtrsim\frac{\alpha}{R^{\eta^3/2}N(s)}}
 |f_\tau^j|^{p_n}.
\end{equation}
Thus \eqref{eq:curveunpruneiteration} is bounded by the sum
over \(\rho\) of the right hand side of \eqref{eq:curvescalarize}.

The Taylor rescaling taking each \(\tau\) to the unit cap
sends the caps at scale \(R_j^{-1/n}\) inside \(\tau\) to caps of \(\mathcal M^n(\sigma^nR_j)\) at scale \(R_j^{-1/n}/\sigma\), where
\(1<\sigma^nR_j\le R^\eta\).
Applying \eqref{eq:curveSone} to each rescaled \(f_\tau^j\),
with parameter \(\sigma^nR_j\), constant \(S_1(R^\eta)\),
and amplitude \(\alpha/(R^{C\eta^3}N(s))\), then returning
to the original coordinates and summing over \(\rho,\tau\),
bounds \eqref{eq:curveunpruneiteration} by
\begin{equation}\label{eq:curveunprunerescaled}
 \bigl(R^{C\eta^3}S_1(R^\eta)\bigr)^{p_n}
 \sum_{R_j^{-1/n}\le u\le\sigma}
 \sum_{\rho'\in\Theta(u)}
 \sum_{\substack{V\parallel V_{\rho',R_j}\\
 V\in\mathcal G_{\rho'}(
 \frac{\alpha}{R^{C\eta^3}S_1(R^\eta)})}}
 |V|\biggl\|
 \Bigl(\sum_{\substack{\tau_j\subset\rho'\\
 |\tau_j|=R_j^{-1/n}}}
 |f_{\tau_j}^j|^2*w_{\tau_j}\Bigr)^{1/2}
 \biggr\|_{L^{p_n/2}_\#(V)}^{p_n},
\end{equation}
where we used
\(\frac{\alpha}{R^{C\eta^3}S_1(R^\eta)N(s)N_\tau(u)}
\gtrsim_\eta
\frac{\alpha}{R^{C\eta^3}S_1(R^\eta)N(u)}\)
by Lemma~\ref{lem:regular}.

We next estimate \eqref{eq:curveunprunerescaled}: for
\(j=N\) we return to the original functions \(f_\theta\), while for
\(j<N\) we first increase the current pruning index.

\subsubsection{Case 1.1: \texorpdfstring{$j=N$}{j=N}}

The sums over \(V\parallel V_{\rho',R}\) in
\eqref{eq:curveunprunerescaled} already involve the fine pieces \(f_\theta^N\).
Using \(|f_\theta^N|\le|f_\theta|\) and \(u\le\sigma\le s\)
therefore bounds \eqref{eq:curveunpruneiteration} by
\(\bigl(R^{C\eta^3}S_1(R^\eta)\bigr)^{p_n}\) times
\eqref{eq:curveunpruneterminal}, with
\(\alpha\mapsto\alpha/(R^{C\eta^3}S_1(R^\eta))\),
and the iteration stops.

\subsubsection{Case 1.2: \texorpdfstring{$j<N$}{j<N}}

After using \(|f_{\tau_j}^j|\le|f_{\tau_j}^{j+1}|\) in
\eqref{eq:curveunprunerescaled}, dyadic pigeonholing reduces
the sum over \(u\) in \eqref{eq:curveunprunerescaled}
to one scale \(u\in[R_j^{-1/n},\sigma]\), with
the logarithmic loss absorbed into \(R^{C\eta^3}\).
Since \(V_{\rho',R_j}=(uR_j^{1/n})^2(\rho')^*\) and
\(uR_j^{1/n}\ge1\), we apply Lemma~\ref{lem:envelopechange}
to the sum over \(V\) for each \(\rho'\), with
\(F=(\sum_{\tau_j\subset\rho'}
|f_{\tau_j}^{j+1}|^2*w_{\tau_j})^{p_n/4}\),
\(B_0=V_{\rho',R_j}\), and \(U_0=(\rho')^*\),
to replace it by a sum over \(Q\parallel(\rho')^*\).
Combining this with \eqref{eq:curveunprunerescaled},
we bound \eqref{eq:curveunpruneiteration} by
\begin{equation}\label{eq:curveunprunenearQ}
 \bigl(R^{C\eta^3}S_1(R^\eta)\bigr)^{p_n}
 \sum_{\rho'\in\Theta(u)}
 \sum_{\substack{Q\parallel(\rho')^*\\
 Q\in\mathcal G_{\rho'}(
 \frac{\alpha}{R^{C\eta^3}S_1(R^\eta)})}}
 |Q|\biggl\|
 \Bigl(\sum_{\substack{\tau_j\subset\rho'\\
 |\tau_j|=R_j^{-1/n}}}
 |f_{\tau_j}^{j+1}|^2*w_{\tau_j}\Bigr)^{1/2}
 \biggr\|_{L^{p_n/2}_\#(Q)}^{p_n}.
\end{equation}
We now apply Lemma~\ref{lem:curveproducer} to
\eqref{eq:curveunprunenearQ}, with
\((j,\sigma,s,\alpha)\) replaced by
\((j+1,R_j^{-1/n},u,
\alpha/(R^{C\eta^3}S_1(R^\eta)))\).
If the first sum over \(V\) on the RHS of
Lemma~\ref{lem:curveproducer} dominates, absorbing the
lemma's \(R^\epsilon\) loss into \(R^{C\eta^3}\) and
using \(u\le\sigma\le s\) bounds
\eqref{eq:curveunpruneiteration} by
\(\bigl(R^{C\eta^3}S_1(R^\eta)\bigr)^{p_n}\) times
\eqref{eq:curveunpruneterminal}, with
\(\alpha\mapsto\alpha/(R^{C\eta^3}S_1(R^\eta))\),
and the iteration stops.

If instead the second sum over \(B\) on the RHS of
Lemma~\ref{lem:curveproducer} dominates, let \(\sigma'\)
and \(j'\) be its cap width and pruning index, and
use dyadic pigeonholing to fix \(s'\in[\sigma',u]\).
Then \eqref{eq:curveunpruneiteration} is bounded by
\begin{equation}\label{eq:curveunprunereturn}
 \bigl(R^{C\eta^3}S_1(R^\eta)\bigr)^{p_n}
 \sum_{\rho''\in\Theta(s')}
 \sum_{\substack{B\parallel(s'/\sigma')V_{\rho'',(\sigma')^{-n}}\\
 B\in\mathcal G_{\rho''}(
 \frac{\alpha}{R^{C\eta^3}S_1(R^\eta)})}}
 |B|\biggl\|
 \Bigl(\sum_{\substack{\tau'\subset\rho''\\|\tau'|=\sigma'}}
 |f_{\tau'}^{j'}|^2*w_{\tau'}\Bigr)^{1/2}
 \biggr\|_{L^{p_n/2}_\#(B)}^{p_n}.
\end{equation}
The sum over \(B\) in \eqref{eq:curveunprunereturn}
has exactly the form \eqref{eq:curveunpruneiteration},
with \((j,\sigma,s,\alpha)\) replaced by
\((j',\sigma',s',\alpha/(R^{C\eta^3}S_1(R^\eta)))\).
Since \(j'\ge j+1\), \(\sigma'\le\sigma\), and \(s'\le s\),
we return to \eqref{eq:curveunpruneiteration} with these
new parameters and repeat the algorithm.

\subsection{Case 2: \texorpdfstring{$s>R^{\eta^3}\sigma$}{s > R^(eta^3) sigma}}
\leavevmode\par\nobreak
\noindent
Applying Lemma~\ref{lem:curvehigh} to
\eqref{eq:curveunpruneiteration}, with the current parameters
\(j,\sigma,s,\alpha\), bounds
\eqref{eq:curveunpruneiteration} by
\begin{equation}\label{eq:curveunprunecase2bound}
\begin{aligned}
 &R^\epsilon
 \sum_{\rho\in\Theta(s)}
 \sum_{\substack{Q\parallel\rho^*\\
 Q\in\mathcal G_\rho(\frac{\alpha}{R^\epsilon})}}
 |Q|\biggl\|
 \Bigl(\sum_{\substack{\tau'\subset\rho\\|\tau'|=\sigma'}}
 |f_{\tau'}^j|^2*w_{\tau'}\Bigr)^{1/2}
 \biggr\|_{L^{p_n/2}_\#(Q)}^{p_n}\\
 &\quad+R^\epsilon
 \sum_{\rho'\in\Theta(s')}
 \sum_{\substack{Q\parallel(\rho')^*\\
 Q\in\mathcal G_{\rho'}(\frac{\alpha}{R^\epsilon})}}
 |Q|\biggl\|
 \Bigl(\sum_{\substack{\tau\subset\rho'\\|\tau|=\sigma}}
 |f_\tau^j|^2*w_\tau\Bigr)^{1/2}
 \biggr\|_{L^{p_n/2}_\#(Q)}^{p_n},
\end{aligned}
\end{equation}
where \(\sigma'=\max\{R_j^{-1/n},\sigma^2/s\}\) and
\(s'=s(\sigma/s)^{1/(n-1)}\).
We distinguish two cases according to which sum in
\eqref{eq:curveunprunecase2bound} dominates.
Case 2.1 treats the second sum, at scales \((\sigma,s')\), while Case 2.2 treats the first sum, at scales \((\sigma',s)\).

\subsubsection{Case 2.1: Smaller parents}

Suppose that the second sum in
\eqref{eq:curveunprunecase2bound} dominates.
We apply Lemma~\ref{lem:curveproducer} to this sum, with
\((j,\sigma,s,\alpha)\) replaced by
\((j,\sigma,s',\alpha/R^\epsilon)\).
If the first sum over \(V\) on the RHS of
Lemma~\ref{lem:curveproducer} dominates, this bounds
\eqref{eq:curveunpruneiteration} by \(R^\epsilon\) times
\eqref{eq:curveunpruneterminal}, with
\(\alpha\mapsto\alpha/R^\epsilon\), and the iteration stops.

If instead the second sum over \(B\) on the RHS of
Lemma~\ref{lem:curveproducer} dominates, let \(\sigma''\)
and \(j'\) be its cap width and pruning index, and apply
dyadic pigeonholing to select one \(u\in[\sigma'',s']\).
This bounds \eqref{eq:curveunpruneiteration} by
\(R^\epsilon\) times \eqref{eq:curveunpruneiteration},
with \((j,\sigma,s,\alpha)\) replaced by
\((j',\sigma'',u,\alpha/R^\epsilon)\).
Since \(j'\ge j\), \(\sigma''\le\sigma\), and
\(u\le s'\le sR^{-\eta^3/(n-1)}<s\), we return to
\eqref{eq:curveunpruneiteration} with these new parameters
and repeat the algorithm.

\subsubsection{Case 2.2: Finer caps}

Now suppose that the first sum in
\eqref{eq:curveunprunecase2bound} dominates.

First suppose that \(\sigma'>R_j^{-1/n}\).
We apply Lemma~\ref{lem:curveproducer} to the first sum in
\eqref{eq:curveunprunecase2bound}, with
\((j,\sigma,s,\alpha)\) replaced by
\((j,\sigma',s,\alpha/R^\epsilon)\).
If the first sum over \(V\) on the RHS of
Lemma~\ref{lem:curveproducer} dominates, this bounds
\eqref{eq:curveunpruneiteration} by \(R^\epsilon\) times
\eqref{eq:curveunpruneterminal}, with
\(\alpha\mapsto\alpha/R^\epsilon\), and the iteration stops.

If instead the second sum over \(B\) on the RHS of
Lemma~\ref{lem:curveproducer} dominates, let \(\sigma''\)
and \(j'\) be its cap width and pruning index, and apply
dyadic pigeonholing to select one \(u\in[\sigma'',s]\).
This bounds \eqref{eq:curveunpruneiteration} by
\(R^\epsilon\) times \eqref{eq:curveunpruneiteration},
with \((j,\sigma,s,\alpha)\) replaced by
\((j',\sigma'',u,\alpha/R^\epsilon)\).
Since \(j'\ge j\),
\(\sigma''\le\sigma'=\sigma^2/s\le\sigma R^{-\eta^3}<\sigma\),
and \(u\le s\), we return to
\eqref{eq:curveunpruneiteration} with these new parameters
and repeat the algorithm.

Next suppose that \(\sigma'=R_j^{-1/n}\) and \(j<N\).
Using \(|f_{\tau_j}^j|\le|f_{\tau_j}^{j+1}|\) in the first
sum of \eqref{eq:curveunprunecase2bound}, we apply
Lemma~\ref{lem:curveproducer} with
\((j,\sigma,s,\alpha)\) replaced by
\((j+1,R_j^{-1/n},s,\alpha/R^\epsilon)\).
If the first sum over \(V\) on the RHS of
Lemma~\ref{lem:curveproducer} dominates, this bounds
\eqref{eq:curveunpruneiteration} by \(R^\epsilon\) times
\eqref{eq:curveunpruneterminal}, with
\(\alpha\mapsto\alpha/R^\epsilon\), and the iteration stops.

If instead the second sum over \(B\) on the RHS of
Lemma~\ref{lem:curveproducer} dominates, let \(\sigma''\)
and \(j'\) be its cap width and pruning index, and apply
dyadic pigeonholing to select one \(u\in[\sigma'',s]\).
This bounds \eqref{eq:curveunpruneiteration} by
\(R^\epsilon\) times \eqref{eq:curveunpruneiteration},
with \((j,\sigma,s,\alpha)\) replaced by
\((j',\sigma'',u,\alpha/R^\epsilon)\).
Since \(j'\ge j+1\), \(\sigma''\le\sigma\), and \(u\le s\),
we return to \eqref{eq:curveunpruneiteration} with these
new parameters and repeat the algorithm.

Finally, suppose that \(\sigma'=R_j^{-1/n}=R^{-1/n}\).
We use \(|f_\theta^N|\le|f_\theta|\) in the first sum of
\eqref{eq:curveunprunecase2bound}.
For each \(\rho\), \eqref{eq:selaverage} then bounds its
sum over \(Q\) by \eqref{eq:curveunprunefinetail}, with
\(\alpha\mapsto\alpha/R^\epsilon\).
Repeating the Taylor rescaling and application of
Lemma~\ref{lem:curvefine} used to estimate
\eqref{eq:curveunprunefinetail} bounds
\eqref{eq:curveunpruneiteration} by \(R^\epsilon\) times
\eqref{eq:curveunpruneterminal}, with
\(\alpha\mapsto\alpha/R^\epsilon\), and the iteration stops.

\subsection{Termination of the iteration}
At every continued step, \(s\) and \(\sigma\) do not increase
and \(j\) does not decrease. Case 2.1 decreases \(s\) by a
fixed power of \(R^{\eta^3}\), while Case 2.2 decreases
\(\sigma\) by such a power or increases \(j\).
Case 1 increases \(j\) after applying \(S_1(R^\eta)\), or
terminates. Since both widths lie in \([R^{-1/n},1]\), there
are \(O_n(\eta^{-3})\) steps in all and at most \(N-j+1\)
applications of \(S_1(R^\eta)\), where \(j\) is the initial
index in the statement.
The factors \(R^\e\) from the refinement steps and
\(R^{C\eta^3}\) from the \(S_1\) steps are bounded together
by \(R^{C\eta}\), after choosing \(\e\) sufficiently small
in terms of \(\eta\) and using \(N=\eta^{-1}\).
Thus the accumulated coefficient is at most
\(\bigl(R^{C\eta}S_1(R^\eta)^{N-j+C}\bigr)^{p_n}\),
and the amplitude in \eqref{eq:curveunpruneterminal} is at least
\(\alpha/(R^{C\eta}S_1(R^\eta)^{N-j+C})\), with the initial
\(\alpha\). This proves \eqref{eq:curveunprunecost}.

\section{\texorpdfstring{\except{toc}{From cones over moment curves to moment curves: }}{From cones over moment curves to moment curves: }Proof of Proposition~\ref{lem:curveSoneRecurrence}.}\label{sec:curvefirstrecurrence}

As outlined in Section~\ref{sec:curve}, we first prove Proposition~\ref{prop:curvebroadpieces}.
We then combine this with Proposition~\ref{prop:curveunpruning} and apply a standard broad narrow analysis to prove Proposition~\ref{lem:curveSoneRecurrence}.

\subsection{Proof of Proposition~\ref{prop:curvebroadpieces}}
\begin{proof}[\unskip\nopunct]
We split \(\operatorname{Br}_{\alpha,\beta}^{K}\) into its
intersections with \(\Omega_k\), \(N_0\le k<N\), and the
remaining part contained in \(L\). On each \(\Omega_k\), we
use the multilinear estimate and the amplitude dependent wave envelope type
estimate for the cone over the moment curve in \(\R^n\) (Theorem~\ref{thm:main-cone}); on \(L\),
we use Lemma~\ref{lem:curvefine}.

\subsubsection*{Case 1: \texorpdfstring{Bounding \(\operatorname{Br}_{\alpha,\beta}^{K}\cap\Omega_k\)}{Omega k}}

Fix \(N_0\le k<N\) and estimate
\(\operatorname{Br}_{\alpha,\beta}^{K}\cap\Omega_k\).
By the two defining inequalities of
\(\operatorname{Br}_{\alpha,\beta}^{K}\), we have
\(|f_{\tau_i}(x)|\ge\alpha/(100K^{100n^2})\),
\(1\le i\le n\).
Since \(K^{-1}\ge R_1^{-1/n}\), after choosing \(A\)
sufficiently large, \eqref{eq:curvepruneerror} gives
\(|f_{\tau_i}^{k+1}|\ge|f_{\tau_i}|/2\) on
\(\operatorname{Br}_{\alpha,\beta}^{K}\cap\Omega_k\),
and hence
\(|f_{\tau_1}^{k+1}\cdots f_{\tau_n}^{k+1}|^{1/n}
\ge\alpha/(200K^{100n})\) there.

Let \(B\) be a ball of radius \(R_k^{1/n}\), and fix
\(x_B\in\operatorname{Br}_{\alpha,\beta}^{K}
\cap\Omega_k\cap B\).
By the multilinear restriction estimate used in
\cite[proof of Proposition~6.1]{GMSF}, we obtain
\begin{equation}\label{eq:curvemultilinear}
 \alpha^{2n}
 |\operatorname{Br}_{\alpha,\beta}^{K}\cap\Omega_k\cap B|
 \lesssim_\eta R^{C\eta}|B|
 \biggl[
 \Bigl(
 \sum_{\tau_k}|f_{\tau_k}^{k+1}|^2*w_{\tau_k}
 \Bigr)(x_B)
 \biggr]^n.
\end{equation}
On \(\Omega_k\), the high frequency estimate
\eqref{eq:curvelowcontrol} from the pruning construction gives
\[
 \biggl|
 \Bigl(
 \sum_{\tau_k}|f_{\tau_k}^{k+1}|^2*w_{\tau_k}
 \Bigr)*\check\eta_{>R_{k+1}^{-1/n}}
 \biggr|
 \gtrsim
 \sum_{\tau_k}|f_{\tau_k}^{k+1}|^2*w_{\tau_k}
 \ge A^{N-k}\beta.
\]
Summing \eqref{eq:curvemultilinear} over the balls \(B\)
intersecting
\(\operatorname{Br}_{\alpha,\beta}^{K}\cap\Omega_k\)
and using the locally constant property, we obtain
\begin{equation}\label{eq:curvebroadhigh}
\begin{aligned}
 \alpha^{2n}
 |\operatorname{Br}_{\alpha,\beta}^{K}\cap\Omega_k|
 &\lesssim_\eta R^{C\eta}
 \sum_B|B|
 \biggl[
 \Bigl(
 \sum_{\tau_k}|f_{\tau_k}^{k+1}|^2*w_{\tau_k}
 \Bigr)(x_B)
 \biggr]^n\\
 &\lesssim_\eta R^{C\eta}
 (A^{N-k}\beta)^{n-p_{n-1}}
 \sum_B|B|
 \biggl[
 \fint_B
 \biggl|
 \Bigl(
 \sum_{\tau_k}|f_{\tau_k}^{k+1}|^2*w_{\tau_k}
 \Bigr)*\check\eta_{>R_{k+1}^{-1/n}}
 \biggr|^{p_{n-1}/2}
 \biggr]^2\\
 &\lesssim_\eta R^{C\eta}
 (A^{N-k}\beta)^{n-p_{n-1}}
 \int_{\left|
 (\sum_{\tau_k}|f_{\tau_k}^{k+1}|^2*w_{\tau_k})
 *\check\eta_{>R_{k+1}^{-1/n}}
 \right|\gtrsim A^{N-k}\beta}
 \biggl|
 \Bigl(
 \sum_{\tau_k}|f_{\tau_k}^{k+1}|^2*w_{\tau_k}
 \Bigr)*\check\eta_{>R_{k+1}^{-1/n}}
 \biggr|^{p_{n-1}},
\end{aligned}
\end{equation}
where we used \(p_{n-1}\ge n\), and the last inequality
follows from \eqref{eq:selaverage}.

We now estimate the last integral in
\eqref{eq:curvebroadhigh} by decomposing its high frequency
cutoff into the dyadic cutoffs \(\eta_\sigma\).
Dyadic pigeonholing gives a fixed width
\(R_{k+1}^{-1/n}<\sigma\le R_k^{-1/n}\) for which the RHS
of \eqref{eq:curvebroadhigh} is bounded by
\begin{equation}\label{eq:curvebroadshell}
 R^{C\eta}(A^{N-k}\beta)^{n-p_{n-1}}
 \int_{|\sum_\tau h_\tau|>
 \frac{A^{N-k}\beta}{R^{C\eta}}}
 \Bigl|\sum_\tau h_\tau\Bigr|^{p_{n-1}},
\end{equation}
where, for \(\tau\in\Theta(\sigma)\),
\[
 h_\tau=
 \Bigl(
 \sum_{\substack{
 \widetilde\tau\subset\tau_k,\,
 |\widetilde\tau|=\sigma\\
 |t_\tau-t_{\widetilde\tau}|\lesssim\sigma}}
 f_\tau^{k+1}\overline{f_{\widetilde\tau}^{k+1}}
 \Bigr)*w_{\tau_k}*\check\eta_\sigma,
\]
and \(\tau_k\) is the unique cap in \(\Theta(R_k^{-1/n})\) containing \(\tau\).

As in the proof of Lemma~\ref{lem:curveproducer}, dilation
by \(\sigma^{-1}\) sends the Fourier supports of \(h_\tau\)
into caps of aperture \(\sigma\) in \(\Gamma^n(\sigma^{-(n-1)})\).
We therefore apply the amplitude dependent wave envelope type estimate for the
cone over the moment curve in \(\R^n\) (Theorem~\ref{thm:main-cone}) to
\(\sum_\tau h_\tau(\,\cdot/\sigma)\), with parameter
\(\sigma^{-(n-1)}\), exponent \(p_{n-1}\), and amplitude
\(A^{N-k}\beta/R^{C\eta}\).
Returning to the original coordinates gives
\begin{equation}\label{eq:curveconeoutput}
\begin{aligned}
 &\int_{|\sum_\tau h_\tau|>
 \frac{A^{N-k}\beta}{R^{C\eta}}}
 \Bigl|\sum_\tau h_\tau\Bigr|^{p_{n-1}}\lesssim_\eta R^{C\eta}
 \sum_{\sigma\le s\le1}
 \sum_{\rho\in\Theta(s)}
 \sum_{\substack{
 B\parallel(s/\sigma)V_{\rho,\sigma^{-n}}\\
 \fint_B(\sum_{\tau\subset\rho}|h_\tau|^2)^{p_{n-1}/4}\\
 \ge[\frac{A^{N-k}\beta}{R^{C\eta}N(s)}]^{p_{n-1}/2}}}
 |B|
 \biggl[
 \fint_B
 \Bigl(\sum_{\tau\subset\rho}|h_\tau|^2\Bigr)^{p_{n-1}/4}
 \biggr]^2,
\end{aligned}
\end{equation}
where \(s=|\rho|\) is the dyadic scale in the estimate for the cone over the moment curve, and the averaging boxes return to
\(B\parallel\sigma^{-1}U_{\rho,\sigma^{-(n-1)}}
=(s/\sigma)V_{\rho,\sigma^{-n}}\), with dimensions
\(\sigma^{-3}s^2\times\sigma^{-3}s\times\cdots
\times\sigma^{-3}s^{3-n}\).

To bound the RHS of \eqref{eq:curveconeoutput} in terms
of \(f_\tau^{k+1}\), we use the pruning \(L^\infty\) bound
\eqref{eq:curveprunebounds}.
Each \(f_\tau^{k+1}\) is a sum of \(O(R^{\eta/n})\)
\(R_{k+1}^{-1/n}\) children, so
\(\|f_\tau^{k+1}\|_\infty
\lesssim_\eta R^{C\eta}A^{N-k}\beta/\alpha\).
Since the cutoff \(\eta_\sigma\) leaves only products
\(f_\tau^{k+1}\overline{f_{\widetilde\tau}^{k+1}}\) with
\(|t_\tau-t_{\widetilde\tau}|\lesssim\sigma\le s\),
the terms with \(\tau\subset\rho\),
\(\rho\in\Theta(s)\), involve only \(\widetilde\tau\)
within distance \(O(\sigma)\) of \(\rho\).
Writing \(\widetilde\rho\sim\rho\) for the \(O(1)\)
neighbouring parents \(\widetilde\rho\in\Theta(s)\)
with \(|t_{\widetilde\rho}-t_\rho|\lesssim s\), we have
\begin{equation}\label{eq:curvebroadpointwise}
\begin{aligned}
 \sum_{\tau\subset\rho}|h_\tau|^2
 &\lesssim_\eta R^{C\eta}
 \left(\frac{A^{N-k}\beta}{\alpha}\right)^2
 \Bigl[
 |\check\eta_\sigma|*
 \sum_{\widetilde\rho\sim\rho}
 \sum_{\tau\subset\widetilde\rho}
 |f_\tau^{k+1}|^2*w_\tau
 \Bigr],\\
 \sum_{\tau\subset\rho}|h_\tau|^2
 &\lesssim
 \Bigl[
 |\check\eta_\sigma|*
 \sum_{\widetilde\rho\sim\rho}
 \sum_{\tau\subset\widetilde\rho}
 |f_\tau^{k+1}|^2*w_\tau
 \Bigr]^2,
\end{aligned}
\end{equation}
where we used
\(|f_\tau^{k+1}|^2*w_{\tau_k}
\lesssim |f_\tau^{k+1}|^2*w_\tau\).
Raising the first inequality in
\eqref{eq:curvebroadpointwise} to the power
\((p_{n-1}-n)/4\) and the second to \(n/4\),
then multiplying, gives
\begin{equation}\label{eq:curvebroadinterp}
\begin{aligned}
 \Bigl(\sum_{\tau\subset\rho}|h_\tau|^2\Bigr)^{p_{n-1}/4}
 \lesssim_\eta {}&
 R^{C\eta}
 \left(\frac{A^{N-k}\beta}{\alpha}\right)^{(p_{n-1}-n)/2}
 \Bigl[
 |\check\eta_\sigma|*
 \sum_{\widetilde\rho\sim\rho}
 \sum_{\tau\subset\widetilde\rho}
 |f_\tau^{k+1}|^2*w_\tau
 \Bigr]^{p_n/4},
\end{aligned}
\end{equation}
where we used \(p_{n-1}+n=p_n\).
The first inequality in \eqref{eq:curvebroadpointwise}
and H\"older also show that each \(B\) on the RHS of
\eqref{eq:curveconeoutput} satisfies
\begin{equation}\label{eq:curvebroadselector}
 \biggl\|
 \Bigl(
 |\check\eta_\sigma|*
 \sum_{\widetilde\rho\sim\rho}
 \sum_{\tau\subset\widetilde\rho}
 |f_\tau^{k+1}|^2*w_\tau
 \Bigr)^{1/2}
 \biggr\|_{L_\#^{p_n/2}(B)}
 \ge \frac{\alpha}{R^{C\eta}N(s)}.
\end{equation}

Combining \eqref{eq:curvebroadhigh}--\eqref{eq:curvebroadselector}
and using bounded overlap of the neighbouring parents,
we obtain
\begin{equation}\label{eq:curvebroadselected}
\begin{aligned}
 \alpha^{p_n}
 |\operatorname{Br}_{\alpha,\beta}^{K}\cap\Omega_k|
 &\lesssim_\eta R^{C\eta}
 \sum_{\sigma\le s\le1}
 \sum_{\rho\in\Theta(s)}
 \sum_{\substack{
 B\parallel(s/\sigma)V_{\rho,\sigma^{-n}}\\
 B\in\mathcal G_\rho(\frac{\alpha}{R^{C\eta}})}}
 |B|
 \biggl\|
 \Bigl(
 |\check\eta_\sigma|*
 \sum_{\tau\subset\rho}|f_\tau^{k+1}|^2*w_\tau
 \Bigr)^{1/2}
 \biggr\|_{L_\#^{p_n/2}(B)}^{p_n},
\end{aligned}
\end{equation}
where we used
\(\alpha^{p_n-2n}(A^{N-k}\beta)^{n-p_{n-1}}
(A^{N-k}\beta/\alpha)^{p_{n-1}-n}=1\).

To remove the remaining convolution with
\(|\check\eta_\sigma|\) from the RHS of
\eqref{eq:curvebroadselected}, we apply H\"older and
\eqref{eq:commonreturn}, with
\(F=(\sum_{\tau\subset\rho}|f_\tau^{k+1}|^2*w_\tau)^{p_n/4}\)
and \(K=|\check\eta_\sigma|\). This gives
\begin{equation}\label{eq:curvebroadreturned}
\begin{aligned}
 \alpha^{p_n}
 |\operatorname{Br}_{\alpha,\beta}^{K}\cap\Omega_k|
 &\lesssim_\eta R^{C\eta}
 \sum_{\sigma\le s\le1}
 \sum_{\rho\in\Theta(s)}
 \sum_{\substack{
 B\parallel(s/\sigma)V_{\rho,\sigma^{-n}}\\
 B\in\mathcal G_\rho(\frac{\alpha}{R^{C\eta}})}}
 |B|
 \biggl\|
 \Bigl(
 \sum_{\tau\subset\rho}|f_\tau^{k+1}|^2*w_\tau
 \Bigr)^{1/2}
 \biggr\|_{L_\#^{p_n/2}(B)}^{p_n}.
\end{aligned}
\end{equation}
This proves the first inequality in
Proposition~\ref{prop:curvebroadpieces}.

\subsubsection*{Case 2: Bounding \(L\)}

The pruning error estimate \eqref{eq:curvepruneerror},
followed by Cauchy--Schwarz over
\(\tau_{N_0}\in\Theta(R_{N_0}^{-1/n})\) and
\eqref{eq:curvehighsets}, gives
\[
 \alpha^2
 \lesssim R_{N_0}^{1/n}
 \Bigl(
 \sum_{\tau_{N_0}}
 |f_{\tau_{N_0}}^{N_0+1}|^2*w_{\tau_{N_0}}
 \Bigr)
 \lesssim_\eta R^{C\eta N_0}\beta
 \qquad\text{on }L.
\]
Since \(L\subset E_{\alpha,\beta}\), we have
\(\beta\le\sum_\theta|f_\theta|^2*w_\theta\) on \(L\),
and hence
\[
 \alpha^{p_n}|L|
 \lesssim_\eta R^{C\eta N_0}
 \int_{\sum_\theta|f_\theta|^2*w_\theta>
 (\frac{\alpha}{R^{C\eta N_0}})^2}
 \Bigl(\sum_\theta|f_\theta|^2*w_\theta\Bigr)^{p_n/2}.
\]
Applying Lemma~\ref{lem:curvefine} to
\(\{|f_\theta|^2*w_\theta\}_{\theta\in\Theta(R^{-1/n})}\),
with parameter \(R\), amplitude
\(\alpha/R^{C\eta N_0}\), and loss \(\eta\),
bounds this integral and gives
\begin{equation}\label{eq:curveLfinal}
 \alpha^{p_n}|L|
 \lesssim_\eta R^{C\eta N_0}
 \sum_{R^{-1/n}\le s\le1}
 \sum_{\rho\in\Theta(s)}
 \sum_{\substack{V\parallel V_{\rho,R}\\
 V\in\mathcal G_\rho(\frac{\alpha}{R^{C\eta N_0}})}}
 |V|
 \biggl\|
 \Bigl(
 \sum_{\theta\subset\rho}|f_\theta|^2*w_\theta
 \Bigr)^{1/2}
 \biggr\|_{L_\#^{p_n/2}(V)}^{p_n}.
\end{equation}
This proves the second inequality in
\eqref{eq:curvebroadreturned or}, completing the proof of
Proposition~\ref{prop:curvebroadpieces}.
\end{proof}

\begin{lemma}\label{lem:curveweakstrong}
Let \(\Lambda(R)\ge1\) be a constant such that, for every
\(1\le r\le R\), every Schwartz function \(g\) with
\(\supp\widehat g\subset\mathcal M^n(r)\), and every
\(\alpha>0\),
\[
 \alpha^{p_n}|\{|g|>\alpha\}|
 \le\Lambda(R)^{p_n}
 \sum_{r^{-1/n}\le s\le1}\sum_{\rho\in\Theta(s)}
 \sum_{\substack{V\parallel V_{\rho,r}\\
 V\in\mathcal G_\rho(\frac{\alpha}{\Lambda(R)})}}
 |V|\biggl\|\bigl(
 \sum_{\tau\subset\rho}|g_\tau|^2*w_\tau
 \bigr)^{1/2}\biggr\|_{L^{p_n/2}_\#(V)}^{p_n}.
\]
Then the same functions satisfy
\[
 \int_{|g|>\alpha}|g|^{p_n}
 \lesssim_\e R^\e\Lambda(R)^{p_n}
 \sum_{r^{-1/n}\le s\le1}\sum_{\rho\in\Theta(s)}
 \sum_{\substack{V\parallel V_{\rho,r}\\
 V\in\mathcal G_\rho(\frac{\alpha}{R^\e\Lambda(R)})}}
 |V|\biggl\|\bigl(
 \sum_{\tau\subset\rho}|g_\tau|^2*w_\tau
 \bigr)^{1/2}\biggr\|_{L^{p_n/2}_\#(V)}^{p_n}.
\]
\end{lemma}

Since the proof uses standard dyadic pigeonholing,
we postpone it to Appendix~\ref{subsec:curveweakstrong}.

\subsection{Proof of Proposition~\ref{lem:curveSoneRecurrence}}
\begin{proof}[\unskip\nopunct]
We first bound \(\alpha^{p_n}|\{|f|>\alpha\}|\), then
deduce the estimate for \(\int_{|f|>\alpha}|f|^{p_n}\)
in \eqref{eq:curveSone}.
For the region where
\(\sum_\theta|f_\theta|^2*w_\theta>\alpha^2\),
applying Lemma~\ref{lem:curvefine} to
\(\{|f_\theta|^2*w_\theta\}_{\theta\in\Theta(R^{-1/n})}\),
with parameter \(R\), amplitude \(\alpha\), and loss
\(\eta^{10}\), gives
\begin{equation}\label{eq:curvebroadlargesquare}
\begin{aligned}
 \alpha^{p_n}\biggl|
 \Bigl\{|f|>\alpha,\ 
 \sum_\theta|f_\theta|^2*w_\theta>\alpha^2\Bigr\}
 \biggr|
 &\le
 \int_{\sum_\theta|f_\theta|^2*w_\theta>\alpha^2}
 \Bigl(\sum_\theta|f_\theta|^2*w_\theta\Bigr)^{p_n/2}\\
 &\lesssim_\eta R^{C\eta^{10}}
 \sum_{R^{-1/n}\le s\le1}
 \sum_{\rho\in\Theta(s)}
 \sum_{\substack{V\parallel V_{\rho,R}\\
 V\in\mathcal G_\rho(\frac{\alpha}{R^{C\eta^{10}}})}}
 |V|
 \biggl\|
 \Bigl(\sum_{\theta\subset\rho}|f_\theta|^2*w_\theta\Bigr)^{1/2}
 \biggr\|_{L^{p_n/2}_\#(V)}^{p_n}.
\end{aligned}
\end{equation}

On the remaining region
\(\{|f|>\alpha,\ \sum_\theta|f_\theta|^2*w_\theta\le\alpha^2\}\),
the dyadic decomposition leading to \eqref{beta reduction}
reduces us to bounding \(\alpha^{p_n}|E_{\alpha,\beta}|\)
for each dyadic \(\beta\).
Fix such a \(\beta\), and choose a dyadic \(K\sim R^\kappa\),
where \(0<\kappa\ll\eta/(1000n^2)\).
The standard broad narrow decomposition gives
\begin{equation}\label{eq:curvebroad}
 |f|\lesssim
 \max_{\tau\in\Theta(K^{-1})}|f_\tau|
 +K^{100n}
 \max_{\substack{\tau_1,\ldots,\tau_n\\
 K^{-1}\text{ separated}}}
 |f_{\tau_1}\cdots f_{\tau_n}|^{1/n}.
\end{equation}
Together with the definition of
\(\operatorname{Br}_{\alpha,\beta}^{K}\), this gives
\[
 E_{\alpha,\beta}\subset
 \bigcup_{\tau_0\in\Theta(K^{-1})}
 \left\{|f_{\tau_0}|>\frac{\alpha}{100n}\right\}
 \,\cup\!
 \bigcup_{\substack{\tau_1,\ldots,\tau_n\\
 K^{-1}\text{ separated}}}
 \operatorname{Br}_{\alpha,\beta}^{K}.
\]

We estimate \(\alpha^{p_n}|E_{\alpha,\beta}|\) by separating
the narrow part, where
\(\max_{\tau_0\in\Theta(K^{-1})}|f_{\tau_0}|>\alpha/(100n)\),
from its complement, the broad part.
In the narrow part, we apply \eqref{eq:curveSone} to each
\(f_{\tau_0}\) after Taylor rescaling, with constant
\(S_1(R/K^n)\); in the broad part, we apply
\eqref{eq:curvebroadbound} to each
\(\operatorname{Br}_{\alpha,\beta}^{K}\).

\subsubsection*{Case 1: Narrow part}

For each \(\tau_0\in\Theta(K^{-1})\), Taylor rescaling
to the unit cap sends each \(\theta\subset\tau_0\) to a cap of \(\mathcal M^n(R/K^n)\) at scale \(K R^{-1/n}\).
Applying \eqref{eq:curveSone} to the rescaled \(f_{\tau_0}\),
with parameter \(R/K^n\) and amplitude \(\alpha/(100n)\),
then returning to the original coordinates gives
\begin{equation}\label{eq:curvebroadnarrow}
\begin{aligned}
 &\alpha^{p_n}\left|
 \left\{|f_{\tau_0}|>\frac{\alpha}{100n}\right\}
 \right|\\
 &\lesssim S_1(R/K^n)^{p_n}
 \sum_{R^{-1/n}\le s\le K^{-1}}
 \sum_{\substack{\rho\in\Theta(s)\\\rho\subset\tau_0}}
 \sum_{\substack{V\parallel V_{\rho,R}\\
 V\in\mathcal G_\rho(
 \frac{\alpha}{R^{C\eta^{10}}S_1(R/K^n)})}}
 |V|
 \biggl\|
 \Bigl(\sum_{\theta\subset\rho}|f_\theta|^2*w_\theta\Bigr)^{1/2}
 \biggr\|_{L^{p_n/2}_\#(V)}^{p_n}.
\end{aligned}
\end{equation}
Here the averaging box corresponding to \(\rho\subset\tau_0\),
\(|\rho|=s\), returns with \(i\)th width
\(K^i(R/K^n)^{2/n}(Ks)^{2-i}=R^{2/n}s^{2-i}\),
as required for \(V_{\rho,R}\), and we used
\(N_{\tau_0}(s)\le N(s)\) in the selection condition.

Summing \eqref{eq:curvebroadnarrow} over \(\tau_0\) and
using \(S_1(R/K^n)\le S_1(R^{1-n\kappa})\), we obtain
\begin{equation}\label{eq:curvebroadnarrowglobal}
\begin{aligned}
 &\alpha^{p_n}\biggl|
 \Bigl\{\max_{\tau_0\in\Theta(K^{-1})}|f_{\tau_0}|
 >\frac{\alpha}{100n}\Bigr\}
 \biggr|\\
 &\lesssim_\eta
 \bigl(R^{C\eta^{10}}S_1(R^{1-n\kappa})\bigr)^{p_n}
 \sum_{R^{-1/n}\le s\le1}
 \sum_{\rho\in\Theta(s)}
 \sum_{\substack{V\parallel V_{\rho,R}\\
 V\in\mathcal G_\rho(
 \frac{\alpha}{R^{C\eta^{10}}S_1(R^{1-n\kappa})})}}
 |V|
 \biggl\|
 \Bigl(\sum_{\theta\subset\rho}|f_\theta|^2*w_\theta\Bigr)^{1/2}
 \biggr\|_{L^{p_n/2}_\#(V)}^{p_n}.
\end{aligned}
\end{equation}

\subsubsection*{Case 2: Broad part}

Apply \eqref{eq:curvebroadbound} to bound
\(\alpha^{p_n}|\operatorname{Br}_{\alpha,\beta}^{K}|\)
for each separated tuple in \eqref{eq:curvebroad}.
There are at most \(K^n\) tuples and logarithmically many
\(\beta\) in \eqref{beta reduction}.
Summing these bounds over the tuples and dyadic \(\beta\) bounds
\(\alpha^{p_n}|\{|f|>\alpha,\ \sum_\theta|f_\theta|^2*w_\theta
\le\alpha^2,\ \max_{\tau_0}|f_{\tau_0}|\le\alpha/(100n)\}|\)
by
\begin{equation}\label{eq:curvebroadglobal}
\begin{aligned}
 &R^{C\eta N_0}
 \sum_{R^{-1/n}\le s\le1}\sum_{\rho\in\Theta(s)}
 \sum_{\substack{V\parallel V_{\rho,R}\\
 V\in\mathcal G_\rho(\frac{\alpha}{R^{C\eta N_0}})}}
 |V|\biggl\|\bigl(
 \sum_{\theta\subset\rho}|f_\theta|^2*w_\theta
 \bigr)^{1/2}\biggr\|_{L^{p_n/2}_\#(V)}^{p_n}\\
 &\quad+
 \bigl(R^{C\eta}S_1(R^\eta)^{N-N_0+C}\bigr)^{p_n}
 \sum_{R^{-1/n}\le s\le1}\sum_{\rho\in\Theta(s)}
 \sum_{\substack{V\parallel V_{\rho,R}\\
 V\in\mathcal G_\rho(\frac{\alpha}{R^{C\eta}S_1(R^\eta)^{N-N_0+C}})}}
 |V|\biggl\|\bigl(
 \sum_{\theta\subset\rho}|f_\theta|^2*w_\theta
 \bigr)^{1/2}\biggr\|_{L^{p_n/2}_\#(V)}^{p_n}.
\end{aligned}
\end{equation}
Here \(K^n\) is absorbed using \(\kappa\ll\eta\).

\subsubsection*{Completion of the strong type estimate}

To combine \eqref{eq:curvebroadlargesquare},
\eqref{eq:curvebroadnarrowglobal}, and \eqref{eq:curvebroadglobal}, put
\[
 \Lambda(R)=\max\Bigl\{
 R^{C\eta N_0},\ R^{C\eta}S_1(R^\eta)^{N-N_0+C},\
 R^{C\eta^{10}}S_1(R^{1-n\kappa})\Bigr\}.
\]
The uniformization reduction following Lemma~\ref{lem:regular}
gives, for every Schwartz function \(g\) with
\(\supp\widehat g\subset \mathcal M^n(R)\),
\begin{equation}\label{eq:curveweakrecurrence}
 \alpha^{p_n}|\{|g|>\alpha\}|
 \lesssim_\eta\Lambda(R)^{p_n}
 \sum_{R^{-1/n}\le s\le1}\sum_{\rho\in\Theta(s)}
 \sum_{\substack{V\parallel V_{\rho,R}\\
 V\in\mathcal G_\rho(\frac{\alpha}{\Lambda(R)})}}
 |V|\biggl\|\bigl(
 \sum_{\theta\subset\rho}|g_\theta|^2*w_\theta
 \bigr)^{1/2}\biggr\|_{L^{p_n/2}_\#(V)}^{p_n}.
\end{equation}

The same bound holds for every \(r\le R\), with
\(\Lambda(r)\le\Lambda(R)\).
Applying Lemma~\ref{lem:curveweakstrong} to this estimate proves
Proposition~\ref{lem:curveSoneRecurrence}.
\end{proof}

\section{\texorpdfstring{\except{toc}{From cones over moment curves to moment curves: }}{From cones over moment curves to moment curves: }Completion of the proof of Proposition~\ref{prop:conetocurve}}\label{sec:curvesecondrecurrence}

In this section, we complete the proof of Proposition~\ref{prop:conetocurve}
by establishing Proposition~\ref{prop:curveStworecurrence}.
Under the hypothesis of Proposition~\ref{prop:conetocurve} on the cone over the moment curve,
this proposition asserts
\[
 S_2(R)\lessapprox S_1(R^{1/n})S_2(CR^{1-1/n}).
\]
Combined with Corollary~\ref{prop:curveSonebound}, this recurrence yields
\(S_2(R)\lessapprox1\) by iteration.

To prove this recurrence, we first apply \(S_1\) at parameter \(R^{1/n}\).
Throughout this section, write \(\sigma_0=R^{-1/n^2}\),
the width of the caps in this application.
We keep the convolution weights attached to these initial caps
as the caps are subdivided.
Their dual boxes have longest side \(R^{1/n}\), so the locally constant property of the original functions \(f_\theta\)
allows us to remove the weights at the fine scale.
For a cap of width \(\sigma\le\sigma_0\), rescaling gives the smaller
parameter \(\sigma^nR\le R^{1-1/n}\), at which we can apply \(S_2\).

\subsection{Auxiliary estimates}

We establish three auxiliary estimates for the proof of
Proposition~\ref{prop:curveStworecurrence}.
For each cap \(\tau\) with \(|\tau|\le\sigma_0\), let
\(\operatorname{par}_{\sigma_0}(\tau)\) denote the cap in
\(\Theta(\sigma_0)\) containing \(\tau\).
The proofs of the last two lemmas follow those of
Lemmas~\ref{lem:curveproducer} and~\ref{lem:curvehigh};
we give the details needed to preserve the weights attached
to the initial \(\sigma_0\) caps.

\begin{lemma}\label{lem:curveS2unsmooth}
Let \(f\) be a Schwartz function with
\(\supp\widehat f\subset\mathcal M^n(R)\).
For every \(\epsilon>0\), \(\alpha>0\), and
\(\rho\in\Theta(s)\) with \(R^{-1/n}\le s\le1\), we have
\begin{equation}\label{eq:curveunsmooth}
\begin{aligned}
 &\sum_{\substack{V\parallel V_{\rho,R}\\V\in\mathcal G_\rho(\alpha)}}
 |V|\biggl\|\bigl(\sum_{\theta\subset\rho}|f_\theta|^2*
 w_{\operatorname{par}_{\sigma_0}(\theta)}\bigr)^{1/2}
 \biggr\|_{L^{p_n/2}_\#(V)}^{p_n}\lesssim_\e R^\e
 \sum_{\substack{V\parallel V_{\rho,R}\\
 V\in\mathcal G_\rho(\frac{\alpha}{R^\e})}}
 |V|\biggl\|\bigl(\sum_{\theta\subset\rho}|f_\theta|^2\bigr)^{1/2}
 \biggr\|_{L^{p_n/2}_\#(V)}^{p_n}.
\end{aligned}
\end{equation}
\end{lemma}
\begin{proof}
Since \(\supp\widehat{|f_\theta|^2}\subset B(0,CR^{-1/n})\)
and the dual boxes of the initial \(\sigma_0\) caps have
longest side \(R^{1/n}\), the locally constant property gives
\(\sum_{\theta\subset\rho}|f_\theta|^2*
w_{\operatorname{par}_{\sigma_0}(\theta)}
\lesssim w_{B_{R^{1/n}}}*
\sum_{\theta\subset\rho}|f_\theta|^2\).
H\"older followed by \eqref{eq:commonreturn}, with kernel
\(w_{B_{R^{1/n}}}\),
\(F=(\sum_{\theta\subset\rho}|f_\theta|^2)^{p_n/4}\)
and \(t=(\alpha/N(s))^{p_n/2}\), gives
\eqref{eq:curveunsmooth}.
\end{proof}

\begin{lemma}\label{lem:curveS2producer}
Let \(f\) be a Schwartz function with
\(\supp\widehat f\subset\mathcal M^n(R)\), and let
\(R^{-1/n}\le\sigma\le\sigma_0\).
For every \(s\in[\sigma,1]\), \(\epsilon>0\), and \(\alpha>0\),
there exists a dyadic \(\sigma'\in[R^{-1/n},\sigma]\) such that
\begin{equation}\label{eq:curveS2producer}
\begin{aligned}
 &\sum_{\rho\in\Theta(s)}
 \sum_{\substack{Q\parallel\rho^*\\Q\in\mathcal G_\rho(\alpha)}}
 |Q|\biggl\|\Bigl(\sum_{\substack{\tau\subset\rho\\|\tau|=\sigma}}
 |f_\tau|^2*w_{\operatorname{par}_{\sigma_0}(\tau)}\Bigr)^{1/2}
 \biggr\|_{L^{p_n/2}_\#(Q)}^{p_n}
 \\
 &\lesssim_\e R^\e
 \sum_{R^{-1/n}\le u\le s}\sum_{\rho'\in\Theta(u)}
 \sum_{\substack{V\parallel V_{\rho',R}\\
 V\in\mathcal G_{\rho'}(\frac{\alpha}{R^\e})}}
 |V|\biggl\|\bigl(\sum_{\theta\subset\rho'}|f_\theta|^2*
 w_{\operatorname{par}_{\sigma_0}(\theta)}\bigr)^{1/2}
 \biggr\|_{L^{p_n/2}_\#(V)}^{p_n}\\
 &\quad+R^\e
 \sum_{\sigma'\le u\le s}\sum_{\rho'\in\Theta(u)}
 \sum_{\substack{B\parallel(u/\sigma')V_{\rho',(\sigma')^{-n}}\\
 B\in\mathcal G_{\rho'}(\frac{\alpha}{R^\e})}}
 |B|\biggl\|\Bigl(\sum_{\substack{\tau'\subset\rho'\\|\tau'|=\sigma'}}
 |f_{\tau'}|^2*w_{\operatorname{par}_{\sigma_0}(\tau')}\Bigr)^{1/2}
 \biggr\|_{L^{p_n/2}_\#(B)}^{p_n}.
\end{aligned}
\end{equation}
\end{lemma}
\begin{proof}
As in the proof of Lemma~\ref{lem:curveproducer}, for each
\(\rho\in\Theta(s)\), decompose
\(\sum_{\tau\subset\rho,\,|\tau|=\sigma}|f_\tau|^2*
w_{\operatorname{par}_{\sigma_0}(\tau)}\)
using \(\eta_{\le R^{-1/n}}\) and \(\eta_{\sigma'}\),
\(R^{-1/n}<\sigma'\le\sigma\).
We will bound the high frequency terms by the second sum
on the RHS of \eqref{eq:curveS2producer}, and the
low frequency term by the first.

\medskip
\noindent\textit{Case 1: The high frequency part dominates.}
After summing the high frequency terms over \(\rho\),
dyadic pigeonholing gives a single width \(\sigma'\)
whose contribution bounds the total high frequency
contribution up to an \(R^\epsilon\) loss.
We follow Case~1 of the proof of
Lemma~\ref{lem:curveproducer}, with \(f_\tau\) and
\(w_{\operatorname{par}_{\sigma_0}(\tau)}\)
in place of \(f_\tau^j\) and \(w_\tau\).

Convolution with this weight does not enlarge Fourier support.
Moreover, neighbouring \(\sigma'\) caps
\(\tau',\tau''\subset\tau\) satisfy
\(\operatorname{par}_{\sigma_0}(\tau')
=\operatorname{par}_{\sigma_0}(\tau'')
=\operatorname{par}_{\sigma_0}(\tau)\),
so Cauchy--Schwarz gives
\[
 \Bigl|(f_{\tau'}\overline{f_{\tau''}})
 *w_{\operatorname{par}_{\sigma_0}(\tau)}
 *\check\eta_{\sigma'}\Bigr|
 \lesssim
 |\check\eta_{\sigma'}|*
 \Bigl[
 (|f_{\tau'}|^2+|f_{\tau''}|^2)
 *w_{\operatorname{par}_{\sigma_0}(\tau)}
 \Bigr].
\]
Thus the rescaling and subsequent estimates in
Case~1 of the proof of Lemma~\ref{lem:curveproducer}
apply with these weights.
Removing \(|\check\eta_{\sigma'}|\) by H\"older and
\eqref{eq:commonreturn}, as in that case, gives the second
sum on the RHS of \eqref{eq:curveS2producer}.

\medskip
\noindent\textit{Case 2: The low frequency part dominates.}
Expand \(f_\tau=\sum_{\theta\subset\tau}f_\theta\).
The cutoff \(\eta_{\le R^{-1/n}}\) selects only products
of neighbouring fine pieces.
Since
\(w_{\operatorname{par}_{\sigma_0}(\theta)}
=w_{\operatorname{par}_{\sigma_0}(\tau)}\)
for \(\theta\subset\tau\), we obtain
\[
 \biggl|
 \Bigl(
 \sum_{\substack{\tau\subset\rho\\|\tau|=\sigma}}
 |f_\tau|^2*w_{\operatorname{par}_{\sigma_0}(\tau)}
 \Bigr)*\check\eta_{\le R^{-1/n}}
 \biggr|
 \lesssim
 |\check\eta_{\le R^{-1/n}}|*
 \sum_{\theta\subset\rho}|f_\theta|^2*
 w_{\operatorname{par}_{\sigma_0}(\theta)}.
\]
Applying \eqref{eq:selaverage} to the sum over \(Q\)
on the LHS of \eqref{eq:curveS2producer}, followed by
\eqref{eq:tailYoung} with
\(F=\sum_{\theta\subset\rho}|f_\theta|^2*
w_{\operatorname{par}_{\sigma_0}(\theta)}\),
kernel \(|\check\eta_{\le R^{-1/n}}|\), and exponent
\(p_n/2\), gives
\begin{equation}\label{eq:curveS2finetail}
 \text{LHS of \eqref{eq:curveS2producer}}
  \lessapprox
 \sum_{\rho\in\Theta(s)}
 \int_{\sum_{\theta\subset\rho}|f_\theta|^2*
 w_{\operatorname{par}_{\sigma_0}(\theta)}
 \,\gtrsim[\alpha/(R^\epsilon N(s))]^2}
 \Bigl(
 \sum_{\theta\subset\rho}|f_\theta|^2*
 w_{\operatorname{par}_{\sigma_0}(\theta)}
 \Bigr)^{p_n/2}.
\end{equation}

To estimate each integral in \eqref{eq:curveS2finetail},
we repeat the Taylor rescaling and application of
Lemma~\ref{lem:curvefine} following
\eqref{eq:curveunprunefinetail}.
For each \(\rho\), let \(t_\rho\) be the left endpoint of the interval indexing \(\rho\).
The rescaled family
\(\{(|f_\theta|^2*
w_{\operatorname{par}_{\sigma_0}(\theta)})
(L_{t_\rho,s}^{-T}\,\cdot)\}_{\theta\subset\rho}\)
satisfies the Fourier support hypothesis of
Lemma~\ref{lem:curvefine} at parameter \(s^nR\).
Applying that lemma with amplitude
\(\alpha/(R^\epsilon N(s))\), returning to the original
coordinates, and summing over \(\rho\) bounds
\eqref{eq:curveS2finetail} by the first sum on the RHS of
\eqref{eq:curveS2producer}, where we used
\(N(s)N_\rho(u)\lesssim_\epsilon N(u)\)
for \(\rho'\subset\rho\), \(|\rho'|=u\).
This proves \eqref{eq:curveS2producer}.
\end{proof}

\begin{lemma}\label{lem:curveS2highreturn}
Let \(f\) be a Schwartz function with
\(\supp\widehat f\subset\mathcal M^n(R)\), and let
\(R^{-1/n}\le\sigma\le\sigma_0\).
For every \(s\in[\sigma,1]\), \(\epsilon>0\), and \(\alpha>0\), we have
\begin{equation}\label{eq:curveS2highreturn}
\begin{aligned}
 &\sum_{\rho\in\Theta(s)}
 \sum_{\substack{B\parallel(s/\sigma)V_{\rho,\sigma^{-n}}\\
 B\in\mathcal G_\rho(\alpha)}}
 |B|\biggl\|\Bigl(\sum_{\substack{\tau\subset\rho\\|\tau|=\sigma}}
 |f_\tau|^2*w_{\operatorname{par}_{\sigma_0}(\tau)}\Bigr)^{1/2}
 \biggr\|_{L^{p_n/2}_\#(B)}^{p_n}
 \\
 &\lesssim_\e R^\e
 \sum_{\rho\in\Theta(s)}
 \sum_{\substack{Q\parallel\rho^*\\
 Q\in\mathcal G_\rho(\frac{\alpha}{R^\e})}}
 |Q|\biggl\|\Bigl(\sum_{\substack{\tau'\subset\rho\\|\tau'|=\sigma'}}
 |f_{\tau'}|^2*w_{\operatorname{par}_{\sigma_0}(\tau')}\Bigr)^{1/2}
 \biggr\|_{L^{p_n/2}_\#(Q)}^{p_n}\\
 &\quad+R^\e
 \sum_{\rho'\in\Theta(s')}
 \sum_{\substack{Q\parallel(\rho')^*\\
 Q\in\mathcal G_{\rho'}(\frac{\alpha}{R^\e})}}
 |Q|\biggl\|\Bigl(\sum_{\substack{\tau\subset\rho'\\|\tau|=\sigma}}
 |f_\tau|^2*w_{\operatorname{par}_{\sigma_0}(\tau)}\Bigr)^{1/2}
 \biggr\|_{L^{p_n/2}_\#(Q)}^{p_n},
\end{aligned}
\end{equation}
where \(\sigma'=\max\{R^{-1/n},\sigma^2/s\}\) and
\(s'=s(\sigma/s)^{1/(n-1)}\in[\sigma,s]\).
\end{lemma}
\begin{proof}
As in the proof of Lemma~\ref{lem:curvehigh}, for each
\(\rho\in\Theta(s)\), decompose
\(\sum_{\tau\subset\rho,\,|\tau|=\sigma}|f_\tau|^2*
w_{\operatorname{par}_{\sigma_0}(\tau)}\)
using \(\eta_{\le\sigma'}\) and \(\eta_\lambda\), where
\(\sigma'<\lambda\lesssim\sigma\), and separate the boxes
\(B\in\mathcal G_\rho(\alpha)\) according to whether their
averages are controlled by the low frequency part or by a
high frequency term.
We will bound these contributions by the first and second
sums on the RHS of \eqref{eq:curveS2highreturn}, respectively.

\medskip
\noindent\textit{Case 1: The low frequency part dominates.}
As in Case~2 of the proof of Lemma~\ref{lem:curveS2producer},
expand each \(f_\tau\) into its \(\sigma'\) pieces.
Since
\(\operatorname{par}_{\sigma_0}(\tau')
=\operatorname{par}_{\sigma_0}(\tau)\)
for \(\tau'\subset\tau\), we obtain
\[
 \biggl|
 \Bigl(
 \sum_{\substack{\tau\subset\rho\\|\tau|=\sigma}}
 |f_\tau|^2*w_{\operatorname{par}_{\sigma_0}(\tau)}
 \Bigr)*\check\eta_{\le\sigma'}
 \biggr|
 \lesssim
 |\check\eta_{\le\sigma'}|*
 \sum_{\substack{\tau'\subset\rho\\|\tau'|=\sigma'}}
 |f_{\tau'}|^2*w_{\operatorname{par}_{\sigma_0}(\tau')}.
\]
Applying H\"older and \eqref{eq:commonreturn} to the sum
over these \(B\), with kernel \(|\check\eta_{\le\sigma'}|\)
and
\(F=(\sum_{\tau'\subset\rho,\,|\tau'|=\sigma'}
|f_{\tau'}|^2*w_{\operatorname{par}_{\sigma_0}(\tau')})^{p_n/4}\),
removes the convolution with \(|\check\eta_{\le\sigma'}|\).
Since
\(B_0=(s/\sigma)V_{\rho,\sigma^{-n}}=(s/\sigma)^3\rho^*\)
and \(s\ge\sigma\), Lemma~\ref{lem:envelopechange},
applied with this \(F\), \(B_0\), and \(U_0=\rho^*\),
bounds the contribution of these boxes by the first sum
on the RHS of \eqref{eq:curveS2highreturn}.

\medskip
\noindent\textit{Case 2: The high frequency part dominates.}
For each remaining \(B\), dyadic pigeonholing gives a width
\(\lambda\), \(\sigma'<\lambda\lesssim\sigma\), for which
the term with \(\eta_\lambda\) controls its average.
Fix \(\lambda\) and consider the boxes with this choice.

We follow Case~2 of the proof of Lemma~\ref{lem:curvehigh},
replacing \(|f_\tau^j|^2*w_\tau\) throughout by
\(|f_\tau|^2*w_{\operatorname{par}_{\sigma_0}(\tau)}\).
Convolution with the inherited weight does not enlarge
Fourier support, so the rescaling and the
Fourier support calculation for multiplication by \(\chi_B\)
remain valid.
We use exponent \(p_n/4\) when \(n\ge4\) and \(7/3\) when
\(n=3\), as in that proof.
The calculation leading to \eqref{eq:curvehighboxtail}
therefore gives that inequality with the same replacement.

Multiplying this inequality by \(|B|\) and summing over
the boxes with this \(\lambda\), we apply
\eqref{eq:tailYoung} with
\(F=\sum_{\tau\subset\rho',\,|\tau|=\sigma}
|f_\tau|^2*w_{\operatorname{par}_{\sigma_0}(\tau)}\),
kernel \(|\check\eta_\lambda|\), and exponent \(p_n/2\).
The locally constant property on \(Q\parallel(\rho')^*\), together with
\(N(s')\ge N(s)\), then bounds the contribution of these
boxes to the LHS of \eqref{eq:curveS2highreturn} by
\[
 R^\epsilon
 \sum_{\substack{\rho'\subset\rho\\|\rho'|=s'}}
 \sum_{\substack{Q\parallel(\rho')^*\\
 Q\in\mathcal G_{\rho'}(\frac{\alpha}{R^\epsilon})}}
 |Q|
 \biggl\|
 \Bigl(
 \sum_{\substack{\tau\subset\rho'\\|\tau|=\sigma}}
 |f_\tau|^2*w_{\operatorname{par}_{\sigma_0}(\tau)}
 \Bigr)^{1/2}
 \biggr\|_{L_\#^{p_n/2}(Q)}^{p_n}.
\]
Summing over \(\rho\) and \(\lambda\) gives the second sum
on the RHS of \eqref{eq:curveS2highreturn}.
Together with Case 1, this proves \eqref{eq:curveS2highreturn}.
\end{proof}

\subsection{Proof of Proposition~\ref{prop:curveStworecurrence}}
\begin{proof}[\unskip\nopunct]
We combine the three preceding lemmas in an iteration analogous
to that in Section~\ref{sec:curveunpruning}.
Fix \(\eta>0\) sufficiently small and choose \(0<\epsilon\ll\eta\).
Let \(f\) be a Schwartz function with
\(\supp\widehat f\subset\mathcal M^n(R)\), and let
\(\alpha>0\). Applying \eqref{eq:curveSone} to \(f\)
at parameter \(R^{1/n}\), with fine pieces
\(\{f_\tau\}_{\tau\in\Theta(\sigma_0)}\), and then
Lemma~\ref{lem:envelopechange} to the sums over
\(V\parallel V_{\rho,R^{1/n}}\), gives
\begin{equation}\label{eq:curveS2initial}
 \int_{|f|>\alpha}|f|^{p_n}
 \lesssim S_1(R^{1/n})^{p_n}
 \sum_{\sigma_0\le s\le1}
 \sum_{\rho\in\Theta(s)}
 \sum_{\substack{Q\parallel\rho^*\\
 Q\in\mathcal G_\rho(\frac{\alpha}{S_1(R^{1/n})})}}
 |Q|\biggl\|
 \Bigl(
 \sum_{\substack{\tau\subset\rho\\|\tau|=\sigma_0}}
 |f_\tau|^2*w_\tau
 \Bigr)^{1/2}
 \biggr\|_{L_\#^{p_n/2}(Q)}^{p_n}.
\end{equation}
Here we used
\(V_{\rho,R^{1/n}}=(s/\sigma_0)^2\rho^*\)
and \(s/\sigma_0\ge1\).

Dyadic pigeonholing reduces the sum over \(s\) in
\eqref{eq:curveS2initial} to one scale
\(s\in[\sigma_0,1]\). For this \(s\), we apply
Lemma~\ref{lem:curveS2producer} to the corresponding sum
over \(\rho\) and \(Q\), with \(\sigma=\sigma_0\) and
amplitude \(\alpha/S_1(R^{1/n})\).

If the first sum over \(V\) on the RHS of
Lemma~\ref{lem:curveS2producer} dominates, applying
Lemma~\ref{lem:curveS2unsmooth} to that sum gives
\begin{equation}\label{eq:curveS2terminal}
 \int_{|f|>\alpha}|f|^{p_n}
 \lesssim_\epsilon R^{C\epsilon} S_1(R^{1/n})^{p_n}
 \sum_{R^{-1/n}\le u\le1}
 \sum_{\rho'\in\Theta(u)}
 \sum_{\substack{V\parallel V_{\rho',R}\\
 V\in\mathcal G_{\rho'}(
 \frac{\alpha}{R^{C\epsilon} S_1(R^{1/n})})}}
 |V|\biggl\|
 \Bigl(\sum_{\theta\subset\rho'}|f_\theta|^2\Bigr)^{1/2}
 \biggr\|_{L_\#^{p_n/2}(V)}^{p_n}.
\end{equation}
Since \(S_2(R^{1-1/n})\ge1\), this gives the bound
required for Proposition~\ref{prop:curveStworecurrence},
and the iteration stops.

If instead the second sum over \(B\) on the RHS of
Lemma~\ref{lem:curveS2producer} dominates, apply dyadic
pigeonholing to fix \(s\). Writing \(\sigma\) for the width of the caps grouped at this scale, we obtain
\begin{equation}\label{eq:curveS2iterate}
 \int_{|f|>\alpha}|f|^{p_n}
 \lesssim_\epsilon R^{C\epsilon} S_1(R^{1/n})^{p_n}
 \sum_{\rho\in\Theta(s)}
 \sum_{\substack{B\parallel(s/\sigma)V_{\rho,\sigma^{-n}}\\
 B\in\mathcal G_\rho(
 \frac{\alpha}{R^{C\epsilon} S_1(R^{1/n})})}}
 |B|\biggl\|
 \Bigl(
 \sum_{\substack{\tau\subset\rho\\|\tau|=\sigma}}
 |f_\tau|^2*w_{\operatorname{par}_{\sigma_0}(\tau)}
 \Bigr)^{1/2}
 \biggr\|_{L_\#^{p_n/2}(B)}^{p_n},
\end{equation}
where \(R^{-1/n}\le\sigma\le\sigma_0\) and
\(\sigma\le s\le1\).

As in Section~\ref{sec:curveunpruning}, caps of width \(\sigma\) are grouped inside caps of width \(s\).
Using Lemmas~\ref{lem:curveS2producer}
and~\ref{lem:curveS2highreturn}, together with the defining
estimate for \(S_2\) after Taylor rescaling, we either
replace the sum over \(B\) in \eqref{eq:curveS2iterate}
by sums over \(V\parallel V_{\rho',R}\) involving the
original functions \(f_\theta\), as in
\eqref{eq:curveS2terminal}, and stop, or return to a bound
of the form \eqref{eq:curveS2iterate} with
\((\sigma,s)\longmapsto(\sigma',s')\),
\(\sigma'\le\sigma\), and \(s'\le s\),
where \(\sigma'\ll\sigma\) or \(s'\ll s\).
At each step, we also update the factor outside the sum
and the denominator in the selection condition.
We will bound the accumulated factor by
\(\bigl(R^{C\eta}S_1(R^{1/n})S_2(R^{1-1/n})\bigr)^{p_n}\)
and keep the final amplitude at least
\(\alpha/\bigl(R^{C\eta}S_1(R^{1/n})S_2(R^{1-1/n})\bigr)\),
with the initial value of \(\alpha\).

We first deal with the case \(\sigma=R^{-1/n}\).
In this case, the derivation of \eqref{eq:curveunpruneterminal} from
\eqref{eq:curveunpruneiteration}
in Section~\ref{sec:curveunpruning}, together with
Lemma~\ref{lem:curveS2unsmooth}, bounds the RHS of
\eqref{eq:curveS2iterate} by the RHS of
\eqref{eq:curveS2terminal}, and the iteration stops.

We may therefore assume \(\sigma>R^{-1/n}\).
We distinguish the cases \(s\le R^{\eta^3}\sigma\) and \(s>R^{\eta^3}\sigma\).

\subsubsection*{Case 1: \texorpdfstring{$s\le R^{\eta^3}\sigma$}{s <= R^(eta^3) sigma}}
\leavevmode\par\nobreak
\noindent
Each \(\rho\) contains at most \(R^{\eta^3}\) caps of
width \(\sigma\). Applying the argument leading to
\eqref{eq:curvescalarize} to the RHS of
\eqref{eq:curveS2iterate}, with \(f_\tau\) and
\(w_{\operatorname{par}_{\sigma_0}(\tau)}\) in place of
\(f_\tau^j\) and \(w_\tau\), gives
\begin{equation}\label{eq:curveS2scalar}
 \int_{|f|>\alpha}|f|^{p_n}
 \lesssim_\eta
 R^{C\eta^3}S_1(R^{1/n})^{p_n}
 \sum_{\rho\in\Theta(s)}
 \sum_{\substack{\tau\subset\rho\\|\tau|=\sigma}}
 \int_{|f_\tau|\gtrsim
 \frac{\alpha}{R^{C\eta^3}S_1(R^{1/n})N(s)}}
 |f_\tau|^{p_n}.
\end{equation}

The Taylor rescaling taking each \(\tau\) to the unit cap
sends each \(\theta\subset\tau\) to a cap of \(\mathcal M^n(\sigma^nR)\) at scale \(R^{-1/n}/\sigma\), where
\(1<\sigma^nR\le R^{1-1/n}\).
To bound each integral on the RHS of
\eqref{eq:curveS2scalar}, we apply the defining estimate
for \(S_2\) to the rescaled \(f_\tau\), with parameter
\(\sigma^nR\) and amplitude
\(\alpha/[R^{C\eta^3}S_1(R^{1/n})N(s)]\).
Returning to the original coordinates as in the derivation
of \eqref{eq:curveunprunerescaled}, and summing over
\(\rho,\tau\), bounds the RHS of \eqref{eq:curveS2scalar} by
\[
 \lesssim_\eta
 R^{C\eta^3}
 \bigl(S_1(R^{1/n})S_2(\sigma^nR)\bigr)^{p_n}
 \sum_{R^{-1/n}\le u\le\sigma}
 \sum_{\rho'\in\Theta(u)}
 \sum_{\substack{V\parallel V_{\rho',R}\\
 V\in\mathcal G_{\rho'}\left(
 \frac{\alpha}{R^{C\eta^3}
 S_1(R^{1/n})S_2(\sigma^nR)}\right)}}
 |V|\biggl\|
 \Bigl(\sum_{\theta\subset\rho'}|f_\theta|^2\Bigr)^{1/2}
 \biggr\|_{L_\#^{p_n/2}(V)}^{p_n},
\]
where we used
\(N(s)N_\tau(u)\le N(\sigma)N_\tau(u)\lesssim_\eta N(u)\)
by Lemma~\ref{lem:regular}.

By monotonicity, we may replace \(S_2(\sigma^nR)\) by
\(S_2(R^{1-1/n})\) in both the coefficient and the amplitude
denominator. Extending the sum over \(u\) to
\(R^{-1/n}\le u\le1\) gives the sums over the original
functions \(f_\theta\) in \eqref{eq:curveS2terminal}, and the iteration stops.

\subsubsection*{Case 2: \texorpdfstring{$s>R^{\eta^3}\sigma$}{s > R^(eta^3) sigma}}
\leavevmode\par\nobreak
\noindent
Applying Lemma~\ref{lem:curveS2highreturn} to the sum over
\(\rho\) and \(B\) on the RHS of \eqref{eq:curveS2iterate},
with amplitude \(\alpha/[R^{C\epsilon} S_1(R^{1/n})]\), gives
\begin{equation}\label{eq:curveS2case2bound}
\begin{aligned}
 \int_{|f|>\alpha}|f|^{p_n}
 &\lesssim_\epsilon
 R^{C\epsilon}S_1(R^{1/n})^{p_n}
 \sum_{\rho\in\Theta(s)}
 \sum_{\substack{Q\parallel\rho^*\\
 Q\in\mathcal G_\rho(
 \frac{\alpha}{R^{C\epsilon}S_1(R^{1/n})})}}
 |Q|\biggl\|
 \Bigl(
 \sum_{\substack{\tau'\subset\rho\\|\tau'|=\sigma'}}
 |f_{\tau'}|^2*w_{\operatorname{par}_{\sigma_0}(\tau')}
 \Bigr)^{1/2}
 \biggr\|_{L_\#^{p_n/2}(Q)}^{p_n}\\
 &\quad+
 R^{C\epsilon}S_1(R^{1/n})^{p_n}
 \sum_{\rho'\in\Theta(s')}
 \sum_{\substack{Q\parallel(\rho')^*\\
 Q\in\mathcal G_{\rho'}(
 \frac{\alpha}{R^{C\epsilon}S_1(R^{1/n})})}}
 |Q|\biggl\|
 \Bigl(
 \sum_{\substack{\tau\subset\rho'\\|\tau|=\sigma}}
 |f_\tau|^2*w_{\operatorname{par}_{\sigma_0}(\tau)}
 \Bigr)^{1/2}
 \biggr\|_{L_\#^{p_n/2}(Q)}^{p_n},
\end{aligned}
\end{equation}
where \(\sigma'=\max\{R^{-1/n},\sigma^2/s\}\) and
\(s'=s(\sigma/s)^{1/(n-1)}\).
We distinguish two cases according to which sum in
\eqref{eq:curveS2case2bound} dominates.
Case~2.1 treats the second sum, where \(s\) decreases to \(s'\), and Case~2.2 treats the first sum, where \(\sigma\) decreases to \(\sigma'\).

\medskip
\noindent\textit{Case 2.1.}
Suppose that the second sum in
\eqref{eq:curveS2case2bound} dominates.
We apply Lemma~\ref{lem:curveS2producer} to its sum over
\(\rho'\) and \(Q\), at scales \((\sigma,s')\) and amplitude
\(\alpha/[R^{C\epsilon}S_1(R^{1/n})]\).

If the first sum over \(V\) on the RHS of
Lemma~\ref{lem:curveS2producer} dominates, applying
Lemma~\ref{lem:curveS2unsmooth} to that sum gives
\eqref{eq:curveS2terminal}, and the iteration stops.

If instead the second sum over \(B\) on the RHS of
Lemma~\ref{lem:curveS2producer} dominates, let
\(\sigma''\) be its cap width. Dyadic pigeonholing reduces
the sum over \(u\) to one scale \(u\in[\sigma'',s']\).
We obtain \eqref{eq:curveS2iterate} with \((\sigma,s)\)
replaced by \((\sigma'',u)\).
Since \(\sigma''\le\sigma\) and
\(u\le s'\le sR^{-\eta^3/(n-1)}\), we return to
\eqref{eq:curveS2iterate} with these new parameters and
repeat the algorithm.

\medskip
\noindent\textit{Case 2.2.}
Suppose that the first sum in
\eqref{eq:curveS2case2bound} dominates.

If \(\sigma'>R^{-1/n}\), we apply
Lemma~\ref{lem:curveS2producer} to its sum over \(\rho\)
and \(Q\), at scales \((\sigma',s)\) and amplitude \(\alpha/[R^{C\epsilon}S_1(R^{1/n})]\).

If the first sum over \(V\) on the RHS of
Lemma~\ref{lem:curveS2producer} dominates, applying
Lemma~\ref{lem:curveS2unsmooth} to that sum gives
\eqref{eq:curveS2terminal}, and the iteration stops.

If instead the second sum over \(B\) on the RHS of
Lemma~\ref{lem:curveS2producer} dominates, let
\(\sigma''\) be its cap width. Dyadic pigeonholing reduces
the sum over \(u\) to one scale \(u\in[\sigma'',s]\).
We obtain \eqref{eq:curveS2iterate} with \((\sigma,s)\)
replaced by \((\sigma'',u)\).
Since \(\sigma''\le\sigma'=\sigma^2/s\le\sigma R^{-\eta^3}\)
and \(u\le s\), we return to \eqref{eq:curveS2iterate}
with these new parameters and repeat the algorithm.

If \(\sigma'=R^{-1/n}\), the first sum over \(Q\) in
\eqref{eq:curveS2case2bound} involves the original
functions \(f_\theta\).
The derivation of \eqref{eq:curveunpruneterminal} from
\eqref{eq:curveunpruneiteration} for \(\sigma=R^{-1/n}\)
in Section~\ref{sec:curveunpruning}, together with
Lemma~\ref{lem:curveS2unsmooth}, applies to this sum and
gives \eqref{eq:curveS2terminal}, and the iteration stops.

\subsubsection*{Termination criterion}

At each return to \eqref{eq:curveS2iterate}, neither
\(\sigma\) nor \(s\) increases, and at least one decreases
by a factor of at least \(R^{\eta^3/(n-1)}\).
Since both lie in \([R^{-1/n},1]\), there are
\(O(\eta^{-3})\) steps. The factor \(S_2(\sigma^nR)\)
occurs at most once, when Case~1 terminates the iteration.

The \(R^{C\epsilon}\) losses in the coefficient and the
amplitude denominator accumulate to at most
\(R^{C\epsilon\eta^{-3}}\).
Taking \(\epsilon\) sufficiently small in terms of \(\eta\),
we absorb these losses into \(R^{C\eta}\).
Combining the estimates above, we bound
\(\int_{|f|>\alpha}|f|^{p_n}\) by
\[
 \bigl(R^{C\eta}S_1(R^{1/n})S_2(R^{1-1/n})\bigr)^{p_n}
 \sum_{R^{-1/n}\le u\le1}
 \sum_{\rho'\in\Theta(u)}
 \sum_{\substack{V\parallel V_{\rho',R}\\
 V\in\mathcal G_{\rho'}\left(
 \frac{\alpha}{R^{C\eta}S_1(R^{1/n})S_2(R^{1-1/n})}
 \right)}}
 |V|\biggl\|
 \Bigl(\sum_{\theta\subset\rho'}|f_\theta|^2\Bigr)^{1/2}
 \biggr\|_{L_\#^{p_n/2}(V)}^{p_n}.
\]
The same proof with \(R\) replaced by any \(1\le r\le R\),
together with the monotonicity of \(S_1,S_2\), gives
\(S_2(R)\lesssim_\eta
R^{C\eta}S_1(R^{1/n})S_2(R^{1-1/n})\)
by the definition of \(S_2\).
Since \(\eta\) can be chosen arbitrarily small, this proves
Proposition~\ref{prop:curveStworecurrence}.
\end{proof}

\addtocontents{toc}{\protect\tocpartheading{II}{From curves to higher dimensional cones over moment curves}\protect\tocsectionindent{3pc}}
\section{\texorpdfstring{\except{toc}{From curves to higher dimensional cones over moment curves: }}{From curves to higher dimensional cones over moment curves: }Overview of the proof}
\label{chap:cone-induction}\label{sec:highercone}

We prove the amplitude dependent wave envelope type estimate for the cone over the moment curve
in \(\R^{d+1}\) (Theorem~\ref{thm:main-cone}), assuming the corresponding
estimates for moment curves in \(\R^k\), \(2\le k\le d\)
(Theorem~\ref{thm:main}). Here is the statement.

\begin{proposition}\label{prop:highercone}
Let \(d\ge3\). Assume that the amplitude dependent wave envelope type estimate
in Theorem~\ref{thm:main} holds for the moment curve in \(\R^k\)
at exponent \(p_k\) for every \(2\le k\le d\).
Then the corresponding estimate in Theorem~\ref{thm:main-cone} holds
for the cone over the moment curve in \(\R^{d+1}\) at exponent \(p_d\), with the averaging boxes defined in
\eqref{eq:conebox}.
\end{proposition}

Recall that \(p_k=k(k+1)/2+1\),
\(\phi_d(t)=(1,t,t^2/2!,\ldots,t^d/d!)\), and
\[
 \Gamma^{d+1}(R)=\Bigl\{\sum_{i=0}^d\lambda_i\phi_d^{(i)}(t):
 0\le t\le1,\quad 1\le\lambda_0\le2,\quad
 |\lambda_i|\le R^{-i/d}\ (1\le i\le d)\Bigr\}.
\]
For a Schwartz function \(f:\R^{d+1}\to\mathbb C\) with
\(\supp\widehat f\subset\Gamma^{d+1}(R)\), we count the caps that contain at least one
\(\theta\in\Xi(R^{-1/d})\) with \(f_\theta\not\equiv0\):
\[
 N(s)=\#\{\tau\in\Xi(s):f_\theta\not\equiv0
                   \text{ for some }\theta\subset\tau\},\qquad
 N_\rho(u)=\#\{\tau\in\Xi(u):\tau\subset\rho,
                   \ f_\theta\not\equiv0\text{ for some }\theta\subset\tau\}.
\]
Let \(t_\rho\) be the left endpoint of the interval corresponding to \(\rho\).
Recall from \eqref{eq:curvebox} and \eqref{eq:conebox} that
\[
\begin{aligned}
 U_{\rho,R}
 &=\bigl\{x\in\R^{d+1}: |x\cdot\phi_d^{(j)}(t_\rho)|\le R^{2/d}s^{2-j}
    \text{ for }0\le j\le d\bigr\}, &&\rho\in\Xi(s),\\
 V_{\rho,R}
 &=\bigl\{x\in\R^k: |x\cdot\gamma_k^{(j)}(t_\rho)|\le R^{2/k}s^{2-j}
    \text{ for }1\le j\le k\bigr\}, &&\rho\in\Theta(s),
\end{aligned}
\]
where \(\phi_d^{(0)}=\phi_d\) and \(2\le k\le d\).

\medskip

\begin{convention}
For a parent \(\rho\) over an interval of length \(s\),
in a sum over boxes of the form
\[
 \sum_{\substack{B\parallel B_0\\B\in\mathcal G_\rho(\alpha)}}
 |B|\|F_\rho\|_{L^{p/2}_\#(B)}^p,
\]
the selection uses the same function, exponent and tiling as the summand:
\[
 \mathcal G_\rho(\alpha)=
 \Bigl\{B\parallel B_0:\|F_\rho\|_{L^{p/2}_\#(B)}
                      \ge\frac{\alpha}{N(s)}\Bigr\}.
\]
\end{convention}

We now write the hypothesis and conclusion of
Proposition~\ref{prop:highercone} in this notation.
For curves in \(\R^k\), \(2\le k\le d\), recall that
\(\gamma_k(t)=(t,t^2,\ldots,t^k)\) and
\[
 \mathcal M^k(R)=\Bigl\{\gamma_k(t)+\sum_{i=2}^k
 \frac{a_i}{i!}\gamma_k^{(i)}(t):
 0\le t\le1,\quad |a_i|\le R^{-i/k}\Bigr\}.
\]
For each \(2\le k\le d\), we apply the amplitude dependent wave envelope type
estimate for the moment curve in \(\R^k\), assumed in
Proposition~\ref{prop:highercone}. For every Schwartz function
\(g:\R^k\to\mathbb C\) with \(\supp\widehat g\subset \mathcal M^k(R)\)
and every \(\alpha>0\), this gives
\begin{equation}\label{eq:hccurvehypothesis}
 \int_{|g|>\alpha}|g|^{p_k}
 \lesssim_\e R^\e
 \sum_{R^{-1/k}\le s\le1}\sum_{\rho\in\Theta(s)}
 \sum_{\substack{V\parallel V_{\rho,R}\\
 V\in\mathcal G_\rho(\frac{\alpha}{R^\e})}}
 |V|\biggl\|\bigl(\sum_{\theta\subset\rho}|g_\theta|^2\bigr)^{1/2}
 \biggr\|_{L^{p_k/2}_\#(V)}^{p_k}.
\end{equation}
The conclusion is that every Schwartz function
\(f:\R^{d+1}\to\mathbb C\) with
\(\supp\widehat f\subset\Gamma^{d+1}(R)\) satisfies
\[
 \int_{|f|>\alpha}|f|^{p_d}
 \lesssim_\e R^\e
 \sum_{R^{-1/d}\le s\le1}\sum_{\rho\in\Xi(s)}
 \sum_{\substack{U\parallel U_{\rho,R}\\
 U\in\mathcal G_\rho(\frac{\alpha}{R^\e})}}
 |U|\biggl\|\bigl(\sum_{\theta\subset\rho}|f_\theta|^2\bigr)^{1/2}
 \biggr\|_{L^{p_d/2}_\#(U)}^{p_d}.
\]

Throughout Sections~\ref{sec:highercone}--\ref{sec:hcgeneral},
we assume the hypothesis of Proposition~\ref{prop:highercone} on moment curves.
The high frequency decomposition produces Fourier supports that
need not satisfy the condition on the cone over the moment curve \(|\lambda_0|\sim1\),
even after rescaling. Following \cite{GMSF}, we therefore introduce
the following neighborhoods of Taylor cones.
For \(R\ge1\) and \(0\le m\le d-2\), define
\begin{equation}\label{eq:hcoverviewtaylorcones}
\begin{aligned}
\Gamma_m^{d+1}(R)
&=
\left\{
\sum_{i=0}^d\lambda_i\phi_d^{(i)}(t):
\begin{array}{l}
t\in[0,1],\quad
|\lambda_i|\le R^{-i/d}\quad(0\le i\le d),\\[2pt]
\displaystyle\max_{0\le i\le m}
2^{m+1-i}R^{i/d}|\lambda_i|\ge1
\end{array}
\right\}=\bigcup_{\theta\in\Xi_m^{d+1}(R)}\theta,
\end{aligned}
\end{equation}
where \(\Xi_m^{d+1}(R)\) is the family of caps indexed by
the standard partition of \([0,1]\) into intervals \(I\)
of length \(R^{-1/d}\), with
\begin{equation}\label{eq:hcoverviewtaylorcaps}
\theta_I=
\left\{
\sum_{i=0}^d\lambda_i\phi_d^{(i)}(t):
\begin{array}{l}
t\in I,\quad
|\lambda_i|\le R^{-i/d}\quad(0\le i\le d),\\[2pt]
\displaystyle\max_{0\le i\le m}
2^{m+1-i}R^{i/d}|\lambda_i|\ge1
\end{array}
\right\}.
\end{equation}

The maximum in \eqref{eq:hcoverviewtaylorcones} requires at least
one of \(\lambda_0,R^{1/d}\lambda_1,\ldots,R^{m/d}\lambda_m\)
to stay away from zero. Section~\ref{subsec:hcoverviewpositive}
explains why the high frequency decomposition raises this order. For \(f\) with Fourier support in \(\Gamma_m^{d+1}(R)\), we use a smooth partition of unity adapted to these caps to decompose \(f\) into pieces \(f_{\theta_I}\), each with Fourier support in a fixed enlargement of \(\theta_I\).

For dyadic \(0<s\le1\), let \(\mathcal J(s)\) denote the standard
partition of \([0,1]\) into intervals of length \(s\).
We identify each interval in \(\mathcal J(s)\) with the cap over it and use \(\mathcal J(s)\) for both families.

Fix \(\e>0\).
We choose \(A=A(\e)\) large, then \(K=K(A,\e)\) large,
and take \(R\) sufficiently large.
To apply the hypothesis on moment curves \eqref{eq:hccurvehypothesis},
we first restrict the lower Taylor coefficients to short intervals.
For \(0\le j\le m\) and \(I\in\mathcal J(R^{-1/d})\), let
\(t_I\) be the left endpoint of \(I\). Define
\begin{equation}\label{eq:hcstableband}
\Gamma^{d+1}_{j,1/K}(R):=
\bigcup_{I\in\mathcal J(R^{-1/d})}
\left\{\sum_{i=0}^d\lambda_i\phi_d^{(i)}(t_I):
\begin{array}{ll}
 |\lambda_i|\le K^{-1}R^{-j/d},&0\le i<j,\\
 |\lambda_j-R^{-j/d}|\le K^{-1}R^{-j/d},&\\
 |\lambda_i|\lesssim R^{-i/d},&j<i\le d
\end{array}\right\}.
\end{equation}
 Geometrically, \(\Gamma^{d+1}_{j,1/K}(R)\) is obtained by
thickening the curve \(R^{-j/d}\phi_d^{(j)}(t)\),
\(0\le t\le1\), by \(K^{-1}R^{-j/d}\) in the directions
\(\phi_d(t),\ldots,\phi_d^{(j)}(t)\), and by \(R^{-i/d}\)
in the directions \(\phi_d^{(i)}(t)\), \(j<i\le d\).

In particular, \(\Gamma^{d+1}_{0,1/K}(R)\) can be viewed
as the anisotropic neighborhood of the moment curve
\((1,t,t^2/2!,\ldots,t^d/d!)\), \(0\le t\le1\),
thickened by \(K^{-1}\) in the \(\phi_d(t)\) direction.
Note that, up to fixed enlargement,
\(\Gamma^{d+1}_{j,1/K}(R)\subset\Gamma_m^{d+1}(R)\),
\(0\le j\le m\).
As before, we write \(f=\sum_I f_{\theta_I}\), allowing any decomposition
with the cap supports specified above. When
\(\operatorname{supp}\widehat f\subset\Gamma^{d+1}_{j,1/K}(R)\),
we further require each \(\widehat{f_{\theta_I}}\) to be supported in
a fixed enlargement of the \(I\)th set in \eqref{eq:hcstableband}.
We now introduce the constants needed to prove
Proposition~\ref{prop:highercone} by descending induction on \(m\).

\begin{definition}\label{def:hcconstants}
Let \(C^m(R)\) be the least constant such that, for every
Schwartz function \(f\) with
\(\supp\widehat f\subset\Gamma_m^{d+1}(R)\) and every
\(\alpha>0\),
\begin{equation}\label{eq:hcconstantestimate}
 \int_{|f|>\alpha}|f|^{p_{d-m}}
 \le C^m(R)
 \sum_{R^{-1/d}\le s\le1}\sum_{\rho\in\mathcal J(s)}
 \sum_{\substack{U\parallel(sR^{1/d})^mU_{\rho,R}\\
 U\in\mathcal G_\rho(\frac{\alpha}{C^m(R)})}}
 |U|\biggl\|\Bigl(\sum_{\theta\subset\rho}|f_\theta|^2\Bigr)^{1/2}
 \biggr\|_{L^{p_{d-m}/2}_\#(U)}^{p_{d-m}}.
\end{equation}
\end{definition}

\begin{definition}\label{def:hcrestrictedconstants}
Let \(C_K^m(R)\) be the least constant such that, for every
\(0\le j\le m\), every Schwartz function \(f\) with
\(\supp\widehat f\subset\Gamma^{d+1}_{j,1/K}(R)\), and every
\(\alpha>0\),
\begin{equation}\label{eq:hcrestrictedconstantestimate}
 \int_{|f|>\alpha}|f|^{p_{d-m}}
 \le C_K^m(R)
 \sum_{R^{-1/d}\le s\le1}\sum_{\rho\in\mathcal J(s)}
 \sum_{\substack{U\parallel(sR^{1/d})^mU_{\rho,R}\\
 U\in\mathcal G_\rho(\frac{\alpha}{C_K^m(R)})}}
 |U|\biggl\|\Bigl(\sum_{\theta\subset\rho}|f_\theta|^2\Bigr)^{1/2}
 \biggr\|_{L^{p_{d-m}/2}_\#(U)}^{p_{d-m}}.
\end{equation}
\end{definition}

Our goal is \(C^0(R)\lessapprox1\). We establish the bounds
\(C_K^m(R)\lessapprox1\) and \(C^m(R)\lessapprox1\)
successively in the order
\[
\begin{aligned}
 C_K^{d-2}(R)
 &\longrightarrow C^{d-2}(R)
 \longrightarrow C_K^{d-3}(R)\longrightarrow C^{d-3}(R)
 \longrightarrow\cdots
 \longrightarrow C_K^0(R)
 \longrightarrow C^0(R).
\end{aligned}
\]
We proceed in decreasing \(m\) because the argument for the high frequency part
uses estimates at exponent \(p_{d-m-1}\) to obtain
estimates at exponent \(p_{d-m}\).
For technical reasons, applying the induction hypothesis to the
high frequency terms requires keeping the convolutions between
successive powers of the square functions. We now formulate these estimates.

For \(m<d-2\), consider angular widths
\(R^{-1/d}\le q_1\le\cdots\le q_{d-m-2}\le u\le1\)
and caps
\(\theta\subset\tau_1\subset\cdots\subset\tau_{d-m-2}
\subset\rho\), where \(|\tau_l|=q_l\) and \(|\rho|=u\).
Let \(B_{\tau_l}\) consist of all \(x\in\R^{d+1}\) satisfying
\begin{equation}\label{eq:hcnestedwidths}
 |x\cdot\phi_d^{(i)}(t_{\tau_l})|\le
 \frac{R^{(m+2)/d}}{q_1\cdots q_{l-1}}
 \begin{cases}
 q_l^{m+l+1-i},&0\le i\le m+l+1,\\
 u^{m+l+1-i},&m+l+1<i\le d,
 \end{cases}
 \qquad 1\le l\le d-m-2.
\end{equation}
The boxes in \eqref{eq:hcnestedwidths} are chosen so that
(i) their form is preserved under Taylor rescaling, and
(ii) the high frequency terms produced by the cutoffs in
Section~\ref{subsec:hcoverviewpositive} can be bounded by
convolutions with \(w_{B_{\tau_l}}\), using the locally constant property
and H\"older.
The role of this choice will become explicit in
Section~\ref{sec:hcpositive}; see \eqref{eq:hchighpieces}
and \eqref{eq:hcalllayerrescaling}.
For notational convenience, we suppress the dependence of
\(B_{\tau_l}\) on \(q_1,\ldots,q_{l-1}\) and on the width
\(u\) of \(\rho\).

In proving \(C^m(R)\lessapprox1\), we will prove the following stronger estimate that for every
Schwartz function \(f=\sum_{\theta\in\Xi_m^{d+1}(R)}f_\theta\) with
\(\supp\widehat f\subset\Gamma_m^{d+1}(R)\),
\begin{equation}\label{eq:hcgeneralconvolution}
\begin{aligned}
 \int_{|f|>\alpha}|f|^{p_{d-m}}
 &\lesssim_\e R^\e
 \sum_{R^{-1/d}\le q_1\le\cdots\le q_{d-m-2}\le u\le1}
 \sum_{\rho\in\mathcal J(u)}
 \sum_{\substack{
 B\parallel
 \frac{u^{d-m-2}}{q_1\cdots q_{d-m-2}}
 (uR^{1/d})^mU_{\rho,R}\\
 B\in\mathcal G_\rho(\frac{\alpha}{R^\e})}}
 |B|\\
 &\quad\times
 \biggl[
 \fint_B
 \sum_{\tau_{d-m-2}\subset\rho}
 \biggl[
 \biggl(
 \cdots
 \sum_{\tau_1\subset\tau_2}
 \Bigl[
 \bigl(\sum_{\theta\subset\tau_1}|f_\theta|^2\bigr)
 *w_{B_{\tau_1}}
 \Bigr]^{p_{d-m}/p_{d-m-1}}
 \cdots
 \biggr)*w_{B_{\tau_{d-m-2}}}
 \biggr]^{7/4}
 \biggr]^2,
\end{aligned}
\end{equation}
where \(B\in\mathcal G_\rho(\alpha/R^\e)\) means that the
weighted average over \(B\) inside the outermost square brackets
is at least \(\bigl[\alpha/(R^\e N(u))\bigr]^{p_{d-m}/2}\).
The boxes \(B\parallel \frac{u^{d-m-2}}{q_1\cdots q_{d-m-2}}(uR^{1/d})^mU_{\rho,R}\) arise by returning the averaging boxes obtained from successive applications of the induction hypothesis to the high frequency terms to the original coordinates.
For example, \eqref{eq:hcgeneralconvolution} with \(m=d-4\)
states that every Schwartz function \(f\) with
\(\supp\widehat f\subset\Gamma_{d-4}^{d+1}(R)\) satisfies
\[
\begin{aligned}
 \int_{|f|>\alpha}|f|^{11}
 &\lesssim_\e R^\e
 \sum_{R^{-1/d}\le q_1\le q_2\le u\le1}
 \sum_{\rho\in\mathcal J(u)}
 \sum_{\substack{
 B\parallel\frac{u^2}{q_1q_2}(uR^{1/d})^{d-4}U_{\rho,R}\\
 B\in\mathcal G_\rho(\frac{\alpha}{R^\e})}}
 |B|\\
 &\quad\times
 \biggl[
 \fint_B
 \sum_{\substack{\tau_2\subset\rho\\|\tau_2|=q_2}}
 \biggl[
 \biggl(
 \sum_{\substack{\tau_1\subset\tau_2\\|\tau_1|=q_1}}
 \Bigl[
 \bigl(\sum_{\theta\subset\tau_1}|f_\theta|^2\bigr)
 *w_{B_{\tau_1}}
 \Bigr]^{11/7}
 \biggr)*w_{B_{\tau_2}}
 \biggr]^{7/4}
 \biggr]^2.
\end{aligned}
\]

In proving \(C_K^m(R)\lessapprox1\), we will also prove the corresponding stronger estimate that every
Schwartz function \(f\) with
\(\supp\widehat f\subset\Gamma_{j,1/K}^{d+1}(R)\),
\begin{equation}\label{eq:hcrestrictedconvolution}
 \int_{|f|>\alpha}|f|^{p_{d-m}}
 \lesssim_\e
 \text{RHS of \eqref{eq:hcgeneralconvolution}},
\end{equation}
uniformly in \(0\le j\le m\).
Indeed, the following lemma shows that
\eqref{eq:hcgeneralconvolution} and
\eqref{eq:hcrestrictedconvolution} imply the estimates
\eqref{eq:hcconstantestimate} and
\eqref{eq:hcrestrictedconstantestimate}, respectively,
with \(C^m(R)\lessapprox1\) and \(C_K^m(R)\lessapprox1\).

\begin{lemma}\label{lem:hcstep2averagecomparison}
Let \(0\le m<d-2\) and \(R\ge1\), and let \(f\) be a Schwartz
function with
\(\supp\widehat f\subset\Gamma_m^{d+1}(R)\).
For every \(\epsilon>0\) and \(\alpha>0\),
\[
 \text{RHS of \eqref{eq:hcgeneralconvolution}}
 \lessapprox
 \sum_{R^{-1/d}\le u\le1}\sum_{\rho\in\mathcal J(u)}
 \sum_{\substack{U\parallel(uR^{1/d})^mU_{\rho,R}\\
 U\in\mathcal G_\rho(\frac{\alpha}{R^\epsilon})}}
 |U|\biggl\|\Bigl(\sum_{\theta\subset\rho}|f_\theta|^2\Bigr)^{1/2}
 \biggr\|_{L^{p_{d-m}/2}_\#(U)}^{p_{d-m}}.
\]
\end{lemma}

We defer the proof to the end of this section.
The induction is organized by the following propositions.

\begin{proposition}[Base case]\label{prop:hcbase}
We have \(C_K^{d-2}(R)\lessapprox1\) and
\(C^{d-2}(R)\lessapprox1\).
\end{proposition}

\begin{proposition}\label{prop:hcstablefunction}
Let \(0\le m\le d-3\). Assume \(C^{m+1}(R)\lessapprox1\)
and that \eqref{eq:hcgeneralconvolution} holds with \(m+1\)
in place of \(m\). Then \(C^m(R)\lessapprox1\), and
\eqref{eq:hcgeneralconvolution} holds for this \(m\). When \(m=d-3\), the assumption on
\eqref{eq:hcgeneralconvolution} is omitted.
\end{proposition}

We prove Proposition~\ref{prop:hcstablefunction} in two steps.

\begin{itemize}
\item \textbf{Step~1.} Under the hypotheses of
Proposition~\ref{prop:hcstablefunction}, we establish
\(C_K^m(R)\lessapprox1\) together with
\eqref{eq:hcrestrictedconvolution} for this \(m\).

\item \textbf{Step~2.} Keeping the hypotheses of
Proposition~\ref{prop:hcstablefunction}, we use
\(C_K^m(R)\lessapprox1\) and \eqref{eq:hcrestrictedconvolution}
from Step~1 to prove \(C^m(R)\lessapprox1\) and
\eqref{eq:hcgeneralconvolution} for this \(m\).
\end{itemize}

We carry out Steps~1 and~2 in Sections~\ref{sec:hcstable}
and~\ref{sec:hcgeneral}, respectively.
Section~\ref{sec:hcbase} proves Proposition~\ref{prop:hcbase},
the base case \(m=d-2\) of the induction on \(m\).

\subsection{Roadmap: Proof of Proposition~\ref{prop:hcstablefunction}, Step 1}
\label{subsec:hcroadmap}\label{sec:hcinitial}

To make the induction easier to follow, we first outline
the proof of \eqref{eq:hcrestrictedconstantestimate} with
\(C_K^m(R)\lessapprox1\), rather than the stronger estimate
involving the convolutions in \eqref{eq:hcgeneralconvolution}.
The argument for the stronger estimate  \eqref{eq:hcrestrictedconvolution} is
discussed in the remark at the end of
Section~\ref{subsec:hcoverviewiteration}.

We first use the hypothesis on moment curves to bound the integral of
\(|f|^{p_{d-m}}\) over \(\{|f|>\alpha\}\)
in terms of the pieces \(f_\tau\), \(|\tau|=A^{-1}\).

By Lemma~\ref{lem:regular}, we split the active caps into \(J\lesssim_\kappa R^\kappa\) regular classes and apply Lemma~\ref{lem:hcseed} to each class at amplitude \(\alpha/J\). The loss from summing over the classes is absorbed into \(R^\epsilon\). In Proposition~\ref{prop:hcstable} and Corollary~\ref{cor:hccoarse}, all counts refer to a fixed class, so \(N(s)N_\rho(v)\lesssim_\kappa N(v)\) for \(v\le s=|\rho|\). For Proposition~\ref{prop:hcstep2positive}, we apply Lemma~\ref{lem:regular} to the family \(\{f_\theta^{j,\nu,\ell}\}_\theta\) after fixing \(j,\nu,\ell\) in Section~\ref{subsec:hcoverviewcoefficients}.

\begin{lemma}\label{lem:hcseed}
Let \(f\) be a Schwartz function with
\(\supp\widehat f\subset\Gamma^{d+1}_{j,1/K}(R)\), \(0\le j\le m\).
Then
\[
 \int_{|f|>\alpha}|f|^{p_{d-m}}
 \lesssim_\e A^\e
 \sum_{A^{-1}\le s\le1}\sum_{\rho\in\mathcal J(s)}
 \sum_{\substack{Q\parallel\rho^*\\
 Q\in\mathcal G_\rho(\frac{\alpha}{A^\e})}}
 |Q|\biggl\|\Bigl(\sum_{\substack{\tau\in\mathcal J(A^{-1})\\\tau\subset\rho}}
 |f_\tau|^2\Bigr)^{1/2}\biggr\|_{L^{p_{d-m}/2}_\#(Q)}^{p_{d-m}}.
\]
\end{lemma}

\begin{proof}
Frequency dilation by \(R^{j/d}\) sends the curve
\(R^{-j/d}\phi_d^{(j)}(t)\) in \eqref{eq:hcstableband}
to \(\phi_d^{(j)}(t)\). Projection onto coordinates
\(j+1,\ldots,j+d-m\), followed by a linear change, sends
\(\phi_d^{(j)}(t)\) to
\(\gamma_{d-m}(t)=(t,t^2,\ldots,t^{d-m})\).
Taking \(K\) and then \(R\) sufficiently large, the Fourier
supports of \(f_\tau(R^{j/d}\,\cdot)\),
\(\tau\in\mathcal J(A^{-1})\), project into the caps of
\(\mathcal M^{d-m}(A^{d-m})\) in these coordinates.
Apply \eqref{eq:hccurvehypothesis} to \(f(R^{j/d}x)\)
in the physical variables \(x_{j+1},\ldots,x_{j+d-m}\),
with parameter \(A^{d-m}\).

Apply \eqref{eq:selaverage} to the sum over \(V\) on the
right hand side of \eqref{eq:hccurvehypothesis}, with
\(F=(\sum_{\tau\in\mathcal J(A^{-1}),\,\tau\subset\rho}
|f_\tau(R^{j/d}x)|^2)^{p_{d-m}/4}\) and
\(t=[\alpha/(A^\e N(s))]^{p_{d-m}/2}\).
Integrating over the remaining physical variables and returning
to the original coordinates gives the inequality in
Lemma~\ref{lem:hcseed}.
\end{proof}

The RHS of Lemma~\ref{lem:hcseed} still involves
\(f_\tau\), \(|\tau|=A^{-1}\).
To obtain \(C_K^m(R)\lessapprox1\), we use the following proposition
to replace these pieces by the original \(f_\theta\). The factor \((\sigma s)^\epsilon\) in the following estimate
provides a gain when either \(\sigma\) or \(s\) decreases, allowing us to absorb the refinement
losses and close the induction.

\begin{proposition}\label{prop:hcstable}
Assume the hypotheses of Proposition~\ref{prop:hcstablefunction}
when \(m<d-2\). Let \(f\) be a Schwartz function with
\(\supp\widehat f\subset\Gamma^{d+1}_{j,1/K}(R)\), \(0\le j\le m\).
If \(R^{-1/d}\le\sigma\le A^{-1}\) and \(\sigma\le s\le1\), then
\[
\begin{split}
 &\sum_{\rho\in\mathcal J(s)}
 \sum_{\substack{Q\parallel\rho^*\\Q\in\mathcal G_\rho(\alpha)}}
 |Q|\biggl\|\Bigl(\sum_{\substack{\tau\subset\rho\\|\tau|=\sigma}}
 |f_\tau|^2\Bigr)^{1/2}\biggr\|_{L^{p_{d-m}/2}_\#(Q)}^{p_{d-m}}\\
 &\lesssim_\e R^\e(\sigma s)^\e
 \sum_{R^{-1/d}\le v\le1}\sum_{\rho'\in\mathcal J(v)}
 \sum_{\substack{U\parallel(vR^{1/d})^mU_{\rho',R}\\
 U\in\mathcal G_{\rho'}(\frac{\alpha}{R^\e(\sigma s)^\e})}}
 |U|\biggl\|\Bigl(\sum_{\theta\subset\rho'}|f_\theta|^2\Bigr)^{1/2}
 \biggr\|_{L^{p_{d-m}/2}_\#(U)}^{p_{d-m}}.
\end{split}
\]
\end{proposition}

For any parameter \(R\) at which Proposition~\ref{prop:hcstable} holds
and any \(f\) satisfying its support hypothesis, apply Lemma~\ref{lem:hcseed} with \(\epsilon/2\) in place of \(\epsilon\),
and then Proposition~\ref{prop:hcstable} with \(\sigma=A^{-1}\)
to the sum over \(\rho\) and \(Q\) on the right hand side of
Lemma~\ref{lem:hcseed}, for each fixed \(s\).
The resulting sum over \(U\) has outer coefficient and amplitude
denominator bounded by \(R^\epsilon A^{-\epsilon/2}s^\epsilon\).
Since \(\sum_{A^{-1}\le s\le1}s^\epsilon\lesssim_\epsilon1\),
this gives
\begin{equation}\label{eq:hcstablefunction}
 \int_{|f|>\alpha}|f|^{p_{d-m}}\lesssim_\e R^\e
 \sum_{R^{-1/d}\le v\le1}\sum_{\rho\in\mathcal J(v)}
 \sum_{\substack{U\parallel(vR^{1/d})^mU_{\rho,R}\\
 U\in\mathcal G_\rho(\frac{\alpha}{R^\e})}}
 |U|\biggl\|\Bigl(\sum_{\theta\subset\rho}|f_\theta|^2\Bigr)^{1/2}
 \biggr\|_{L^{p_{d-m}/2}_\#(U)}^{p_{d-m}}.
\end{equation}

We will prove Proposition~\ref{prop:hcstable} by induction on \(R\)
and, for each fixed \(R\), on \(\sigma s\).
Assume that the proposition holds at every parameter \(R'<R/2\).
The deduction above then gives that \eqref{eq:hcstablefunction}
holds at every such \(R'\).
With \(R\) fixed, we now prove Proposition~\ref{prop:hcstable}
by induction on \(\sigma s\). For the pair \((\sigma,s)\),
assume that Proposition~\ref{prop:hcstable} holds for every
pair \((\sigma',s')\) satisfying
\(R^{-1/d}\le\sigma'\le A^{-1}\), \(\sigma'\le s'\le1\),
and \(\sigma's'<\sigma s\).
In this argument, Taylor rescaling
of \(f_\tau\), \(|\tau|=\sigma\), gives the parameter
\(R'=R\sigma^d\le RA^{-d}\le R/2\), so
\eqref{eq:hcstablefunction} is available at \(R'\).
The following subsection describes this induction on \(\sigma s\).

Once Proposition~\ref{prop:hcstable} is proved at \(R\), applying
Lemma~\ref{lem:hcseed} and summing over \(s\) as above gives
\eqref{eq:hcstablefunction} at \(R\).
In particular, this gives \(C_K^m(R)\lessapprox1\).

\subsection{Roadmap continued: induction on $\sigma s$}
\label{subsec:hcoverviewiteration}

We now show how the induction on \(\sigma s\) at fixed \(R\)
reduces the proof of Proposition~\ref{prop:hcstable} to
Lemmas~\ref{lem:hcproducer} and~\ref{lem:hcoverviewrefinement} below.
For the argument below, we will choose the parameters such that
\begin{itemize}
    \item \(0<\delta\ll\eta\ll\epsilon\).
\end{itemize}
Specifically, we will take \(\delta\ll\epsilon\eta^5\).
Let \(f\) be a Schwartz function with
\(\supp\widehat f\subset\Gamma^{d+1}_{j,1/K}(R)\), \(0\le j\le m\).
To estimate the left hand side of Proposition~\ref{prop:hcstable},
we first dilate the Fourier supports of \(f_\tau\), \(|\tau|=\sigma\),
into the caps of a Taylor cone.
For \(\tau\in\mathcal J(\sigma)\), write
\(\xi=\sum_{i=0}^d\lambda_i\phi_d^{(i)}(t_\tau)\)
for \(\xi\in\supp\widehat{f_\tau}\).
Taylor expansion and \eqref{eq:hcstableband} give
\[
 |\lambda_i|\lesssim
 \begin{cases}
 R^{-j/d},&0\le i\le j,\\
 R^{-j/d}\sigma^{i-j},&j<i\le d,
 \end{cases}
 \qquad |\lambda_j|\sim R^{-j/d}.
\]
Here we used \(R^{-1/d}\le\sigma\); the last comparison follows
from \(\lambda_j=R^{-j/d}(1+O(K^{-1}+\sigma))\)
and \(\sigma\le A^{-1}\).

Set \(g(x)=f(R^{j/d}\sigma^j x)\). This dilation multiplies each coefficient \(\lambda_i\)
by \(R^{j/d}\sigma^j\), so
\[
 |R^{j/d}\sigma^j\lambda_i|\lesssim\sigma^i
 \quad(0\le i\le d),\qquad
 |R^{j/d}\sigma^j\lambda_j|\sim\sigma^j.
\]
Since \(\sigma^i=(\sigma^{-d})^{-i/d}\) and \(j\le m\), this gives
\[
 \supp\widehat{f_\tau(R^{j/d}\sigma^j\,\cdot)}
 \subset\theta_\tau,
 \qquad \theta_\tau\in\Xi_m^{d+1}(\sigma^{-d}).
\]
We therefore use the following lemma with \(g\) and
\(r=\sigma^{-d}\).

\begin{lemma}\label{lem:hcproducer}\label{prop:hcpositivecanonical}\label{prop:hcmomentpositivecanonical}
Assume the hypotheses of Proposition~\ref{prop:hcstablefunction}
for a fixed \(0\le m\le d-3\). Let \(g\) be a Schwartz function with
\(\supp\widehat g\subset\Gamma_m^{d+1}(r)\). Then
\[
 \int_{\sum_\theta|g_\theta|^2\gtrsim\alpha^2}
       \bigl(\sum_\theta|g_\theta|^2\bigr)^{p_{d-m}/2}
 \lesssim_\e r^\e
 \sum_{r^{-1/d}\le u\le1}\sum_{\rho\in\mathcal J(u)}
 \sum_{\substack{U\parallel(ur^{1/d})^mU_{\rho,r}\\
 U\in\mathcal G_\rho(\frac{\alpha}{r^\e})}}
 |U|\biggl\|\bigl(\sum_{\theta\subset\rho}|g_\theta|^2\bigr)^{1/2}
 \biggr\|_{L^{p_{d-m}/2}_\#(U)}^{p_{d-m}}.
\]
\end{lemma}

Apply Lemma~\ref{lem:hcproducer} to
\(g(x)=f(R^{j/d}\sigma^j x)\), with \(r=\sigma^{-d}\)
and pieces \(f_\tau(R^{j/d}\sigma^j x)\), \(|\tau|=\sigma\).
On returning to the original coordinates, an averaging box of
angular width \(u\) has dimensions
\[
 R^{j/d}\sigma^{j-m-2}u^{m+2}\times
 R^{j/d}\sigma^{j-m-2}u^{m+1}\times\cdots\times
 R^{j/d}\sigma^{j-m-2}u^{m+2-d}.
\]
These boxes are \((u/\sigma)^{m-j}\times\frac{(uR^{1/d})^j}{(\sigma R^{1/d})^2}U_{\rho,R}\).
Since \(u\ge\sigma\) and \(j\le m\), we can apply
Lemma~\ref{lem:envelopechange}, with
\(F=(\sum_{\tau\subset\rho}|f_\tau|^2)^{p_{d-m}/4}\),
from these larger boxes to
\(B\parallel\frac{(uR^{1/d})^j}{(\sigma R^{1/d})^2}U_{\rho,R}\).
The output, with an arbitrarily small loss \(\delta>0\), is
\begin{equation}\label{eq:hcoverviewfullparent}
\begin{split}
 &\int_{\sum_{|\tau|=\sigma}|f_\tau|^2\gtrsim\alpha^2}
       \Bigl(\sum_{|\tau|=\sigma}|f_\tau|^2\Bigr)^{p_{d-m}/2}\\
 &\lesssim_\delta\sigma^{-\delta}
 \sum_{\sigma\le u\le1}\sum_{\rho\in\mathcal J(u)}
 \sum_{\substack{B\parallel\frac{(uR^{1/d})^j}{(\sigma R^{1/d})^2}U_{\rho,R}\\
 B\in\mathcal G_\rho(\frac{\alpha}{\sigma^{-\delta}})}}
 |B|\biggl\|\Bigl(\sum_{\tau\subset\rho,\ |\tau|=\sigma}|f_\tau|^2\Bigr)^{1/2}
 \biggr\|_{L^{p_{d-m}/2}_\#(B)}^{p_{d-m}}.
\end{split}
\end{equation}

For a parent \(\rho=[t_\rho,t_\rho+s]\) on the left hand side of
Proposition~\ref{prop:hcstable}, use \eqref{eq:hcoverviewfullparent}
when \(s\gg\sigma\), and the locally constant property when \(s\sim\sigma\).
In the first case, apply \eqref{eq:hcoverviewfullparent} to
\(f_\rho(L_{t_\rho,s}^{-T}x)\), with pieces
\(f_\tau(L_{t_\rho,s}^{-T}x)\), \(\tau\subset\rho\), where
\(L_{t_\rho,s}=(\phi_d(t_\rho),s\phi_d'(t_\rho),\ldots,s^d\phi_d^{(d)}(t_\rho))\).
The parameter is \(Rs^d\), the cap width is \(\sigma/s\), and
the amplitude is \(\alpha/N(s)\). When \(s\sim\sigma\), the locally constant property on the dual of the Fourier support box, whose dimensions are
\(R^{-j/d}\sigma^{-j}\times R^{-j/d}\sigma^{1-j}
\times\cdots\times R^{-j/d}\sigma^{d-j}\), gives directly
\begin{equation}\label{eq:hcoverviewcomparable}
\begin{split}
 &\sum_{\substack{Q\parallel\rho^*\\Q\in\mathcal G_\rho(\alpha)}}
 |Q|\biggl\|\Bigl(\sum_{\tau\subset\rho,\ |\tau|=\sigma}
 |f_\tau|^2\Bigr)^{1/2}\biggr\|_{L^{p_{d-m}/2}_\#(Q)}^{p_{d-m}}\lesssim\sum_{\substack{B\parallel\frac{(sR^{1/d})^j}{(\sigma R^{1/d})^2}U_{\rho,R}\\
 B\in\mathcal G_\rho(\alpha)}}
 |B|\biggl\|\Bigl(\sum_{\tau\subset\rho,\ |\tau|=\sigma}|f_\tau|^2\Bigr)^{1/2}
 \biggr\|_{L^{p_{d-m}/2}_\#(B)}^{p_{d-m}}.
\end{split}
\end{equation}
Returning to the original coordinates in \eqref{eq:hcoverviewfullparent},
we combine it with \eqref{eq:hcoverviewcomparable} and use
\(N(s)N_\rho(u)\lesssim N(u)\) from Lemma~\ref{lem:regular} to bound
the left hand side of Proposition~\ref{prop:hcstable} by
\begin{equation}\label{eq:hcoverviewcoarse}
 \lesssim_\delta\sigma^{-\delta}
 \sum_{\sigma\le u\le s}\sum_{\rho\in\mathcal J(u)}
 \sum_{\substack{B\parallel
       \frac{(uR^{1/d})^j}{(\sigma R^{1/d})^2}U_{\rho,R}\\
 B\in\mathcal G_\rho(\frac{\alpha}{\sigma^{-\delta}})}}
 |B|\biggl\|\Bigl(\sum_{\substack{\tau\subset\rho\\|\tau|=\sigma}}
                         |f_\tau|^2\Bigr)^{1/2}
       \biggr\|_{L^{p_{d-m}/2}_\#(B)}^{p_{d-m}}.
\end{equation}
The averages in
\eqref{eq:hcoverviewcoarse} still contain \(f_\tau\), whereas
Proposition~\ref{prop:hcstable} requires the original \(f_\theta\).
We make this replacement by rescaling when a parent contains few
\(\sigma\) intervals, and by refining when it contains many.
Fix \(u\). We estimate the corresponding term on the RHS of
\eqref{eq:hcoverviewcoarse}.

\medskip
\noindent\textit{Case 1: \(u<\sigma^{1-\eta}\).}
There are at most \(\sigma^{-\eta}\) intervals \(\tau\subset\rho\)
of width \(\sigma\). Apply \eqref{eq:selaverage} to the
fixed \(u\) sum in \eqref{eq:hcoverviewcoarse}, with
\(F=(\sum_{\tau\subset\rho,\ |\tau|=\sigma}|f_\tau|^2)^{p_{d-m}/4}\)
and \(t=[\alpha\sigma^\delta/N(u)]^{p_{d-m}/2}\).
The fixed \(u\) term is therefore bounded by
\[
 \lesssim\sigma^{-\delta}\sum_{\rho\in\mathcal J(u)}
 \int_{\sum_{\tau\subset\rho,\ |\tau|=\sigma}|f_\tau|^2
        \gtrsim[\alpha\sigma^\delta/N(u)]^2}
 \Bigl(\sum_{\tau\subset\rho,\ |\tau|=\sigma}|f_\tau|^2\Bigr)^{p_{d-m}/2}.
\]
In this integral, discard the terms
\(|f_\tau|\ll\alpha\sigma^\delta/(N(u)N_\rho(\sigma))\) and apply H\"older.
The fixed \(u\) term is bounded by
\begin{equation}\label{eq:hcoverviewnarrow}
 \lesssim\sigma^{-C\eta}
 \sum_{\rho\in\mathcal J(u)}
 \sum_{\substack{\tau\subset\rho\\|\tau|=\sigma}}
 \int_{|f_\tau|\gtrsim\alpha\sigma^\delta/(N(u)N_\rho(\sigma))}
 |f_\tau|^{p_{d-m}}.
\end{equation}
We use \eqref{eq:hcstablefunction} inductively at parameters
smaller than \(R\). For \(\tau=[t_\tau,t_\tau+\sigma]\),
the function \(f_\tau(L_{t_\tau,\sigma}^{-T}x)\) has parameter
\(R\sigma^d<R\). The left endpoint \(t_I\) is sent to
\((t_I-t_\tau)/\sigma\), and its coefficients become
\(\sigma^{-i}\lambda_i\).
For \(i<j\), \(|\sigma^{-i}\lambda_i|\le K^{-1}(R\sigma^d)^{-j/d}\),
while
\(|\sigma^{-j}\lambda_j-(R\sigma^d)^{-j/d}|
\le K^{-1}(R\sigma^d)^{-j/d}\).
Thus the Fourier support of
\(f_\tau(L_{t_\tau,\sigma}^{-T}\,\cdot)\) lies in
\(\Gamma_{j,1/K}^{d+1}(R\sigma^d)\) with the same \(K\).
Apply \eqref{eq:hcstablefunction} to this function with
amplitude \(\alpha\sigma^\delta/(N(u)N_\rho(\sigma))\).
The averaging boxes return to translates of
\((vR^{1/d})^mU_{\rho',R}\), where \(\rho'\subset\tau\)
and \(|\rho'|=v\), since
\begin{equation}\label{eq:hcscaleidentities}
 \sigma^{-i}(R\sigma^d)^{(m+2)/d}(v/\sigma)^{m+2-i}
 =R^{(m+2)/d}v^{m+2-i},\qquad 0\le i\le d.
\end{equation}
Moreover, \(N(u)N_\rho(\sigma)N_\tau(v)\lesssim N(v)\)
by Lemma~\ref{lem:regular} for the uniformized family.
Consequently \eqref{eq:hcoverviewnarrow} is bounded by
\begin{equation}\label{eq:hcoverviewnarrowreturn}
 \lesssim_\e\sigma^{-C\eta}(R\sigma^d)^\e
 \sum_{R^{-1/d}\le v\le\sigma}\sum_{\rho'\in\mathcal J(v)}
 \sum_{\substack{U\parallel(vR^{1/d})^mU_{\rho',R}\\
 U\in\mathcal G_{\rho'}(
       \frac{\alpha}{\sigma^{-C\eta}(R\sigma^d)^\e})}}
 |U|\biggl\|\Bigl(\sum_{\theta\subset\rho'}|f_\theta|^2\Bigr)^{1/2}
 \biggr\|_{L^{p_{d-m}/2}_\#(U)}^{p_{d-m}}.
\end{equation}
Since \(s\ge\sigma\),
\[
 \frac{\sigma^{-C\eta}(R\sigma^d)^\epsilon}
      {R^\epsilon(\sigma s)^\epsilon}
 \le \sigma^{(d-2)\epsilon-C\eta}.
\]
Since \eqref{eq:hcoverviewnarrowreturn} already contains
the original \(f_\theta\), the iteration stops in this case.
We keep this factor for \eqref{eq:hccontractions}, where
the bounds from both cases are combined.

\medskip
\noindent\textit{Case 2: \(u\ge\sigma^{1-\eta}\).}
For each \(\rho\in\mathcal J(u)\), Taylor rescaling sends
the boxes \(B\) in \eqref{eq:hcoverviewcoarse} to cubes of
side length \((Ru^d)^{j/d}(\sigma/u)^{-2}\), with parameter
\(Ru^d\) and cap width \(\sigma/u\).
We therefore use the following lemma to reduce either \(\sigma\) or \(u\), allowing us to continue the induction
on \(\sigma s\).

\begin{lemma}\label{lem:hcoverviewrefinement}
Let \(f\) be a Schwartz function with
\(\supp\widehat f\subset\Gamma_{j,1/K}^{d+1}(R)\).
For \(\eta>0\) sufficiently small and
\(R^{-1/d}<\sigma\) sufficiently small in terms of \(\eta\),
there are dyadic
\(\sigma',s'\) with
\(R^{-1/d}\le\sigma'\le\max\{R^{-1/d},\sigma^{3/2}\}\),
and \(\sigma\le s'\lesssim\sigma^{\eta^4}\) such that
\[
\begin{split}
 \sum_{\substack{B\parallel[0,R^{j/d}\sigma^{-2}]^{d+1}\\
 B\in\mathcal G_{[0,1]}(\alpha)}}
 |B|\biggl\|\Bigl(\sum_{|\tau|=\sigma}|f_\tau|^2\Bigr)^{1/2}
 \biggr\|_{L^{p_{d-m}/2}_\#(B)}^{p_{d-m}}
 &\lesssim_{\delta,\eta}\sigma^{-C\eta^4}
 \sum_{\substack{Q\parallel[0,1]^{d+1}\\
 Q\in\mathcal G_{[0,1]}(\frac{\alpha}{\sigma^{-C\eta^4}})}}
 |Q|\biggl\|\Bigl(\sum_{|\tau'|=\sigma'}|f_{\tau'}|^2\Bigr)^{1/2}
 \biggr\|_{L^{p_{d-m}/2}_\#(Q)}^{p_{d-m}}\\
 &\quad+\sigma^{-\delta}\sum_{\rho\in\mathcal J(s')}
 \sum_{\substack{Q\parallel\rho^*\\
 Q\in\mathcal G_\rho(\frac{\alpha}{\sigma^{-\delta}})}}
 |Q|\biggl\|\Bigl(\sum_{\substack{\tau\subset\rho\\|\tau|=\sigma}}
 |f_\tau|^2\Bigr)^{1/2}\biggr\|_{L^{p_{d-m}/2}_\#(Q)}^{p_{d-m}}.
\end{split}
\]
For \(p_{d-m}=4\), only the first sum over \(Q\) is needed.
\end{lemma}

Assuming Lemma~\ref{lem:hcoverviewrefinement}, we now complete Case~2.
For \(\rho=[t_\rho,t_\rho+u]\), the function
\(f_\rho(L_{t_\rho,u}^{-T}\,\cdot)\) has Fourier support in
\(\Gamma_{j,1/K}^{d+1}(Ru^d)\), by the calculation in Case~1
with \(u\) in place of \(\sigma\).
Since \(\sigma/u\le\sigma^\eta\le A^{-\eta}\), taking \(A\)
sufficiently large in terms of \(\eta\) allows us to apply
Lemma~\ref{lem:hcoverviewrefinement} to this function with
parameter \(Ru^d\), cap width \(\sigma/u\), and amplitude
\(\alpha\sigma^\delta/N(u)\).
Sum over \(\rho\) before choosing common dyadic scales;
the losses from these choices are absorbed into
\(\sigma^{-C\eta^4}\) in the first term and
\(\sigma^{-\delta}\) in the second.
Returning to the original coordinates and using
\(\sigma/u\le\sigma^\eta\) and \(\delta\ll\eta^4\),
we bound the fixed \(u\) term on the RHS of
\eqref{eq:hcoverviewcoarse} by
\begin{equation}\label{eq:hcoverviewrefinement}
\begin{split}
 &\lesssim\sigma^{-C\eta^4}
 \sum_{\rho\in\mathcal J(u)}
 \sum_{\substack{Q\parallel\rho^*\\
 Q\in\mathcal G_\rho(\frac{\alpha}{\sigma^{-C\eta^4}})}}
 |Q|\biggl\|\Bigl(\sum_{\substack{\tau'\subset\rho\\|\tau'|=\sigma'}}
 |f_{\tau'}|^2\Bigr)^{1/2}\biggr\|_{L^{p_{d-m}/2}_\#(Q)}^{p_{d-m}}\\
 &\quad+\sigma^{-2\delta}\sum_{\rho'\in\mathcal J(u')}
 \sum_{\substack{Q\parallel(\rho')^*\\
 Q\in\mathcal G_{\rho'}(\frac{\alpha}{\sigma^{-2\delta}})}}
 |Q|\biggl\|\Bigl(\sum_{\substack{\tau\subset\rho'\\|\tau|=\sigma}}
 |f_\tau|^2\Bigr)^{1/2}\biggr\|_{L^{p_{d-m}/2}_\#(Q)}^{p_{d-m}},
\end{split}
\end{equation}
where \(\sigma'\le\max\{R^{-1/d},\sigma^{1+\eta/4}\}\)
and \(u'\le u\sigma^{\eta^5/2}\). The first sum over \(Q\) is
comparable to the left hand side of Proposition~\ref{prop:hcstable} at
\((\sigma',u)\); the second is comparable to that left hand side at
\((\sigma,u')\), with the amplitudes in \eqref{eq:hcoverviewrefinement}.
In both cases the new product is smaller than \(\sigma u\).

Suppose first that \(\sigma'>R^{-1/d}\).
Applying Proposition~\ref{prop:hcstable} inductively at
\((\sigma',u)\), with amplitude
\(\alpha/\sigma^{-C\eta^4}\), bounds the first term
on the right hand side of \eqref{eq:hcoverviewrefinement} by
\[
 \lesssim_\epsilon
 \sigma^{-C\eta^4}R^\epsilon(\sigma'u)^\epsilon
 \sum_{R^{-1/d}\le v\le1}\sum_{\rho'\in\mathcal J(v)}
 \sum_{\substack{
 U\parallel(vR^{1/d})^mU_{\rho',R}\\
 U\in\mathcal G_{\rho'}(
 \frac{\alpha}{\sigma^{-C\eta^4}R^\epsilon(\sigma'u)^\epsilon})}}
 |U|\biggl\|\Bigl(
 \sum_{\theta\subset\rho'}|f_\theta|^2
 \Bigr)^{1/2}\biggr\|_{L^{p_{d-m}/2}_\#(U)}^{p_{d-m}}.
\]
Since \(\sigma'\le\sigma^{1+\eta/4}\) and \(u\le s\),
\[
 \frac{\sigma^{-C\eta^4}R^\epsilon(\sigma'u)^\epsilon}
      {R^\epsilon(\sigma s)^\epsilon}
 =
 \sigma^{-C\eta^4}
 \Bigl(\frac{\sigma'}{\sigma}\Bigr)^\epsilon
 \Bigl(\frac us\Bigr)^\epsilon
 \le \sigma^{\epsilon\eta/4-C\eta^4}.
\]

If instead \(\sigma'=R^{-1/d}\), the first sum over \(Q\) on the
right hand side of \eqref{eq:hcoverviewrefinement} already
contains \(f_\theta\).
Apply Lemma~\ref{lem:hcproducer} with \(\eta\) in place of
\(\epsilon\) to each function \(f_\rho(L_{t_\rho,u}^{-T}x)\),
\(\rho\in\mathcal J(u)\), with parameter \(Ru^d\) and
amplitude \(\alpha/(\sigma^{-C\eta^4}N(u))\).
For each returned parent \(\rho'\subset\rho\) of width \(v\),
Lemma~\ref{lem:regular} gives \(N(u)N_\rho(v)\lesssim N(v)\).
Using this bound in the selection condition and
\eqref{eq:hcscaleidentities} for the averaging boxes, we sum over
\(\rho\) to obtain the right hand side of
\eqref{eq:hcstablefunction} with \(R^\epsilon\) replaced,
both outside the sum and in \(\mathcal G_{\rho'}\), by
\(\sigma^{-C\eta^4}R^\eta\le R^{C\eta}\).

For the second term on the right hand side of
\eqref{eq:hcoverviewrefinement}, regardless of \(\sigma'\),
applying Proposition~\ref{prop:hcstable} inductively at
\((\sigma,u')\), with amplitude \(\alpha/\sigma^{-2\delta}\),
gives
\[
 \lesssim_\epsilon
 \sigma^{-2\delta}R^\epsilon(\sigma u')^\epsilon
 \sum_{R^{-1/d}\le v\le1}\sum_{\rho'\in\mathcal J(v)}
 \sum_{\substack{
 U\parallel(vR^{1/d})^mU_{\rho',R}\\
 U\in\mathcal G_{\rho'}(
 \frac{\alpha}{\sigma^{-2\delta}R^\epsilon(\sigma u')^\epsilon})}}
 |U|\biggl\|\Bigl(
 \sum_{\theta\subset\rho'}|f_\theta|^2
 \Bigr)^{1/2}\biggr\|_{L^{p_{d-m}/2}_\#(U)}^{p_{d-m}}.
\]
Since \(u'\le u\sigma^{\eta^5/2}\) and \(u\le s\),
\[
 \frac{\sigma^{-2\delta}R^\epsilon(\sigma u')^\epsilon}
      {R^\epsilon(\sigma s)^\epsilon}
 =
 \sigma^{-2\delta}\Bigl(\frac{u'}s\Bigr)^\epsilon
 \le \sigma^{\epsilon\eta^5/2-2\delta}.
\]

For the initial case \(\sigma=R^{-1/d}\), apply
Lemma~\ref{lem:hcproducer} with \(\eta\) in place of \(\epsilon\) to
\(f_\rho(L_{t_\rho,s}^{-T}x)\), \(\rho\in\mathcal J(s)\),
with parameter \(Rs^d\) and amplitude \(\alpha/N(s)\).
Its output is the right hand side of
\eqref{eq:hcstablefunction}, with both occurrences of
\(R^\epsilon\) replaced by \(R^\eta\).

The estimates for the two terms in
\eqref{eq:hcoverviewrefinement}
(with \(\sigma'>R^{-1/d}\) for the first term)
and the bound \eqref{eq:hcoverviewnarrowreturn}
give, respectively, the following bounds for the ratios
of the outer coefficients and amplitude denominators
to \(R^\epsilon(\sigma s)^\epsilon\):
\begin{equation}\label{eq:hccontractions}
 \sigma^{\epsilon\eta/4-C\eta^4},\qquad
 \sigma^{\epsilon\eta^5/2-2\delta},\qquad
 \sigma^{(d-2)\epsilon-C\eta}.
\end{equation}
Our choices of \(\eta,\delta\) make all three exponents in
\eqref{eq:hccontractions} positive.
Since \(\sigma\le A^{-1}\), taking \(A\) sufficiently large
therefore allows these factors to absorb the implicit constants
and the dyadic sum over \(u\) in \eqref{eq:hcoverviewcoarse}.
For \(\sigma'=R^{-1/d}\) or \(\sigma=R^{-1/d}\),
the inequalities \(s\ge\sigma\ge R^{-1/d}\) give
\(\frac{R^{C\eta}}{R^\epsilon(\sigma s)^\epsilon}
\le R^{-\epsilon/C}\), which also absorbs the constants
and the dyadic sum over \(u\) for sufficiently large \(R\).

Summing the estimates for the terms in
\eqref{eq:hcoverviewcoarse} therefore proves
Proposition~\ref{prop:hcstable} at the current \(R\).
Applying Lemma~\ref{lem:hcseed} as in the derivation of
\eqref{eq:hcstablefunction} then gives
\eqref{eq:hcstablefunction} at the same \(R\),
completing the induction on \(R\).

We prove Lemmas~\ref{lem:hcproducer} and~\ref{lem:hcoverviewrefinement}
for \(m<d-2\) in Section~\ref{sec:hcpositive}.

\begin{remark}
To complete Step~1 of Proposition~\ref{prop:hcstablefunction}, we
also need \eqref{eq:hcrestrictedconvolution} for the current \(m\).
In Section~\ref{sec:hcpositive}, we first bound the left hand side
of Lemma~\ref{lem:hcproducer} by the RHS of
\eqref{eq:hcgeneralconvolution} with \(R=r\) and \(f=g\).
Lemma~\ref{lem:hcstep2averagecomparison} then gives
Lemma~\ref{lem:hcproducer} as stated.

Following the induction in Section~\ref{subsec:hcoverviewiteration},
we use this strengthened form of Lemma~\ref{lem:hcproducer} together
with Lemma~\ref{lem:hcoverviewrefinement} to strengthen
Proposition~\ref{prop:hcstable}; we bound its left hand side by
the RHS of \eqref{eq:hcgeneralconvolution} for the current \(m\),
with both occurrences of \(R^\epsilon\) replaced by
\(R^\epsilon(\sigma s)^\epsilon\).
Applying this bound at \(\sigma=A^{-1}\) to the sum in
Lemma~\ref{lem:hcseed} gives \eqref{eq:hcrestrictedconvolution}
at the current \(R\), as in Section~\ref{subsec:hcroadmap}.
Together with \(C_K^m(R)\lessapprox1\) from
\eqref{eq:hcstablefunction}, this completes Step~1.
\end{remark}

\subsection{Why the Taylor order induction is needed}
\label{subsec:hcoverviewpositive}

The proof of Proposition~\ref{prop:hcstable} has been reduced to
Lemmas~\ref{lem:hcproducer} and~\ref{lem:hcoverviewrefinement}.
We explain the Taylor order induction in Lemma~\ref{lem:hcproducer},
first for \(m=0\). Let \(g\) be a Schwartz function with
\(\supp\widehat g\subset\Gamma^{d+1}(r)\). The integral to be
bounded is
\begin{equation}\label{eq:hcoverviewpositiveinput}
 \int_{\sum_\theta|g_\theta|^2\gtrsim\alpha^2}
       \bigl(\sum_\theta|g_\theta|^2\bigr)^{p_d/2},
 \qquad \theta\in\Xi(r^{-1/d}).
\end{equation}
We split the Fourier transform of \(\sum_\theta|g_\theta|^2\)
at scales adapted to the averaging boxes on the right hand side of
Lemma~\ref{lem:hcproducer}. The lowest frequencies will give
its \(u=1\) term. To see which coordinates must be cut off,
compare the dimensions
\[
 \theta^*:1\times r^{1/d}\times\cdots\times r,
 \qquad U_{[0,1],r}:r^{2/d}\times\cdots\times r^{2/d}.
\]
The locally constant property already gives averaging on scale \(r^{2/d}\)
in directions \(2,\ldots,d\). To obtain the \(U_{[0,1],r}\)
average, we therefore restrict the first two Taylor coordinates
to width \(r^{-2/d}\).

For an interval \(\tau\) of length \(u\), the corresponding
averaging box is
\[
 U_{\tau,r}:r^{2/d}u^2\times r^{2/d}u\times\cdots\times r^{2/d}u^{2-d}.
\]
Write
\(\xi=\sum_i\lambda_i\phi_d^{(i)}(t_\theta)
      =\sum_i\widetilde\lambda_i\phi_d^{(i)}(t_\tau)\).
For \(\theta\subset\tau\), put \(h=t_\theta-t_\tau\).
On \(\supp\widehat{|g_\theta|^2}\),
\(|\lambda_i|\lesssim r^{-i/d}\), and Taylor expansion gives
\[
 \widetilde\lambda_i-\lambda_i
 =\frac{h^i}{i!}\lambda_0+
   \frac{h^{i-1}}{(i-1)!}\lambda_1
   +O(r^{-2/d}u^{i-2}),\qquad 2\le i\le d.
\]
Thus the restrictions
\(|\lambda_0|\lesssim r^{-2/d}u^{-2}\) and
\(|\lambda_1|\lesssim r^{-2/d}u^{-1}\) keep the change in
the \(i\)th coefficient within the \(i\)th width of
\(U_{\tau,r}^*\). Put \(\sigma=r^{-1/d}/u\).
Choose \(\chi\in C_c^\infty(\mathbb R)\), \(0\le\chi\le1\),
with \(\chi=1\) on \([-C,C]\), and set
\[
 \eta^2_{\theta,\le\sigma}(\xi)
 =\chi\!\Bigl(\frac{\lambda_0}{\sigma^2}\Bigr)
  \chi\!\Bigl(\frac{\lambda_1}{r^{-1/d}\sigma}\Bigr)
  \prod_{i=2}^d\chi(r^{i/d}\lambda_i),
 \qquad
 \eta^2_{\theta,\sigma}
 =\eta^2_{\theta,\le\sigma}-\eta^2_{\theta,\le\sigma/2}.
\]

On \(\operatorname{supp}\widehat{|g_\theta|^2}\), we have
\(\prod_{i=2}^d\chi(r^{i/d}\lambda_i)=1\).
Thus the factors with \(2\le i\le d\) do not change
\(|g_\theta|^2*\check\eta^2_{\theta,\le\sigma}\).
At \(\sigma=r^{-1/d}\), the cutoff widths in directions
\(0,1\) are both \(r^{-2/d}\),
and \(\eta^2_{\theta,\le1}=1\) on
\(\supp\widehat{|g_\theta|^2}\).

Dyadic pigeonholing now gives \(r^{-1/d}<\sigma\le1\) such that,
with
\[
 z_{\mathrm{low}}=\sum_\theta|g_\theta|^2*
       \check\eta^2_{\theta,\le r^{-1/d}},\qquad
 z_{\mathrm{high}}=\sum_\theta|g_\theta|^2*
       \check\eta^2_{\theta,\sigma},
\]
the integral \eqref{eq:hcoverviewpositiveinput} is bounded by
\begin{equation}\label{eq:hcoverviewdecomposition}
 \lessapprox\int_{|z_{\mathrm{low}}|\gtrapprox\alpha^2}|z_{\mathrm{low}}|^{p_d/2}
 +\int_{|z_{\mathrm{high}}|\gtrapprox\alpha^2}|z_{\mathrm{high}}|^{p_d/2}.
\end{equation}
For the first integral, the locally constant property, H\"older's inequality, and
\eqref{eq:weightedYoung} give
\begin{equation}\label{eq:hcoverviewlowtail}
 \int_{|z_{\mathrm{low}}|\gtrsim\alpha^2}|z_{\mathrm{low}}|^{p_d/2}
 \lesssim\sum_{\substack{U\parallel U_{[0,1],r}\\
                         U\in\mathcal G_{[0,1]}(\alpha)}}
 |U|\biggl\|\bigl(\sum_\theta|g_\theta|^2\bigr)^{1/2}
      \biggr\|_{L_\#^{p_d/2}(U)}^{p_d}.
\end{equation}
This is the \(u=1\) term in Lemma~\ref{lem:hcproducer}.

We next show why estimating the \(z_{\mathrm{high}}\) term in
\eqref{eq:hcoverviewdecomposition} leads to the Taylor cone of order one. On
\(\supp\widehat{|g_\theta|^2*\check\eta^2_{\theta,\sigma}}\),
the coefficients satisfy
\begin{equation}\label{eq:hcoverviewfirstshell}
 |\lambda_0|\lesssim\sigma^2,\quad
 |\lambda_1|\lesssim r^{-1/d}\sigma,\quad
 \max\{\sigma^{-2}|\lambda_0|,
        r^{1/d}\sigma^{-1}|\lambda_1|\}\sim1.
\end{equation}
Group the \(\theta\)'s in intervals \(\tau\) of length
\(u=r^{-1/d}\sigma^{-1}\), so that
\(z_{\mathrm{high}}=\sum_\tau\sum_{\theta\subset\tau}
|g_\theta|^2*\check\eta^2_{\theta,\sigma}\).
By \eqref{eq:hcoverviewfirstshell}, frequency dilation by
\(\sigma^{-2}\) sends each inner sum to a cap with coefficients
\begin{equation}\label{eq:hcoverviewrescaledshell}
 |\lambda_i|\lesssim(r\sigma^d)^{-i/d}\ (0\le i\le d),
 \qquad \max\{|\lambda_0|,(r\sigma^d)^{1/d}|\lambda_1|\}\sim1.
\end{equation}
On the cone over the moment curve in Proposition~\ref{prop:highercone}, the radial
coefficient satisfies \(|\lambda_0|\sim1\). But in
\eqref{eq:hcoverviewrescaledshell}, \(|\lambda_0|\) can be small,
in which case \((r\sigma^d)^{1/d}|\lambda_1|\sim1\).
This is precisely the additional
support allowed by \(\Gamma_1^{d+1}(r\sigma^d)\) in
\eqref{eq:hcoverviewtaylorcones}. More explicitly,
\begin{equation}\label{eq:hcoverviewfirsttaylorcap}
 \sigma^{-2}\supp\widehat{\Bigl(
 \sum_{\theta\subset\tau}|g_\theta|^2*
                         \check\eta^2_{\theta,\sigma}\Bigr)}
 \subset\omega_\tau\in\Xi_1^{d+1}(r\sigma^d).
\end{equation}
where \(\omega_\tau\) is the cap corresponding to \(\tau\).

In general, for \(\supp\widehat g\subset\Gamma_m^{d+1}(r)\), the low frequency
integral must give the average over \(r^{m/d}U_{[0,1],r}\),
while, for \(m<d-2\), the high frequency support must lie at
Taylor order \(m+1\). We therefore define
\[
 \eta^{m+2}_{\theta,\le\sigma}(\xi)
 =\prod_{i=0}^{m+1}
   \chi\!\Bigl(\frac{\lambda_i}{r^{-i/d}\sigma^{m+2-i}}\Bigr)
  \prod_{i=m+2}^d\chi(r^{i/d}\lambda_i),
 \qquad
 \eta^{m+2}_{\theta,\sigma}
 =\eta^{m+2}_{\theta,\le\sigma}-\eta^{m+2}_{\theta,\le\sigma/2}.
\]
On \(\operatorname{supp}\widehat{|g_\theta|^2}\), we have
\(\prod_{i=m+2}^d\chi(r^{i/d}\lambda_i)=1\).
Thus the factors with \(m+2\le i\le d\) do not change
\(|g_\theta|^2*\check\eta^{m+2}_{\theta,\le\sigma}\).
For dyadic \(r^{-1/d}<\sigma\le1\), set
\[
 z_{\mathrm{low}}^m=\sum_{\theta\in\Xi_m^{d+1}(r)}|g_\theta|^2*
                         \check\eta^{m+2}_{\theta,\le r^{-1/d}},
 \qquad
 z_{\mathrm{high}}^m=\sum_{\theta\in\Xi_m^{d+1}(r)}|g_\theta|^2*
                         \check\eta^{m+2}_{\theta,\sigma}.
\]
As in \eqref{eq:hcoverviewdecomposition}, dyadic pigeonholing reduces
the estimate in Lemma~\ref{lem:hcproducer} to two integrals with
exponent \(p_{d-m}/2\), one for \(z_{\mathrm{low}}^m\) and one for
\(z_{\mathrm{high}}^m\) at a single dyadic scale \(\sigma\).

For the low frequency integral, the cutoff at \(\sigma=r^{-1/d}\)
has width \(r^{-(m+2)/d}\) in directions \(0,\ldots,m+1\).
By the locally constant property at scale \(r^{(m+2)/d}\) in the remaining
directions, we can apply the argument of \eqref{eq:hcoverviewlowtail}
to \(z_{\mathrm{low}}^m\), with \(p_{d-m}\) in place of \(p_d\)
and boxes \(U\parallel r^{m/d}U_{[0,1],r}\).
This gives the \(u=1\) term on the right hand side of
Lemma~\ref{lem:hcproducer}.

For the high frequency integral, the same Taylor calculation
as in \eqref{eq:hcoverviewfirsttaylorcap} gives
\begin{equation}\label{eq:highfrequencyhigherTaylororder}
 \sigma^{-(m+2)}\supp\widehat{\Bigl(
 \sum_{\theta\subset\tau}|g_\theta|^2*
                       \check\eta^{m+2}_{\theta,\sigma}\Bigr)}
 \subset\omega_\tau\in\Xi_{m+1}^{d+1}(r\sigma^d),
 \qquad 0\le m\le d-3.
\end{equation}
Thus \eqref{eq:hcgeneralconvolution}, with \(m+1\) in place of
\(m\), applies to the sum of
\(\bigl(\sum_{\theta\subset\tau}|g_\theta|^2*
\check\eta^{m+2}_{\theta,\sigma}\bigr)(\sigma^{-(m+2)}x)\),
with parameter \(r\sigma^d\).
Its exponent is \(p_{d-m-1}\), whereas the high frequency
integral has exponent \(p_{d-m}/2\). Section~\ref{sec:hcpositive}
performs the wave packet interpolation between these exponents and
first bounds the integral in Lemma~\ref{lem:hcproducer} by
the right hand side of \eqref{eq:hcgeneralconvolution} for the
function \(g\) at parameter \(r\), before passing to the
sum over \(U\) in Lemma~\ref{lem:hcproducer}.
At \(m=d-2\), the exponent of the high frequency integral is
\(p_2/2=2\), and its proof uses \(L^2\) orthogonality instead.

The averaging box \((ur^{1/d})^mU_{\tau,r}\) also follows from
rescaling the \(u=1\) averaging box. For \(\tau=[t_\tau,t_\tau+u]\),
rescale \(g_\tau\) to \(g_\tau(L_{t_\tau,u}^{-T}x)\), with
parameter \(ru^d\). Returning the box
\((ru^d)^{m/d}U_{[0,1],ru^d}\) multiplies its \(i\)th width
by \(u^{-i}\), giving
\[
 r^{(m+2)/d}u^{m+2}\times r^{(m+2)/d}u^{m+1}
 \times\cdots\times r^{(m+2)/d}u^{m+2-d}.
\]
These are the dimensions of \((ur^{1/d})^mU_{\tau,r}\).

\subsection{Roadmap: Proof of Proposition~\ref{prop:hcstablefunction}, Step 2}
\label{subsec:hcoverviewcoefficients}

In this subsection, we reduce Step~2 of Proposition~\ref{prop:hcstablefunction}
to Proposition~\ref{prop:hcstep2positive}, which we prove in
Section~\ref{sec:hcgeneral} by adapting the induction on \(\sigma s\)
from Step~1.

Fix \(0\le m<d-2\). Under the hypotheses of
Proposition~\ref{prop:hcstablefunction}, Step~1 gives
\(C_K^m(R)\lessapprox1\) and \eqref{eq:hcrestrictedconvolution}
for this \(m\). We will prove that for \(R\ge1\),
\(\epsilon>0\), and \(\alpha>0\), every Schwartz function \(f\)
with \(\operatorname{supp}\widehat f\subset\Gamma_j^{d+1}(R)\)
satisfies
\begin{equation}\label{eq:hcstep2inductiongoal}
 \int_{|f|>\alpha}|f|^{p_{d-m}}
 \lesssim_\epsilon
 \text{RHS of \eqref{eq:hcgeneralconvolution}}
 \lessapprox
 \text{RHS of \eqref{eq:hcstablefunction}}.
\end{equation}
The last inequality follows from
Lemma~\ref{lem:hcstep2averagecomparison}, so it suffices to prove
the first inequality. We do this successively in the order
\(j=0,\ldots,m\), keeping the exponent \(p_{d-m}\).
The case \(j=0\) follows directly from Step~1 after the
coefficient decomposition and dilation below.
Taking \(j=m\) will then give \(C^m(R)\lessapprox1\) and
\eqref{eq:hcgeneralconvolution} for the current \(m\).

Let \(f\) be a Schwartz function with
\(\operatorname{supp}\widehat f\subset\Gamma_m^{d+1}(R)\).
For \(\xi\in\operatorname{supp}\widehat f_{\theta_I}\), write
\(\xi=\sum_{i=0}^d\lambda_i\phi_d^{(i)}(t_I)\).
The condition defining \(\Gamma_m^{d+1}(R)\) gives
\(|\lambda_j|\sim R^{-j/d}\) for some \(j\le m\).
Using a smooth partition of unity in frequency space, we obtain
\[
 f_\theta
 =\sum_{j=0}^m
  \sum_{\substack{\nu\in\mathbb Z\\|\nu|\sim K}}
  \sum_{R^{-1/d}\le\ell\le1}f_\theta^{j,\nu,\ell},
\]
where \(\ell\) is dyadic and
\[
 \operatorname{supp}\widehat{f_{\theta_I}^{j,\nu,\ell}}
 \subset
 \left\{
 \sum_{i=0}^d\lambda_i\phi_d^{(i)}(t_I):
 \begin{array}{l}
 |\lambda_i|\lesssim R^{-i/d}\quad(0\le i\le d),\\
 |R^{j/d}\lambda_j-\nu/K|\lesssim K^{-1},\\
 \displaystyle
 \max\bigl\{1,(R^{j/d}|\lambda_i|)^{1/(j-i)}:
             0\le i<j\bigr\}\sim\ell^{-1}
 \end{array}
 \right\}.
\]

By Lemma~\ref{lem:regular}, we may assume that each family
\(\{f_\theta^{j,\nu,\ell}\}_{\theta\in\Xi_m^{d+1}(R)}\)
satisfies the counting inequalities \eqref{eq:regularcounts}.
Dyadic pigeonholing then gives fixed \(j,\nu,\ell\) such that,
putting \(\widetilde f_\theta=f_\theta^{j,\nu,\ell}\) and
\(\widetilde f=\sum_{\theta\in\Xi_m^{d+1}(R)}
\widetilde f_\theta\),
\begin{equation}\label{eq:hcstep2tailreduction}
 \int_{|f|>\alpha}|f|^{p_{d-m}}
 \lessapprox
 \int_{|\widetilde f|\gtrapprox\alpha}
 |\widetilde f|^{p_{d-m}}.
\end{equation}

Since \(|\nu|\sim K\), replacing each \(\widetilde f_\theta(x)\)
by \(\widetilde f_\theta((K/\nu)x)\), we may assume \(\nu=K\).
Thus \(\operatorname{supp}\widehat{\widetilde f}
\subset\Gamma_{j,1/K,\ell}^{d+1}(R)\), where
\begin{equation}\label{eq:hcstep2support}
 \Gamma_{j,1/K,\ell}^{d+1}(R)
 :=
 \bigcup_{I\in\mathcal J(R^{-1/d})}
 \left\{
 \sum_{i=0}^d\lambda_i\phi_d^{(i)}(t_I):
 \begin{array}{l}
 |\lambda_i|\lesssim R^{-i/d}\quad(0\le i\le d),\\
 |\lambda_j-R^{-j/d}|\lesssim K^{-1}R^{-j/d},\\
 \displaystyle
 \max\bigl\{1,(R^{j/d}|\lambda_i|)^{1/(j-i)}:
             0\le i<j\bigr\}\sim\ell^{-1}
 \end{array}
 \right\}.
\end{equation}

We will use induction for this \(j\) since, after dilation,
\(\operatorname{supp}
\widehat{\widetilde f(R^{j/d}\ell^j\,\cdot)}
\subset\Gamma_{j-1}^{d+1}(\ell^{-d})\), as we verify below.

By the following lemma and \eqref{eq:hcstep2tailreduction},
to prove \eqref{eq:hcstep2inductiongoal} for the original \(f\),
it suffices to show that, for \(\epsilon>0\) and \(\alpha>0\),
\[
 \int_{|\widetilde f|>\alpha}|\widetilde f|^{p_{d-m}}
 \lesssim_\epsilon
 \text{RHS of \eqref{eq:hcgeneralconvolution} with }
 \widetilde f\text{ in place of }f.
\]

\begin{lemma}\label{lem:hcstep2cutoffcomparison}
Let \(0\le m<d-2\), \(R\ge1\), and let \(f\) be a Schwartz
function with
\(\operatorname{supp}\widehat f\subset\Gamma_m^{d+1}(R)\).
Let \(\widetilde f\) be obtained from \(f\) by the coefficient
decomposition, uniformization, and dilation above.
For \(\epsilon>0\) and \(\alpha>0\),
\begin{equation}\label{eq:hcstep2cutoffreturn}
 \text{RHS of \eqref{eq:hcgeneralconvolution} for }\widetilde f
 \lesssim_\epsilon
 \text{RHS of \eqref{eq:hcgeneralconvolution} for }f.
\end{equation}
\end{lemma}

We defer the proof to the end of this section.

In summary, assuming \eqref{eq:hcstep2inductiongoal}
for every \(j'<j\), it suffices to prove
the following proposition.

\begin{proposition}\label{prop:hcstep2reducedtail}
Fix \(0\le m<d-2\) and \(0\le j\le m\).
Assume the hypotheses of Proposition~\ref{prop:hcstablefunction},
together with \(C_K^m(R)\lessapprox1\) and
\eqref{eq:hcrestrictedconvolution} for this \(m\),
and suppose that \eqref{eq:hcstep2inductiongoal} holds
for every \(j'<j\).
Let \(R\ge1\), let \(R^{-1/d}\le\ell\le1\) be dyadic,
and let \(f\) be a Schwartz function with
\(\operatorname{supp}\widehat f
\subset\Gamma_{j,1/K,\ell}^{d+1}(R)\).
Then, for \(\epsilon>0\) and \(\alpha>0\),
\[
 \int_{|f|>\alpha}|f|^{p_{d-m}}
 \lesssim_\epsilon
 \text{RHS of \eqref{eq:hcgeneralconvolution}}.
\]
\end{proposition}

We now prove Proposition~\ref{prop:hcstep2reducedtail}.
Let \(f\) be a Schwartz function with
\(\operatorname{supp}\widehat f
\subset\Gamma_{j,1/K,\ell}^{d+1}(R)\).
We first consider \(\ell\sim1\).
In this case, the lower coefficients satisfy
\(|\lambda_i|\lesssim R^{-j/d}\) for \(i<j\), whereas the
support hypothesis in \eqref{eq:hcrestrictedconvolution}
requires an additional factor \(K^{-1}\).
We obtain this factor by Taylor rescaling on sufficiently short caps.

Choose \(a\le A^{-1}\) sufficiently small in terms of \(K\).
For \(\tau\in\mathcal J(a)\), let \(t_\tau\) denote the left endpoint of \(\tau\).
Taylor rescaling by \(L_{t_\tau,a}^{-T}\), which sends \(\tau\)
to a cap over the unit interval, gives parameter \(Ra^d\) and
\[
 \frac{|a^{-i}\lambda_i|}{(Ra^d)^{-j/d}}
 =a^{j-i}R^{j/d}|\lambda_i|
 \lesssim a^{j-i}\le K^{-1},
 \qquad i<j.
\]
Thus
\(\operatorname{supp}\widehat{f_\tau(L_{t_\tau,a}^{-T}\,\cdot)}
\subset\Gamma_{j,1/K}^{d+1}(Ra^d)\).
H\"older and \eqref{eq:hcrestrictedconvolution}, applied to
these functions at parameter \(Ra^d\) and amplitude
\(\alpha/N(a)\), give, on returning to the original coordinates,
\begin{equation}\label{eq:hcstep2boundedtail}
 \int_{|f|>\alpha}|f|^{p_{d-m}}
 \lesssim
 \sum_{\tau\in\mathcal J(a)}
 \int_{|f_\tau|\gtrsim\alpha/N(a)}|f_\tau|^{p_{d-m}}
 \lesssim_\epsilon
 \text{RHS of \eqref{eq:hcgeneralconvolution} for }f.
\end{equation}
This also treats the initial index \(j=0\), since only
\(\ell\sim1\) occurs in that case.

Suppose now that \(\ell\ll1\).
For \(\tau\in\mathcal J(\ell)\), Taylor expansion at \(t_\tau\)
and the conditions defining \(\Gamma_{j,1/K,\ell}^{d+1}(R)\) give
\begin{equation}\label{eq:hcstep2lowerorder}
 \operatorname{supp}
 \widehat{f_\tau(R^{j/d}\ell^j\,\cdot)}
 \subset
 \left\{
 \sum_{i=0}^{d}\lambda_i\phi_d^{(i)}(t_\tau):
 \begin{array}{l}
 |\lambda_i|\lesssim\ell^i,\quad 0\le i\le d,\\
 \displaystyle\max_{i<j}\ell^{-i}|\lambda_i|\gtrsim1
 \end{array}
 \right\}.
\end{equation}
Thus
\(\operatorname{supp}\widehat{f(R^{j/d}\ell^j\,\cdot)}
\subset\Gamma_{j-1}^{d+1}(\ell^{-d})\).
Apply \eqref{eq:hcstablefunction}, available by the induction
hypothesis at \(j-1\), to \(f(R^{j/d}\ell^j x)\), with
parameter \(\ell^{-d}\) and amplitude \(\alpha\).
Returning to the original coordinates and applying
\eqref{eq:selaverage} to the sums over averaging boxes on the RHS of
\eqref{eq:hcstablefunction}, with
\(F=(\sum_{\tau\subset\rho,\,|\tau|=\ell}
|f_\tau|^2)^{p_{d-m}/4}\) and
\(t=[\alpha/(R^\epsilon N(s))]^{p_{d-m}/2}\),
followed by the locally constant property on \(Q\parallel\rho^*\), gives
\begin{equation}\label{eq:hcstep2coarsetail}
 \int_{|f|>\alpha}|f|^{p_{d-m}}
 \lesssim_\epsilon R^\epsilon
 \sum_{\ell\le s\le1}
 \sum_{\rho\in\mathcal J(s)}
 \sum_{\substack{Q\parallel\rho^*\\
 Q\in\mathcal G_\rho(\alpha/R^\epsilon)}}
 |Q|
 \biggl\|
 \Bigl(
 \sum_{\substack{\tau\subset\rho\\|\tau|=\ell}}
 |f_\tau|^2
 \Bigr)^{1/2}
 \biggr\|_{L_\#^{p_{d-m}/2}(Q)}^{p_{d-m}}.
\end{equation}

We next subdivide the \(\ell\) caps so that either the pieces
are the original \(f_\theta\), or the hypothesis
\eqref{eq:hcrestrictedconvolution} for the current \(m\)
applies after Taylor rescaling.
For this purpose, set \(\sigma=\max\{R^{-1/d},a\ell\}\).
When \(\sigma>R^{-1/d}\), the support condition needed to apply
\eqref{eq:hcrestrictedconvolution} is
\begin{equation}\label{eq:hcstableexit}
 \operatorname{supp}
 \widehat{f_\tau(L_{t_\tau,\sigma}^{-T}\,\cdot)}
 \subset\Gamma_{j,1/K}^{d+1}(R\sigma^d),
 \qquad \tau\in\mathcal J(\sigma).
\end{equation}
This follows from the \(j\)th coefficient condition in the
definition of \(\Gamma_{j,1/K,\ell}^{d+1}(R)\) and
\(\frac{|\sigma^{-i}\lambda_i|}{(R\sigma^d)^{-j/d}}
\lesssim(\sigma/\ell)^{j-i}\le K^{-1}\) for \(i<j\).

In either case, each \(\tau'\in\mathcal J(\ell)\) contains
\(\ell/\sigma\le a^{-1}\) caps of width \(\sigma\).
Cauchy--Schwarz therefore gives
\(\lvert f_{\tau'}\rvert^2\lesssim_a
\sum_{\tau\subset\tau',\,|\tau|=\sigma}|f_\tau|^2\),
so the RHS of \eqref{eq:hcstep2coarsetail} is bounded by
\begin{equation}\label{eq:hcoverviewshelltail}
 \lesssim_\epsilon R^\epsilon
 \sum_{\ell\le s\le1}
 \sum_{\rho\in\mathcal J(s)}
 \sum_{\substack{Q\parallel\rho^*\\
 Q\in\mathcal G_\rho(\alpha/R^\epsilon)}}
 |Q|
 \biggl\|
 \Bigl(
 \sum_{\substack{\tau\subset\rho\\|\tau|=\sigma}}
 |f_\tau|^2
 \Bigr)^{1/2}
 \biggr\|_{L_\#^{p_{d-m}/2}(Q)}^{p_{d-m}}.
\end{equation}
We estimate the sum over \(\rho\) and \(Q\) in
\eqref{eq:hcoverviewshelltail}, for each fixed \(s\), by the
following extension of Proposition~\ref{prop:hcstable}.
In its proof in Section~\ref{sec:hcgeneral}, we use
\eqref{eq:hcrestrictedconvolution} for the current \(m\),
since the range \(\sigma\le a\ell\) gives the support inclusion
\eqref{eq:hcstableexit}.
We state a stronger estimate, with the additional factor
\((\sigma s)^\epsilon\), to close the induction on \(\sigma s\)
in its proof.

\begin{proposition}\label{prop:hcstep2positive}
\label{cor:hcoverviewshellpositive}
\label{cor:hccurvaturepositive}
Assume the hypotheses of Proposition~\ref{prop:hcstablefunction}
for a fixed \(0\le m<d-2\), together with
\(C_K^m(R)\lessapprox1\) and
\eqref{eq:hcrestrictedconvolution} for this \(m\).
Let \(1\le j\le m\), \(R^{-1/d}\le\ell\ll1\), and let
\(f\) be a Schwartz function with
\(\supp\widehat f\subset\Gamma_{j,1/K,\ell}^{d+1}(R)\)
For
\(
 R^{-1/d}\le\sigma
 \le\max\{R^{-1/d},a\ell\},
 ~ \sigma\le s\le1,
\)
and \(\alpha>0\), we have
\[
\begin{aligned}
 &\sum_{\rho\in\mathcal J(s)}
 \sum_{\substack{Q\parallel\rho^*\\Q\in\mathcal G_\rho(\alpha)}}
 |Q|\biggl\|\Bigl(
 \sum_{\substack{\tau\subset\rho\\|\tau|=\sigma}}
 |f_\tau|^2\Bigr)^{1/2}
 \biggr\|_{L^{p_{d-m}/2}_\#(Q)}^{p_{d-m}}\lesssim_\epsilon
 \text{RHS of
 \eqref{eq:hcgeneralconvolution}},
\end{aligned}
\]
with both occurrences of \(R^\epsilon\) replaced by
\(R^\epsilon(\sigma s)^\epsilon\).
\end{proposition}

Assuming Proposition~\ref{prop:hcstep2positive}, use
\eqref{eq:hcoverviewshelltail} with \(\epsilon/2\) in place of
\(\epsilon\). For each fixed \(s\), apply the proposition to the sum
over \(\rho\) and \(Q\) in \eqref{eq:hcoverviewshelltail}, with
\(\epsilon/2\) in place of \(\epsilon\) and amplitude
\(\alpha/R^{\epsilon/2}\).
The outer coefficient and amplitude denominator are bounded by
\(R^{\epsilon/2}R^{\epsilon/2}(\sigma s)^{\epsilon/2}
\le R^\epsilon\ell^{\epsilon/2}s^{\epsilon/2}\le R^\epsilon\).
Summing over \(s\) gives
\begin{equation}\label{eq:hcstep2shellconvolution}
 \int_{|f|>\alpha}|f|^{p_{d-m}}
 \lesssim_\epsilon
 \ell^{\epsilon/2}\sum_{\ell\le s\le1}s^{\epsilon/2}
 \times\bigl[\text{RHS of \eqref{eq:hcgeneralconvolution}}\bigr]
 \lesssim_\epsilon
 \text{RHS of \eqref{eq:hcgeneralconvolution}}.
\end{equation}
Together with \eqref{eq:hcstep2boundedtail} for \(\ell\sim1\),
this shows that Proposition~\ref{prop:hcstep2positive} implies
Proposition~\ref{prop:hcstep2reducedtail}.

The reduction \eqref{eq:hcstep2tailreduction} and comparison
\eqref{eq:hcstep2cutoffreturn} then give
\eqref{eq:hcgeneralconvolution} for every Schwartz function
\(f\) with
\(\operatorname{supp}\widehat f\subset\Gamma_j^{d+1}(R)\).
Lemma~\ref{lem:hcstep2averagecomparison} gives
\eqref{eq:hcstablefunction}, so both estimates are established
at index \(j\). Proceeding successively through \(j=0,\ldots,m\)
completes Step~2.

\begin{remark}
Section~\ref{sec:hcgeneral} proves
Proposition~\ref{prop:hcstep2positive} by induction on \(\sigma s\),
with \(R,j,\ell\) fixed.

Section~\ref{sec:hcbase} proves
Proposition~\ref{prop:hcbase}, the base case \(m=d-2\)
of the induction on \(m\). The proof follows Steps~1 and~2
in the proof of Proposition~\ref{prop:hcstablefunction}, with
\(p_{d-m}=4\).
\end{remark}

We close this section with the deferred proofs of
Lemmas~\ref{lem:hcstep2cutoffcomparison} and~\ref{lem:hcstep2averagecomparison}.

\begin{proof}[Proof of Lemma~\ref{lem:hcstep2cutoffcomparison}]
For \(\xi=\sum_{i=0}^d\lambda_i\phi_d^{(i)}(t_\theta)
\in\operatorname{supp}\widehat{\widetilde f_\theta}\),
the bounds in \eqref{eq:hcstep2support} give
\[
 |\lambda_i|\lesssim
 \begin{cases}
 R^{-j/d}\ell^{i-j},&0\le i\le j,\\
 R^{-j/d}\ell^{-j},&j<i\le d.
 \end{cases}
\]
For \(i\le k\le d\), the coefficient of \(x\cdot\phi_d^{(k)}(t_\theta)\) in the Taylor expansion of
\(x\cdot\phi_d^{(i)}(t_\tau)\) is \((t_\tau-t_\theta)^{k-i}/(k-i)!\).
Since \(\theta\subset\tau\) and \(|\tau|=\ell\), we have \(|t_\tau-t_\theta|\lesssim\ell\), so
\[
\begin{aligned}
 R^{j/d}\ell^{j-k}|t_\theta-t_\tau|^{k-i}
 &\lesssim R^{j/d}\ell^{j-i},
 &&0\le i\le k\le j,\\
 R^{j/d}\ell^j|t_\theta-t_\tau|^{k-i}
 &\lesssim R^{j/d}\ell^j,
 &&0\le i\le k,\quad j<k\le d.
\end{aligned}
\]
Thus, for \(\tau\in\mathcal J(\ell)\), let \(D_\tau\) consist of all
\(x\in\R^{d+1}\) satisfying
\[
 |x\cdot\phi_d^{(i)}(t_\tau)|\le
 \begin{cases}
 R^{j/d}\ell^{j-i},&0\le i\le j,\\
 R^{j/d}\ell^j,&j<i\le d.
 \end{cases}
\]
The locally constant property then gives
\(|\widetilde f_\theta|^2\lesssim |f_\theta|^2*w_{D_\tau}\)
for \(\theta\subset\tau\).

Fix the scales \(q_1\le\cdots\le q_{d-m-2}\le u\) and
\(\rho\in\mathcal J(u)\) on the RHS of
\eqref{eq:hcgeneralconvolution}.
Its weighted average over \(B\) contains sums over
\(\theta\subset\tau_1\subset\cdots\subset\tau_{d-m-2}
\subset\rho\), where \(|\tau_i|=q_i\), with convolution
weights \(w_{B_{\tau_i}}\) defined in
\eqref{eq:hcnestedwidths}.
We divide the argument into three cases according to the size
of \(\ell\) relative to \(q_{d-m-2}\) and \(u\).

\noindent\emph{Case 1: \(\ell\le q_{d-m-2}\).}
Let \(l\) be the smallest index such that \(q_l\ge\ell\). If \(l>1\), then each
\(\tau_1\) lies in one \(\tau\in\mathcal J(\ell)\),
and the preceding pointwise bound and H\"older give
\begin{equation}
    \biggl[
 \Bigl(\sum_{\theta\subset\tau_1}|\widetilde f_\theta|^2\Bigr)
 *w_{B_{\tau_1}}
 \biggr]^{p_{d-m}/p_{d-m-1}}
 \lesssim
 \biggl[
 \Bigl(\sum_{\theta\subset\tau_1}|f_\theta|^2\Bigr)
 *w_{B_{\tau_1}}
 \biggr]^{p_{d-m}/p_{d-m-1}}*w_{D_\tau}.
\end{equation}
The same application of H\"older is valid at each subsequent
index \(i<l\), since \(\tau_i\) is still contained in one
\(\ell\) cap.

At the \(l\)th convolution, we show
\(w_{D_\tau}*w_{B_{\tau_l}}\lesssim w_{B_{\tau_l}}\)
for \(\tau\subset\tau_l\).
Once this is proved, for \(\theta\subset\tau\) we have
\(
 |\widetilde f_\theta|^2*w_{B_{\tau_l}}
 \lesssim |f_\theta|^2*(w_{D_\tau}*w_{B_{\tau_l}})
 \lesssim |f_\theta|^2*w_{B_{\tau_l}}.
\)
Together with the preceding applications of H\"older,
this bounds the contribution of these fixed scales to the
RHS of \eqref{eq:hcgeneralconvolution} for \(\widetilde f\)
by the corresponding contribution for \(f\).
To prove the kernel bound, we verify that
\(D_\tau\subset C B_{\tau_l}\).

For \(x\in D_\tau\) and \(0\le i\le d\),
\[
 |x\cdot\phi_d^{(i)}(t_\tau)|\le
\begin{cases}
R^{j/d}\ell^{j-i},&0\le i\le j,\\
R^{j/d}\ell^j,&j<i\le d
\end{cases}
\lesssim
\frac{R^{(m+2)/d}}{q_1\cdots q_{l-1}}
\begin{cases}
q_l^{m+l+1-i},&0\le i\le m+l+1,\\
u^{m+l+1-i},&m+l+1<i\le d.
\end{cases}
\]
where we used \(q_1\cdots q_{l-1}\le q_l^{l-1}\),
\(R^{-1/d}\le\ell\le q_l\le u\le1\), and \(j\le m\).
To pass from \(t_\tau\) to \(t_{\tau_l}\), use the identity
\[
 x\cdot\phi_d^{(i)}(b)=\sum_{k=i}^d
 \frac{(b-a)^{k-i}}{(k-i)!}\,x\cdot\phi_d^{(k)}(a).
\]
Since \(|t_{\tau_l}-t_\tau|\lesssim q_l\le u\), the same upper bound holds for
\(|x\cdot\phi_d^{(i)}(t_{\tau_l})|\). By \eqref{eq:hcnestedwidths},
\(D_\tau\subset C B_{\tau_l}\), and hence
\(w_{D_\tau}*w_{B_{\tau_l}}\lesssim w_{B_{\tau_l}}\).

\noindent\emph{Case 2: \(q_{d-m-2}<\ell\le u\).}
After applying H\"older at every index in
\eqref{eq:hcgeneralconvolution}, we absorb \(w_{D_\tau}\)
into the weight \(W_B\) of the outer weighted average over \(B\).
To do this, we show \(w_{D_\tau}*W_B\lesssim W_B\),
so that \eqref{eq:weightedYoung} bounds the weighted average over \(B\)
for \(\widetilde f\) by that for \(f\).

For \(x\in D_\tau\) and \(0\le i\le d\),
\[
 |x\cdot\phi_d^{(i)}(t_\tau)|\le
\begin{cases}
R^{j/d}\ell^{j-i},&0\le i\le j,\\
R^{j/d}\ell^j,&j<i\le d
\end{cases}
\lesssim
\frac{R^{(m+2)/d}u^{d-i}}{q_1\cdots q_{d-m-2}}.
\]
where we used \(q_1\cdots q_{d-m-2}\le u^{d-m-2}\),
\(R^{-1/d}\le\ell\le u\le1\), and \(j\le m\).
Since \(|t_\rho-t_\tau|\lesssim u\), Taylor's formula gives the same upper bound for
\(|x\cdot\phi_d^{(i)}(t_\rho)|\). Thus \(D_\tau\subset C(B-c_B)\), where
\(c_B\) is the center of \(B\). Hence \(w_{D_\tau}*W_B\lesssim W_B\), and
\eqref{eq:weightedYoung} bounds the weighted average over \(B\) in
\eqref{eq:hcgeneralconvolution} for \(\widetilde f\)
by that for \(f\).

\noindent\emph{Case 3: \(u<\ell\).}
Then \(\rho\) is contained in one \(\tau\in\mathcal J(\ell)\).
In this case, \(w_{D_\tau}*W_B\lesssim W_B\) need not hold,
so we use \eqref{eq:commonreturn} to remove \(w_{D_\tau}\)
from the sum over \(B\) in \eqref{eq:hcgeneralconvolution}.

To apply \eqref{eq:commonreturn}, let \(F\) be the integrand
inside the weighted average over \(B\) in \eqref{eq:hcgeneralconvolution}
for \(f\), and let \(\widetilde F\) be the same expression
with \(f_\theta\) replaced by \(\widetilde f_\theta\).
Since every \(\theta\subset\rho\) lies in this single
\(\tau\), applying H\"older as in the first case gives
\(\widetilde F\lesssim F*w_{D_\tau}\).
With \(t=[\alpha/(R^\epsilon N(u))]^{p_{d-m}/2}\),
where \(N(u)\) is computed from \(\{\widetilde f_\theta\}\),
we therefore obtain
\[
\begin{aligned}
 \sum_{\|\widetilde F\|_{L_\#^1(B)}\gtrsim t}
 |B|\,\|\widetilde F\|_{L_\#^1(B)}^2
 &\lesssim
 \sum_{\|F*w_{D_\tau}\|_{L_\#^1(B)}\gtrsim t}
 |B|\,\|F*w_{D_\tau}\|_{L_\#^1(B)}^2\lesssim
 \sum_{\|F\|_{L_\#^1(B)}\gtrsim t}
 |B|\,\|F\|_{L_\#^1(B)}^2,
\end{aligned}
\]
where all sums use the \(B\) tiling in
\eqref{eq:hcgeneralconvolution}, and the last inequality
follows from \eqref{eq:commonreturn}.
The final sum is bounded by the sum over \(B\) on the RHS of
\eqref{eq:hcgeneralconvolution} for \(f\).

The counts \(N(u)\) computed from \(\{\widetilde f_\theta\}\)
are no larger than those from \(\{f_\theta\}\), so using
the latter only enlarges the sums obtained in all three
cases. This proves the lemma.
\end{proof}

\begin{proof}[Proof of Lemma~\ref{lem:hcstep2averagecomparison}]
Starting with the innermost convolution in
\eqref{eq:hcgeneralconvolution}, we move each power inside
its convolution by H\"older and absorb the kernel into
the next convolution. At the last layer, we absorb the
kernel into the outer weight \(W_B\).

To absorb \(w_{B_{\tau_l}}\) into the next convolution,
we show that \(B_{\tau_l}\subset C B_{\tau_{l+1}}\).
For \(x\in B_{\tau_l}\) and \(0\le i\le d\), \eqref{eq:hcnestedwidths} gives
\[
 |x\cdot\phi_d^{(i)}(t_{\tau_l})|\le
 \frac{R^{(m+2)/d}}{q_1\cdots q_{l-1}}
 \begin{cases}
 q_l^{m+l+1-i},&0\le i\le m+l+1,\\
 u^{m+l+1-i},&m+l+1<i\le d
 \end{cases}
 \lesssim
 \frac{R^{(m+2)/d}}{q_1\cdots q_l}
 \begin{cases}
 q_{l+1}^{m+l+2-i},&0\le i\le m+l+2,\\
 u^{m+l+2-i},&m+l+2<i\le d,
 \end{cases}
\]
where we used \(q_l\le q_{l+1}\le u\).
Since \(|t_{\tau_{l+1}}-t_{\tau_l}|\lesssim q_{l+1}\le u\), Taylor's formula gives the same
upper bound for \(|x\cdot\phi_d^{(i)}(t_{\tau_{l+1}})|\). Hence
\(B_{\tau_l}\subset C B_{\tau_{l+1}}\).

At the last convolution, we compare \(B_{\tau_{d-m-2}}\)
with the outer box \(B\).
For \(x\in B_{\tau_{d-m-2}}\) and \(0\le i\le d\), \eqref{eq:hcnestedwidths} gives
\[
 |x\cdot\phi_d^{(i)}(t_{\tau_{d-m-2}})|\le
 \frac{R^{(m+2)/d}}{q_1\cdots q_{d-m-3}}
 \begin{cases}
 q_{d-m-2}^{d-1-i},&0\le i\le d-1,\\
 u^{-1},&i=d
 \end{cases}
 \lesssim
 \frac{R^{(m+2)/d}u^{d-i}}{q_1\cdots q_{d-m-2}},
\]
where we used \(q_{d-m-2}\le u\). Since \(|t_{\tau_{d-m-2}}-t_\rho|\lesssim u\),
Taylor's formula gives the same upper bound for \(|x\cdot\phi_d^{(i)}(t_\rho)|\).
Thus \(B_{\tau_{d-m-2}}\subset C(B-c_B)\), where \(c_B\) is the center of \(B\).
These comparisons give
\(
 w_{B_{\tau_l}}*w_{B_{\tau_{l+1}}}
 \lesssim w_{B_{\tau_{l+1}}},~~
 w_{B_{\tau_{d-m-2}}}*W_B\lesssim W_B.
\)

For example, when \(m=d-4\), H\"older and the first kernel
bound give
\[
 \biggl(\sum_{\tau_1\subset\tau_2}
 \Bigl[\bigl(\sum_{\theta\subset\tau_1}|f_\theta|^2\bigr)
 *w_{B_{\tau_1}}\Bigr]^{11/7}\biggr)*w_{B_{\tau_2}}
 \lesssim
 \Bigl(\sum_{\theta\subset\tau_2}|f_\theta|^2\Bigr)^{11/7}
 *w_{B_{\tau_2}}.
\]
Applying \eqref{eq:weightedYoung} with exponent \(7/4\)
and summing over \(\tau_2\subset\rho\) gives
\[
\begin{aligned}
 &\fint_B\sum_{\tau_2\subset\rho}
 \biggl[\biggl(\sum_{\tau_1\subset\tau_2}
 \Bigl[\bigl(\sum_{\theta\subset\tau_1}|f_\theta|^2\bigr)
                 *w_{B_{\tau_1}}\Bigr]^{11/7}\biggr)
                 *w_{B_{\tau_2}}\biggr]^{7/4}\lesssim\fint_B
        \bigl(\sum_{\theta\subset\rho}|f_\theta|^2\bigr)^{11/4}.
\end{aligned}
\]
The powers in \eqref{eq:hcgeneralconvolution} have product
\(p_{d-m}/4\). Repeating the same calculation bounds its
weighted average over \(B\) by
\(\fint_B(\sum_{\theta\subset\rho}|f_\theta|^2)^{p_{d-m}/4}\).

For fixed \(q_1,\ldots,q_{d-m-2},u,\rho\), put
\(F=(\sum_{\theta\subset\rho}|f_\theta|^2)^{p_{d-m}/4}\)
and \(t=[\alpha/(R^\e N(u))]^{p_{d-m}/2}\).
The comparison of the weighted averages over \(B\) shows that every
\(B\) counted in \eqref{eq:hcgeneralconvolution}
satisfies \(\fint_BF\gtrsim t\).
Since \(u^{d-m-2}/(q_1\cdots q_{d-m-2})\ge1\),
Lemma~\ref{lem:envelopechange} bounds this sum over \(B\) by
the sum over \(U\parallel(uR^{1/d})^mU_{\rho,R}\), giving
\begin{equation}\label{eq:hcreturnedaverage}
 \sum_{\substack{U\parallel(uR^{1/d})^mU_{\rho,R}\\
 U\in\mathcal G_\rho(\frac{\alpha}{R^\e})}}
 |U|\biggl[\fint_U
        \bigl(\sum_{\theta\subset\rho}|f_\theta|^2\bigr)^{p_{d-m}/4}
        \biggr]^2.
\end{equation}
Summing
\eqref{eq:hcreturnedaverage} over \(u,\rho\) gives the
right hand side of \eqref{eq:hcstablefunction}, with
the original function \(f\). Summing over the dyadic scales \(q_1,\ldots,q_{d-m-2}\) proves the lemma.
\end{proof}

\section{\texorpdfstring{\except{toc}{From curves to higher dimensional cones over moment curves: }}{From curves to higher dimensional cones over moment curves: }Proof of Proposition~\ref{prop:hcstablefunction}, Step 1}
\label{sec:hcpositive}\label{sec:hcstable}

Assume the hypotheses of Proposition~\ref{prop:hcstablefunction}
for \(0\le m\le d-3\). We first prove Lemma~\ref{lem:hcproducer}
and then Lemma~\ref{lem:hcoverviewrefinement}. These give the
four branches as in Section~\ref{subsec:hcoverviewiteration}.
We then use those branches to prove
Proposition~\ref{prop:hcstable} together with
\eqref{eq:hcrestrictedconvolution}, completing Step~1 in the proof of
Proposition~\ref{prop:hcstablefunction}.
The case \(m=d-2\) is treated in Section~\ref{sec:hcbase}.

\subsection{Proof of Lemma~\ref{lem:hcproducer}}
\label{subsec:hcpositiveaverages}

The argument for the high frequency part applies the estimate at exponent
\(p_{d-m-1}\) to a sum of quadratic expressions. We first
record the resulting estimate at exponent \(p_{d-m}/2\),
keeping the convolutions needed for the induction.

\begin{lemma}\label{lem:hcpackettail}
Assume the hypotheses of Proposition~\ref{prop:hcstablefunction}
for a fixed \(0\le m\le d-3\).
Let
\(
 h=\sum_{\theta\in\Xi_{m+1}^{d+1}(r)}h_\theta,
\)
where each \(h_\theta\) is a Schwartz function with
\(\supp\widehat h_\theta\subset\theta\). Then
\[
\begin{aligned}
 &\int_{|h|>\alpha^2}|h|^{p_{d-m}/2}\\
 &\lesssim_\e r^\e
 \sum_{r^{-1/d}\le q_1\le\cdots\le q_{d-m-3}\le u\le1}
 \sum_{\rho\in\mathcal J(u)}
 \sum_{\substack{
 B\parallel
 \frac{u^{d-m-3}}{q_1\cdots q_{d-m-3}}
 (ur^{1/d})^{m+1}U_{\rho,r}\\
 B\in\mathcal G_\rho(\frac{\alpha}{r^\e})}}
 |B|\\
 &\quad\times
 \biggl[
 \fint_B
 \sum_{\tau_{d-m-3}\subset\rho}
 \biggl[
 \biggl(
 \cdots
 \sum_{\tau_1\subset\tau_2}
 \biggl[
 \Bigl(
 \sum_{\theta\subset\tau_1}
 |h_\theta|^{p_{d-m}/p_{d-m-1}}
 \Bigr)*w_{B_{\tau_1}}
 \biggr]^{p_{d-m-1}/p_{d-m-2}}
 \cdots
 \biggr)*w_{B_{\tau_{d-m-3}}}
 \biggr]^{7/4}
 \biggr]^2\\
 &\quad+
 r^\e\sum_\theta
 \sum_{\substack{
 Q\parallel\theta^*\\
 Q\in\mathcal G_\theta(\frac{\alpha}{r^\e})}}
 |Q|\biggl[\fint_Q|h_\theta|^{p_{d-m}/4}\biggr]^2.
\end{aligned}
\]
In both sums, the condition on \(B\) or \(Q\) means that the average
inside the outer square is at least
\(
 \bigl[\frac{\alpha}{r^\e N(u)}\bigr]^{p_{d-m}/2},
\)
with \(u=r^{-1/d}\) in the second sum.
\end{lemma}
Since the proof of the lemma is rather standard,
we postpone the proof to Appendix~\ref{app:hcpackets}.

\medskip
\medskip

In the following lemma, we take \(q_1=\cdots=q_{d-m-2}=u\)
and \(\tau_1=\cdots=\tau_{d-m-2}=\rho\) in
\eqref{eq:hcnestedwidths}, with parameter \(r\).
Then \(B_{\tau_l}=(ur^{1/d})^mU_{\rho,r}\) for
\(1\le l\le d-m-2\). The convolutions and powers follow
the order in \eqref{eq:hcgeneralconvolution}.

\begin{lemma}\label{lem:hcendpointaverages}
Let \(0\le m\le d-3\), and let \(g\) be a Schwartz function
with \(\supp\widehat g\subset\Gamma_m^{d+1}(r)\).
For \(r^{-1/d}\le u\le1\), \(\rho\in\mathcal J(u)\), and
\(B\parallel B_0=(ur^{1/d})^mU_{\rho,r}\), we have
\[
\begin{aligned}
 &\fint_B
 \biggl[
 \Bigl(\sum_{\theta\subset\rho}|g_\theta|^2\Bigr)*w_{B_0}
 \biggr]^{p_{d-m}/4}\\
 &\quad\lesssim
 \fint_B
 \biggl[
 \biggl(
 \cdots
 \biggl[
 \biggl[
 \Bigl(\sum_{\theta\subset\rho}|g_\theta|^2\Bigr)*w_{B_0}
 \biggr]^{p_{d-m}/p_{d-m-1}}
 *w_{B_0}
 \biggr]^{p_{d-m-1}/p_{d-m-2}}
 \cdots
 \biggr)
 *w_{B_0}
 \biggr]^{7/4}.
\end{aligned}
\]
When \(m=d-3\), the two averages agree.
\end{lemma}

\begin{proof}
Since \(w_{B_0}\lesssim w_{B_0}*w_{B_0}\),
H\"older's inequality gives, for \(\beta\ge1\),
\begin{equation}\label{eq:hcendpointpadding}
 \biggl[
 \Bigl(\sum_{\theta\subset\rho}|g_\theta|^2\Bigr)*w_{B_0}
 \biggr]^\beta
 \lesssim_\beta
 \biggl[
 \Bigl(\sum_{\theta\subset\rho}|g_\theta|^2\Bigr)*w_{B_0}
 \biggr]^\beta*w_{B_0}.
\end{equation}
Apply \eqref{eq:hcendpointpadding} successively for
\(l=2,\ldots,d-m-2\), with
\(\beta=p_{d-m}/p_{d-m-l+1}\), and take the power
\(p_{d-m-l+1}/p_{d-m-l}\) after each convolution.
The successive powers have product \(p_{d-m}/4\),
so averaging over \(B\) proves the lemma.
\end{proof}

\begin{proof}[Proof of Lemma~\ref{lem:hcproducer} for \(m<d-2\)]

Let \(g\) be a Schwartz function with
\(\supp\widehat g\subset\Gamma_m^{d+1}(r)\). 
We prove the stronger estimate
\begin{equation}\label{eq:hcpositiveconvolution}
 \int_{\sum_\theta|g_\theta|^2\gtrsim\alpha^2}
       \bigl(\sum_\theta|g_\theta|^2\bigr)^{p_{d-m}/2}
 \lesssim_\e
 \text{RHS of \eqref{eq:hcgeneralconvolution}
 with \(R=r\), \(f=g\)}.
\end{equation}
Lemma~\ref{lem:hcstep2averagecomparison} will then give
Lemma~\ref{lem:hcproducer}. The estimate
\eqref{eq:hcpositiveconvolution} is also the input at the fine
endpoint in the proof of Proposition~\ref{prop:hcstable} below.

We use the cutoffs \(\eta^{m+2}_{\theta,\le\sigma}\) and
\(\eta^{m+2}_{\theta,\sigma}\) from
Section~\ref{subsec:hcoverviewpositive}.
For
\(\xi=\sum_{i=0}^d\lambda_i\phi_d^{(i)}(t_\theta)
\in\operatorname{supp}\eta^{m+2}_{\theta,\sigma}\),
we have
\begin{equation}\label{eq:hcpositiveshell}
 |\lambda_i|\lesssim r^{-i/d}\sigma^{m+2-i}
 \quad(0\le i\le m+1),\qquad
 \max_{0\le i\le m+1}
 r^{i/d}\sigma^{i-m-2}|\lambda_i|\sim1.
\end{equation}
Since \(\eta^{m+2}_{\theta,\le1}=1\) on
\(\supp\widehat{|g_\theta|^2}\), dyadic decomposition gives
\begin{equation}\label{eq:hcpositiveparts}
\begin{aligned}
 \int_{\sum_\theta|g_\theta|^2\gtrsim\alpha^2}
 \bigl(\sum_\theta|g_\theta|^2\bigr)^{p_{d-m}/2}
 &\lessapprox
 \int_{|\sum_\theta|g_\theta|^2*
             \check\eta^{m+2}_{\theta,\le r^{-1/d}}|
                  \gtrapprox\alpha^2}
 \Bigl|\sum_\theta|g_\theta|^2*
             \check\eta^{m+2}_{\theta,\le r^{-1/d}}\Bigr|^{p_{d-m}/2}\\
 &\quad+\sum_{r^{-1/d}<\sigma\le1}
 \int_{|\sum_\theta|g_\theta|^2*
             \check\eta^{m+2}_{\theta,\sigma}|\gtrapprox\alpha^2}
 \Bigl|\sum_\theta|g_\theta|^2*
             \check\eta^{m+2}_{\theta,\sigma}\Bigr|^{p_{d-m}/2}.
\end{aligned}
\end{equation}
We will bound the first integral on the RHS of
\eqref{eq:hcpositiveparts} by the terms with
\(q_1=\cdots=q_{d-m-2}=u=1\) in
\eqref{eq:hcgeneralconvolution}.
For each \(\sigma>r^{-1/d}\), we will group the
\(\theta\) caps and apply Lemma~\ref{lem:hcpackettail}
after dilation to bound the corresponding integral in
the second term of \eqref{eq:hcpositiveparts}.

For the first integral on the RHS of
\eqref{eq:hcpositiveparts}, take
\(B_0=r^{m/d}U_{[0,1],r}\), as in
Lemma~\ref{lem:hcendpointaverages} with \(u=1\) and
\(\rho=[0,1]\).
This box has width \(r^{(m+2)/d}\) in every dual Taylor
direction.

At \(\sigma=r^{-1/d}\), the first \(m+2\) factors of
\(\eta^{m+2}_{\theta,\le r^{-1/d}}\) cut off the frequency
coordinates at widths

\(r^{-i/d}(r^{-1/d})^{m+2-i} = r^{-(m+2)/d}\),
\(0\le i\le m+1\).
The remaining factors equal \(1\) on
\(\operatorname{supp}\widehat{|g_\theta|^2}\),
and the locally constant property allows averaging over width
\(r^{(m+2)/d}\) in directions \(m+2,\ldots,d\),
since \(r^{(m+2)/d}\le r^{i/d}\).

Since \(\theta\subset[0,1]\), Taylor expansion shows that
boxes with width \(r^{(m+2)/d}\) in every dual Taylor
direction at \(t_\theta\) are comparable to \(B_0\).
The locally constant property then gives
\[
 \Bigl|\sum_\theta |g_\theta|^2*
 \check\eta^{m+2}_{\theta,\le r^{-1/d}}\Bigr|
 \lesssim
 \Bigl(\sum_\theta|g_\theta|^2\Bigr)*w_{B_0}.
\]
and bound the first
integral on the RHS of \eqref{eq:hcpositiveparts} by
\begin{equation}\label{eq:hclowaverage}
 \sum_{\substack{U\parallel r^{m/d}U_{[0,1],r}\\
 U\in\mathcal G_{[0,1]}(\frac{\alpha}{r^\e})}}
 |U|\biggl[\fint_U
 \Bigl(\bigl(\sum_\theta|g_\theta|^2\bigr)*
          w_{B_0}\Bigr)^{p_{d-m}/4}\biggr]^2.
\end{equation}
Applying Lemma~\ref{lem:hcendpointaverages} to \(g\), with
\(u=1\) and \(\rho=[0,1]\), bounds \eqref{eq:hclowaverage}
by the terms of \eqref{eq:hcgeneralconvolution} with
\(q_1=\cdots=q_{d-m-2}=u=1\), including the amplitude condition.

We next estimate the integral indexed by \(\sigma>r^{-1/d}\)
in the second term of \eqref{eq:hcpositiveparts}.
Put \(q_1=r^{-1/d}\sigma^{-1}\), and group the
\(\theta\) caps into \(\tau_1\in\mathcal J(q_1)\).
With \(B_{\tau_1}\) as in \eqref{eq:hcnestedwidths},
where \(R=r\) and \(l=1\), we have
\begin{equation}\label{eq:hchighpieces}
 h_{\tau_1}
 =\sum_{\theta\subset\tau_1}|g_\theta|^2*
 \check\eta^{m+2}_{\theta,\sigma},
 \qquad
 |h_{\tau_1}|\lesssim
 \Bigl(\sum_{\theta\subset\tau_1}|g_\theta|^2\Bigr)
 *w_{B_{\tau_1}},
\end{equation}
uniformly for \(q_1\le u\le1\).

Since \(\theta\subset\tau_1\), we have \(|t_\theta-t_{\tau_1}|\lesssim q_1\le u\).
By Taylor expansion, \(B_{\tau_1}\) agrees up to a fixed dilation with the box defined by
\[
 |x\cdot\phi_d^{(i)}(t_\theta)|\le r^{(m+2)/d}
 \begin{cases}
 q_1^{m+2-i},&0\le i\le m+2,\\
 u^{m+2-i},&m+2<i\le d.
 \end{cases}
\]
In directions \(0\le i\le m+1\), the convolution kernel
\(\check\eta_{\theta,\sigma}^{m+2}\) decays at spatial scales
\(r^{i/d}\sigma^{i-m-2}
=r^{(m+2)/d}q_1^{m+2-i}\), which are the corresponding
widths of \(B_{\tau_1}\).
In directions \(m+2\le i\le d\), the cutoff factors equal one
on \(\operatorname{supp}\widehat{|g_\theta|^2}\), and the
widths of \(B_{\tau_1}\) satisfy
\(r^{(m+2)/d}u^{m+2-i}\le r^{i/d}\).
The locally constant property of \(g_\theta\) therefore gives
\(\bigl||g_\theta|^2*\check\eta_{\theta,\sigma}^{m+2}\bigr|
\lesssim |g_\theta|^2*w_{B_{\tau_1}}\) for each
\(\theta\subset\tau_1\). Summing over \(\theta\) proves the
pointwise bound in \eqref{eq:hchighpieces}.

Taylor expansion as in \eqref{eq:highfrequencyhigherTaylororder} also
gives \(\sigma^{-(m+2)}\operatorname{supp}\widehat{h_{\tau_1}}
\subset\omega_{\tau_1}\) for a cap
\(\omega_{\tau_1}\in\Xi_{m+1}^{d+1}(r\sigma^d)\).
We therefore apply Lemma~\ref{lem:hcpackettail} to
\(h(x)=\sum_{\tau_1}h_{\tau_1}(\sigma^{-(m+2)}x)\),
with parameter \(r\sigma^d\), pieces
\(h_{\omega_{\tau_1}}(x)
=h_{\tau_1}(\sigma^{-(m+2)}x)\),
and amplitude \(\alpha/r^\epsilon\).
Returning to the original coordinates by
\(y=\sigma^{-(m+2)}x\) multiplies all spatial widths
by \(\sigma^{-(m+2)}\).

In the first sum of Lemma~\ref{lem:hcpackettail} at parameter
\(r\sigma^d\), denote the intermediate caps by \(\tau_l'\),
the convolution weights by \(w_{B_{\tau_l'}}\), and the
averaging boxes by \(B'\).
Since \eqref{eq:hchighpieces} gives the first convolution
with \(w_{B_{\tau_1}}\), returning by
\(y=\sigma^{-(m+2)}x\), we write
\(\tau_l=\sigma^{m+2}\tau_{l-1}'\) and
\(B=\sigma^{-(m+2)}B'\), where
\(2\le l\le d-m-2\) and \(|\tau_{l-1}'|=q_l\).
The \(i\)th widths of
\(\sigma^{-(m+2)}B_{\tau_{l-1}'}\) satisfy
\[
\begin{aligned}
 \sigma^{-(m+2)}
 \frac{(r\sigma^d)^{(m+3)/d}q_l^{m+l+1-i}}
      {q_2\cdots q_{l-1}}
 &=
 \frac{r^{(m+2)/d}q_l^{m+l+1-i}}
      {q_1\cdots q_{l-1}},
 &&0\le i\le m+l+1,\\
 \sigma^{-(m+2)}
 \frac{(r\sigma^d)^{(m+3)/d}u^{m+l+1-i}}
      {q_2\cdots q_{l-1}}
 &=
 \frac{r^{(m+2)/d}u^{m+l+1-i}}
      {q_1\cdots q_{l-1}},
 &&m+l+1<i\le d.
\end{aligned}
\]
These are exactly the widths of \(B_{\tau_l}\) in
\eqref{eq:hcnestedwidths}, with \(R=r\).

The \(i\)th widths of the averaging boxes
\(B=\sigma^{-(m+2)}B'\) satisfy
\[
 \sigma^{-(m+2)}
 \frac{(r\sigma^d)^{(m+3)/d}u^{d-i}}
      {q_2\cdots q_{d-m-2}}
 =
 \frac{r^{(m+2)/d}u^{d-i}}
      {q_1\cdots q_{d-m-2}},
 \qquad 0\le i\le d.
\]
Hence
\(B\parallel
\frac{u^{d-m-2}}{q_1\cdots q_{d-m-2}}
(ur^{1/d})^mU_{\rho,r}\),
exactly as in \eqref{eq:hcgeneralconvolution} with \(R=r\).

Using \eqref{eq:hchighpieces} with this \(u\), the innermost
expressions satisfy
\[
 |h_{\tau_1}|^{p_{d-m}/p_{d-m-1}}
 \lesssim
 \biggl[
 \Bigl(\sum_{\theta\subset\tau_1}|g_\theta|^2\Bigr)
 *w_{B_{\tau_1}}
 \biggr]^{p_{d-m}/p_{d-m-1}}.
\]
Substituting this bound into the first sum on the RHS of
Lemma~\ref{lem:hcpackettail}, after returning to the original
coordinates, gives the convolution weights \(w_{B_{\tau_l}}\)
and successive powers
\(p_{d-m}/p_{d-m-1},\ldots,p_3/p_2\)
appearing in \eqref{eq:hcgeneralconvolution}.
The counts computed from \(\{h_{\tau_1}\}\) are no larger
than those from \(\{g_\theta\}\), so using the latter
\(N(u)\) only lowers the selection thresholds.
Hence the first sum in Lemma~\ref{lem:hcpackettail}, after
returning to the original coordinates, is bounded by the
terms on the RHS of \eqref{eq:hcgeneralconvolution} with
\(R=r\), \(f=g\), and \(q_1=r^{-1/d}\sigma^{-1}\).

For the second sum in Lemma~\ref{lem:hcpackettail}, use
\eqref{eq:hchighpieces} with \(u=q_1\).
After the same change of coordinates, writing
\(B=\sigma^{-(m+2)}Q\), this sum is bounded by
\begin{equation}\label{eq:hcpacketaverage}
 \lessapprox
 \sum_{\tau_1\in\mathcal J(q_1)}
 \sum_{\substack{
 B\parallel\sigma^{-(m+2)}\omega_{\tau_1}^*\\
 B\in\mathcal G_{\tau_1}(\alpha/r^\epsilon)}}
 |B|\biggl[\fint_B
 \Bigl[
 \bigl(\sum_{\theta\subset\tau_1}|g_\theta|^2\bigr)
 *w_{B_{\tau_1}}
 \Bigr]^{p_{d-m}/4}\biggr]^2.
\end{equation}
Both \(B_{\tau_1}\) and the averaging boxes in
\eqref{eq:hcpacketaverage} have dimensions
\(\sigma^{-m-2}\times r^{1/d}\sigma^{-m-1}
\times\cdots\times r\sigma^{d-m-2}\).
Thus these boxes are parallel to
\((ur^{1/d})^mU_{\tau_1,r}\), where
\(u=q_1=r^{-1/d}\sigma^{-1}\).
Lemma~\ref{lem:hcendpointaverages}, applied to \(g\) with
this \(u\) and \(\rho=\tau_1\), bounds
\eqref{eq:hcpacketaverage} by the terms of
\eqref{eq:hcgeneralconvolution} with
\(q_1=\cdots=q_{d-m-2}=u=r^{-1/d}\sigma^{-1}\),
including the amplitude condition.

Summing these bounds over \(r^{-1/d}<\sigma\le1\) and
using the bound for \eqref{eq:hclowaverage} in
\eqref{eq:hcpositiveparts} proves
\eqref{eq:hcpositiveconvolution}.
Lemma~\ref{lem:hcstep2averagecomparison}, applied to \(g\)
with parameter \(r\), bounds the RHS of
\eqref{eq:hcpositiveconvolution} by the RHS of
Lemma~\ref{lem:hcproducer} and completes the proof.
\end{proof}

The following consequence of Lemma~\ref{lem:hcproducer}
will be used in both Steps~1 and~2.
We first introduce a common notation for the Fourier support
conditions in these two steps.
For \(R\ge1\), \(0\le j\le d-2\), and
\(R^{-1/d}\le\ell\le1\), define
\begin{equation}\label{eq:hccurvaturesupport}
\widetilde\Gamma_{j,\ell}^{d+1}(R)
:=
\bigcup_{I\in\mathcal J(R^{-1/d})}
\left\{
\sum_{i=0}^{d}\lambda_i\phi_d^{(i)}(t_I):
\begin{array}{ll}
|\lambda_i|\lesssim R^{-j/d}\ell^{i-j},
   &0\le i<j,\\
|\lambda_j|\sim R^{-j/d},&\\
|\lambda_i|\lesssim R^{-i/d},
   &j<i\le d
\end{array}
\right\}.
\end{equation}
The supports \(\Gamma_{j,1/K}^{d+1}(R)\) in Step~1 and
\(\Gamma_{j,1/K,\ell}^{d+1}(R)\) in Step~2 are contained in
\(\widetilde\Gamma_{j,1}^{d+1}(R)\) and
\(\widetilde\Gamma_{j,\ell}^{d+1}(R)\), respectively.

\begin{corollary}\label{cor:hccoarse}
Assume the hypotheses of
Proposition~\ref{prop:hcstablefunction}
for a fixed \(0\le m\le d-3\).
Let \(0\le j\le m\), \(R^{-1/d}\le\ell\le1\), and let
\(f\) be a Schwartz function with
\(\operatorname{supp}\widehat f
\subset\widetilde\Gamma_{j,\ell}^{d+1}(R)\).
For \(R^{-1/d}\le\sigma\ll\ell\), \(\sigma\le s\le1\),
\(\delta>0\), and \(\alpha>0\), we have
\[
\begin{aligned}
&\sum_{\rho\in\mathcal J(s)}
 \sum_{\substack{
 Q\parallel\rho^*\\
 Q\in\mathcal G_\rho(\alpha)}}
 |Q|
 \biggl\|
 \Bigl(
 \sum_{\substack{\tau\subset\rho\\|\tau|=\sigma}}
 |f_\tau|^2
 \Bigr)^{1/2}
 \biggr\|_{L_\#^{p_{d-m}/2}(Q)}^{p_{d-m}}
\\
&\lesssim_\delta
 \sigma^{-\delta}
 \sum_{\sigma\le u\le s}
 \sum_{\rho'\in\mathcal J(u)}
 \sum_{\substack{
 B\parallel
 \frac{(uR^{1/d})^j}{(\sigma R^{1/d})^2}U_{\rho',R}\\
 B\in\mathcal G_{\rho'}(\frac{\alpha}{\sigma^{-\delta}})}}
 |B|
 \biggl\|
 \Bigl(
 \sum_{\substack{\tau\subset\rho'\\|\tau|=\sigma}}
 |f_\tau|^2
 \Bigr)^{1/2}
 \biggr\|_{L_\#^{p_{d-m}/2}(B)}^{p_{d-m}}.
\end{aligned}
\]
\end{corollary}

\begin{proof}
We first prove the integral estimate in
\eqref{eq:hcoverviewfullparent} under the support hypothesis
of this corollary, that is,
\begin{equation}\label{eq:hccoarsefullparent}
\begin{aligned}
&\int_{\sum_{|\tau|=\sigma}|f_\tau|^2\gtrsim\alpha^2}
 \Bigl(\sum_{|\tau|=\sigma}|f_\tau|^2\Bigr)^{p_{d-m}/2}
\lesssim_\delta \sigma^{-\delta}
 \sum_{\sigma\le u\le1}
 \sum_{\rho\in\mathcal J(u)}
 \sum_{\substack{
 B\parallel\frac{(uR^{1/d})^j}{(\sigma R^{1/d})^2}U_{\rho,R}\\
 B\in\mathcal G_\rho(\alpha/\sigma^{-\delta})}}
 |B|
 \biggl\|
 \Bigl(
 \sum_{\substack{\tau\subset\rho\\|\tau|=\sigma}}
 |f_\tau|^2
 \Bigr)^{1/2}
 \biggr\|_{L_\#^{p_{d-m}/2}(B)}^{p_{d-m}}.
\end{aligned}
\end{equation}
For \(\tau\in\mathcal J(\sigma)\), write
\(\xi=\sum_{i=0}^d\lambda_i\phi_d^{(i)}(t_\tau)
\in\operatorname{supp}\widehat f_\tau\).
Since \(f_\tau=\sum_{\theta\subset\tau}f_\theta\) and
\(|t_\theta-t_\tau|\lesssim\sigma\), Taylor expansion and
\eqref{eq:hccurvaturesupport} give
\[
 |R^{j/d}\sigma^j\lambda_i|\lesssim\sigma^i
 \quad(0\le i\le d),\qquad
 |R^{j/d}\sigma^j\lambda_j|\sim\sigma^j.
\]
The \(j\)th coefficient differs from its value at
\(t_\theta\) by \(O(R^{-j/d}\sigma/\ell)\), so
\(\sigma\ll\ell\) gives the last comparison. Thus
\[
 \operatorname{supp}
 \widehat{f_\tau(R^{j/d}\sigma^j\,\cdot)}
 \subset\theta_\tau,
 \qquad
 \theta_\tau\in\Xi_m^{d+1}(\sigma^{-d}).
\]
We therefore apply Lemma~\ref{lem:hcproducer} to
\(g(x)=f(R^{j/d}\sigma^j x)\), with parameter
\(\sigma^{-d}\), pieces
\(f_\tau(R^{j/d}\sigma^j x)\), and amplitude \(\alpha\).

As in the derivation of \eqref{eq:hcoverviewfullparent},
returning to the original coordinates gives boxes parallel to
\((u/\sigma)^{m-j}
\frac{(uR^{1/d})^j}{(\sigma R^{1/d})^2}U_{\rho,R}\).
Since \((u/\sigma)^{m-j}\ge1\),
Lemma~\ref{lem:envelopechange}, applied with
\(F=(\sum_{\tau\subset\rho,\,|\tau|=\sigma}
|f_\tau|^2)^{p_{d-m}/4}\), replaces these boxes by
\(B\parallel
\frac{(uR^{1/d})^j}{(\sigma R^{1/d})^2}U_{\rho,R}\)
while preserving the amplitude condition.
This proves \eqref{eq:hccoarsefullparent}.

To prove the corollary, we first use the locally constant property when
\(\sigma\sim s\) and then apply \eqref{eq:hccoarsefullparent}
when \(\sigma\ll s\).
For each \(\rho\in\mathcal J(s)\) and \(\tau\subset\rho\) of
width \(\sigma\), the Fourier support of \(f_\tau\) lies in a box
of dimensions \(R^{-j/d}\sigma^{-j}\times
R^{-j/d}\sigma^{1-j}\times\cdots\times R^{-j/d}\sigma^{d-j}\).
As in \eqref{eq:hcoverviewcomparable}, the locally constant property bounds
the left hand side of the corollary by the terms with \(u=s\)
on its right hand side.

Now suppose that \(\sigma\ll s\). By \eqref{eq:selaverage} with
\(F=(\sum_{\tau\subset\rho,\,|\tau|=\sigma}|f_\tau|^2)^{p_{d-m}/4}\)
and \(t=[\alpha/N(s)]^{p_{d-m}/2}\), it suffices to bound
\[
 \sum_{\rho\in\mathcal J(s)}
 \int_{\sum_{\tau\subset\rho,\,|\tau|=\sigma}|f_\tau|^2
       \gtrsim[\alpha/N(s)]^2}
 \Bigl(\sum_{\substack{\tau\subset\rho\\|\tau|=\sigma}}
 |f_\tau|^2\Bigr)^{p_{d-m}/2}.
\]
For \(\rho=[t_\rho,t_\rho+s]\in\mathcal J(s)\) and each cap \(\theta\subset\rho\) at scale \(R^{-1/d}\),
the Fourier change of variables is
\[
 \xi=\sum_{i=0}^d\lambda_i\phi_d^{(i)}(t_\theta)
 \quad\Longrightarrow\quad L_{t_\rho,s}^{-1}\xi
 =\sum_{i=0}^d s^{-i}\lambda_i\phi_d^{(i)}\!\left(\frac{t_\theta-t_\rho}{s}\right).
\]
Hence
\(\operatorname{supp}\widehat{f_\rho(L_{t_\rho,s}^{-T}\,\cdot)}
\subset\widetilde\Gamma_{j,\min\{\ell/s,1\}}^{d+1}(Rs^d)\).
Since \(\sigma/s\ll\min\{\ell/s,1\}\), we apply
\eqref{eq:hccoarsefullparent} to this function with parameter
\(Rs^d\), pieces \(f_\tau(L_{t_\rho,s}^{-T}\,\cdot)\) of width
\(\sigma/s\), and amplitude \(\alpha/N(s)\).

In the resulting sum, denote the parent of width \(u/s\) by
\(\rho''\) and the averaging boxes by \(B'\).
Returning by \(x=L_{t_\rho,s}^{-T}y\), write
\(\rho'=L_{t_\rho,s}\rho''\subset\rho\), with \(|\rho'|=u\),
and \(B=L_{t_\rho,s}^{-T}B'\).
The \(i\)th width of \(B\) is
\(s^{-i}(Rs^d)^{j/d}(\sigma/s)^{-2}(u/s)^{j+2-i}
=R^{j/d}\sigma^{-2}u^{j+2-i}\) for \(0\le i\le d\).
Thus \(B\parallel
\frac{(uR^{1/d})^j}{(\sigma R^{1/d})^2}U_{\rho',R}\),
as in the corollary.

The active parent count at rescaled width \(u/s\) is \(N_\rho(u)\).
Using \(N(s)N_\rho(u)\lesssim N(u)\) from Lemma~\ref{lem:regular}
and \((\sigma/s)^{-\delta}\le\sigma^{-\delta}\), we obtain
\(\frac{\alpha}{(\sigma/s)^{-\delta}N(s)N_\rho(u)}
\gtrsim\frac{\alpha}{\sigma^{-\delta}N(u)}\)
in the selection condition.
Summing over \(\rho\) proves the corollary.
\end{proof}

\subsection{Proof of Lemma~\ref{lem:hcoverviewrefinement}}

We now prove the estimate in Lemma~\ref{lem:hcoverviewrefinement}.
We allow Fourier support in \(\widetilde\Gamma_{j,\ell}^{d+1}(R)\),
so that the same result applies to the coefficient shells in
Step~2. Taking \(\ell=1\) gives
Lemma~\ref{lem:hcoverviewrefinement}.
The proof uses the coefficient bounds and the following
square function estimate for Taylor cones in \(\mathbb R^{d+1}\).

\begin{proposition}[{\cite[Theorem~1, Proposition~1.1]{GMSF}}]
\label{prop:hctaylorsf}
Let \(0\le k\le d-2\) and \(2\le p\le p_{d-k}\).
For a Schwartz function \(g\) with
\(\supp\widehat g\subset\Gamma_k^{d+1}(r)\),
\[
 \int|g|^p\lesssim_\e r^\e
 \int\Bigl(\sum_{\theta\in\Xi_k^{d+1}(r)}
 |g_\theta|^2\Bigr)^{p/2}.
\]
\end{proposition}

\begin{lemma}\label{lem:hcrefinement}
Let \(0\le j\le m\le d-3\) and \(R^{-1/d}\le\ell\le1\). Let \(f\) be a Schwartz function with \(\operatorname{supp}\widehat f\subset\widetilde\Gamma_{j,\ell}^{d+1}(R)\). Let \(R^{-1/d}<\sigma\ll\ell\) and \(\alpha,\delta>0\). For \(\eta>0\) sufficiently small and \(\sigma\) sufficiently small in terms of \(\eta\), there are dyadic scales \(\sigma',s'\) with
\(R^{-1/d}\le\sigma'\le\max\{R^{-1/d},\sigma^{3/2}\}\) and \(\sigma\le s'\lesssim\sigma^{\eta^4}\) such that
\[
\begin{aligned}
 \sum_{\substack{
 B\parallel[0,R^{j/d}\sigma^{-2}]^{d+1}\\
 B\in\mathcal G_{[0,1]}(\alpha)}}
 |B|
 \biggl\|\Bigl(\sum_{|\tau|=\sigma}|f_\tau|^2\Bigr)^{1/2}
 \biggr\|_{L_\#^{p_{d-m}/2}(B)}^{p_{d-m}}
 &\lesssim_{\delta,\eta}
 \sigma^{-C\eta^4}
 \sum_{\substack{
 Q\parallel[0,1]^{d+1}\\
 Q\in\mathcal G_{[0,1]}
       (\frac{\alpha}{\sigma^{-C\eta^4}})}}
 |Q|
 \biggl\|\Bigl(\sum_{|\tau'|=\sigma'}|f_{\tau'}|^2\Bigr)^{1/2}
 \biggr\|_{L_\#^{p_{d-m}/2}(Q)}^{p_{d-m}}
\\
 &\quad+
 \sigma^{-\delta}
 \sum_{\rho\in\mathcal J(s')}
 \sum_{\substack{
 Q\parallel\rho^*\\
 Q\in\mathcal G_\rho
       (\frac{\alpha}{\sigma^{-\delta}})}}
 |Q|
 \biggl\|\Bigl(
 \sum_{\substack{\tau\subset\rho\\|\tau|=\sigma}}
 |f_\tau|^2\Bigr)^{1/2}
 \biggr\|_{L_\#^{p_{d-m}/2}(Q)}^{p_{d-m}}.
\end{aligned}
\]

\end{lemma}

\begin{proof}
For \(\tau\in\mathcal J(\sigma)\), write
\(\xi=\sum_{i=0}^d\lambda_i\phi_d^{(i)}(t_\tau)
\in\operatorname{supp}\widehat f_\tau\).
Taylor expansion using the coefficient bounds in
\eqref{eq:hccurvaturesupport}, together with
\(\sigma\ll\ell\), gives
\[
 |\lambda_i|\lesssim
 \begin{cases}
 R^{-j/d}\ell^{i-j},&0\le i<j,\\
 R^{-j/d}\sigma^{i-j},&j\le i\le d,
 \end{cases}
 \qquad |\lambda_j|\sim R^{-j/d}.
\]
For \(\sigma\le u\le1\), choose a nonnegative, \(L^1\) normalized kernel \(w_{\tau,u}\)
adapted to the box of all \(x\in\R^{d+1}\) satisfying
\[
 |x\cdot\phi_d^{(i)}(t_\tau)|\le
 \begin{cases}
 R^{j/d}\sigma^{-2}u^{d-i},&0\le i\le j+1,\\
 R^{j/d}\sigma^{-2}u^d,&j+2\le i\le d.
 \end{cases}
\]
Then \(w_{\tau,u}\lesssim w_{\tau,2u}\) for dyadic
\(\sigma\le u<1\), and
\[
 \operatorname{supp}\widehat w_{\tau,u}
 \subset
 \left\{
 \sum_{i=0}^d\lambda_i\phi_d^{(i)}(t_\tau):
 |\lambda_i|\lesssim
 \begin{cases}
 R^{-j/d}\sigma^2u^{i-d},&0\le i\le j+1,\\
 R^{-j/d}\sigma^2u^{-d},&j+2\le i\le d
 \end{cases}
 \right\}.
\]
At \(u=\sigma\), the kernel widths satisfy
\(R^{j/d}\sigma^{d-i-2}\le R^{j/d}\sigma^{j-i}\)
for \(0\le i\le j+1\), and
\(R^{j/d}\sigma^{d-2}\le R^{j/d}\sigma^{j-i}\)
for \(j+2\le i\le d\). The locally constant property therefore gives
\begin{equation}\label{eq:hcrefinementinitialaverage}
 \sum_{|\tau|=\sigma}|f_\tau|^2
 \lesssim
 \sum_{|\tau|=\sigma}|f_\tau|^2*w_{\tau,\sigma}.
\end{equation}

We divide the boxes \(B\in\mathcal G_{[0,1]}(\alpha)\),
\(B\parallel[0,R^{j/d}\sigma^{-2}]^{d+1}\), according to whether
\begin{equation}\label{eq:hcrefinementlowinput}
 \Bigl\|\sum_{|\tau|=\sigma}|f_\tau|^2
 \Bigr\|_{L_\#^{p_{d-m}/4}(B)}
 \le \sigma^{-\eta^4}
 \Bigl\|\sum_{|\tau|=\sigma}|f_\tau|^2*w_{\tau,u}
 \Bigr\|_{L_\#^{p_{d-m}/4}(B)}
\end{equation}
holds for some dyadic \(\sigma^{\eta^4}\le u\le1\).
If this holds, we will bound the RHS using finer pieces
\(f_{\tau'}\), with \(|\tau'|=\sigma'<\sigma\), to obtain
the first term in Lemma~\ref{lem:hcrefinement}.
For the remaining boxes, we will use a further high frequency
decomposition and Proposition~\ref{prop:hctaylorsf}
to obtain the second term, with parents of width
\(s'\lesssim\sigma^{\eta^4}\).
After summing over \(B\), dyadic pigeonholing will give
common scales \(\sigma'\) and \(s'\).

\medskip
\noindent\textit{Case 1.}
Suppose that \eqref{eq:hcrefinementlowinput} holds for some
dyadic \(u\in[\sigma^{\eta^4},1]\).
We first determine which products
\(f_\theta\overline{f_{\theta'}}\), \(\theta,\theta'\subset\tau\),
survive convolution with \(w_{\tau,u}\).
Take \(\xi\in\operatorname{supp}\widehat f_\theta\) and
\(\xi'\in\operatorname{supp}\widehat f_{\theta'}\) with
\(\xi-\xi'\in\operatorname{supp}\widehat w_{\tau,u}\).
Let \(t_\theta,t_{\theta'}\) be the left endpoints of the intervals indexing \(\theta,\theta'\), respectively, and write
\[
 \xi=\sum_{i=0}^d\lambda_i(t_\theta)\phi_d^{(i)}(t_\theta),\qquad
 \xi'=\sum_{i=0}^d\lambda_i(t_{\theta'})\phi_d^{(i)}(t_{\theta'}).
\]
Since \(|t_{\theta'}-t_\tau|\le\sigma\le u\), applying the support bounds for \(\widehat w_{\tau,u}\)
to the coefficient of \(\phi_d^{(j+1)}(t_{\theta'})\) in the Taylor expansion of \(\xi-\xi'\) gives
\[
 \biggl|
 \lambda_{j+1}(t_\theta)-\lambda_{j+1}(t_{\theta'})
 +\lambda_j(t_\theta)(t_\theta-t_{\theta'})
 +\sum_{i<j}\lambda_i(t_\theta)
       \frac{(t_\theta-t_{\theta'})^{j+1-i}}{(j+1-i)!}
 \biggr|
 \lesssim R^{-j/d}\sigma^2u^{j+1-d}.
\]
By \eqref{eq:hccurvaturesupport} and
\(|t_\theta-t_{\theta'}|\le\sigma\),
\[
 \biggl|
 \sum_{i<j}\lambda_i(t_\theta)
       \frac{(t_\theta-t_{\theta'})^{j+1-i}}{(j+1-i)!}
 \biggr|
 \lesssim
 \frac{\sigma}{\ell}R^{-j/d}|t_\theta-t_{\theta'}|
 \ll R^{-j/d}|t_\theta-t_{\theta'}|.
\]
Together with \(|\lambda_j(t_\theta)|\sim R^{-j/d}\) and
\(|\lambda_{j+1}(t_\theta)|+|\lambda_{j+1}(t_{\theta'})|
\lesssim R^{-(j+1)/d}\), this gives
\begin{equation}\label{eq:hcpairseparation}
 |t_\theta-t_{\theta'}|\lesssim
 \max\{R^{-1/d},\sigma^2u^{-(d-j-1)}\}.
\end{equation}
Since \(u\ge\sigma^{\eta^4}\) and \(\eta\) is sufficiently
small, we may choose a dyadic
\(\sigma'\sim\max\{R^{-1/d},\sigma^2u^{-(d-j-1)}\}\)
with \(R^{-1/d}\le\sigma'\le
\max\{R^{-1/d},\sigma^{3/2}\}<\sigma\).
Only neighbouring \(\sigma'\) intervals interact, so
Cauchy--Schwarz gives
\[
 \sum_{|\tau|=\sigma}|f_\tau|^2*w_{\tau,u}
 \lesssim
 \sum_{|\tau|=\sigma}
 \sum_{\substack{\tau'\subset\tau\\|\tau'|=\sigma'}}
 |f_{\tau'}|^2*w_{\tau,u}.
\]

To combine these terms into one convolution, let \(w\) be a
nonnegative, \(L^1\) normalized radial kernel with sufficiently
rapid decay at scale \(R^{j/d}\sigma^{-2}u^{d-j-1}\), so that
every width of \(w_{\tau,u}\) is at most this radius.
The ratio of the ball volume to the kernel box volume is at most
\(u^{-C}\), so choosing \(w\) to decay more slowly than
\(w_{\tau,u}\) gives
\(w_{\tau,u}\lesssim u^{-C}w\lesssim\sigma^{-C\eta^4}w\).
Substituting into \eqref{eq:hcrefinementlowinput} gives
\begin{equation}\label{eq:hcrefinementlowaverage}
 \Bigl\|\sum_{|\tau|=\sigma}|f_\tau|^2
 \Bigr\|_{L_\#^{p_{d-m}/4}(B)}
 \lesssim\sigma^{-C\eta^4}
 \biggl\|\Bigl(\sum_{|\tau'|=\sigma'}|f_{\tau'}|^2\Bigr)*w
 \biggr\|_{L_\#^{p_{d-m}/4}(B)}.
\end{equation}

Fix \(u,\sigma'\), raise \eqref{eq:hcrefinementlowaverage}
to the power \(p_{d-m}/2\), multiply by \(|B|\), and sum over
the boxes \(B\) satisfying \eqref{eq:hcrefinementlowinput}
with this \(u\).
Applying \eqref{eq:selaverage} with
\(F=[(\sum_{|\tau'|=\sigma'}|f_{\tau'}|^2)*w]^{p_{d-m}/4}\)
and \(t=[\alpha/\sigma^{-C\eta^4}]^{p_{d-m}/2}\)
bounds their contribution to the LHS of
Lemma~\ref{lem:hcrefinement} by
\begin{equation}\label{eq:hcrefinementlowconvolvedtail}
 \lesssim_\eta\sigma^{-C\eta^4}
 \int_{(\sum_{|\tau'|=\sigma'}|f_{\tau'}|^2)*w
       \,\gtrsim\,[\alpha/\sigma^{-C\eta^4}]^2}
 \biggl[\Bigl(\sum_{|\tau'|=\sigma'}|f_{\tau'}|^2\Bigr)*w
 \biggr]^{p_{d-m}/2}.
\end{equation}
Apply \eqref{eq:tailYoung} with
\(F=\sum_{|\tau'|=\sigma'}|f_{\tau'}|^2\), \(K=w\),
and exponent \(p_{d-m}/2\).
The locally constant property on unit cubes then bounds
\eqref{eq:hcrefinementlowconvolvedtail} by
\begin{equation}\label{eq:hcrefinementlowout}
 \lesssim_\eta\sigma^{-C\eta^4}
 \sum_{\substack{Q\parallel[0,1]^{d+1}\\
 Q\in\mathcal G_{[0,1]}(\frac{\alpha}{\sigma^{-C\eta^4}})}}
 |Q|\biggl\|\Bigl(\sum_{|\tau'|=\sigma'}|f_{\tau'}|^2\Bigr)^{1/2}
 \biggr\|_{L_\#^{p_{d-m}/2}(Q)}^{p_{d-m}}.
\end{equation}
This gives the first sum in Lemma~\ref{lem:hcrefinement}.

\medskip
\noindent\textit{Case 2.}
Suppose that \eqref{eq:hcrefinementlowinput} fails for every
dyadic \(u\in[\sigma^{\eta^4},1]\), and fix such a box \(B\).
Starting at \(u=\sigma^{\eta^4}\), decrease \(u\) dyadically
and let \(u_*\) be the first scale for which
\[
 \Bigl\|\sum_{|\tau|=\sigma}|f_\tau|^2*w_{\tau,u_*}
 \Bigr\|_{L_\#^{p_{d-m}/4}(B)}
 \gtrsim
 \Bigl\|\sum_{|\tau|=\sigma}|f_\tau|^2
 \Bigr\|_{L_\#^{p_{d-m}/4}(B)}.
\]
At \(u=\sigma\), the locally constant property gives \eqref{eq:hcrefinementinitialaverage}, so such a scale exists.
The failure of \eqref{eq:hcrefinementlowinput} gives
\(u_*<\sigma^{\eta^4}\).
By the choice of \(u_*\) and comparison of the kernels
at comparable scales,
\begin{equation}\label{eq:hcstoppingupper}
 \fint_B\Bigl(\sum_{|\tau|=\sigma}|f_\tau|^2*w_{\tau,u}
 \Bigr)^{p_{d-m}/4}
 \lesssim
 \fint_B\Bigl(\sum_{|\tau|=\sigma}|f_\tau|^2\Bigr)^{p_{d-m}/4},
 \qquad
 \sigma\le u\le\sigma^{\eta^4},\quad u_*\lesssim u.
\end{equation}

We now decompose \(\sum_{|\tau|=\sigma}|f_\tau|^2*w_{\tau,u_*}\)
into low and high frequencies.
As in Section~\ref{subsec:hcoverviewpositive}, choose
\(\chi\in C_c^\infty(\mathbb R)\), \(0\le\chi\le1\),
with \(\chi=1\) on \([-C,C]\) and
\(\operatorname{supp}\chi\subset[-2C,2C]\).
For \(\xi=\sum_i\lambda_i\phi_d^{(i)}(t_\tau)\), put
\[
 \chi_{\tau,\le u}(\xi)
 =\prod_{i=0}^{j+1}
 \chi\!\bigl(R^{j/d}\sigma^{-2}u^{d-i}\lambda_i\bigr)
 \prod_{i=j+2}^d
 \chi\!\bigl(R^{j/d}\sigma^{-2}u^d\lambda_i\bigr),
 \qquad
 \chi_{\tau,u}=\chi_{\tau,\le u}-\chi_{\tau,\le2u}.
\]
On \(\operatorname{supp}\widehat{|f_\tau|^2}\), we have
\(|\lambda_i|\lesssim R^{-j/d}\sigma^{i-j}\).
Thus \(\chi_{\tau,\le\sigma}=1\) there, and the factors with
\(i\ge j+2\) equal one for \(\sigma\le u\le1\).
Consequently,
\[
\begin{aligned}
 \sum_{|\tau|=\sigma}|f_\tau|^2*w_{\tau,u_*}
 &=\sum_{|\tau|=\sigma}|f_\tau|^2*w_{\tau,u_*}
       *\check\chi_{\tau,\le\sigma^{\eta^4}}+\sum_{\sigma\le s'<\sigma^{\eta^4}}
       \sum_{|\tau|=\sigma}|f_\tau|^2*w_{\tau,u_*}
       *\check\chi_{\tau,s'}.
\end{aligned}
\]
The first term is small: since \(u_*\le\sigma^{\eta^4}\),
the widths of \(w_{\tau,u_*}\) are at most the corresponding
widths of \(\check\chi_{\tau,\le\sigma^{\eta^4}}\), so
\(w_{\tau,u_*}*|\check\chi_{\tau,\le\sigma^{\eta^4}}|
\lesssim w_{\tau,\sigma^{\eta^4}}\).
The failure of \eqref{eq:hcrefinementlowinput} at
\(u=\sigma^{\eta^4}\) therefore gives
\[
 \Bigl\|\sum_{|\tau|=\sigma}|f_\tau|^2*w_{\tau,u_*}
       *\check\chi_{\tau,\le\sigma^{\eta^4}}
 \Bigr\|_{L_\#^{p_{d-m}/4}(B)}
 \lesssim\sigma^{\eta^4}
 \Bigl\|\sum_{|\tau|=\sigma}|f_\tau|^2
 \Bigr\|_{L_\#^{p_{d-m}/4}(B)}.
\]
Moreover, the Fourier support bounds give
\(\widehat w_{\tau,u_*}\chi_{\tau,s'}=0\)
unless \(s'\gtrsim u_*\).
Dyadic pigeonholing therefore gives
\(\sigma\le s'<\sigma^{\eta^4}\), with \(u_*\lesssim s'\),
such that the functions
\[
 g_\rho=
 \sum_{\substack{\tau\subset\rho\\|\tau|=\sigma}}
 |f_\tau|^2*w_{\tau,u_*}*\check\chi_{\tau,s'},
 \qquad \rho\in\mathcal J(s').
\]
satisfy
\begin{equation}\label{eq:hcrefinementdominance}
 \fint_B\Bigl(\sum_{|\tau|=\sigma}|f_\tau|^2\Bigr)^{p_{d-m}/4}
 \lesssim_\delta\sigma^{-\delta}
 \fint_B\Bigl|\sum_{\rho\in\mathcal J(s')}g_\rho
 \Bigr|^{p_{d-m}/4}.
\end{equation}

To estimate the RHS of \eqref{eq:hcrefinementdominance},
we first bound each \(g_\rho\) in terms of the original
pieces \(f_\tau\).
Since \(u_*\lesssim s'\), comparing the kernel widths gives
\(w_{\tau,u_*}*|\check\chi_{\tau,s'}|\lesssim w_{\tau,s'}\),
and hence
\begin{equation}\label{eq:hcrefinementfieldmajorant}
 |g_\rho|\lesssim
 \sum_{\substack{\tau\subset\rho\\|\tau|=\sigma}}
 |f_\tau|^2*w_{\tau,s'}.
\end{equation}

We now show that dilation places the Fourier support of each
\(g_\rho\) in a cap of \(\Gamma_{j+1}^{d+1}((s')^{-d})\).
This will allow us to apply Proposition~\ref{prop:hctaylorsf}.
Take \(\xi=\sum_i\lambda_i\phi_d^{(i)}(t_\tau)\)
in the Fourier support of the \(\tau\) summand defining
\(g_\rho\). The cutoff \(\chi_{\tau,s'}\) gives
\[
 |\lambda_i|\lesssim R^{-j/d}\sigma^2(s')^{i-d}
 \quad(0\le i\le j+1),\qquad
 \max_{0\le i\le j+1}
 R^{j/d}\sigma^{-2}(s')^{d-i}|\lambda_i|\sim1.
\]
For \(j+2\le i\le d\), the bound
\(|\lambda_i|\lesssim R^{-j/d}\sigma^{i-j}\) on
\(\operatorname{supp}\widehat{|f_\tau|^2}\) gives
\begin{equation}\label{eq:hcrefinementsupport}
 \frac{|\lambda_i|}{R^{-j/d}\sigma^2(s')^{i-d}}
 \lesssim\sigma^{i-j-2}(s')^{d-i}\le1.
\end{equation}
Thus dilation by \(R^{j/d}\sigma^{-2}(s')^d\) gives
\[
 |R^{j/d}\sigma^{-2}(s')^d\lambda_i|\lesssim(s')^i
 \quad(0\le i\le d),\qquad
 \max_{0\le i\le j+1}(s')^{-i}
 |R^{j/d}\sigma^{-2}(s')^d\lambda_i|\sim1.
\]
Hence
\[
 \operatorname{supp}
 \widehat{g_\rho(R^{j/d}\sigma^{-2}(s')^d\,\cdot)}
 \subset\theta_\rho,\qquad
 \theta_\rho\in\Xi_{j+1}^{d+1}((s')^{-d}).
\]

We now estimate the RHS of \eqref{eq:hcrefinementdominance}
using Proposition~\ref{prop:hctaylorsf}.
As in the proof of Lemma~\ref{lem:curvehigh}, we use exponent
\(p_{d-m}/4\) when \(p_{d-m}\ge11\), and first apply
H\"older to pass from exponent \(7/4\) to \(7/3\) when
\(p_{d-m}=7\).

First suppose \(p_{d-m}\ge11\).
Choose \(\chi_B\sim W_B^{4/p_{d-m}}\) with
\(\operatorname{supp}\widehat{\chi_B}\subset B^*\),
as in the proof of Lemma~\ref{lem:curvehigh}. Then
\[
 \fint_B
 \Bigl|\sum_{\rho\in\mathcal J(s')}g_\rho\Bigr|^{p_{d-m}/4}
 \sim \frac1{|B|}\int
 \Bigl|\chi_B\sum_{\rho\in\mathcal J(s')}g_\rho
 \Bigr|^{p_{d-m}/4}.
\]
Since \(B\) has side length \(R^{j/d}\sigma^{-2}\),
multiplication by \(\chi_B\) adds Fourier widths
\(R^{-j/d}\sigma^2\).
Their ratios to the coefficient widths
\(R^{-j/d}\sigma^2(s')^{i-d}\) are
\((s')^{d-i}\le1\), and for \(i\le j+1\) they are
at most \((s')^{d-j-1}\ll1\). Hence
\[
 \operatorname{supp}
 \widehat{(\chi_Bg_\rho)
 (R^{j/d}\sigma^{-2}(s')^d\,\cdot)}
 \subset\theta_\rho,\qquad
 \theta_\rho\in\Xi_{j+1}^{d+1}((s')^{-d}).
\]
Apply Proposition~\ref{prop:hctaylorsf} at exponent
\(p_{d-m}/4\in[2,p_{d-j-1}]\) and parameter \((s')^{-d}\) to
\(\bigl(\chi_B\sum_{\rho\in\mathcal J(s')}g_\rho\bigr)
(R^{j/d}\sigma^{-2}(s')^d\,\cdot)\), with fine pieces
\((\chi_Bg_\rho)(R^{j/d}\sigma^{-2}(s')^d\,\cdot)\).
Returning to the original coordinates and using
\eqref{eq:hcrefinementdominance} and
\eqref{eq:hcrefinementfieldmajorant}, we obtain
\begin{equation}\label{eq:hcrefinementordinary}
 \fint_B
 \Bigl(\sum_{|\tau|=\sigma}|f_\tau|^2\Bigr)^{p_{d-m}/4}
 \lesssim_\delta\sigma^{-\delta}
 \fint_B
 \biggl[
 \sum_{\rho\in\mathcal J(s')}
 \Bigl(
 \sum_{\substack{\tau\subset\rho\\|\tau|=\sigma}}
 |f_\tau|^2*w_{\tau,s'}
 \Bigr)^2
 \biggr]^{p_{d-m}/8}.
\end{equation}

For \(p_{d-m}=7\), choose
\(\chi_B^{7/3}\sim W_B\) with
\(\operatorname{supp}\widehat{\chi_B}\subset B^*\).
As in the \(n=3\) case of the proof of
Lemma~\ref{lem:curvehigh}, H\"older gives
\[
 \biggl(
 \fint_B
 \Bigl|\sum_{\rho\in\mathcal J(s')}g_\rho\Bigr|^{7/4}
 \biggr)^{4/3}
 \lesssim
 \frac1{|B|}\int
 \Bigl|\chi_B\sum_{\rho\in\mathcal J(s')}g_\rho
 \Bigr|^{7/3}.
\]
With this choice of \(\chi_B\), apply Proposition~\ref{prop:hctaylorsf}
to \(\sum_{\rho\in\mathcal J(s')}
(\chi_Bg_\rho)(R^{j/d}\sigma^{-2}(s')^d\,\cdot)\)
at exponent \(7/3\in[2,p_{d-j-1}]\) and parameter \((s')^{-d}\).
Returning to the original coordinates and using
\eqref{eq:hcrefinementdominance} and
\eqref{eq:hcrefinementfieldmajorant} gives
\begin{equation}\label{eq:hcrefinementseven}
 \biggl[
 \fint_B
 \Bigl(\sum_{|\tau|=\sigma}|f_\tau|^2\Bigr)^{7/4}
 \biggr]^{4/3}
 \lesssim_\delta\sigma^{-\delta}
 \fint_B
 \biggl[
 \sum_{\rho\in\mathcal J(s')}
 \Bigl(
 \sum_{\substack{\tau\subset\rho\\|\tau|=\sigma}}
 |f_\tau|^2*w_{\tau,s'}
 \Bigr)^2
 \biggr]^{7/6}.
\end{equation}

To obtain the second sum in Lemma~\ref{lem:hcrefinement},
we first bound \eqref{eq:hcrefinementordinary} and
\eqref{eq:hcrefinementseven} by integrals over the set where
\(\sum_{\tau\subset\rho,\,|\tau|=\sigma}
|f_\tau|^2*w_{\tau,s'}\)
exceeds the threshold \([\alpha/\sigma^{-\delta}]^2\).
Since \(B\in\mathcal G_{[0,1]}(\alpha)\), we have
\(\fint_B(\sum_{|\tau|=\sigma}|f_\tau|^2)^{p_{d-m}/4}
\ge\alpha^{p_{d-m}/2}\).
Together with \eqref{eq:hcstoppingupper} at \(u=s'\),
this allows us to repeat the argument leading to
\eqref{eq:curvehighboxtail}, with
\(\sum_{\tau\subset\rho,\,|\tau|=\sigma}
|f_\tau|^2*w_{\tau,s'}\)
in place of \(a_{\rho'}\).
Applying that argument to \eqref{eq:hcrefinementordinary}
and \eqref{eq:hcrefinementseven} gives
\begin{equation}\label{eq:hcrefinementextract}
\begin{aligned}
 &\biggl[
 \fint_B
 \Bigl(\sum_{|\tau|=\sigma}|f_\tau|^2\Bigr)^{p_{d-m}/4}
 \biggr]^2\lesssim_\delta
 \frac{\sigma^{-\delta}}{|B|}
 \sum_{\rho\in\mathcal J(s')}
 \int_{\sum_{\tau\subset\rho,\,|\tau|=\sigma}
 |f_\tau|^2*w_{\tau,s'}>[\alpha/\sigma^{-\delta}]^2}
 \Bigl(
 \sum_{\substack{\tau\subset\rho\\|\tau|=\sigma}}
 |f_\tau|^2*w_{\tau,s'}
 \Bigr)^{p_{d-m}/2}W_B.
\end{aligned}
\end{equation}

To remove the convolutions in \eqref{eq:hcrefinementextract}
after summing over \(B\), we first compare the kernels
\(w_{\tau,s'}\), \(\tau\subset\rho\), with a single kernel
for each \(\rho\).
For \(w_{\rho,s'}\), take a nonnegative, \(L^1\) normalized kernel with sufficiently rapid decay,
adapted to the box in the definition of \(w_{\tau,u}\) with \(u=s'\) and
\(|x\cdot\phi_d^{(i)}(t_\rho)|\) in place of \(|x\cdot\phi_d^{(i)}(t_\tau)|\).

Indeed, since \(\tau\subset\rho\), we have
\(|t_\tau-t_\rho|\lesssim s'\).
Taylor expansion gives
\[
\begin{aligned}
 R^{j/d}\sigma^{-2}(s')^{d-k}
 |t_\tau-t_\rho|^{k-i}
 &\lesssim R^{j/d}\sigma^{-2}(s')^{d-i},
 &&0\le i\le k\le j+1,\\
 R^{j/d}\sigma^{-2}(s')^d
 |t_\tau-t_\rho|^{k-i}
 &\lesssim R^{j/d}\sigma^{-2}(s')^d,
 &&0\le i\le k,\quad j+2\le k\le d.
\end{aligned}
\]
These inequalities also hold with \(t_\tau\) and \(t_\rho\) interchanged.
Thus the bounds on \(|x\cdot\phi_d^{(i)}(t_\tau)|\) and on \(|x\cdot\phi_d^{(i)}(t_\rho)|\)
define boxes comparable up to a fixed dilation, uniformly in \(\tau\subset\rho\).
We may therefore choose \(w_{\rho,s'}\) so that \(w_{\tau,s'}\lesssim w_{\rho,s'}\) for all such \(\tau\), giving
\[
 \sum_{\substack{\tau\subset\rho\\|\tau|=\sigma}}
 |f_\tau|^2*w_{\tau,s'}
 \lesssim
 \Bigl(
 \sum_{\substack{\tau\subset\rho\\|\tau|=\sigma}}
 |f_\tau|^2
 \Bigr)*w_{\rho,s'}.
\]

Fix \(s'\), multiply \eqref{eq:hcrefinementextract} by
\(|B|\), and sum over the boxes \(B\) in Case 2 for which
this \(s'\) was chosen.
Using \(\sum_BW_B\lesssim1\) and lowering the threshold to
\([\alpha/(\sigma^{-\delta}N(s'))]^2\), we bound the
contribution of these boxes to the LHS of
Lemma~\ref{lem:hcrefinement} by
\[
 \lesssim_\delta\sigma^{-\delta}
 \sum_{\rho\in\mathcal J(s')}
 \int_{(\sum_{\tau\subset\rho,\,|\tau|=\sigma}|f_\tau|^2)*w_{\rho,s'}
       \,\gtrsim\,[\alpha/(\sigma^{-\delta}N(s'))]^2}
 \biggl[
 \Bigl(
 \sum_{\substack{\tau\subset\rho\\|\tau|=\sigma}}
 |f_\tau|^2
 \Bigr)*w_{\rho,s'}
 \biggr]^{p_{d-m}/2}.
\]
We now apply \eqref{eq:tailYoung} to each integral, with
\(F=\sum_{\tau\subset\rho,\,|\tau|=\sigma}|f_\tau|^2\),
\(K=w_{\rho,s'}\), exponent \(p_{d-m}/2\), and threshold
\([\alpha/(\sigma^{-\delta}N(s'))]^2\).
This removes the convolution, and the locally constant property on
\(\rho^*\) then bounds these integrals by
\begin{equation}\label{eq:hcrefinementhighout}
 \lesssim_\delta\sigma^{-\delta}
 \sum_{\rho\in\mathcal J(s')}
 \sum_{\substack{Q\parallel\rho^*\\
 Q\in\mathcal G_\rho(\frac{\alpha}{\sigma^{-\delta}})}}
 |Q|
 \biggl\|
 \Bigl(
 \sum_{\substack{\tau\subset\rho\\|\tau|=\sigma}}
 |f_\tau|^2
 \Bigr)^{1/2}
 \biggr\|_{L_\#^{p_{d-m}/2}(Q)}^{p_{d-m}}.
\end{equation}

Summing \eqref{eq:hcrefinementlowout} over \((u,\sigma')\)
from Case 1 and \eqref{eq:hcrefinementhighout} over \(s'\)
from Case 2 covers every \(B\in\mathcal G_{[0,1]}(\alpha)\).
Dyadic pigeonholing gives one common \(\sigma'\) in the first
sum and one common \(s'\) in the second.
The logarithmic losses are absorbed into
\(\sigma^{-C\eta^4}\) and \(\sigma^{-\delta}\), respectively,
both outside the sums and in the selection conditions.
This proves Lemma~\ref{lem:hcrefinement}.
\end{proof}

\begin{proof}[Proof of Lemma~\ref{lem:hcoverviewrefinement}]
By definition, \(\Gamma_{j,1/K}^{d+1}(R)
\subset\widetilde\Gamma_{j,1}^{d+1}(R)\) for \(0\le j\le m\).
Applying Lemma~\ref{lem:hcrefinement} with \(\ell=1\)
therefore proves the claim.
\end{proof}

\subsection{Proof of Proposition~\ref{prop:hcstable} and \eqref{eq:hcrestrictedconvolution}}
\label{subsec:hcgroupedinduction}

\begin{proof}[Proof of Proposition~\ref{prop:hcstable} and
\eqref{eq:hcrestrictedconvolution}]
Under the hypotheses of Proposition~\ref{prop:hcstablefunction},
fix \(m<d-2\) and let \(f\) be a Schwartz function with
\(\operatorname{supp}\widehat f\subset\Gamma_{j,1/K}^{d+1}(R)\),
\(0\le j\le m\).
We choose small parameters \(0<\delta\ll\eta\ll\epsilon\)
as in Section~\ref{subsec:hcoverviewiteration}.
For \(R^{-1/d}\le\sigma\le A^{-1}\) and \(\sigma\le s\le1\),
we prove the stronger bound
\begin{equation}\label{eq:hcgroupedconvolution}
 \sum_{\rho\in\mathcal J(s)}
 \sum_{\substack{Q\parallel\rho^*\\Q\in\mathcal G_\rho(\alpha)}}
 |Q|\biggl\|\Bigl(
 \sum_{\substack{\tau\subset\rho\\|\tau|=\sigma}}
 |f_\tau|^2\Bigr)^{1/2}
 \biggr\|_{L_\#^{p_{d-m}/2}(Q)}^{p_{d-m}}
 \lesssim_\epsilon
 \text{RHS of \eqref{eq:hcgeneralconvolution}},
\end{equation}
where both occurrences of \(R^\epsilon\) on the right are
replaced by \(R^\epsilon(\sigma s)^\epsilon\).
We prove \eqref{eq:hcgroupedconvolution} and
\eqref{eq:hcrestrictedconvolution} together by induction on
\(R\) and, at fixed \(R\), on the dyadic parameter \(\sigma s\).
For bounded \(R\), both bounds follow from H\"older and the locally constant property, followed by Lemma~\ref{lem:hcendpointaverages}
applied to \(f\) with parameter \(R\) and \(u=R^{-1/d}\),
using the terms \(q_1=\cdots=q_{d-m-2}=u=R^{-1/d}\) in
\eqref{eq:hcgeneralconvolution}.

Assume that \eqref{eq:hcgroupedconvolution} and
\eqref{eq:hcrestrictedconvolution} hold at every parameter
\(R'<R/2\).
With \(R\) fixed, assume also that
\eqref{eq:hcgroupedconvolution} holds for pairs
\((\sigma',s')\) satisfying \(\sigma's'<\sigma s/2\).

We first deal with the case \(\sigma=R^{-1/d}\).
By \eqref{eq:selaverage},
\[
 \text{LHS of \eqref{eq:hcgroupedconvolution}}
 \lesssim
 \sum_{\rho\in\mathcal J(s)}
 \int_{\sum_{\theta\subset\rho}|f_\theta|^2
                    \gtrsim[\alpha/N(s)]^2}
 \Bigl(\sum_{\theta\subset\rho}|f_\theta|^2
 \Bigr)^{p_{d-m}/2}.
\]
For each interval \(\rho=[t_\rho,t_\rho+s]\in\mathcal J(s)\), set
\(g(y)=f_\rho(L_{t_\rho,s}^{-T}y)\).
Then
\(\operatorname{supp}\widehat g\subset\Gamma_m^{d+1}(Rs^d)\),
with fine pieces \(f_\theta(L_{t_\rho,s}^{-T}y)\),
\(\theta\subset\rho\), of angular width \((Rs^d)^{-1/d}\).

We therefore apply \eqref{eq:hcpositiveconvolution} to \(g\)
with parameter \(Rs^d\), amplitude \(\alpha/N(s)\),
and loss \(\eta\).
This further bounds the integral for this \(\rho\) in the
sum above, after rescaling, by the RHS of
\eqref{eq:hcgeneralconvolution} for \(g\) at parameter \(Rs^d\).
In this sum, denote the intermediate caps by \(\tau_l'\),
with \(|\tau_l'|=q_l/s\), their parent by \(\rho'\),
with \(|\rho'|=u/s\), the convolution weights by
\(w_{B_{\tau_l'}}\), and the averaging boxes by \(B'\).
Returning by \(x=L_{t_\rho,s}^{-T}y\), write
\(\tau_l=L_{t_\rho,s}\tau_l'\subset\rho\),
with \(|\tau_l|=q_l\), and \(B=L_{t_\rho,s}^{-T}B'\).
The \(i\)th widths of \(L_{t_\rho,s}^{-T}B_{\tau_l'}\)
and \(B\) satisfy
\begin{equation}\label{eq:hcalllayerrescaling}
\begin{aligned}
 L_{t_\rho,s}^{-T}B_{\tau_l'}:\quad&
 s^{-i}\frac{(Rs^d)^{(m+2)/d}(q_l/s)^{m+l+1-i}}
                   {\prod_{b<l}(q_b/s)}
 =\frac{R^{(m+2)/d}q_l^{m+l+1-i}}{\prod_{b<l}q_b},
 &&0\le i\le m+l+1,\\
 &
 s^{-i}\frac{(Rs^d)^{(m+2)/d}(u/s)^{m+l+1-i}}
                   {\prod_{b<l}(q_b/s)}
 =\frac{R^{(m+2)/d}u^{m+l+1-i}}{\prod_{b<l}q_b},
 &&m+l+1<i\le d,\\
 B:\quad&
 s^{-i}\frac{(Rs^d)^{(m+2)/d}(u/s)^{d-i}}
                   {\prod_{b=1}^{d-m-2}(q_b/s)}
 =\frac{R^{(m+2)/d}u^{d-i}}{\prod_{b=1}^{d-m-2}q_b},
 &&0\le i\le d.
\end{aligned}
\end{equation}
These are exactly the widths of \(B_{\tau_l}\) in
\eqref{eq:hcnestedwidths} and of the averaging boxes in
\eqref{eq:hcgeneralconvolution}, with parent
\(L_{t_\rho,s}\rho'\subset\rho\) of width \(u\).
Here
\(R^{-1/d}\le q_1\le\cdots\le q_{d-m-2}\le u\le s\).

The active parent count at rescaled width \(u/s\) is
\(N_\rho(u)\).
Since \((Rs^d)^\eta\le R^\eta\), the bound
\(N(s)N_\rho(u)\lesssim N(u)\) from Lemma~\ref{lem:regular} gives
\(\frac{\alpha}{(Rs^d)^\eta N(s)N_\rho(u)}
\gtrsim\frac{\alpha}{R^\eta N(u)}\).
Consequently, summing over \(\rho\) and extending \(u\le s\)
to \(u\le1\) bounds the LHS of \eqref{eq:hcgroupedconvolution}
by the RHS of \eqref{eq:hcgeneralconvolution}, with both
occurrences of \(R^\epsilon\) replaced by \(R^\eta\).
Since \(\sigma=R^{-1/d}\) and \(s\ge\sigma\),
\[
 \frac{R^\eta}{R^\epsilon(\sigma s)^\epsilon}
 \le R^{-\epsilon/C}
\]
for \(\eta\) sufficiently small in terms of \(\epsilon\).
This proves \eqref{eq:hcgroupedconvolution} at
\(\sigma=R^{-1/d}\).

We now consider the case \(\sigma>R^{-1/d}\).
Corollary~\ref{cor:hccoarse}, applied to \(f\) with \(\ell=1\)
and loss \(\delta\), bounds the LHS of
\eqref{eq:hcgroupedconvolution} by the RHS of
\eqref{eq:hcoverviewcoarse}.
Fix \(\sigma\le u\le s\).
We estimate the corresponding term on the RHS of
\eqref{eq:hcoverviewcoarse}, namely
\[
 \sigma^{-\delta}\sum_{\rho\in\mathcal J(u)}
 \sum_{\substack{B\parallel
       \frac{(uR^{1/d})^j}{(\sigma R^{1/d})^2}U_{\rho,R}\\
 B\in\mathcal G_\rho(\alpha\sigma^\delta)}}
 |B|\biggl\|\Bigl(
 \sum_{\substack{\tau\subset\rho\\|\tau|=\sigma}}
 |f_\tau|^2\Bigr)^{1/2}
 \biggr\|_{L^{p_{d-m}/2}_\#(B)}^{p_{d-m}}.
\]
For \(u<\sigma^{1-\eta}\), we use
\eqref{eq:hcoverviewnarrow} and then apply
\eqref{eq:hcrestrictedconvolution} at parameter
\(R\sigma^d<R/2\) after Taylor rescaling \(f_\tau\).
For \(u\ge\sigma^{1-\eta}\),
Lemma~\ref{lem:hcoverviewrefinement} gives the scales \((\sigma',u)\) with \(\sigma'<\sigma\), or
\((\sigma,u')\) with \(u'<u\).

\medskip
\noindent\textit{Case 1: \(u<\sigma^{1-\eta}\).}
As in \eqref{eq:hcoverviewnarrow}, we can bound the fixed \(u\)
term on the RHS of \eqref{eq:hcoverviewcoarse} by
\[
 \lesssim \sigma^{-C\eta}
 \sum_{\rho\in\mathcal J(u)}
 \sum_{\substack{\tau\subset\rho\\|\tau|=\sigma}}
 \int_{|f_\tau|\gtrsim\alpha\sigma^\delta/[N(u)N_\rho(\sigma)]}
 |f_\tau|^{p_{d-m}}.
\]
Since \(\operatorname{supp}
\widehat{f_\tau(L_{t_\tau,\sigma}^{-T}\,\cdot)}
\subset\Gamma_{j,1/K}^{d+1}(R\sigma^d)\) for
\(\tau\in\mathcal J(\sigma)\), we can apply
\eqref{eq:hcrestrictedconvolution} inductively to this function
with parameter \(R\sigma^d<R/2\) and amplitude
\(\alpha\sigma^\delta/[N(u)N_\rho(\sigma)]\).

Returning to the original coordinates,
\eqref{eq:hcalllayerrescaling}, with \(\sigma\) in place of \(s\)
and \(v\) in place of \(u\), gives the convolution weights
\(w_{B_{\tau_l}}\) and boxes \(B\) in
\eqref{eq:hcgeneralconvolution}, with
\(R^{-1/d}\le q_1\le\cdots\le q_{d-m-2}\le v\le\sigma\).
For a parent \(\rho'\subset\tau\) with \(|\rho'|=v\), applying
Lemma~\ref{lem:regular} twice gives
\(N(u)N_\rho(\sigma)N_\tau(v)\lesssim N(v)\).
Hence the amplitude in the selection condition satisfies
\[
 \frac{\alpha\sigma^\delta}
 {(R\sigma^d)^\epsilon N(u)N_\rho(\sigma)N_\tau(v)}
 \gtrsim
 \frac{\alpha}
 {\sigma^{-C\eta}(R\sigma^d)^\epsilon N(v)}.
\]
Together with the factor \(\sigma^{-C\eta}\) in
\eqref{eq:hcoverviewnarrow}, summing over \(\rho,\tau\) and
extending \(v\le\sigma\) to \(v\le1\) therefore bounds the
fixed \(u\) term on the RHS of \eqref{eq:hcoverviewcoarse}
by the RHS of \eqref{eq:hcgeneralconvolution}, with both
occurrences of \(R^\epsilon\) replaced by
\(\sigma^{-C\eta}(R\sigma^d)^\epsilon\).

Since \(s\ge\sigma\),
\[
 \frac{\sigma^{-C\eta}(R\sigma^d)^\epsilon}
      {R^\epsilon(\sigma s)^\epsilon}
 \le \sigma^{(d-2)\epsilon-C\eta}.
\]
This is the third factor in \eqref{eq:hccontractions} and
finishes the fixed \(u\) estimate in Case~1.

\medskip
\noindent\textit{Case 2: \(u\ge\sigma^{1-\eta}\).}
For each \(\rho\in\mathcal J(u)\), apply
Lemma~\ref{lem:hcoverviewrefinement} to
\(f_\rho(L_{t_\rho,u}^{-T}\,\cdot)\), with parameter \(Ru^d\),
cap width \(\sigma/u\), and amplitude
\(\alpha\sigma^\delta/N(u)\).
After applying the lemma and returning to the original
coordinates, denote the new scales by \((\sigma',u')\).
Summing over \(\rho\) and applying dyadic pigeonholing as in
\eqref{eq:hcoverviewrefinement}, we can bound the fixed \(u\)
term on the RHS of \eqref{eq:hcoverviewcoarse} by
\[
\begin{aligned}
 &\lesssim \sigma^{-C\eta^4}
 \sum_{\rho\in\mathcal J(u)}
 \sum_{\substack{Q\parallel\rho^*\\
 Q\in\mathcal G_\rho(\frac{\alpha}{\sigma^{-C\eta^4}})}}
 |Q|\biggl\|\Bigl(
 \sum_{\substack{\tau'\subset\rho\\|\tau'|=\sigma'}}
 |f_{\tau'}|^2\Bigr)^{1/2}
 \biggr\|_{L^{p_{d-m}/2}_\#(Q)}^{p_{d-m}}
 \\
 &\quad+\sigma^{-2\delta}
 \sum_{\rho'\in\mathcal J(u')}
 \sum_{\substack{Q\parallel(\rho')^*\\
 Q\in\mathcal G_{\rho'}(\frac{\alpha}{\sigma^{-2\delta}})}}
 |Q|\biggl\|\Bigl(
 \sum_{\substack{\tau\subset\rho'\\|\tau|=\sigma}}
 |f_\tau|^2\Bigr)^{1/2}
 \biggr\|_{L^{p_{d-m}/2}_\#(Q)}^{p_{d-m}},
\end{aligned}
\]
where \(R^{-1/d}\le\sigma'\le
\max\{R^{-1/d},\sigma^{1+\eta/4}\}\) and
\(\sigma\le u'\le u\sigma^{\eta^5/2}\).

For the first sum over \(Q\) in \eqref{eq:hcoverviewrefinement},
suppose that \(\sigma'>R^{-1/d}\).
Since \(\sigma\le A^{-1}\), we have
\(\sigma'u\le\sigma^{\eta/4}\sigma s<\sigma s/2\) for
sufficiently large \(A\).
Apply \eqref{eq:hcgroupedconvolution} inductively at the
fixed parameter \(R\), at scales
\((\sigma',u)\) and amplitude \(\alpha/\sigma^{-C\eta^4}\).
This bounds the first term in \eqref{eq:hcoverviewrefinement}
by the RHS of \eqref{eq:hcgeneralconvolution}, with both
occurrences of \(R^\epsilon\) replaced by
\(\sigma^{-C\eta^4}R^\epsilon(\sigma'u)^\epsilon\), and
\[
 \frac{\sigma^{-C\eta^4}R^\epsilon(\sigma'u)^\epsilon}
      {R^\epsilon(\sigma s)^\epsilon}
 \le \sigma^{\epsilon\eta/4-C\eta^4}.
\]
This is the first factor in \eqref{eq:hccontractions}.

If \(\sigma'=R^{-1/d}\), the first sum over \(Q\) in
\eqref{eq:hcoverviewrefinement} contains the original
functions \(f_\theta\).
Use the estimate proved for \(\sigma=R^{-1/d}\) before
Case~1, with \(s=u\) and amplitude
\(\alpha/\sigma^{-C\eta^4}\).
It bounds the first term in \eqref{eq:hcoverviewrefinement}
by the RHS of \eqref{eq:hcgeneralconvolution}, with both
occurrences of \(R^\epsilon\) replaced by
\(\sigma^{-C\eta^4}R^\eta\le R^{C\eta}\).
Since \(s\ge\sigma\ge R^{-1/d}\) and \(C\eta\ll\epsilon\),
we have
\(\frac{R^{C\eta}}{R^\epsilon(\sigma s)^\epsilon}
\le R^{-\epsilon/C}\).

For the second sum over \(Q\) in \eqref{eq:hcoverviewrefinement},
regardless of \(\sigma'\), we have
\(\sigma u'\le\sigma^{\eta^5/2}\sigma s<\sigma s/2\).
Apply \eqref{eq:hcgroupedconvolution} inductively at scales \((\sigma,u')\) and amplitude
\(\alpha/\sigma^{-2\delta}\).
This bounds the second term in \eqref{eq:hcoverviewrefinement}
by the RHS of \eqref{eq:hcgeneralconvolution}, with both
occurrences of \(R^\epsilon\) replaced by
\(\sigma^{-2\delta}R^\epsilon(\sigma u')^\epsilon\), and
\[
 \frac{\sigma^{-2\delta}R^\epsilon(\sigma u')^\epsilon}
      {R^\epsilon(\sigma s)^\epsilon}
 \le \sigma^{\epsilon\eta^5/2-2\delta}.
\]
This is the second factor in \eqref{eq:hccontractions}.

\medskip
We now sum the estimates of Cases~1 and~2 over
\(\sigma\le u\le s\) to prove \eqref{eq:hcgroupedconvolution}
at the current \(\sigma s\).
Our choices of \(\eta,\delta\) make all three exponents in
\eqref{eq:hccontractions} positive.
Since \(\sigma\le A^{-1}\), the positive powers of \(\sigma\)
in \eqref{eq:hccontractions} and the gain \(R^{-\epsilon/C}\)
absorb the dyadic sum over \(u\) and the implicit constants
when \(A\) and then \(R\) are sufficiently large.
The sum of the outer coefficients and each amplitude
denominator are therefore bounded by
\(R^\epsilon(\sigma s)^\epsilon\).
This proves \eqref{eq:hcgroupedconvolution} at the current
\(R\) by induction on \(\sigma s\).

To prove \eqref{eq:hcrestrictedconvolution} at the same \(R\),
apply Lemma~\ref{lem:hcseed} with \(\epsilon/2\) in place
of \(\epsilon\).
For each \(s\), apply \eqref{eq:hcgroupedconvolution} to
the sum over \(\rho,Q\) on the RHS of the lemma, with
\(\sigma=A^{-1}\) and amplitude \(\alpha/A^{\epsilon/2}\).
Including the factor \(A^{\epsilon/2}\) from
Lemma~\ref{lem:hcseed}, this bounds the fixed \(s\) term
by the RHS of \eqref{eq:hcgeneralconvolution}, with both
occurrences of \(R^\epsilon\) replaced by
\[
 A^{\epsilon/2}R^\epsilon(A^{-1}s)^\epsilon
 =R^\epsilon A^{-\epsilon/2}s^\epsilon.
\]
Since \(A^{-\epsilon/2}s^\epsilon\le1\), each averaging box in
these sums satisfies the selection condition in
\eqref{eq:hcrestrictedconvolution}.
Summing over \(s\) and using
\(A^{-\epsilon/2}\sum_{A^{-1}\le s\le1}s^\epsilon
\lesssim_\epsilon A^{-\epsilon/2}\)
therefore proves \eqref{eq:hcrestrictedconvolution} at \(R\).
This completes the induction on \(R\), proving
\eqref{eq:hcgroupedconvolution} and
\eqref{eq:hcrestrictedconvolution} for the current \(m\).

To obtain Proposition~\ref{prop:hcstable} and
\eqref{eq:hcstablefunction}, use
\eqref{eq:hcgroupedconvolution} and
\eqref{eq:hcrestrictedconvolution} with \(\epsilon/2\) in
place of \(\epsilon\), and apply
Lemma~\ref{lem:hcstep2averagecomparison} to their
right hand sides.
The additional loss, including fixed constants, may be
taken to be \(R^\eta\) in both the outer coefficient and
the amplitude denominator.
Since \(R\sigma s\ge R^{1-2/d}\) and \(\eta\ll\epsilon\),
\[
 R^\eta(R\sigma s)^{\epsilon/2}\le(R\sigma s)^\epsilon,
 \qquad R^{\epsilon/2+\eta}\le R^\epsilon.
\]
The first inequality gives Proposition~\ref{prop:hcstable};
the second gives \eqref{eq:hcstablefunction} for
\(\operatorname{supp}\widehat f
\subset\Gamma_{j,1/K}^{d+1}(R)\), \(0\le j\le m\).
Hence \(C_K^m(R)\lessapprox1\), and together with
\eqref{eq:hcrestrictedconvolution} this completes Step~1 in the proof of
Proposition~\ref{prop:hcstablefunction}.
\end{proof}

\section{\texorpdfstring{\except{toc}{From curves to higher dimensional cones over moment curves: }}{From curves to higher dimensional cones over moment curves: }Proof of Proposition~\ref{prop:hcstablefunction}, Step 2}
\label{sec:hcgeneral}

\begin{proof}[Proof of Proposition~\ref{prop:hcstep2positive}]
Let \(f\) be a Schwartz function with
\(\operatorname{supp}\widehat f
\subset\Gamma_{j,1/K,\ell}^{d+1}(R)\).
We choose small parameters \(0<\delta\ll\eta\ll\epsilon\)
as in Section~\ref{subsec:hcoverviewiteration} and induct
on the dyadic product \(\sigma s\), keeping \(R,j,\ell\) fixed.

We first deal with the case \(\sigma=R^{-1/d}\), where the
pieces \(f_\tau\) are the fine pieces \(f_\theta\).
By \eqref{eq:selaverage},
\[
 \text{LHS of Proposition~\ref{prop:hcstep2positive}}
 \lesssim
 \sum_{\rho\in\mathcal J(s)}
 \int_{\sum_{\theta\subset\rho}|f_\theta|^2
                    \gtrsim[\alpha/N(s)]^2}
 \Bigl(\sum_{\theta\subset\rho}|f_\theta|^2\Bigr)^{p_{d-m}/2}.
\]
For each interval \(\rho=[t_\rho,t_\rho+s]\in\mathcal J(s)\), set
\(g(x)=f_\rho(L_{t_\rho,s}^{-T}x)\).
The coefficient bounds in \eqref{eq:hcstep2support} give
\(\operatorname{supp}\widehat g\subset\Gamma_m^{d+1}(Rs^d)\),
up to fixed enlargement, with fine pieces
\(f_\theta(L_{t_\rho,s}^{-T}x)\), \(\theta\subset\rho\).
Apply \eqref{eq:hcpositiveconvolution} to \(g\) with parameter
\(Rs^d\), amplitude \(\alpha/N(s)\), and loss \(\eta\).

Returning to the original coordinates,
\eqref{eq:hcalllayerrescaling} gives the convolution weights
\(w_{B_{\tau_l}}\) and boxes \(B\) in
\eqref{eq:hcgeneralconvolution}, with
\(R^{-1/d}\le q_1\le\cdots\le q_{d-m-2}\le u\le s\).
For a returned parent \(\rho'\subset\rho\) of width \(u\),
Lemma~\ref{lem:regular} gives \(N(s)N_\rho(u)\lesssim N(u)\).
Since \((Rs^d)^\eta\le R^\eta\), it follows that
\(\frac{\alpha}{(Rs^d)^\eta N(s)N_\rho(u)}
\gtrsim\frac{\alpha}{R^\eta N(u)}\).
Thus summing over \(\rho\) and extending \(u\le s\) to \(u\le1\)
bounds the LHS of Proposition~\ref{prop:hcstep2positive}
by the RHS of \eqref{eq:hcgeneralconvolution},
with both occurrences of \(R^\epsilon\) replaced by \(R^\eta\).
Since \(\sigma=R^{-1/d}\) and \(s\ge\sigma\),
\[
 \frac{R^\eta}{R^\epsilon(\sigma s)^\epsilon}
 \le R^{-\epsilon/C}
\]
for \(\eta\) sufficiently small in terms of \(\epsilon\).
This proves Proposition~\ref{prop:hcstep2positive}
at \(\sigma=R^{-1/d}\).

We now consider the case \(\sigma>R^{-1/d}\), so that
\(\sigma\le a\ell\ll\ell\).
Since \(\Gamma_{j,1/K,\ell}^{d+1}(R)
\subset\widetilde\Gamma_{j,\ell}^{d+1}(R)\),
Corollary~\ref{cor:hccoarse}, applied to \(f\) with loss
\(\delta\), bounds the LHS of Proposition~\ref{prop:hcstep2positive}
by the RHS of \eqref{eq:hcoverviewcoarse}.
Fix \(\sigma\le u\le s\).
We estimate the corresponding term on the RHS of
\eqref{eq:hcoverviewcoarse}, namely
\begin{equation}\label{eq:hcstep2fixedu}
 \sigma^{-\delta}\sum_{\rho\in\mathcal J(u)}
 \sum_{\substack{
 B\parallel\frac{(uR^{1/d})^j}{(\sigma R^{1/d})^2}U_{\rho,R}\\
 B\in\mathcal G_\rho(\alpha\sigma^\delta)}}
 |B|\biggl\|\Bigl(
 \sum_{\substack{\tau\subset\rho\\|\tau|=\sigma}}
 |f_\tau|^2\Bigr)^{1/2}
 \biggr\|_{L^{p_{d-m}/2}_\#(B)}^{p_{d-m}}.
\end{equation}
As in the proof of \eqref{eq:hcgroupedconvolution}, we distinguish
\(u<\sigma^{1-\eta}\), where we use \eqref{eq:hcoverviewnarrow}
and then apply \eqref{eq:hcrestrictedconvolution} after Taylor
rescaling \(f_\tau\), from \(u\ge\sigma^{1-\eta}\), where
Lemma~\ref{lem:hcrefinement} gives either \((\sigma',u)\) with \(\sigma'<\sigma\) or \((\sigma,u')\) with \(u'<u\).

\medskip
\noindent\textit{Case 1: \(u<\sigma^{1-\eta}\).}
As in \eqref{eq:hcoverviewnarrow}, we can bound
\eqref{eq:hcstep2fixedu} by
\[
 \lesssim \sigma^{-C\eta}
 \sum_{\rho\in\mathcal J(u)}
 \sum_{\substack{\tau\subset\rho\\|\tau|=\sigma}}
 \int_{|f_\tau|\gtrsim\alpha\sigma^\delta/[N(u)N_\rho(\sigma)]}
 |f_\tau|^{p_{d-m}}.
\]
Since \(\operatorname{supp}
\widehat{f_\tau(L_{t_\tau,\sigma}^{-T}\,\cdot)}
\subset\Gamma_{j,1/K}^{d+1}(R\sigma^d)\) for
\(\tau\in\mathcal J(\sigma)\), we can apply
\eqref{eq:hcrestrictedconvolution}, proved in Step~1 for the
current \(m\), to this function with parameter \(R\sigma^d\)
and amplitude \(\alpha\sigma^\delta/[N(u)N_\rho(\sigma)]\).

Returning to the original coordinates, \eqref{eq:hcalllayerrescaling},
with \(\sigma\) in place of \(s\) and \(v\) in place of \(u\),
gives the convolution weights \(w_{B_{\tau_l}}\) and boxes \(B\)
in \eqref{eq:hcgeneralconvolution}, with
\(R^{-1/d}\le q_1\le\cdots\le q_{d-m-2}\le v\le\sigma\).
For a parent \(\rho'\subset\tau\) with \(|\rho'|=v\), applying
Lemma~\ref{lem:regular} twice gives
\(N(u)N_\rho(\sigma)N_\tau(v)\lesssim N(v)\).
Hence the amplitude in the selection condition satisfies
\[
 \frac{\alpha\sigma^\delta}
 {(R\sigma^d)^\epsilon N(u)N_\rho(\sigma)N_\tau(v)}
 \gtrsim
 \frac{\alpha}{\sigma^{-C\eta}(R\sigma^d)^\epsilon N(v)}.
\]
Summing over \(\rho,\tau\) and extending \(v\le\sigma\) to
\(v\le1\) therefore bounds \eqref{eq:hcstep2fixedu} by the RHS of
\eqref{eq:hcgeneralconvolution}, with both occurrences of
\(R^\epsilon\) replaced by \(\sigma^{-C\eta}(R\sigma^d)^\epsilon\).

Since \(s\ge\sigma\),
\[
 \frac{\sigma^{-C\eta}(R\sigma^d)^\epsilon}
      {R^\epsilon(\sigma s)^\epsilon}
 \le \sigma^{(d-2)\epsilon-C\eta}.
\]
This is the third factor in \eqref{eq:hccontractions} and
finishes the fixed \(u\) estimate in Case~1.

\medskip
\noindent\textit{Case 2: \(u\ge\sigma^{1-\eta}\).}
For each interval \(\rho=[t_\rho,t_\rho+u]\in\mathcal J(u)\), we have
\(\operatorname{supp}
\widehat{f_\rho(L_{t_\rho,u}^{-T}\,\cdot)}
\subset\widetilde\Gamma_{j,\min\{\ell/u,1\}}^{d+1}(Ru^d)\),
as in the proof of Corollary~\ref{cor:hccoarse}.
Since \(\sigma/u\le\sigma^\eta\) and
\(\sigma/u\ll\min\{\ell/u,1\}\), we can apply
Lemma~\ref{lem:hcrefinement} to this function with cap width
\(\sigma/u\) and amplitude \(\alpha\sigma^\delta/N(u)\).
After applying Lemma~\ref{lem:hcrefinement} and returning to
the original coordinates, we denote the new scales by \((\sigma',u')\).
Summing over \(\rho\) and applying dyadic pigeonholing as in
\eqref{eq:hcoverviewrefinement}, we can bound
\eqref{eq:hcstep2fixedu} by
\begin{equation}\label{eq:hcstep2refinement}
\begin{aligned}
 &\lesssim \sigma^{-C\eta^4}
 \sum_{\rho\in\mathcal J(u)}
 \sum_{\substack{Q\parallel\rho^*\\
 Q\in\mathcal G_\rho(\frac{\alpha}{\sigma^{-C\eta^4}})}}
 |Q|\biggl\|\Bigl(
 \sum_{\substack{\tau'\subset\rho\\|\tau'|=\sigma'}}
 |f_{\tau'}|^2\Bigr)^{1/2}
 \biggr\|_{L^{p_{d-m}/2}_\#(Q)}^{p_{d-m}}
 \\
 &\quad+\sigma^{-2\delta}
 \sum_{\rho'\in\mathcal J(u')}
 \sum_{\substack{Q\parallel(\rho')^*\\
 Q\in\mathcal G_{\rho'}(\frac{\alpha}{\sigma^{-2\delta}})}}
 |Q|\biggl\|\Bigl(
 \sum_{\substack{\tau\subset\rho'\\|\tau|=\sigma}}
 |f_\tau|^2\Bigr)^{1/2}
 \biggr\|_{L^{p_{d-m}/2}_\#(Q)}^{p_{d-m}},
\end{aligned}
\end{equation}
where \(R^{-1/d}\le\sigma'\le
\max\{R^{-1/d},\sigma^{1+\eta/4}\}\) and
\(\sigma\le u'\le u\sigma^{\eta^5/2}\).
In particular, \(\sigma'\le\sigma\le a\ell\), so both pairs
\((\sigma',u)\) and \((\sigma,u')\) lie in the scale range of
Proposition~\ref{prop:hcstep2positive}, with \(R,j,\ell\) fixed.

For the first sum over \(Q\) in \eqref{eq:hcstep2refinement},
suppose that \(\sigma'>R^{-1/d}\). Then
\(\sigma'u\le\sigma^{\eta/4}\sigma s<\sigma s\).
Apply the induction hypothesis of Proposition~\ref{prop:hcstep2positive}
at scales \((\sigma',u)\) and amplitude
\(\alpha/\sigma^{-C\eta^4}\).
This bounds the first term in \eqref{eq:hcstep2refinement}
by the RHS of \eqref{eq:hcgeneralconvolution}, with both
occurrences of \(R^\epsilon\) replaced by
\(\sigma^{-C\eta^4}R^\epsilon(\sigma'u)^\epsilon\), and
\[
 \frac{\sigma^{-C\eta^4}R^\epsilon(\sigma'u)^\epsilon}
      {R^\epsilon(\sigma s)^\epsilon}
 \le \sigma^{\epsilon\eta/4-C\eta^4}.
\]
This is the first factor in \eqref{eq:hccontractions}.

If \(\sigma'=R^{-1/d}\), the first sum over \(Q\) in
\eqref{eq:hcstep2refinement} contains the fine pieces \(f_\theta\).
Use the \(\sigma=R^{-1/d}\) calculation at the beginning of
this proof, with \(s=u\) and amplitude
\(\alpha/\sigma^{-C\eta^4}\).
It bounds the first term in \eqref{eq:hcstep2refinement}
by the RHS of \eqref{eq:hcgeneralconvolution}, with both
occurrences of \(R^\epsilon\) replaced by
\(\sigma^{-C\eta^4}R^\eta\le R^{C\eta}\).
Since \(s\ge\sigma\ge R^{-1/d}\) and \(C\eta\ll\epsilon\),
we have
\(\frac{R^{C\eta}}{R^\epsilon(\sigma s)^\epsilon}
\le R^{-\epsilon/C}\).

For the second sum over \(Q\) in \eqref{eq:hcstep2refinement},
regardless of \(\sigma'\), we have
\(\sigma u'\le\sigma^{\eta^5/2}\sigma s<\sigma s\).
Apply the induction hypothesis of Proposition~\ref{prop:hcstep2positive}
at scales \((\sigma,u')\) and amplitude
\(\alpha/\sigma^{-2\delta}\).
This bounds the second term in \eqref{eq:hcstep2refinement}
by the RHS of \eqref{eq:hcgeneralconvolution}, with both
occurrences of \(R^\epsilon\) replaced by
\(\sigma^{-2\delta}R^\epsilon(\sigma u')^\epsilon\), and
\[
 \frac{\sigma^{-2\delta}R^\epsilon(\sigma u')^\epsilon}
      {R^\epsilon(\sigma s)^\epsilon}
 \le \sigma^{\epsilon\eta^5/2-2\delta}.
\]
This is the second factor in \eqref{eq:hccontractions}.

We now sum the estimates of Cases~1 and~2 over
\(\sigma\le u\le s\) to prove Proposition~\ref{prop:hcstep2positive}
at the current \(\sigma s\).
Our choices of \(\eta,\delta\) make all three exponents in
\eqref{eq:hccontractions} positive.
Since \(\sigma\le a\ell\le A^{-1}\), the positive powers of
\(\sigma\) in \eqref{eq:hccontractions} and the gain
\(R^{-\epsilon/C}\) absorb the dyadic sum over \(u\) and the
implicit constants when \(A\) and then \(R\) are sufficiently large.
The amplitude denominators in \eqref{eq:hcgeneralconvolution}
and the sum of the outer coefficients are therefore bounded by
\(R^\epsilon(\sigma s)^\epsilon\).
This closes the induction on \(\sigma s\) and proves
Proposition~\ref{prop:hcstep2positive}.
\end{proof}

\begin{proof}[Completion of Proposition~\ref{prop:hcstablefunction}, Step~2]
Proposition~\ref{prop:hcstep2positive} completes the proof of
Proposition~\ref{prop:hcstep2reducedtail} as described in
Section~\ref{subsec:hcoverviewcoefficients}.
Combining Proposition~\ref{prop:hcstep2reducedtail} with
\eqref{eq:hcstep2tailreduction} and
\eqref{eq:hcstep2cutoffreturn} gives
\eqref{eq:hcgeneralconvolution} for functions with Fourier
support in \(\Gamma_j^{d+1}(R)\).
Lemma~\ref{lem:hcstep2averagecomparison} then gives
\eqref{eq:hcstablefunction}, completing the induction step in \(j\).
Proceeding through \(j=0,\ldots,m\), we obtain
\eqref{eq:hcgeneralconvolution} and \(C^m(R)\lessapprox1\),
proving Step~2.
\end{proof}

\section{The base case \(m=d-2\)}
\label{sec:hcbase}

We prove Proposition~\ref{prop:hcbase}. Here \(m=d-2\), so
\(p_{d-m}=4\).
We first establish the inequalities in Lemmas~\ref{lem:hcproducer}
and~\ref{lem:hcrefinement} for \(m=d-2\), and then follow
Steps~1 and~2 in the proof of Proposition~\ref{prop:hcstablefunction}
as in Sections~\ref{sec:hcstable} and~\ref{sec:hcgeneral}.

\begin{proof}[Proof of Proposition~\ref{prop:hcbase}]
To prove \(C_K^{d-2}(R)\lessapprox1\), it suffices to establish
the estimates of Lemmas~\ref{lem:hcproducer} and~\ref{lem:hcrefinement}
for \(m=d-2\).
Taking \(\ell=1\) in Lemma~\ref{lem:hcrefinement} gives
Lemma~\ref{lem:hcoverviewrefinement}, after which the argument in
Sections~\ref{sec:hcinitial}--\ref{subsec:hcoverviewiteration} applies.

We first prove the estimate of Lemma~\ref{lem:hcproducer}
for \(m=d-2\). Let \(g\) be a Schwartz function with
\(\operatorname{supp}\widehat g\subset\Gamma_{d-2}^{d+1}(r)\).
For each \(r^{-1/d}<\sigma\le1\), put
\(u=r^{-1/d}\sigma^{-1}\) and define
\[
 h_\tau=\sum_{\theta\subset\tau}|g_\theta|^2*
 \check\eta^d_{\theta,\sigma},
 \qquad \tau\in\mathcal J(u).
\]
Then \eqref{eq:hcpositiveparts}, with \(m=d-2\), gives
\[
\begin{aligned}
 \int_{\sum_\theta|g_\theta|^2\gtrsim\alpha^2}
 \Bigl(\sum_\theta|g_\theta|^2\Bigr)^2
 &\lesssim_\epsilon r^\epsilon
 \int_{|\sum_\theta|g_\theta|^2*
 \check\eta^d_{\theta,\le r^{-1/d}}|
 \gtrsim(\alpha/r^\epsilon)^2}
 \Bigl|\sum_\theta|g_\theta|^2*
 \check\eta^d_{\theta,\le r^{-1/d}}\Bigr|^2\\
 &\quad+r^\epsilon
 \sum_{r^{-1/d}<\sigma\le1}
 \int_{|\sum_{\tau\in\mathcal J(u)}h_\tau|
 \gtrsim(\alpha/r^\epsilon)^2}
 \Bigl|\sum_{\tau\in\mathcal J(u)}h_\tau\Bigr|^2.
\end{aligned}
\]
As in \eqref{eq:hclowaverage}, the first term is bounded by
the \(u=1\) term on the RHS of Lemma~\ref{lem:hcproducer}.
We therefore estimate each integral in the \(\sigma\) sum.

Fix \(\sigma\). As in the proof of \eqref{eq:hchighpieces},
\(\lvert h_\tau\rvert
\lesssim(\sum_{\theta\subset\tau}|g_\theta|^2)*w_{B_\tau}\),
where, since \(r^{i/d}\sigma^{i-d}=ru^{d-i}\),
\[
 B_\tau=\bigl\{x\in\R^{d+1}: |x\cdot\phi_d^{(i)}(t_\tau)|\le ru^{d-i}
 \text{ for }0\le i\le d\bigr\}.
\]
To apply \(L^2\) orthogonality to
\(\sum_{\tau\in\mathcal J(u)}h_\tau\), we first show that
the Fourier supports of the \(h_\tau\) have bounded overlap.

Write \(\xi=\sum_{i=0}^d\lambda_i(t;\xi)\phi_d^{(i)}(t)\).
Taylor expansion at \(t_\tau\), using
\eqref{eq:hcpositiveshell}, gives, for
\(\xi\in\operatorname{supp}\widehat h_\tau\),
\[
 |\lambda_i(t_\tau;\xi)|\lesssim r^{-1}u^{i-d}
 \quad(0\le i\le d),\qquad
 \max_{0\le i<d}ru^{d-i}|\lambda_i(t_\tau;\xi)|\sim1.
\]
For fixed \(\xi\), the coefficient \(\lambda_d(t;\xi)\)
is a polynomial of degree at most \(d\), whose derivatives
are the preceding Taylor coefficients up to signs.
The set \(\{t:|\lambda_d(t;\xi)|\lesssim r^{-1}\}\)
has \(O(1)\) components.
If a component \(J\) contains \(t_\tau\) with
\(\xi\in\operatorname{supp}\widehat h_\tau\),
the polynomial derivative bounds give
\(|\lambda_i(t_\tau;\xi)|\lesssim r^{-1}|J|^{i-d}\), \(i<d\).
The annular condition therefore forces \(|J|\lesssim u\).
Each component contains only \(O(1)\) left endpoints \(t_\tau\),
so the Fourier supports of the \(h_\tau\) have bounded overlap.

To obtain the condition selecting the averaging boxes in
Lemma~\ref{lem:hcproducer}, we decompose \(h_\tau\) using
the cutoffs \(\psi_{\tau,T}\) constructed before
\eqref{eq:packetBernstein}:
\[
 h_\tau=
 \sum_{T\parallel(ur^{1/d})^{d-2}U_{\tau,r}}
 \psi_{\tau,T}h_\tau,
 \qquad
 \sum_T\psi_{\tau,T}=1.
\]
Discard the wave packets satisfying
\(\|\psi_{\tau,T}^{1/2}h_\tau\|_\infty
\ll(\alpha/r^\epsilon)^2/N(u)\), and denote the boxes of
the remaining wave packets by \(\mathbb T_\tau\).
There are at most \(N(u)\) nonzero functions \(h_\tau\), so
\(\sum_T\psi_{\tau,T}^{1/2}\lesssim1\) gives
\[
 \Bigl|\sum_\tau\sum_{T\notin\mathbb T_\tau}
 \psi_{\tau,T}h_\tau\Bigr|
 \le
 \sum_\tau\sum_{T\notin\mathbb T_\tau}
 \psi_{\tau,T}^{1/2}
 \|\psi_{\tau,T}^{1/2}h_\tau\|_\infty
 \le\frac12(\alpha/r^\epsilon)^2.
\]
Since \(T\) has width \(ru^{d-i}\) in direction \(i\),
the wave packet cutoff \(\psi_{\tau,T}\) has Fourier width
\(\ll r^{-1}u^{i-d}\) in that direction.
The Fourier supports of
\(\sum_{T\in\mathbb T_\tau}\psi_{\tau,T}h_\tau\)
therefore also have bounded overlap as \(\tau\) varies.
By this bounded overlapping property and Cauchy--Schwarz
in \(T\), we obtain
\[
\begin{aligned}
 \int_{|\sum_\tau h_\tau|\gtrsim(\alpha/r^\epsilon)^2}
 \Bigl|\sum_\tau h_\tau\Bigr|^2
 &\lesssim
 \sum_\tau\int
 \Bigl|\sum_{T\in\mathbb T_\tau}
 \psi_{\tau,T}h_\tau\Bigr|^2\\
 &\lesssim
 \sum_\tau\sum_{T\in\mathbb T_\tau}
 \int\psi_{\tau,T}^{3/2}|h_\tau|^2
 \lesssim
 \sum_\tau\sum_{T\in\mathbb T_\tau}
 |T|\|\psi_{\tau,T}^{1/2}h_\tau\|_\infty^2.
\end{aligned}
\]
Here Cauchy--Schwarz uses
\(\sum_T\psi_{\tau,T}^{1/2}\lesssim1\), and the last
inequality uses
\(\int\psi_{\tau,T}^{1/2}\lesssim|T|\).

Applying \eqref{eq:packetBernstein} to \(h_\tau\) with
exponent one and using the pointwise bound for \(h_\tau\),
we obtain
\[
 \|\psi_{\tau,T}^{1/2}h_\tau\|_\infty
 \lesssim
 \biggl\|\Bigl(\sum_{\theta\subset\tau}|g_\theta|^2\Bigr)
 *w_{B_\tau}\biggr\|_{L^1_\#(T)}
 \lesssim
 \biggl\|\Bigl(\sum_{\theta\subset\tau}|g_\theta|^2\Bigr)^{1/2}
 \biggr\|_{L^2_\#(T)}^2.
\]
Since \(T\) is a translate of \(B_\tau\), we have \(w_{B_\tau}*W_T\lesssim W_T\).
The last inequality follows from \eqref{eq:weightedYoung}.

Since \((\alpha/r^\epsilon)^2/N(u)
\ge[\alpha/(r^\epsilon N(u))]^2\), every
\(T\in\mathbb T_\tau\) belongs to
\(\mathcal G_\tau(\alpha/r^\epsilon)\).
Substituting the wave packet bound into the preceding \(L^2\)
estimate gives
\begin{equation}\label{eq:hcbasepacket}
 \int_{|\sum_\tau h_\tau|\gtrsim(\alpha/r^\epsilon)^2}
 \Bigl|\sum_\tau h_\tau\Bigr|^2
 \lesssim
 \sum_{\tau\in\mathcal J(u)}
 \sum_{\substack{
 T\parallel(ur^{1/d})^{d-2}U_{\tau,r}\\
 T\in\mathcal G_\tau(\alpha/r^\epsilon)}}
 |T|
 \biggl\|\Bigl(\sum_{\theta\subset\tau}|g_\theta|^2\Bigr)^{1/2}
 \biggr\|_{L^2_\#(T)}^4.
\end{equation}
Summing \eqref{eq:hcbasepacket} over \(\sigma\), with
\(u=r^{-1/d}\sigma^{-1}\), and combining this with the
\(u=1\) bound for the first term in
\eqref{eq:hcpositiveparts} proves the estimate of
Lemma~\ref{lem:hcproducer} for \(m=d-2\).

We next prove the estimate of Lemma~\ref{lem:hcrefinement}
for \(m=d-2\).
Fix \(0\le j\le d-2\) and \(R^{-1/d}<\sigma\ll\ell\le1\).
Let \(f\) be a Schwartz function with
\(\operatorname{supp}\widehat f
\subset\widetilde\Gamma_{j,\ell}^{d+1}(R)\).
Since \(p_{d-m}/4=1\), we now refine the caps directly in the
weighted average over \(B\) of \(\sum_{|\tau|=\sigma}|f_\tau|^2\).

For \(B\parallel[0,R^{j/d}\sigma^{-2}]^{d+1}\), replace \(W_B\)
by a comparable weight whose Fourier transform is supported
in \(CB^*\), and expand the weighted average over \(B\) of
\(\sum_{|\tau|=\sigma}|f_\tau|^2\).
The calculation leading to \eqref{eq:hcpairseparation},
with \(u=1\), shows that only pairs of caps at scale \(R^{-1/d}\) whose
intervals have left endpoints satisfying \(|t-t'|\lesssim R^{-1/d}+\sigma^2\)
contribute to the average.
Choose a dyadic \(\sigma'\sim\max\{R^{-1/d},\sigma^2\}\).
Only neighbouring \(\sigma'\) intervals interact, so
Cauchy--Schwarz gives
\begin{equation}\label{eq:hcbaseaveragerefinement}
 \biggl\|\Bigl(\sum_{|\tau|=\sigma}|f_\tau|^2\Bigr)^{1/2}
 \biggr\|_{L^2_\#(B)}^2
 \lesssim
 \biggl\|\Bigl(\sum_{|\tau'|=\sigma'}|f_{\tau'}|^2\Bigr)^{1/2}
 \biggr\|_{L^2_\#(B)}^2.
\end{equation}

For \(B\in\mathcal G_{[0,1]}(\alpha)\), the LHS of
\eqref{eq:hcbaseaveragerefinement} is at least \(\alpha^2\).
Apply \eqref{eq:selaverage} with
\(F=\sum_{|\tau'|=\sigma'}|f_{\tau'}|^2\) and \(t=\alpha^2\).
Together with \eqref{eq:hcbaseaveragerefinement} and the locally constant property on unit cubes, this gives
\begin{equation}\label{eq:hcbasequarticrefinement}
 \sum_{\substack{
 B\parallel[0,R^{j/d}\sigma^{-2}]^{d+1}\\
 B\in\mathcal G_{[0,1]}(\alpha)}}
 |B|
 \biggl\|\Bigl(\sum_{|\tau|=\sigma}|f_\tau|^2\Bigr)^{1/2}
 \biggr\|_{L^2_\#(B)}^4
 \lesssim
 \sum_{\substack{
 Q\parallel[0,1]^{d+1}\\
 Q\in\mathcal G_{[0,1]}(\alpha)}}
 |Q|
 \biggl\|\Bigl(\sum_{|\tau'|=\sigma'}|f_{\tau'}|^2\Bigr)^{1/2}
 \biggr\|_{L^2_\#(Q)}^4.
\end{equation}
This proves the estimate of Lemma~\ref{lem:hcrefinement}
for \(m=d-2\).
In particular, for
\(\operatorname{supp}\widehat f\subset\Gamma_{j,1/K}^{d+1}(R)\),
taking \(\ell=1\) gives Lemma~\ref{lem:hcoverviewrefinement}
for \(m=d-2\).

The argument in
Sections~\ref{sec:hcinitial}--\ref{subsec:hcoverviewiteration}
therefore gives Proposition~\ref{prop:hcstable} and
\eqref{eq:hcstablefunction} for \(m=d-2\), and hence
\(C_K^{d-2}(R)\lessapprox1\), completing Step~1.

We now prove \(C^{d-2}(R)\lessapprox1\).
As in the reduction in Section~\ref{subsec:hcoverviewcoefficients},
it suffices to establish Lemma~\ref{lem:hcstep2cutoffcomparison}
and Proposition~\ref{prop:hcstep2reducedtail} for \(m=d-2\),
replacing the RHS of \eqref{eq:hcgeneralconvolution} by the
RHS of \eqref{eq:hcstablefunction} in both statements.
Lemma~\ref{lem:hcstep2cutoffcomparison} holds also for \(m=d-2\),
since its proof uses no inductive hypothesis.

In proving Proposition~\ref{prop:hcstep2reducedtail} for \(m=d-2\),
we use \eqref{eq:hcstablefunction} in place of the hypothesis
\eqref{eq:hcrestrictedconvolution} for functions with
\(\operatorname{supp}\widehat f
\subset\Gamma_{j,1/K}^{d+1}(R)\).
This follows from the estimates of Lemmas~\ref{lem:hcproducer}
and~\ref{lem:hcrefinement} proved in Step~1, by the argument in
Sections~\ref{sec:hcinitial}--\ref{subsec:hcoverviewiteration}.
Following the proof of Proposition~\ref{prop:hcstep2reducedtail},
it therefore remains to prove the estimate of
Proposition~\ref{prop:hcstep2positive} for \(m=d-2\), with the RHS of
\eqref{eq:hcgeneralconvolution} replaced by the RHS of
\eqref{eq:hcstablefunction}.

Using \eqref{eq:hcstablefunction} at \(m=d-2\) in place of the
hypothesis \eqref{eq:hcrestrictedconvolution}, we follow the
proof in Section~\ref{sec:hcgeneral} with the RHS of
\eqref{eq:hcgeneralconvolution} replaced throughout by the
RHS of \eqref{eq:hcstablefunction}.
The estimates of Lemmas~\ref{lem:hcproducer}
and~\ref{lem:hcrefinement} established above, together with the
rescaling and counting arguments in Section~\ref{sec:hcgeneral},
give the same bounds for the outer coefficients and amplitude
denominators.
The comparisons involving the convolution weights in
\eqref{eq:hcgeneralconvolution} are therefore unnecessary,
and the induction on \(\sigma s\) closes as before.
This proves the estimate of Proposition~\ref{prop:hcstep2positive}
for \(m=d-2\).
The reduction in Section~\ref{subsec:hcoverviewcoefficients}
then gives \(C^{d-2}(R)\lessapprox1\), completing the proof of
Proposition~\ref{prop:hcbase}.
\end{proof}

Starting with Proposition~\ref{prop:hcbase},
Proposition~\ref{prop:hcstablefunction} proves
\eqref{eq:hcgeneralconvolution} and \(C^m(R)\lessapprox1\)
successively for \(m=d-3,\ldots,0\). Taking \(m=0\) proves
Proposition~\ref{prop:highercone}.
The induction on dimension starts with the cone over the moment curve in \(\mathbb R^3\), the \(d=2\) case of Theorem~\ref{thm:main-cone}. Proposition~\ref{prop:conetocurve}
then gives the corresponding estimate for the moment curve in \(\mathbb R^3\).
Now suppose the estimates for moment curves hold in dimensions \(2\le k\le d\).
Proposition~\ref{prop:highercone} gives the estimate for the cone over the moment curve in
\(\mathbb R^{d+1}\), and Proposition~\ref{prop:conetocurve} gives the estimate
for the moment curve in the same dimension. This completes the induction and
proves Theorems~\ref{thm:main} and~\ref{thm:main-cone}.

\addtocontents{toc}{\protect\tocpartheading{III}{Applications}\protect\tocsectionindent{3pc}}
\section{Small cap decoupling inequalities}\label{chap:smallcap}\label{sec:commontools}

In this section, we prove Theorem~\ref{thm:SC} using the amplitude dependent
wave envelope type estimate of Theorem~\ref{thm:main}.
Throughout this section, $\gamma$ denotes a small cap in the
decomposition of $\mathcal M^n(R^\beta,R)$ from the introduction,
and $\supp\widehat f_\gamma\subset\gamma$.
For each dyadic scale $R^{-\beta}\le s\le1$, we use $\nu$ to
denote the cap obtained by adding the vertical segment
$\{0\}^{n-1}\times[-R^{-1},R^{-1}]$ to the caps in $\Theta(s)$.
Write $|\nu|=s$ in such a case, and $t_\nu$
for the left endpoint of the interval corresponding to the cap. The cap is contained in a rectangle with dimension
\[
 s\times s^2\times\cdots\times s^{n-1}
 \times\max\{s^n,R^{-1}\}.
\]
Let us state the main proposition in this section.

\begin{proposition}\label{prop:densityweak}
Let $n\ge4$, $R\ge1$, $\beta\ge1/n$, and $p=n(n-1)+2/\beta$.
If $\|f_\gamma\|_\infty\le1$ for every small cap $\gamma$, then
for every $\e>0$ and $\alpha>0$,
\begin{equation}
 \alpha^p\biggl|\biggl\{\Bigl|\sum_\gamma f_\gamma\Bigr|>
 \alpha\biggr\}\biggr|
 \lesssim_\e
 R^{1+\beta(n(n-1)/2-1)+\beta\e}
 \sum_\gamma\|f_\gamma\|_2^2.
 \label{eq:smallcapweak}
\end{equation}
\end{proposition}

After the standard reduction, one may see that this proposition implies Theorem~\ref{thm:SC}. For the completeness, we give the details at the end of this section.

To prove Proposition \ref{prop:densityweak}, we consider two cases. The first case $n\ge5$ or $n=4$, $\beta\le1$ is proved in
Subsection~\ref{sec:smallcaps}, and the remaining case $n=4$,
$\beta>1$ is proved in Subsection~\ref{sec:quarticlarge}.

\subsection{Key inequalities}
We state three lemmas to prove the proposition.
Recall from Subsection~\ref{subsec:mainresults} that $\Theta(s)$
denotes the family of caps in $\R^n$ of dimensions
$s\times s^2\times\cdots\times s^n$, and $\Xi(s)$ the family
of caps in $\R^{d+1}$ of dimensions
$1\times s\times\cdots\times s^d$.

\begin{lemma}[Canonical decoupling]\label{lem:canonicaldecoupling} Let $n,d \geq 2$.
 For $2 \leq p \leq n(n+1)$ we have
\begin{equation}
 \Bigl\|\sum_{\tau\in\Theta(s)} f_\tau\Bigr\|_p
 \lesssim_\e s^{-\e}
 \bigl(\sum_{\tau\in\Theta(s)}\|f_\tau\|_p^2\bigr)^{1/2}
 \label{eq:classicalcurve}
\end{equation}
for any function $f_{\gamma}$ with $\supp\widehat f_\tau\subset\tau$. For $2 \leq p \leq d(d+1)$,
\begin{equation}
 \Bigl\|\sum_{\tau\in\Xi(s)} f_\tau\Bigr\|_p
 \lesssim_\e s^{-\e}
 \bigl(\sum_{\tau\in\Xi(s)}\|f_\tau\|_p^2\bigr)^{1/2}
 \label{eq:classicalcone}
\end{equation}
for any function $f_{\tau}$ with $\supp\widehat f_\tau\subset\tau$.
\end{lemma}

The first estimate is the decoupling of Bourgain-Demeter-Guth
\cite{BDG}; decoupling for the cone over the moment curve in \(\mathbb R^{d+1}\) follows from
\cite[Corollary 2.13]{KLO}. 
Recall that $p_n=\frac{n(n+1)}2+1$.
The following is a consequence of Theorem~\ref{thm:main}.

\begin{lemma}\label{lem:bodyAD}
Let $n\ge3$ and let $\supp\widehat f_\theta\subset\theta$ for
$\theta\in\Theta(R^{-1/n})$. For every $\e>0$ and $\alpha>0$,
\begin{equation}
\begin{aligned}
 &\alpha^{p_n}
 \biggl|\biggl\{\Bigl|\sum_\theta f_\theta\Bigr|>\alpha\biggr\}\biggr| \le C_\e R^\e
 \sum_{R^{-1/n}\le s\le1}\sum_{\tau\in\Theta(s)}
 \int_{\{\sum_{\theta\subset\tau}|f_\theta|^2>
 C_\e^{-1} R^{-\e}\alpha^2s^2\}}
 \bigl(\sum_{\theta\subset\tau}|f_\theta|^2\bigr)^{p_n/2}.
\end{aligned}
 \label{eq:bodyAD}
\end{equation}
\end{lemma}

\begin{proof}
Since $\#\tau\lesssim s^{-1}$, every good averaging box in
Theorem~\ref{thm:main}, applied with loss $\e/2$, satisfies
\[
 \frac1{|V|}\int
 \bigl(\sum_{\theta\subset\tau}|f_\theta|^2\bigr)^{p_n/4}W_V
 \gtrsim_\e (R^{-\e}\alpha^2s^2)^{p_n/4}.
\]
For $C_\e$ sufficiently large, the part where
$\sum_{\theta\subset\tau}|f_\theta|^2
\le C_\e^{-1}R^{-\e}\alpha^2s^2$ contributes at most half
this average. In other words,
\[
 \frac1{|V|}\int_{\{\sum_{\theta\subset\tau}|f_\theta|^2 \leq 
 C_\e^{-1} R^{-\e}\alpha^2s^2\}}
 \bigl(\sum_{\theta\subset\tau}|f_\theta|^2\bigr)^{p_n/4}W_V
 \leq C_{\e}^{-p_n/4} (R^{-\e}\alpha^2s^2)^{p_n/4}.
\]
So we can ignore this term. For the good averaging box,
Cauchy--Schwarz on the remaining set therefore gives
\[
 \frac1{|V|}\biggl[\int
 \bigl(\sum_{\theta\subset\tau}|f_\theta|^2\bigr)^{p_n/4}W_V\biggr]^2
 \le 4
 \int_{\{\sum_{\theta\subset\tau}|f_\theta|^2>
 C_\e^{-1}R^{-\e}\alpha^2s^2\}}
 \bigl(\sum_{\theta\subset\tau}|f_\theta|^2\bigr)^{p_n/2}W_V.
\]
Summing over the translates of $V_{\tau,R}$ and using
Theorem~\ref{thm:main} proves the lemma.
\end{proof}

\begin{lemma}[Key lemma]\label{lem:smallcapdensity}
Let $n\ge4$, $R^{-\beta}\le\sigma\le s\le1$, and  $\nu$ satisfy $|\nu|=s$. Suppose $\|f_\gamma\|_\infty\le1$.
For  $\nu'\subset\nu$ with $|\nu'|=\sigma$, we have the following estimates.
\begin{enumerate}
\renewcommand{\labelenumi}{\textup{(\roman{enumi})}}
\item \textup{Pointwise bound.}
\begin{equation}
 \sum_{\substack{\nu'\subset \nu\\|\nu'|=\sigma}}\Bigl|\sum_{\gamma\subset\nu'}f_\gamma\Bigr|^2
 \lesssim R^{2\beta}\sigma s.
 \label{eq:densitypointwise}
\end{equation}
\item \textup{Moment bound.}
\begin{equation}
 \int\Bigl(\sum_{\substack{\nu'\subset \nu\\|\nu'|=\sigma}}\Bigl|\sum_{\gamma\subset\nu'}f_\gamma\Bigr|^2
 \Bigr)^{n(n-1)/2}
 \lesssim_\e R^{\beta\e}(R^\beta s)^{n(n-1)/2-1}
 \sum_{\gamma\subset \nu}\|f_\gamma\|_2^2.
 \label{eq:densitycanonical}
\end{equation}
\item \textup{Restricted moment bound.} For a sufficiently large constant $C$,
\begin{equation}
\begin{aligned}
 &\int_{\{\sum_{\nu'\subset \nu,\ |\nu'|=\sigma}\lvert\sum_{\gamma\subset\nu'}f_\gamma\rvert^2>CR^\beta s\}}
 \Bigl(\sum_{\substack{\nu'\subset \nu\\|\nu'|=\sigma}}\Bigl|\sum_{\gamma\subset\nu'}f_\gamma\Bigr|^2
 \Bigr)^{n(n-1)/2}\\
 &\quad\quad\quad\quad\quad\quad\quad\lesssim_\e
 R^{\beta\e}(R^\beta s)^{n(n-1)/2-1} \Bigl(\frac{\sigma}{s}\Bigr)^{n(n-1)/4}
 \sum_{\gamma\subset \nu}\|f_\gamma\|_2^2.
\end{aligned}
 \label{eq:densitybroad}
\end{equation}
\end{enumerate}
\end{lemma}
The first two items are straightforward, while the third is the key to Proposition~\ref{prop:densityweak}.
We implement the iterated broad narrow decomposition of each
$|\sum_{\gamma\subset\nu'}f_\gamma|^2$ through the  Whitney decomposition.
The separated pairs form the broad terms, while $E_{\nu'}$ below is the
terminal narrow remainder at scale $R^{-\beta}$.
For each fixed $x_n$, a suitable change of variables in the first $n-1$ coordinates places the Fourier supports of the broad terms in a small neighborhood of a cone generated by a moment curve.
We estimate these terms by decoupling for the cone generated by \(\phi_{n-2}\) in \(\mathbb R^{n-1}\), as in \cite[Lemma~7.1]{Glidewell}.
Related bilinear changes of variables appear in
\cite{BourgainZeta,BDSurfaces,BDCurves} and \cite[Section~2]{OhBilinear}.
The terminal narrow contribution, summed over $\nu'$, is $O(R^\beta s)$ and can be absorbed
on the restricted set in \eqref{eq:densitybroad}.
Combining these estimates proves the third item.

\begin{proof}
For item~\textup{(i)}, the triangle inequality and
$\|f_\gamma\|_\infty\le1$ give
$|\sum_{\gamma\subset\nu'}f_\gamma|\le\sum_{\gamma\subset\nu'}|f_\gamma|\lesssim R^\beta\sigma$.
There are $O(s/\sigma)$ such caps $\nu'\subset\nu$, so
\[
 \sum_{\substack{\nu'\subset\nu\\|\nu'|=\sigma}}
 \Bigl|\sum_{\gamma\subset\nu'}f_\gamma\Bigr|^2
 \lesssim \frac{s}{\sigma}(R^\beta\sigma)^2
 =R^{2\beta}\sigma s.
\]

To prove item~\textup{(ii)}, let $r_{\nu'}$ be independent random
signs taking the values $\pm1$ with equal probability.
Khintchine's inequality gives
\[
 \int\Bigl(\sum_{\substack{\nu'\subset\nu\\|\nu'|=\sigma}}
 \Bigl|\sum_{\gamma\subset\nu'}f_\gamma\Bigr|^2
 \Bigr)^{n(n-1)/2}
 \lesssim \mathbb E\Bigl\|
 \sum_{\substack{\nu'\subset\nu\\|\nu'|=\sigma}}
 r_{\nu'}\sum_{\gamma\subset\nu'}f_\gamma
 \Bigr\|_{n(n-1)}^{n(n-1)}.
\]
We next apply Lemma~\ref{lem:canonicaldecoupling}
to the first $n-1$ coordinates,
after rescaling the  cap $\nu$. Then applying H\"older's inequality gives
\[
\begin{aligned}
 \mathbb E\Bigl\|
 \sum_{\substack{\nu'\subset\nu\\|\nu'|=\sigma}}
 r_{\nu'}\sum_{\gamma\subset\nu'}f_\gamma
 \Bigr\|_{n(n-1)}^{n(n-1)}
 &\lesssim_\e R^{\beta\e}
 \bigl(\sum_{\gamma\subset\nu}\|f_\gamma\|_{n(n-1)}^2
 \bigr)^{n(n-1)/2}\\
 &\lesssim R^{\beta\e}(R^\beta s)^{n(n-1)/2-1}
 \sum_{\gamma\subset\nu}\|f_\gamma\|_{n(n-1)}^{n(n-1)}.
\end{aligned}
\]
Since $\|f_\gamma\|_\infty\le1$, this proves
\eqref{eq:densitycanonical}. \medskip

To prove item~\textup{(iii)}, use a disjoint Whitney decomposition, and let $\mathcal W_\omega(\nu')$ denote its pairs of subcaps of $\nu'$ satisfying $|\nu_1|=|\nu_2|=\omega$ and $t_{\nu_1}-t_{\nu_2}\sim\omega$. For a sufficiently large fixed constant $C_0$, summing over dyadic $\omega$ gives
\[
 \Bigl|\sum_{\gamma\subset\nu'}f_\gamma\Bigr|^2
 =2\operatorname{Re}\sum_{C_0R^{-\beta}<\omega\le\sigma}
 \sum_{(\nu_1,\nu_2)\in\mathcal W_\omega(\nu')}
 \bigl(\sum_{\gamma\subset\nu_1}f_\gamma\bigr)
 \overline{\sum_{\gamma\subset\nu_2}f_\gamma}+E_{\nu'}.
\]
Here the remainder is defined by
\[
 E_{\nu'}:=\sum_{\gamma\subset\nu'}|f_\gamma|^2
 +2\operatorname{Re}\sum_{R^{-\beta}\le\omega\le C_0R^{-\beta}}
 \sum_{(\nu_1,\nu_2)\in\mathcal W_\omega(\nu')}
 \bigl(\sum_{\gamma\subset\nu_1}f_\gamma\bigr)
 \overline{\sum_{\gamma\subset\nu_2}f_\gamma}.
\]
At the finest scale, $\mathcal W_{R^{-\beta}}(\nu')$ also includes the remaining pairs of distinct neighbouring small caps, each counted once with $t_{\nu_1}>t_{\nu_2}$. Writing $t_\gamma$ for the left endpoint of the interval indexing $\gamma$, each $\gamma$ has only $O(C_0)$ such neighbours, so
\[
 \sum_{\substack{\nu'\subset\nu\\|\nu'|=\sigma}}|E_{\nu'}|
 \lesssim\sum_{\substack{\gamma,\gamma'\subset\nu\\
 |t_\gamma-t_{\gamma'}|\lesssim C_0R^{-\beta}}}
 |f_\gamma||f_{\gamma'}|
 \lesssim\sum_{\gamma\subset\nu}|f_\gamma|^2
 \lesssim R^\beta s.
\]
Write $\mathcal W_\omega(\nu)$ for the union of $\mathcal W_\omega(\nu')$ over all $\nu'\subset\nu$ with $|\nu'|=\sigma$. On the set in
\eqref{eq:densitybroad}, we consequently have
\[
 \sum_{\substack{\nu'\subset \nu\\|\nu'|=\sigma}}\Bigl|\sum_{\gamma\subset\nu'}f_\gamma\Bigr|^2
 \lesssim
 \Bigl|\sum_{C_0R^{-\beta}<\omega\le\sigma}
 \sum_{(\nu_1,\nu_2)\in\mathcal W_\omega(\nu)}
 \bigl(\sum_{\gamma\subset\nu_1}f_\gamma\bigr)\overline{\sum_{\gamma\subset\nu_2}f_\gamma}\Bigr|.
\]
Applying this bound to the left hand side of \eqref{eq:densitybroad}, and then using the triangle inequality, we obtain
\[
\begin{aligned}
&\int_{\{\sum_{\nu'\subset\nu,\,|\nu'|=\sigma}|\sum_{\gamma\subset\nu'}f_\gamma|^2>CR^\beta s\}}
\Bigl(\sum_{\substack{\nu'\subset\nu\\|\nu'|=\sigma}}
 \Bigl|\sum_{\gamma\subset\nu'}f_\gamma\Bigr|^2\Bigr)^{n(n-1)/2}\\
&\lesssim\biggl(\sum_{C_0R^{-\beta}<\omega\le\sigma}
 \Bigl\|\sum_{(\nu_1,\nu_2)\in\mathcal W_\omega(\nu)}
 \bigl(\sum_{\gamma\subset\nu_1}f_\gamma\bigr)
 \overline{\sum_{\gamma\subset\nu_2}f_\gamma}\Bigr\|_{n(n-1)/2}\biggr)^{n(n-1)/2}.
\end{aligned}
\]
We now estimate each summand. At each $\omega$, we may split these pairs into finitely many
collections in which the points $(t_{\nu_1}+t_{\nu_2})/2$ are $\omega$ separated. Each cap occurs
in boundedly many pairs.

We check the Fourier support needed to apply decoupling for the cone generated by \(\phi_{n-2}\) in \(\mathbb R^{n-1}\). If
$t,t'$ are the left endpoints of the intervals indexing small caps grouped in
$\nu_1,\nu_2$, respectively, and $t_0=(t+t')/2$, then
\begin{equation}
 \gamma_{n-1}(t)-\gamma_{n-1}(t')
 =(t-t')\gamma_{n-1}'(t_0)
 +\sum_{k\ge1}\frac{(t-t')^{2k+1}}{2^{2k}(2k+1)!}
 \gamma_{n-1}^{(2k+1)}(t_0).
 \label{eq:smallcapsecant}
\end{equation}
Since $t-t'\sim\omega$, division by $\omega$ gives a radial
coefficient comparable to $\sim1$. In the $(\ell+1)$st derivative coordinate, the errors from the widths of the two small caps satisfy
\[
 \omega^{-1}\sum_{j=1}^{\ell+1}
 R^{-j\beta}(C\omega)^{\ell+1-j}\lesssim\omega^\ell,
 \qquad 0\le\ell\le n-2.
\]
In the radial coordinate this error is $O((R^\beta\omega)^{-1})$,
so our choice of $C_0$ keeps the radial coefficient comparable to $\sim1$.
Recentering at $(t_{\nu_1}+t_{\nu_2})/2$ preserves the widths
$1,\omega,\ldots,\omega^{n-2}$, and note that
\[
 \gamma_{n-1}'(t)
 =\operatorname{diag}(1!,2!,\ldots,(n-1)!)\phi_{n-2}(t).
\]
Thus, after this fixed diagonal change of coordinates and division
by $\omega$, the projected Fourier supports of the products
\[
 \bigl(\sum_{\gamma\subset\nu_1}f_\gamma\bigr)
 \overline{\sum_{\gamma\subset\nu_2}f_\gamma}
\]
lie in caps of aperture comparable to $\omega$ in
$\mathbb R^{n-1}$.

The exponent $n(n-1)/2$ is admissible in \eqref{eq:classicalcone},
since $n(n-1)/2\le(n-2)(n-1)$ for $n\ge4$.
Applying the cylindrical version of the decoupling estimate for the cone generated by $\phi_{n-2}$ in $\mathbb R^{n-1}$ (see \eqref{eq:classicalcone}) and then H\"older's
inequality over the $O(s/\omega)$ pairs $(\nu_1,\nu_2)\in\mathcal W_\omega(\nu)$ gives
\[
\begin{aligned}
 &\Bigl\|\sum_{(\nu_1,\nu_2)\in\mathcal W_\omega(\nu)}
 \bigl(\sum_{\gamma\subset\nu_1}f_\gamma\bigr)\overline{\sum_{\gamma\subset\nu_2}f_\gamma}\Bigr\|_{n(n-1)/2}^{n(n-1)/2}\\
 &\quad\lesssim_\e
 R^{\beta\e}(s/\omega)^{n(n-1)/4-1}
 \sum_{(\nu_1,\nu_2)\in\mathcal W_\omega(\nu)}
 \Bigl\|\bigl(\sum_{\gamma\subset\nu_1}f_\gamma\bigr)\overline{\sum_{\gamma\subset\nu_2}f_\gamma}\Bigr\|_{n(n-1)/2}^{n(n-1)/2}.
\end{aligned}
\]
The cylindrical version of the decoupling estimate for the moment curve $\gamma_{n-1}$ (see equation~\ref{eq:classicalcurve}), now applied at exponent $n(n-1)$ on each
cap of parameter length $\omega$, also gives
\[
 \Bigl\|\bigl(\sum_{\gamma\subset\nu_1}f_\gamma\bigr)\overline{\sum_{\gamma\subset\nu_2}f_\gamma}\Bigr\|_{n(n-1)/2}^{n(n-1)/2}
 \lesssim_\e R^{\beta\e}(R^\beta\omega)^{n(n-1)/2-1}
 \sum_{\gamma\subset \nu_1\cup \nu_2}\|f_\gamma\|_2^2.
\]
Each small cap occurs in boundedly many pairs. Combining the last two estimates, taking roots, and summing in $\omega$ gives
\[
\begin{aligned}
&\sum_{C_0R^{-\beta}<\omega\le\sigma}
\Bigl\|
\sum_{(\nu_1,\nu_2)\in\mathcal W_\omega(\nu)}
\bigl(\sum_{\gamma\subset\nu_1}f_\gamma\bigr)
\overline{\sum_{\gamma\subset\nu_2}f_\gamma}
\Bigr\|_{n(n-1)/2}\\
&\qquad\lesssim_\e
R^{\frac{\beta(n(n-1)-2+2\e)}{n(n-1)}}
s^{\frac12-\frac{2}{n(n-1)}}\sigma^{1/2}
\bigl(\sum_{\gamma\subset\nu}\|f_\gamma\|_2^2
\bigr)^{\frac{2}{n(n-1)}},
\end{aligned}
\]
where we used
$\sum_{\omega\le\sigma}\omega^{1/2}\lesssim\sigma^{1/2}$.
By the reduction above, raising this bound to the power $n(n-1)/2$ proves \eqref{eq:densitybroad}.
\end{proof}
Let us now prove Proposition \ref{prop:densityweak}. Consider two cases:
\begin{equation}
    n \geq 5, \; \beta \ge 1/n, \; \textrm{or} \; n=4, \; 1/4 \le \beta \le 1, \qquad \textrm{or} \qquad n=4, \; \beta>1. 
\end{equation}

\subsection{Case 1: \texorpdfstring{$n\ge5$, $\beta\ge1/n$, or $n=4$, $1/4\le\beta\le1$}{n >= 5, beta >= 1/n, or n = 4, 1/4 <= beta <= 1}}\label{sec:smallcaps}

We prove Proposition~\ref{prop:densityweak} for such $n,\beta$.   
Fix the final loss $\e>0$. In the estimates below we use an auxiliary
parameter $\delta=\delta(\e)>0$, to be chosen sufficiently small at the end.
Group the small caps into $\theta\in\Theta(R^{-1/n})$, writing
$f_\theta=\sum_{\gamma\subset\theta}f_\gamma$.
Applying Lemma~\ref{lem:bodyAD} with $\delta/n$, and using $\beta\ge1/n$,
\begin{equation}
\begin{aligned}
 \alpha^p\biggl|\biggl\{\Bigl|\sum_\gamma f_\gamma\Bigr|>
 \alpha\biggr\}\biggr|\lesssim_\delta R^{\beta\delta}\alpha^{p-p_n}
 \sum_{R^{-1/n}\le s\le1}\sum_{\tau\in\Theta(s)}
 \int_{\{\sum_{\theta\subset\tau}|f_\theta|^2>
 C_\delta^{-1} R^{-\beta\delta}\alpha^2s^2\}}
 \bigl(\sum_{\theta\subset\tau}|f_\theta|^2\bigr)^{p_n/2}.
\end{aligned}
 \label{eq:smallcapADreduction}
\end{equation}
Since $p_n/2<n(n-1)/2$, on the indicated set
\[
 \bigl(\sum_{\theta\subset\tau}|f_\theta|^2\bigr)^{p_n/2}
 \lesssim_\delta
 (R^{-\beta\delta}\alpha^2s^2)^{p_n/2-n(n-1)/2}
 \bigl(\sum_{\theta\subset\tau}|f_\theta|^2\bigr)^{n(n-1)/2}.
\]
Substituting this bound into \eqref{eq:smallcapADreduction}, we obtain
\begin{align*}
 \alpha^p\biggl|\biggl\{\Bigl|\sum_\gamma f_\gamma\Bigr|>\alpha\biggr\}\biggr|
 &\lesssim_\delta R^{C\beta\delta}\alpha^{2/\beta}
 \sum_{R^{-1/n}\le s\le1}s^{n+1-n(n-1)/2}\\
 &\qquad\times\sum_{\tau\in\Theta(s)}
 \int_{\{\sum_{\theta\subset\tau}|f_\theta|^2>
 C_\delta^{-1}R^{-\beta\delta}\alpha^2s^2\}}
 \bigl(\sum_{\theta\subset\tau}|f_\theta|^2\bigr)^{n(n-1)/2}.
\end{align*}
We estimate the integral using Lemma~\ref{lem:smallcapdensity}
with $\nu=\tau$, $\nu'=\theta$, and $\sigma=R^{-1/n}$, so that
$\sum_{\gamma\subset\nu'}f_\gamma=f_\theta$, distinguishing two ranges of $s$. Consider two ranges
\begin{equation}
    s\le CR^{\beta(1+\delta)}\alpha^{-2}, \;\;\;\;\;\; s>CR^{\beta(1+\delta)}\alpha^{-2}.
\end{equation}
This split compares the threshold $C_\delta^{-1}R^{-\beta\delta}\alpha^2s^2$
in \eqref{eq:smallcapADreduction} with the threshold $CR^\beta s$
in \eqref{eq:densitybroad}.
We use \eqref{eq:densitycanonical} in the first range and
\eqref{eq:densitybroad} in the second.
We split the entire right hand side of the preceding level set
estimate into the two terms
\begin{align*}
 \mathrm{I}
 &:=R^{C\beta\delta}\alpha^{2/\beta}
 \sum_{\substack{R^{-1/n}\le s\le1\\
 s\le CR^{\beta(1+\delta)}\alpha^{-2}}}
 s^{n+1-n(n-1)/2}\sum_{\tau\in\Theta(s)}
 \int_{\{\sum_{\theta\subset\tau}|f_\theta|^2>
 C_\delta^{-1}R^{-\beta\delta}\alpha^2s^2\}}
 \bigl(\sum_{\theta\subset\tau}|f_\theta|^2\bigr)^{n(n-1)/2},\\
 \mathrm{II}
 &:=R^{C\beta\delta}\alpha^{2/\beta}
 \sum_{\substack{R^{-1/n}\le s\le1\\
 s>CR^{\beta(1+\delta)}\alpha^{-2}}}
 s^{n+1-n(n-1)/2}\sum_{\tau\in\Theta(s)}
 \int_{\{\sum_{\theta\subset\tau}|f_\theta|^2>
 C_\delta^{-1}R^{-\beta\delta}\alpha^2s^2\}}
 \bigl(\sum_{\theta\subset\tau}|f_\theta|^2\bigr)^{n(n-1)/2}.
\end{align*}
It suffices to prove the following two bounds, with $\delta$ chosen
sufficiently small in terms of $\e$:
\begin{align*}
 \mathrm{I}+\mathrm{II}
 &\lesssim_\e R^{1+\beta(n(n-1)/2-1)+\beta\e}
 \sum_\gamma\|f_\gamma\|_2^2,.
\end{align*}
Indeed, the preceding level set estimate gives
\[
 \alpha^p\biggl|\biggl\{\Bigl|\sum_\gamma f_\gamma\Bigr|>\alpha\biggr\}\biggr|
 \lesssim_\delta\mathrm{I}+\mathrm{II}.
\]
Assuming the two bounds above, adding them therefore yields
\[
 \alpha^p\biggl|\biggl\{\Bigl|\sum_\gamma f_\gamma\Bigr|>\alpha\biggr\}\biggr|
 \lesssim_\e R^{1+\beta(n(n-1)/2-1)+\beta\e}
 \sum_\gamma\|f_\gamma\|_2^2,
\]
which is \eqref{eq:smallcapweak}. We now estimate $\mathrm{I}$ and
$\mathrm{II}$ separately.

\subsubsection{Estimate for \texorpdfstring{$\mathrm{I}$}{I}}

Fix a dyadic scale
$s\le CR^{\beta(1+\delta)}\alpha^{-2}$ in its defining sum.
For every $\tau\in\Theta(s)$, enlarge the integration region to
$\mathbb R^n$ and apply \eqref{eq:densitycanonical} with $\nu=\tau$,
$\nu'=\theta$, and $\sigma=R^{-1/n}$ to obtain
\begin{align*}
 &\int_{\{\sum_{\theta\subset\tau}|f_\theta|^2>
 C_\delta^{-1}R^{-\beta\delta}\alpha^2s^2\}}
 \bigl(\sum_{\theta\subset\tau}|f_\theta|^2\bigr)^{n(n-1)/2}\lesssim_\delta
 R^{\beta\delta}(R^\beta s)^{n(n-1)/2-1}
 \sum_{\gamma\subset\tau}\|f_\gamma\|_2^2.
\end{align*}
Substituting this estimate directly into the definition of $\mathrm{I}$,
we obtain
\[
 \mathrm{I}\lesssim_\delta
 R^{\beta(n(n-1)/2-1)+C\beta\delta}
 \sum_{\substack{R^{-1/n}\le s\le1\\
 s\le CR^{\beta(1+\delta)}\alpha^{-2}}}
 \alpha^{2/\beta}s^n
 \sum_\gamma\|f_\gamma\|_2^2.
\]
For every scale in this sum, the condition
$s\le CR^{\beta(1+\delta)}\alpha^{-2}$ gives
$\alpha^2s\lesssim R^{\beta(1+\delta)}$. Hence
\[
 \alpha^{2/\beta}s^n
 =(\alpha^2s)^{1/\beta}s^{n-1/\beta}
 \lesssim R^{1+\delta}s^{n-1/\beta}
 \le R^{1+\delta},
\]
where the last inequality uses $s\le1$ and $\beta\ge1/n$.
Substituting this bound into the preceding estimate for $\mathrm{I}$
and summing over the $O(\log(2R))$ dyadic scales yields
\begin{align*}
 \mathrm{I}
 &\lesssim_\delta
 R^{1+\beta(n(n-1)/2-1)+C\beta\delta}
 \sum_{\substack{R^{-1/n}\le s\le1\\
 s\le CR^{\beta(1+\delta)}\alpha^{-2}}}
 \sum_\gamma\|f_\gamma\|_2^2\\
 &\lesssim_\e
 R^{1+\beta(n(n-1)/2-1)+\beta\e}
 \sum_\gamma\|f_\gamma\|_2^2.
\end{align*}
 This proves the desired bound
for $\mathrm{I}$.

\subsubsection{Estimate for \texorpdfstring{$\mathrm{II}$}{II}}

Fix a dyadic scale
$s>CR^{\beta(1+\delta)}\alpha^{-2}$ in its defining sum.
With the constant in this range condition chosen sufficiently large
compared with the threshold constant in \eqref{eq:densitybroad},
we have
\[
 \Bigl\{\sum_{\theta\subset\tau}|f_\theta|^2>
 C_\delta^{-1}R^{-\beta\delta}\alpha^2s^2\Bigr\}
 \subset
 \Bigl\{\sum_{\theta\subset\tau}|f_\theta|^2>CR^\beta s\Bigr\}.
\]
For every $\tau\in\Theta(s)$, we may therefore apply
\eqref{eq:densitybroad} with $\nu=\tau$, $\nu'=\theta$,
and $\sigma=R^{-1/n}$ to obtain
\begin{align*}
 &\int_{\{\sum_{\theta\subset\tau}|f_\theta|^2>
 C_\delta^{-1}R^{-\beta\delta}\alpha^2s^2\}}
 \bigl(\sum_{\theta\subset\tau}|f_\theta|^2\bigr)^{n(n-1)/2} \lesssim
 R^{\beta(n(n-1)/2-1)+\beta\delta}
 R^{-(n-1)/4}s^{n(n-1)/4-1}
 \sum_{\gamma\subset\tau}\|f_\gamma\|_2^2.
\end{align*}
Discard scales for which all integration regions in $\mathrm{II}$ are empty;
the sums below run only over the remaining scales.
Substituting this estimate directly into the definition of $\mathrm{II}$,
we obtain
\[
 \mathrm{II}\lesssim_\delta
 R^{\beta(n(n-1)/2-1)+C\beta\delta}
 \sum_{\substack{R^{-1/n}\le s\le1\\
 s>CR^{\beta(1+\delta)}\alpha^{-2}}}
 \alpha^{2/\beta}R^{-(n-1)/4}s^{n-n(n-1)/4}
 \sum_\gamma\|f_\gamma\|_2^2.
\]
 The inequality \eqref{eq:densitypointwise} gives
\[
 C_\delta^{-1}R^{-\beta\delta}\alpha^2s^2
 <\sum_{\theta\subset\tau}|f_\theta|^2
 \lesssim R^{2\beta}R^{-1/n}s.
\]
It follows that
\[ \alpha^{2/\beta}\lesssim_\delta
 R^{2+\delta}\Bigl(\frac{R^{-1/n}}s\Bigr)^{1/\beta}.
\]
Consequently, the remaining factor satisfies
\begin{align*}
 &\alpha^{2/\beta}R^{-(n-1)/4}s^{n-n(n-1)/4}\lesssim
 R^{1+\delta}\Bigl(\frac{R^{-1/n}}s\Bigr)^{1/\beta+n(n-5)/4}
 \le R^{1+\delta}.
\end{align*}
For the last inequality, $R^{-1/n}/s\le1$, and the exponent is
nonnegative: this is immediate for $n\ge5$, while for $n=4$ it is
$1/\beta-1\ge0$ because $\beta\le1$.
Summing over the $O(\log(2R))$ dyadic scales gives the required bound
\begin{align*}
 \mathrm{II}
\lesssim_\e R^{1+\beta(n(n-1)/2-1)+\beta\e}
 \sum_\gamma\|f_\gamma\|_2^2.
\end{align*}
Here we choose $\delta>0$ sufficiently small that $C\delta\le\e/2$
and use $\log(2R)\lesssim_\e R^{\beta\e/2}$.
\medskip

\subsection{Case 2: \texorpdfstring{$n=4$, $\beta>1$}{n = 4, beta > 1}}
\label{sec:quarticlarge}

It remains to prove Proposition~\ref{prop:densityweak} for $n=4$
and $\beta>1$. We first establish the following broad estimate,
then pass to the general case by a broad narrow analysis.

\begin{proposition}[Broad estimate]\label{prop:quarticbroad}
Let $n=4$, $\beta\ge1$, and let $K$ be a sufficiently large number. Let $\nu_1,\ldots,\nu_4$ be caps with
$|\nu_i|=K^{-1}$ such that $\mathrm{dist}(\nu_i,\nu_j)>K^{-1}$ for any $i \neq j$. Suppose $\|f_\gamma\|_\infty\le1$ for every
small cap $\gamma$. For $\alpha>0$, define the broad set
\[
 U_\alpha=
 \biggl\{x:\min_{1\le i\le4}
 \Bigl|\sum_{\gamma\subset\nu_i}f_\gamma(x)\Bigr|>\alpha\biggr\}.
\]
Then, for every $\e>0$,
\begin{equation}
 \alpha^{12+2/\beta}|U_\alpha|
 \lesssim_{\e,K}R^{5\beta+1+\beta\e}
 \sum_\gamma\|f_\gamma\|_2^2.
 \label{eq:quarticbroad}
\end{equation}
\end{proposition}

Assuming Proposition~\ref{prop:quarticbroad}, the remaining case of
Proposition~\ref{prop:densityweak} follows by a standard broad narrow analysis.
For completeness, we give this deduction at the end of this subsection.
We prepare the proof of Proposition~\ref{prop:quarticbroad}, keeping
the notation in its statement.
\medskip

Let $B_{s^{-1},R}$ be the box centered at the origin, with sides parallel
to the coordinate axes and dimensions
\[
 s^{-1}\times s^{-1}\times s^{-1}\times\min\{s^{-1},R\}.
\]
Let us start with the observation.

\begin{lemma}\label{lem:quarticendpointaverages}
For every $x\in\R^4$,
\begin{equation}
 \frac1{|B_{R^\beta,R}|}\int
 \sum_\gamma|f_\gamma(y)|^2W_{B_{R^\beta,R}}(x-y)\,dy
 \lesssim R^\beta.
 \label{eq:quarticfinestaverage}
\end{equation}
For every $x\in U_\alpha$,
\[
 \alpha^2\lesssim\frac1{|B_{K,R}|}\int
 \sum_{|\nu'|=K^{-1}}\Bigl|\sum_{\gamma\subset\nu'}f_\gamma(y)\Bigr|^2
 W_{B_{K,R}}(x-y)\,dy.
\]
\end{lemma}

\begin{proof}
The first bound follows from $\|f_\gamma\|_\infty\le1$,
$\#\{\gamma\}\lesssim R^\beta$, and
$\int W_{B_{R^\beta,R}}=|B_{R^\beta,R}|$.
For $R^{-\beta}\le s\le K^{-1}$ and $|\nu'|=s$, we have
\[
 \supp\widehat{\sum_{\gamma\subset\nu'}f_\gamma}
 \subset\gamma_4(t_{\nu'})+CB_{s^{-1},R}^*.
\]
The locally constant property, as in \cite[Lemma~2]{GMAD}, therefore gives
\begin{equation}
 \sum_{|\nu'|=s}\Bigl|\sum_{\gamma\subset\nu'}f_\gamma(x)\Bigr|^2
 \lesssim\frac1{|B_{s^{-1},R}|}\int
 \sum_{|\nu'|=s}\Bigl|\sum_{\gamma\subset\nu'}f_\gamma(y)\Bigr|^2
 W_{B_{s^{-1},R}}(x-y)\,dy.
 \label{eq:quarticmoderate}
\end{equation}
Taking $s=K^{-1}$  proves the
second bound, since
$|\sum_{\gamma\subset\nu_1}f_\gamma(x)|>\alpha$ on $U_\alpha$.
\end{proof}

Choose $C_*$ sufficiently large, and then $C_K$ sufficiently large
depending on $C_*$. Suppose that $\alpha>C_KR^{\beta/2}$.
By Lemma~\ref{lem:quarticendpointaverages}, the average at scale
$K^{-1}$ exceeds $C_*R^\beta$, whereas the average at scale
$R^{-\beta}$ is at most $C_*R^\beta$.
For dyadic $R^{-\beta}\le s<K^{-1}$, define
\[
\begin{split}
 U_{\alpha,s}:=\biggl\{&x\in U_\alpha:\;
 \frac1{|B_{s^{-1},R}|}\int
 \sum_{|\nu'|=s}\Bigl|\sum_{\gamma\subset\nu'}f_\gamma(y)\Bigr|^2
 W_{B_{s^{-1},R}}(x-y)\,dy\le C_*R^\beta,\\
 &\frac1{|B_{(s')^{-1},R}|}\int
 \sum_{|\nu'|=s'}\Bigl|\sum_{\gamma\subset\nu'}f_\gamma(y)\Bigr|^2
 W_{B_{(s')^{-1},R}}(x-y)\,dy>C_*R^\beta \;\; \text{for every dyadic }s<s'\le K^{-1}\biggr\}.
\end{split}
\]
Thus, $U_{\alpha,s}$ consists of the points where the weighted
average first drops to at most $C_*R^\beta$ at scale $s$,
as we pass from coarse to fine scales.
These sets partition $U_\alpha$:
\[
 U_\alpha=\bigsqcup_{\substack{R^{-\beta}\le s<K^{-1}\\s\text{ dyadic}}}U_{\alpha,s}.
\]
Using a fixed tiling by translates $B\parallel B_{s^{-1},R}$, define
\[
 \Lambda_s:=\bigcup_{\substack{B\parallel B_{s^{-1},R}\\B\cap U_{\alpha,s}\ne\varnothing}}B.
\]
For $R^{-\beta}\le s<K^{-1}$ and $x\in\Lambda_s$,
\begin{equation}
\begin{aligned}
 \frac1{|B_{s^{-1},R}|}\int
 \sum_{|\nu'|=s}\Bigl|\sum_{\gamma\subset\nu'}f_\gamma(y)\Bigr|^2
 W_{B_{s^{-1},R}}(x-y)\,dy&\lesssim C_*R^\beta,\\
 \frac1{|B_{(2s)^{-1},R}|}\int
 \sum_{|\nu'|=2s}\Bigl|\sum_{\gamma\subset\nu'}f_\gamma(y)\Bigr|^2
 W_{B_{(2s)^{-1},R}}(x-y)\,dy&\gtrsim C_*R^\beta.
\end{aligned}
 \label{eq:quarticstoppingaverages}
\end{equation}
Indeed, for $x\in\Lambda_s$, choose $x_0\in U_{\alpha,s}$ in the same tile
$B\parallel B_{s^{-1},R}$.
Since $W_{B_{t^{-1},R}}(x-y)\sim W_{B_{t^{-1},R}}(x_0-y)$ for
$t=s,2s$, the defining inequalities at $x_0$ give
\eqref{eq:quarticstoppingaverages} at $x$.

\begin{lemma}\label{lem:quarticstopping}
For $R^{-\beta}\le s<K^{-1}$ and every $\e>0$,
\begin{equation}
 |\Lambda_s|\lesssim_\e R^{\beta(-1+\e)}s^3
 \sum_\gamma\|f_\gamma\|_2^2,
 \qquad
 \alpha^2\lesssim R^\beta s^{-1}
 \quad\text{if }U_{\alpha,s}\ne\varnothing.
 \label{eq:quarticstop}
\end{equation}
\end{lemma}

\begin{lemma}[Multilinear restriction estimate]
\label{lem:quarticfourfold}
Let $0<\delta\le1$, and let $\nu_1,\ldots,\nu_4\subset[0,1]$
be intervals satisfying $\mathrm{dist}(\nu_i,\nu_j)>K^{-1}$ for $i\ne j$.
Suppose that $f_1,\ldots,f_4$ are Schwartz functions with
\[
 \supp\widehat f_i\subset
 \gamma_4(\nu_i)+C[-\delta,\delta]^4.
\]
Then
\begin{equation}
 \int_{\R^4}\prod_{i=1}^4|f_i(x)|^2\,dx
 \lesssim_K\delta^{12}\prod_{i=1}^4\|f_i\|_2^2.
 \label{eq:quarticrestriction}
\end{equation}
\end{lemma}

Assuming these two lemmas, we first prove Proposition~\ref{prop:quarticbroad}, and then return to their proofs.

\begin{proof}[Proof of Proposition~\ref{prop:quarticbroad}]

The case $\alpha\le C_KR^{\beta/2}$ is easy. By $L^2$ orthogonality,
\[
 \alpha^{12+2/\beta}|U_\alpha|
 \lesssim\alpha^{10+2/\beta}\sum_\gamma\|f_\gamma\|_2^2
 \lesssim_K R^{5\beta+1}\sum_\gamma\|f_\gamma\|_2^2.
\]
We may therefore assume $\alpha>C_KR^{\beta/2}$.
Fix $R^{-\beta}\le s<K^{-1}$ and $B\parallel B_{s^{-1},R}$ meeting $U_{\alpha,s}$.
Choose a Schwartz cutoff $\chi_B$ with $|\chi_B|\ge1$ on $B$,
$\supp\widehat\chi_B\subset CB_{s^{-1},R}^*$, and $|\chi_B|^2\lesssim W_B$.
The pieces $\chi_B\sum_{\gamma\subset\nu'}f_\gamma$, $|\nu'|=s$,
have bounded Fourier overlap in the first coordinate. Hence almost
orthogonality and the upper bound in \eqref{eq:quarticstoppingaverages} give
\[
 \Bigl\|\chi_B\sum_{\gamma\subset\nu_i}f_\gamma\Bigr\|_2^2
 \lesssim\sum_{\substack{\nu'\subset\nu_i\\|\nu'|=s}}
 \Bigl\|\chi_B\sum_{\gamma\subset\nu'}f_\gamma\Bigr\|_2^2
 \lesssim R^\beta|B|.
\]
On the Fourier supports of these localized functions,
$|\xi_4-\xi_1^4|\lesssim\max\{s,R^{-1}\}$.
To apply Lemma~\ref{lem:quarticfourfold}, partition
$\xi_4-\xi_1^4$ smoothly into intervals of length $s$ and write
$\chi_B\sum_{\gamma\subset\nu_i}f_\gamma=\sum_j g_{i,j}$.
There are $O(\max\{1,(Rs)^{-1}\})$ pieces, with bounded Fourier overlap, so
\[
 \sum_j\|g_{i,j}\|_2^2
 \lesssim\Bigl\|\chi_B\sum_{\gamma\subset\nu_i}f_\gamma\Bigr\|_2^2.
\]
After modulation in $x_4$, each $g_{i,j}$ has Fourier support in
$\gamma_4(\nu_i)+C[-s,s]^4$.
The triangle inequality, Lemma~\ref{lem:quarticfourfold} with $\delta=s$,
and Cauchy--Schwarz in $j$ therefore give
\begin{align*}
 \Bigl\|\prod_{i=1}^4\bigl(\chi_B\sum_{\gamma\subset\nu_i}f_\gamma\bigr)\Bigr\|_2
 &\le\sum_{j_1,\ldots,j_4}\Bigl\|\prod_{i=1}^4g_{i,j_i}\Bigr\|_2\\
 &\lesssim_K s^6\prod_{i=1}^4\sum_j\|g_{i,j}\|_2\\
 &\lesssim_K s^6\max\{1,(Rs)^{-1}\}^2
 \prod_{i=1}^4\Bigl\|\chi_B\sum_{\gamma\subset\nu_i}f_\gamma\Bigr\|_2.
\end{align*}
Squaring and using the preceding $L^2$ bound together with
$|B|=s^{-3}\max\{s,R^{-1}\}^{-1}$, we obtain
\begin{equation}
\begin{aligned}
 \int_B\prod_{i=1}^4\Bigl|\sum_{\gamma\subset\nu_i}f_\gamma\Bigr|^2
 &\lesssim_K s^{12}\max\{1,(Rs)^{-1}\}^4R^{4\beta}|B|^4\\
 &=R^{4\beta}\max\{1,(Rs)^{-1}\}|B|.
\end{aligned}
 \label{eq:quarticfourfold}
\end{equation}
Summing \eqref{eq:quarticfourfold} over the tiles in $\Lambda_s$
and using \eqref{eq:quarticstop}, we obtain
\begin{align*}
 \alpha^8|U_{\alpha,s}|
 &\lesssim_K R^{4\beta}\max\{1,(Rs)^{-1}\}|\Lambda_s|\\
 &\lesssim_{\e,K}R^{\beta(3+\e)}s^3
 \max\{1,(Rs)^{-1}\}\sum_\gamma\|f_\gamma\|_2^2.
\end{align*}
Multiplying by
$\alpha^{4+2/\beta}\lesssim(R^\beta/s)^{2+1/\beta}$, we have
\[
 \alpha^{12+2/\beta}|U_{\alpha,s}|
 \lesssim_{\e,K}R^{5\beta+1+\beta\e}
 s^{1-1/\beta}\max\{1,(Rs)^{-1}\}
 \sum_\gamma\|f_\gamma\|_2^2.
\]
The scale factor satisfies
\[
 s^{1-1/\beta}\max\{1,(Rs)^{-1}\}
 =\begin{cases}
 s^{1-1/\beta}\le1,&s\ge R^{-1},\\
 R^{-1}s^{-1/\beta}\le R^{-1}(R^{-\beta})^{-1/\beta}=1,
 &R^{-\beta}\le s<R^{-1}.
 \end{cases}
\]
Summing over the dyadic scales and relabeling $\e$ gives
\eqref{eq:quarticbroad}.
\end{proof}

\begin{proof}[Proof of Lemma~\ref{lem:quarticstopping}]
Applying \eqref{eq:densitybroad} with $n=4$, $|\nu|=1$, and
$\sigma=s$ gives
\begin{equation}
 \int_{\{\sum_{|\nu'|=s}\lvert\sum_{\gamma\subset\nu'}f_\gamma\rvert^2>CR^\beta\}}
 \Bigl(\sum_{|\nu'|=s}\Bigl|\sum_{\gamma\subset\nu'}f_\gamma\Bigr|^2\Bigr)^6
 \lesssim_\e R^{\beta(5+\e)}s^3
 \sum_\gamma\|f_\gamma\|_2^2,
 \label{eq:quarticsquare}
\end{equation}
for $R^{-\beta}\le s\le1$.
If $C$ is the constant in \eqref{eq:quarticsquare}, then
\[
 \frac1{|B_{(2s)^{-1},R}|}
 \int_{\{\sum_{|\nu'|=2s}\lvert\sum_{\gamma\subset\nu'}f_\gamma(y)\rvert^2\le CR^\beta\}}
 \sum_{|\nu'|=2s}\Bigl|\sum_{\gamma\subset\nu'}f_\gamma(y)\Bigr|^2
 W_{B_{(2s)^{-1},R}}(x-y)\,dy\le CR^\beta.
\]
By the choice of $C_*$, subtracting this contribution from the
second average in \eqref{eq:quarticstoppingaverages} gives
\[
 R^\beta\lesssim\frac1{|B_{(2s)^{-1},R}|}
 \int_{\{\sum_{|\nu'|=2s}\lvert\sum_{\gamma\subset\nu'}f_\gamma(y)\rvert^2>CR^\beta\}}
 \sum_{|\nu'|=2s}\Bigl|\sum_{\gamma\subset\nu'}f_\gamma(y)\Bigr|^2
 W_{B_{(2s)^{-1},R}}(x-y)\,dy,
 \qquad x\in\Lambda_s.
\]
Integrating over $\Lambda_s$, and then applying Young's inequality
and \eqref{eq:quarticsquare}, yields
\begin{align*}
 R^{6\beta}|\Lambda_s|
 &\lesssim\int_{\Lambda_s}\biggl[\frac1{|B_{(2s)^{-1},R}|}
 \int_{\{\sum_{|\nu'|=2s}\lvert\sum_{\gamma\subset\nu'}f_\gamma(y)\rvert^2>CR^\beta\}}
 \sum_{|\nu'|=2s}\Bigl|\sum_{\gamma\subset\nu'}f_\gamma(y)\Bigr|^2
 W_{B_{(2s)^{-1},R}}(x-y)\,dy\biggr]^6dx\\
 &\lesssim
 \int_{\{\sum_{|\nu'|=2s}\lvert\sum_{\gamma\subset\nu'}f_\gamma\rvert^2>CR^\beta\}}
 \Bigl(\sum_{|\nu'|=2s}\Bigl|\sum_{\gamma\subset\nu'}f_\gamma\Bigr|^2\Bigr)^6 \lesssim_\e R^{\beta(5+\e)}s^3\sum_\gamma\|f_\gamma\|_2^2.
\end{align*}
Finally, for $x\in U_{\alpha,s}$, Cauchy--Schwarz over the subcaps of
$\nu_1$, followed by \eqref{eq:quarticmoderate}, gives
\[
\begin{aligned}
 \alpha^2
 &<\Bigl|\sum_{\substack{\nu'\subset\nu_1\\|\nu'|=s}}
 \sum_{\gamma\subset\nu'}f_\gamma(x)\Bigr|^2\
 \lesssim s^{-1}\sum_{|\nu'|=s}
 \Bigl|\sum_{\gamma\subset\nu'}f_\gamma(x)\Bigr|^2
 \lesssim s^{-1}R^\beta.
\end{aligned}
\]
\end{proof}

\begin{proof}[Proof of Lemma~\ref{lem:quarticfourfold}]
By Bernstein's inequality, we may assume $\delta<c/K$ for a sufficiently
small absolute constant $c$. Write $t=\xi_1$ and
$v_j=\xi_j-t^j$ for $j=2,3,4$, and set
\[
 g_{i,v}(t)=\widehat f_i(t,t^2+v_2,t^3+v_3,t^4+v_4),
 \qquad Eg(x)=\int g(t)e^{2\pi i x\cdot\gamma_4(t)}\,dt.
\]
Then $g_{i,v}$ is supported in $\nu_i+[-C\delta,C\delta]$, these enlarged
intervals remain separated, and the $v$ support has volume $O(\delta^3)$.
Fourier inversion and Plancherel give
\[
 f_i(x)=\int e^{2\pi i(x_2v_2+x_3v_3+x_4v_4)}Eg_{i,v}(x)\,dv,
 \qquad \int\|g_{i,v}\|_2^2\,dv=\|f_i\|_2^2.
\]
For fixed shifts $v^{(1)},\ldots,v^{(4)}$, the change of variables
$\Phi(t_1,\ldots,t_4)=\sum_{i=1}^4\gamma_4(t_i)$ is injective on the product of these four enlarged intervals and satisfies
$|\det D\Phi|=24\prod_{i<j}|t_i-t_j|\gtrsim_K1$.
Thus Plancherel yields
\[
 \Bigl\|\prod_{i=1}^4Eg_{i,v^{(i)}}\Bigr\|_2
 \lesssim_K\prod_{i=1}^4\|g_{i,v^{(i)}}\|_2.
\]
Minkowski and Cauchy--Schwarz in $v$ now imply
\[
 \Bigl\|\prod_{i=1}^4f_i\Bigr\|_2
 \lesssim_K\prod_{i=1}^4\int\|g_{i,v}\|_2\,dv
 \lesssim_K\delta^6\prod_{i=1}^4\|f_i\|_2.
\]
Squaring proves \eqref{eq:quarticrestriction}.
\end{proof}

\begin{proof}[Proof of Proposition~\ref{prop:densityweak} for case 2]
Our proof uses the induction on $R^{\beta}$.
In the standard broad narrow analysis, we apply the rescaling to the narrow case. After applying the rescaling, we use the induction hypothesis with different $\beta$. To do so, we prove the stronger
statement
\begin{equation}
 \alpha^p
 \biggl|\biggl\{\Bigl|\sum_\gamma f_\gamma\Bigr|>\alpha\biggr\}\biggr|
 \le C_\e R^{-1+\beta(p-7)+\beta\e}
 \sum_\gamma\|f_\gamma\|_2^2
 \label{eq:quarticupperbranch}
\end{equation}
for $12\le p\le14$ and $1\le\beta(p-12)/2$, by induction on
$R^\beta$. The bounded scale case follows from $L^2$ orthogonality.
The standard broad narrow decomposition gives
\[
\begin{aligned}
 \biggl\{\Bigl|\sum_\gamma f_\gamma\Bigr|>\alpha\biggr\}
 &\subset\bigcup_{|\nu|=K^{-1}}
 \biggl\{\Bigl|\sum_{\gamma\subset\nu}f_\gamma\Bigr|>c\alpha\biggr\}\cup
 \bigcup_{\substack{\nu_1,\ldots,\nu_4\text{ separated}\\|\nu_i|=K^{-1}}}
 \biggl\{\min_i\Bigl|\sum_{\gamma\subset\nu_i}f_\gamma\Bigr|>\alpha/(CK)\biggr\},
\end{aligned}
\]
where $c>0$ is independent of $K$.
Since $p\ge12+2/\beta$ and
$\|\sum_{\gamma\subset\nu_i}f_\gamma\|_\infty\lesssim R^\beta$,
Proposition~\ref{prop:quarticbroad} bounds the broad contribution by
\[
 C_{\e,K}R^{\beta(p-12-2/\beta)}R^{5\beta+1+\beta\e}
 \sum_\gamma\|f_\gamma\|_2^2
 =C_{\e,K}R^{-1+\beta(p-7)+\beta\e}
 \sum_\gamma\|f_\gamma\|_2^2.
\]

For the narrow contribution, rescale the interval of length $K^{-1}$ corresponding to $\nu$ to $[0,1]$. The small cap dimensions
become
\[
 KR^{-\beta}\times K^2R^{-2\beta}\times K^3R^{-3\beta}
 \times K^4R^{-1}.
\]
If $R/K^4>1$, the inductive range is preserved because
\[
 \frac{R}{K^4}
 \le\frac{R^{\beta(p-12)/2}}{K^4}
 \le\Bigl(\frac{R^\beta}{K}\Bigr)^{(p-12)/2}.
\]
Applying the inductive hypothesis at fine scale $KR^{-\beta}$
and last width $K^4/R$, and returning to the original variables,
gives the coefficient
\[
 (R/K^4)^{-1}(R^\beta/K)^{p-7+\e}
 =K^{11-p-\e}R^{-1+\beta(p-7)+\beta\e}.
\]
The Jacobian occurs on both sides of \eqref{eq:quarticupperbranch}.
Since $p\ge12$, summing over $\nu$ bounds the narrow contribution by
\[
 C_0C_\e K^{-1-\e}R^{-1+\beta(p-7)+\beta\e}
 \sum_\gamma\|f_\gamma\|_2^2,
\]
where $C_0$ is independent of $K$ and $p$.
For the remaining case $R/K^4\le1$, cylindrical decoupling for
the cubic moment curve, followed by H\"older's inequality over the
$R^\beta|\nu|$ small caps and the normalization
$\|f_\gamma\|_\infty\le1$, gives
\begin{equation}
\begin{aligned}
 \Bigl\|\sum_{\gamma\subset\nu}f_\gamma\Bigr\|_{12}^{12}
 &\lesssim_\e (R^\beta|\nu|)^\e
 \bigl(\sum_{\gamma\subset\nu}\|f_\gamma\|_{12}^2\bigr)^6\\
 &\le (R^\beta|\nu|)^{5+\e}
 \sum_{\gamma\subset\nu}\|f_\gamma\|_{12}^{12}\le (R^\beta|\nu|)^{5+\e}
 \sum_{\gamma\subset\nu}\|f_\gamma\|_2^2.
\end{aligned}
 \label{eq:quarticcubic}
\end{equation}
Together with
$\|\sum_{\gamma\subset\nu}f_\gamma\|_\infty\lesssim R^\beta/K$,
this gives the required bound directly:
\[
\begin{aligned}
 \Bigl\|\sum_{\gamma\subset\nu}f_\gamma\Bigr\|_p^p
 &\lesssim_\e (R^\beta/K)^{p-7+\e}
 \sum_{\gamma\subset\nu}\|f_\gamma\|_2^2\le (R/K^4)^{-1}(R^\beta/K)^{p-7+\e}
 \sum_{\gamma\subset\nu}\|f_\gamma\|_2^2.
\end{aligned}
\]
Choose $K$ so that $C_0K^{-1-\e}<1/2$, and then choose $C_\e$
large enough to close the induction. Setting $p=12+2/\beta$ gives
$-1+\beta(p-7)=1+5\beta$, which proves \eqref{eq:smallcapweak}
for $n=4$ and $\beta\ge1$. This completes the proof of
Proposition~\ref{prop:densityweak}.
\end{proof}

\subsection{Proof of Theorem~\ref{thm:SC}}
\label{sec:applications}

Set $p_0=n(n-1)+2/\beta$. Choose smooth cutoffs $\eta_\gamma=1$ on
$\gamma$, supported in $2\gamma$, with $\|\check\eta_\gamma\|_1\lesssim1$,
and define $Tg=\sum_\gamma g_\gamma*\check\eta_\gamma$.
We view the input as a function on the disjoint union of one copy of
$\R^n$ for each $\gamma$. Orthogonality and Young's inequality give
\[
 \|Tg\|_2^2\lesssim\sum_\gamma\|g_\gamma\|_2^2,
 \qquad \|Tg\|_\infty\lesssim R^\beta\sup_\gamma\|g_\gamma\|_\infty.
\]
The proof of Proposition~\ref{prop:densityweak} also applies to fixed
enlargements of the caps. Thus, for $|g_\gamma|\le\mathbf1_{E_\gamma}$
with $|E_\gamma|<\infty$, Young's inequality and that proposition imply
\[
 \alpha^{p_0}|\{|Tg|>\alpha\}|
 \lesssim_\delta R^{\beta(p_0/2-1+\delta)}\sum_\gamma|E_\gamma|.
\]
For $0<t<p_0-2$, the $L^2$ and $L^\infty$ bounds above give restricted
weak type estimates at $p_0-t$ and $p_0+t$, with the common coefficient
$R^{\beta(p_0/2-1+\delta+t)}$: use the $L^2$ bound for $\alpha\le1$
and the preceding estimate for $1<\alpha\lesssim R^\beta$.
Taking $\delta+t<\e$ and applying restricted weak type interpolation
\cite[Theorem~8.5]{TaoInterpolation}, we obtain
\[
 \|Tg\|_{p_0}^{p_0}
 \lesssim_\e R^{\beta(p_0/2-1+\e)}
 \sum_\gamma\|g_\gamma\|_{p_0}^{p_0}.
\]
Interpolating with the $L^2$ bound for $2\le p\le p_0$ and the
$L^\infty$ bound for $p\ge p_0$, and using
$T((f_\gamma)_\gamma)=\sum_\gamma f_\gamma$, gives
\begin{equation}\label{eq:smallcapuniform}
 \Bigl\|\sum_\gamma f_\gamma\Bigr\|_p^p
 \lesssim_{\e,p}R^{\beta\e}
 \bigl(R^{\beta(p/2-1)}+
 R^{\beta(p-n(n-1)/2-1)-1}\bigr)
 \sum_\gamma\|f_\gamma\|_p^p.
\end{equation}
This finishes the proof.

\section{Wave envelope estimate by Gan-Maldague-Oh}

Let \(n\ge2\) and let \(\gamma:[0,1]\to\mathbb R^n\) be a smooth
nondegenerate curve, so that
\[
 \inf_{t\in[0,1]}
 \bigl|\det(\gamma'(t),\ldots,\gamma^{(n)}(t))\bigr|>0.
\]
For \(0<s\le1\), let \(\Theta_\gamma(s)\) be the canonical cap
cover indexed by \(I=[ks,(k+2)s]\cap[0,1]\), \(k\in\mathbb Z\),
\(0\le ks<1\), with
\[
 \tau_I=\Bigl\{\gamma(t)+\sum_{j=1}^n\lambda_j\gamma^{(j)}(t):
              t\in I,\ |\lambda_j|\le s^j\Bigr\}.
\]
Thus a cap at scale \(s\) has dimensions
\(s\times s^2\times\cdots\times s^n\).
For \(R\ge1\), 
\[
 \mathcal M_\gamma(R)=\bigcup_{\theta\in\Theta_\gamma(R^{-1/n})}\theta.
\]
For \(\tau\in\Theta_\gamma(s)\), let \(t_\tau\) be the left endpoint
of the interval \(I\) corresponding to \(\tau\), and define the dual box and the wave envelope by
\[
 \tau^*=\bigl\{x\in\mathbb R^n:
     |x\cdot\gamma^{(j)}(t_\tau)|\le s^{-j},\ 1\le j\le n\bigr\},
 \qquad \widetilde V_{\tau,R}=Rs^n\tau^*.
\]
The wave envelope \(\widetilde V_{\tau,R}\) has coordinate half widths \(Rs^{n-1},Rs^{n-2},\ldots,R\).
For \(\tau\in\Xi(s)\), the corresponding wave envelope \(\widetilde U_{\tau,R}\) of
Gan-Maldague-Oh \cite{GMO} is centered at the origin and has half widths
\(Rs^d,Rs^{d-1},\ldots,R\) in the coordinates of \eqref{eq:conebox}.
We write \(V\parallel\widetilde V_{\tau,R}\) and \(U\parallel\widetilde U_{\tau,R}\)
for the tiles in fixed tilings by translates of these boxes.
Up to a fixed dilation, each envelope contains the centered dual boxes of all
caps \(\theta\subset\tau\) at the finest scale. For \(n=3\), the wave envelope \(\widetilde V_{\tau,R}\)
is comparable to the one used by Maldague \cite{Mal}.

For the moment curve \(\gamma=\gamma_n\) and the cone over the moment curve,
comparison with \eqref{eq:curvebox} and \eqref{eq:conebox} gives
\[
 \widetilde V_{\tau,R}=(sR^{1/n})^{n-2}V_{\tau,R},
 \qquad
 \widetilde U_{\tau,R}=(sR^{1/d})^{d-2}U_{\tau,R}.
\]
Thus \(\widetilde V_{\tau,R}=V_{\tau,R}\) at \(s=R^{-1/n}\)
and \(\widetilde U_{\tau,R}=U_{\tau,R}\) at \(s=R^{-1/d}\). For \(d=2\), the latter equality holds at every scale,
as in \cite{GWZ,GMAD}. At coarser scales, our averaging boxes are smaller for curves
with \(n\ge3\) and cones with \(d\ge3\). At \(s=1\), our boxes are comparable
to balls of radii \(R^{2/n}\) and \(R^{2/d}\), respectively, while both envelopes
of Gan-Maldague-Oh are comparable to balls of radius \(R\).

\begin{theorem}[{\cite[Theorem~2.1]{GMO}}]
\label{gmo8:curve-theorem}
\label{gmo8:curve-implication}
Let \(n\ge2\), \(4\le p\le n(n+1)-2\), and \(\e>0\).
For every Schwartz function \(f\) with
\(\operatorname{supp}\widehat f\subset\mathcal M_\gamma(R)\),
\begin{equation}\label{0924.thm31}
 \int_{\mathbb R^n}|f|^p
 \lesssim R^{\e+\frac{p-4}{2n}}
 \sum_{R^{-1/n}\le s\le1}\sum_{\tau\in\Theta_\gamma(s)}
 \sum_{V\parallel\widetilde V_{\tau,R}}|V|^{-1}
 \biggl(\int_{\mathbb R^n}\sum_{\theta\subset\tau}
                     |f_\theta|^{p/2}W_V\biggr)^2.
\end{equation}
\end{theorem}

\begin{proposition}
\label{prop:ADWE-implies-GMO}
Let \(n\ge3\).  Theorem~\ref{thm:main}
implies Theorem~\ref{gmo8:curve-theorem} for the moment curve \(\gamma(t)=(t,t^2,\ldots,t^n)\).
\end{proposition}

\begin{proof}
It suffices to prove
\begin{equation}
 \int|f|^p\lesssim R^{\e/2+\frac{p-4}{2n}}
          \int\bigl(\sum_\theta|f_\theta|^{p/2}\bigr)^2,
 \qquad 4\le p\le n(n+1)-2.
 \label{gmo3:ell-p-half}
\end{equation}
Indeed, by the essentially constant property, we have
\[
 |f_\theta|^{p/2}\lesssim
 |f_\theta|^{p/2}*\Bigl(\frac{W_{\theta^*}}{|\theta^*|}\Bigr).
\]
Assuming \eqref{gmo3:ell-p-half}, \cite[Proposition~4.6]{GMO},
applied with exponent \(p/2\ge2\), gives
\begin{equation}
\begin{aligned}
 \int|f|^p
 &\lesssim R^{\e/2+\frac{p-4}{2n}}
 \int\Bigl(\sum_\theta |f_\theta|^{p/2}*
                   \frac{W_{\theta^*}}{|\theta^*|}\Bigr)^2\\
 &\lesssim R^{\e+\frac{p-4}{2n}}
 \sum_{R^{-1/n}\le s\le1}\sum_{\tau\in\Theta_\gamma(s)}
 \sum_{V\parallel\widetilde V_{\tau,R}}|V|^{-1}
 \biggl(\int\sum_{\theta\subset\tau}|f_\theta|^{p/2}W_V\biggr)^2.
\end{aligned}
\label{gmo8:curve-weighted}
\end{equation}
This is the desired result.

We now prove \eqref{gmo3:ell-p-half}. The case \(4\le p\le p_n\) is immediate by the square function estimate
\cite{GMSF}. By the estimate and H\"older's inequality,
\[
 \int|f|^p\lesssim R^{\e/2}
       \int\bigl(\sum_\theta|f_\theta|^2\bigr)^{p/2}
 \lesssim R^{\e/2+\frac{p-4}{2n}}
       \int\bigl(\sum_\theta|f_\theta|^{p/2}\bigr)^2.
\]
So we may assume \(p_n<p\le n(n+1)-2\). We'll apply Theorem~\ref{thm:main}.
For the averaging box \(V\) in the theorem, its amplitude condition gives
\begin{equation}\label{0924.00}
 \alpha\lesssim R^{\epsilon}\#\tau
 \biggl[\frac1{|V|}\int
       \bigl(\sum_{\theta\subset\tau}|f_\theta|^2\bigr)^{p_n/4}
       W_V\biggr]^{2/p_n}.
\end{equation}
To apply the inequality, we use the following identity.
\[
 \int|f|^p=p\int_0^\infty\alpha^{p-p_n-1}
             \bigl[\alpha^{p_n}|\{|f|>\alpha\}|\bigr]\,d\alpha.
\]
We now apply Theorem~\ref{thm:main} together with \eqref{0924.00}, and obtain
\begin{equation}
\begin{aligned}
 \int|f|^p
 &\lesssim R^{\e/4}
 \sum_{R^{-1/n}\le s\le1}\sum_{\tau\in\Theta_\gamma(s)}
 (\#\tau)^{p-p_n}\sum_{V\parallel V_{\tau,R}}|V|
 \biggl[\frac1{|V|}\int
       \bigl(\sum_{\theta\subset\tau}|f_\theta|^2\bigr)^{p_n/4}
       W_V\biggr]^{2p/p_n}\\
 &\lesssim R^{\e/4}
 \sum_{R^{-1/n}\le s\le1}\sum_{\tau\in\Theta_\gamma(s)}
 (\#\tau)^{p-p_n}
 \int\bigl(\sum_{\theta\subset\tau}|f_\theta|^2\bigr)^{p/2}.
\end{aligned}
\end{equation}
The second inequality is by H\"{o}lder's inequality. Since \(\#\tau\lesssim s^{-1}\) and
\(\tau\) contains \(O(R^{1/n}s)\) caps at scale \(R^{-1/n}\), H\"older gives
\[
\begin{aligned}
 (\#\tau)^{p-p_n}
 \bigl(\sum_{\theta\subset\tau}|f_\theta|^2\bigr)^{p/2}
 &\lesssim s^{p_n-p}(R^{1/n}s)^{p/2-2}
       \bigl(\sum_{\theta\subset\tau}|f_\theta|^{p/2}\bigr)^2\\
 &=R^{\frac{p-4}{2n}}s^{\frac{n(n+1)-2-p}{2}}
       \bigl(\sum_{\theta\subset\tau}|f_\theta|^{p/2}\bigr)^2.
\end{aligned}
\]
The power of \(s\) is nonnegative. Hence, summing up the parameter $s$ and $\tau$, we obtain the desired result.
\end{proof}

\section{Restriction estimates}
\label{sec:restriction}

Write
\[
 \gamma_n(t)=(t,t^2,\ldots,t^n),\qquad
 Eg(x)=\int_0^1g(t)e^{2\pi i x\cdot\gamma_n(t)}\,dt,
 \qquad p_n=\frac{n(n+1)}2+1.
\]

\begin{proposition}\label{prop:restriction-moment}
Let \(n\ge3\). The amplitude dependent wave envelope type estimate for the moment
curve in \(\R^n\) (Theorem~\ref{thm:main}) implies that, for every \(\e>0\),
every ball \(B_R\subset\R^n\) of radius \(R\ge2\), and every
\(g\in L^{p_n}([0,1])\), we have
\begin{equation}\label{eq:restriction-local}
 \|Eg\|_{L^{p_n}(B_R)}\lesssim_{\e,n} R^\e\|g\|_{L^{p_n}([0,1])}.
\end{equation}
Consequently, for every \(q>p_n\),
\begin{equation}\label{eq:restriction-global}
 \|Eg\|_{L^q(\R^n)}\lesssim_{q,n}\|g\|_{L^{p_n}([0,1])}.
\end{equation}
\end{proposition}

\begin{proof}
\emph{The square function consequence.}
Let \(0<\delta<1\) be dyadic, and apply
Theorem~\ref{thm:main} at scale \(\delta^{-n}\).
Enlarging the collection of selected averaging boxes to the full tiling
and letting \(\alpha\downarrow0\) gives the ordinary
wave envelope type estimate. For each parent \(\tau\),
Cauchy--Schwarz and \(\int W_V=|V|\) give
\[
 |V|\biggl[\frac1{|V|}\int
       \Bigl(\sum_{\theta\subset\tau}|F_\theta|^2\Bigr)^{p_n/4}
       W_V\biggr]^2
 \le\int\Bigl(\sum_{\theta\subset\tau}|F_\theta|^2\Bigr)^{p_n/2}W_V.
\]
Sum over \(V\), using \(\sum_{V\parallel V_{\tau,\delta^{-n}}}W_V
\lesssim_n1\), and then over \(\tau\). The \(O(\log\delta^{-1})\)
angular scales are absorbed into an arbitrarily small power loss.
Thus, for every \(\eta>0\),
\begin{equation}\label{eq:restriction-square}
 \int_{\R^n}\Bigl|\sum_{\theta\in\Theta(\delta)}F_\theta\Bigr|^{p_n}
 \lesssim_{\eta,n}\delta^{-\eta}
 \int_{\R^n}\Bigl(\sum_{\theta\in\Theta(\delta)}|F_\theta|^2\Bigr)^{p_n/2},
 \qquad \supp\widehat F_\theta\subset\theta.
\end{equation}

\emph{Overlap of centered dual boxes.}
Partition \([0,1]\) into intervals \(I\) of length \(\delta\),
and let \(t_I\) be their midpoints. Define
\[
 T_I=\Bigl\{x\in\R^n:
 \bigl|\delta^j x\cdot\gamma_n^{(j)}(t_I)/j!\bigr|\le1
 \ (1\le j\le n)\Bigr\},\qquad
 \mu=\sum_I\mathbf1_{T_I}.
\]
For each fixed \(I_0\), we claim
\begin{equation}\label{eq:restriction-overlap}
 |T_{I_0}\cap\{\mu\ge m\}|
 \lesssim_n\delta^{-n(n+1)/2}m^{-n(n-1)/2}
 \qquad(m\ge1).
\end{equation}
Put \(t_0=t_{I_0}\) and use coordinates
\[
 y_j=\frac{\delta^j}{j!}x\cdot\gamma_n^{(j)}(t_0)
 \quad(1\le j\le n),\qquad
 dx=\delta^{-n(n+1)/2}dy.
\]
Then \(T_{I_0}\) becomes \([-1,1]^n\). If
\(t_I=t_0+k\delta\), the first defining inequality for
\(T_I\) gives
\[
 |P_y(k)|\le1,\qquad P_y(u)=\sum_{j=1}^n j y_j u^{j-1}.
\]
The set \(\{u:|P_y(u)|\le1\}\) has \(O_n(1)\) interval
components. Thus, if \(\mu(x)\ge m\) and \(m\) exceeds
a fixed dimensional constant, one component contains an
interval of length \(\gtrsim_n m\). Fix \(h=c_nm\), where
\(c_n>0\) is a sufficiently small dimensional constant
independent of \(x\). This interval then contains a cell
\([kh,(k+1)h]\) of the fixed \(h\) grid.
Polynomial derivative bounds on this cell give
\[
 z_j=\frac{h^j}{j!}P_y^{(j)}(kh),\qquad
 |z_j|\lesssim_n1\quad(0\le j\le n-1),\qquad
 dy=\frac{h^{-n(n-1)/2}}{n!}\,dz.
\]
We also have \(|P_y(0)|=|y_1|\le1\), since \(x\in T_{I_0}\).
In the \(z\) coordinates this is the slab condition
\[
 \Bigl|\sum_{j=0}^{n-1}(-k)^jz_j\Bigr|\le1.
\]
Inside the bounded \(z\) box, its volume is
\(O_n((1+|k|)^{-(n-1)})\): for \(|k|\ge1\), integrate
first in \(z_{n-1}\). Summing over the grid cells therefore gives
\[
 |T_{I_0}\cap\{\mu\ge m\}|
 \lesssim_n\delta^{-n(n+1)/2}h^{-n(n-1)/2}
 \sum_{k\in\mathbb Z}(1+|k|)^{-(n-1)}
 \lesssim_n\delta^{-n(n+1)/2}m^{-n(n-1)/2}.
\]
The series converges because \(n\ge3\); bounded \(m\) is
immediate. This proves \eqref{eq:restriction-overlap}.
Since
\[
 \frac{n(n-1)}2-\Bigl(\frac{p_n}2-1\Bigr)
 =\frac{(n-1)(n-2)}4>0,
\]
layer cake yields
\(\int_{T_{I_0}}\mu^{p_n/2-1}\lesssim_n\delta^{-n(n+1)/2}\).
Consequently, H\"older's inequality gives, for arbitrary
coefficients \(c_I\),
\begin{equation}\label{eq:restriction-weighted-boxes}
 \int\Bigl(\sum_I|c_I|^2\mathbf1_{T_I}\Bigr)^{p_n/2}
 \le\sum_I|c_I|^{p_n}\int_{T_I}\mu^{p_n/2-1}
 \lesssim_n\delta^{-n(n+1)/2}\sum_I|c_I|^{p_n}.
\end{equation}

\emph{Approximation of the extension operator.}
Choose a fixed Schwartz function \(\varphi\) on \(\R^n\)
with \(\varphi(0)=1\) and Fourier support in a sufficiently
small cube about the origin, and set
\[
 \Phi_I(x)=e^{2\pi i x\cdot\gamma_n(t_I)}
 \varphi\Bigl(\bigl(\delta^j x\cdot\gamma_n^{(j)}(t_I)/j!\bigr)_{j=1}^n\Bigr).
\]
The exact Taylor formula for \(\gamma_n\) shows that
\(\supp\widehat\Phi_I\subset\theta_I\in\Theta(\delta)\),
where \(\theta_I\) is the cap corresponding to \(I\).
Schwartz decay gives
\[
 |\Phi_I|^2\lesssim_n\sum_{\ell\ge0}2^{-10n\ell}
                  \mathbf1_{2^\ell T_I}.
\]
Apply \eqref{eq:restriction-weighted-boxes} to the common
dilations \(2^\ell T_I\). Their \(L^{p_n/2}\) norms acquire
the factor \(2^{2n\ell/p_n}\), so Minkowski and the rapid
decay imply
\begin{equation}\label{eq:restriction-smooth-boxes}
 \int\Bigl(\sum_I|c_I\Phi_I|^2\Bigr)^{p_n/2}
 \lesssim_n\delta^{-n(n+1)/2}\sum_I|c_I|^{p_n}.
\end{equation}

It suffices to consider \(B_R=B(0,R)\); a translated ball is
handled by modulation of \(g\). Choose dyadic
\(\delta\sim R^{-2}\), and put
\[
 c_I=\int_Ig(t)\,dt,\qquad F=\sum_Ic_I\Phi_I.
\]
The functions \(c_I\Phi_I\) are auxiliary Fourier bumps.
On \(B_R\), the smoothness of \(\gamma_n\),
\(\varphi(0)=1\), and \(R\delta\lesssim1\) give
\[
 \Bigl|Eg(x)-\sum_Ic_Ie^{2\pi i x\cdot\gamma_n(t_I)}\Bigr|
 +\Bigl|F(x)-\sum_Ic_Ie^{2\pi i x\cdot\gamma_n(t_I)}\Bigr|
 \lesssim_n R\delta\|g\|_{L^1([0,1])}.
\]
Consequently, since \(p_n>n\),
\[
 \|Eg-F\|_{L^{p_n}(B_R)}
 \lesssim_n R^{1+n/p_n}\delta\|g\|_{L^1([0,1])}
 \lesssim_n R^{-1+n/p_n}\|g\|_{L^{p_n}([0,1])}.
\]
On the other hand, \eqref{eq:restriction-square},
\eqref{eq:restriction-smooth-boxes}, and H\"older on each \(I\)
give
\[
 \|F\|_{p_n}^{p_n}
 \lesssim_{\eta,n}\delta^{-\eta-n(n+1)/2}
       \sum_I\Bigl|\int_Ig\Bigr|^{p_n}
 \lesssim_{\eta,n}\delta^{-\eta}\sum_I\int_I|g|^{p_n}
 \lesssim_{\eta,n}R^{2\eta}\|g\|_{p_n}^{p_n}.
\]
Here the cancellation of powers of \(\delta\) is exactly
\(p_n-1=n(n+1)/2\). Taking \(\eta>0\) sufficiently small
in terms of \(\e\) proves \eqref{eq:restriction-local}.

Finally, a smoothly truncated moment curve in \(\R^n\) has
Fourier decay \(O_n((1+|x|)^{-1/n})\) by van der Corput's lemma.
The local estimate holds on every ball and, by a fixed affine
change of parameter, on any fixed interval slightly larger than $[0,1]$. The epsilon removal theorem
\cite[Theorem~4.1]{GuoOhRestriction}, with both local exponents
equal to \(p_n\), therefore gives \eqref{eq:restriction-global}.
\end{proof}

This argument uses only the ordinary wave envelope type estimate
that follows from Theorem~\ref{thm:main}. The local endpoint estimate
\eqref{eq:restriction-local} has an \(R^\e\) loss; the global
conclusion \eqref{eq:restriction-global} is for \(q>p_n\).

\addtocontents{toc}{\protect\tocsectionindent{1pc}\protect\addvspace{12pt plus2pt}}

\appendix
\renewcommand{\theequation}{\thesection.\arabic{equation}}

\section{Standard interpolation argument}\label{sec:curveauxiliary}

\subsection{The estimate for the cone over the moment curve}\label{subsec:momentconeinterpolation}

\begin{proof}[Proof of Lemma~\ref{lem:curvconetail}]
By the layer cake formula,
\(\int_{|g|>\alpha^2}|g|^{p_n/2}
\lesssim\int_{\alpha^2/2}^\infty
t^{p_n/2-1}|\{|g|>t\}|\,dt\).
For each \(\tau\in\Xi(\sigma)\), take the partition of unity
\(\{\psi_{\tau,B}\}_{B\parallel\tau^*}\) constructed before
\eqref{eq:packetBernstein} and write
\(g_\tau=\sum_{B\parallel\tau^*}\psi_{\tau,B}g_\tau\).
Thus if \(\xi=\sum_{i=0}^{n-1}\lambda_i\phi_{n-1}^{(i)}(t_\tau)
\in\operatorname{supp}\widehat{\psi_{\tau,B}g_\tau}\),
then \(\frac12\le|\lambda_0|\le3\) and
\(|\lambda_i|\lesssim\sigma^i\) for \(1\le i\le n-1\).
For \(t>0\), put
\[g_{\tau,t}=g_\tau
\sum_{\substack{B\parallel\tau^*\\
\|\psi_{\tau,B}^{1/2}g_\tau\|_\infty\le t}}\psi_{\tau,B}.\]
Then \(|g_{\tau,t}|\lesssim\min\{|g_\tau|,t\}\), and the same
Fourier support bounds hold for \(g_{\tau,t}\).
Since \(g_{\tau,t}=0\) whenever \(g_\tau=0\), every active
parent for \(\{g_{\tau,t}\}_\tau\) is also active for
\(\{g_\tau\}_\tau\); hence its count is at most \(N(s)\).

Splitting
\(g=\sum_\tau(g_\tau-g_{\tau,t})+\sum_\tau g_{\tau,t}\),
we will estimate the level sets of the first sum by Chebyshev's
inequality and those of the second by the hypothesis on the cone,
each with threshold \(t/4\).

We first estimate the level set integral for
\(\sum_\tau(g_\tau-g_{\tau,t})\).
The bound
\(\|\psi_{\tau,B}g_\tau\|_1
\lesssim |B|\|\psi_{\tau,B}^{1/2}g_\tau\|_\infty\)
and Chebyshev's inequality give
\[
\begin{aligned}
 \int_{\alpha^2/2}^\infty t^{p_n/2-1}
 \biggl|\biggl\{\Bigl|\sum_\tau(g_\tau-g_{\tau,t})\Bigr|
 >t/4\biggr\}\biggr|\,dt
 &\lesssim
 \int_{\alpha^2/2}^\infty t^{p_n/2-2}
 \sum_\tau
 \sum_{\substack{B\parallel\tau^*\\
 \|\psi_{\tau,B}^{1/2}g_\tau\|_\infty>t}}
 \|\psi_{\tau,B}g_\tau\|_1\,dt\\
 &\lesssim
 \sum_\tau
 \sum_{\substack{B\parallel\tau^*\\
 \|\psi_{\tau,B}^{1/2}g_\tau\|_\infty>\alpha^2/2}}
 |B|\|\psi_{\tau,B}^{1/2}g_\tau\|_\infty^{p_n/2}.
\end{aligned}
\]
Apply \eqref{eq:packetBernstein} to each \(g_\tau\), with
exponent \(p_n/4\). Since the averaging box \(U_{\tau,\sigma^{-(n-1)}}\) at \(s=\sigma\)
is the dual box \(\tau^*\), the last sum is bounded by
\[
 \sum_\tau
 \sum_{\substack{U\parallel\tau^*\\
 U\in\mathcal G_\tau
 \bigl(\frac{\alpha}{C_\epsilon\sigma^{-\epsilon}}\bigr)}}
 |U|\,
 \bigl\||g_\tau|^{1/2}\bigr\|_{L_\#^{p_n/2}(U)}^{p_n}.
\]

We next estimate
\(\int_{\alpha^2/2}^\infty t^{p_n/2-1}
|\{|\sum_\tau g_{\tau,t}|>t/4\}|\,dt\).
Apply the amplitude dependent wave envelope type estimate for the cone over the moment curve
in \(\mathbb R^n\) (Theorem~\ref{thm:main-cone} with \(d=n-1\)) to
\(\sum_\tau g_{\tau,t}\) at scale \(R=\sigma^{-(n-1)}\) and amplitude
\(t/4\), with \(g_{\tau,t}\) as the fine pieces.
Using the identity
\[
 \biggl\|\Bigl(\sum_{\tau\subset\rho}|g_{\tau,t}|^2\Bigr)^{1/2}
 \biggr\|_{L_\#^{p_{n-1}/2}(U)}^{p_{n-1}}
 =
 \biggl\|\Bigl(\sum_{\tau\subset\rho}|g_{\tau,t}|^2
 \Bigr)^{p_{n-1}/(2p_n)}
 \biggr\|_{L_\#^{p_n/2}(U)}^{p_n}
\]
and \(N(s)^{p_{n-1}/p_n}\le N(s)\) in the selection
condition, we obtain
\[
 t^{p_{n-1}}
 \biggl|\biggl\{\Bigl|\sum_\tau g_{\tau,t}\Bigr|>t/4
 \biggr\}\biggr|
 \lesssim_\epsilon \sigma^{-\epsilon}
 \sum_{\sigma\le s\le1}
 \sum_{\rho\in\Xi(s)}
 \sum_{\substack{
 U\parallel U_{\rho,\sigma^{-(n-1)}}\\
 U\in\mathcal G_\rho\Bigl(
 \bigl(\frac{t}{C_\epsilon\sigma^{-\epsilon}}\bigr)^{p_{n-1}/p_n}
 \Bigr)}}
 |U|
 \biggl\|\Bigl(\sum_{\tau\subset\rho}|g_{\tau,t}|^2
 \Bigr)^{p_{n-1}/(2p_n)}
 \biggr\|_{L_\#^{p_n/2}(U)}^{p_n}.
\]

We next verify the selection condition in
Lemma~\ref{lem:curvconetail}.
For \(t\ge\alpha^2/2\), the selection condition in the estimate for the cone over the moment curve
and
\(\min\{a,t\}^2\le
t^{2-p_n/p_{n-1}}a^{p_n/p_{n-1}}\)
give
\[
 \biggl\|\Bigl(\sum_{\tau\subset\rho}|g_\tau|\Bigr)^{1/2}
 \biggr\|_{L_\#^{p_n/2}(U)}
 \ge
 \biggl\|\Bigl(\sum_{\tau\subset\rho}
 |g_\tau|^{p_n/p_{n-1}}\Bigr)^{p_{n-1}/(2p_n)}
 \biggr\|_{L_\#^{p_n/2}(U)}
 \gtrsim_\epsilon
 \frac{\sqrt t}{\sigma^{-\epsilon p_{n-1}/p_n}N(s)}.
\]
Since \(t\ge\alpha^2/2\) and \(p_{n-1}<p_n\), we have
\(U\in\mathcal G_\rho
(\alpha/(C_\epsilon\sigma^{-\epsilon}))\).
Multiplying the estimate for this cone by
\(t^{p_n/2-p_{n-1}-1}\) and integrating over
\(t\ge\alpha^2/2\), Tonelli's theorem gives
\begin{equation}\label{eq:coneweakstrongintegrated}
\begin{aligned}
 &\int_{\alpha^2/2}^\infty t^{p_n/2-1}
 \biggl|\biggl\{\Bigl|\sum_\tau g_{\tau,t}\Bigr|
 >t/4\biggr\}\biggr|\,dt\\
 &\lesssim_\epsilon \sigma^{-\epsilon}
 \sum_{\sigma\le s\le1}
 \sum_{\rho\in\Xi(s)}
 \sum_{\substack{
 U\parallel U_{\rho,\sigma^{-(n-1)}}\\
 U\in\mathcal G_\rho
 \bigl(\frac{\alpha}{C_\epsilon\sigma^{-\epsilon}}\bigr)}}
 |U|
 \int_0^\infty
 \biggl\|\Bigl(\sum_{\tau\subset\rho}|g_{\tau,t}|^2
 \Bigr)^{p_{n-1}/(2p_n)}
 \biggr\|_{L_\#^{p_n/2}(U)}^{p_n}
 t^{p_n/2-p_{n-1}-1}\,dt.
\end{aligned}
\end{equation}

We now bound the \(t\) integral in each summand of
\eqref{eq:coneweakstrongintegrated} by
\(\|(\sum_{\tau\subset\rho}|g_\tau|)^{1/2}
\|_{L_\#^{p_n/2}(U)}^{p_n}\).
For \(a>0\), splitting at \(t=a\) gives
\(\int_0^\infty\min\{a,t\}^{p_{n-1}}
t^{p_n/2-p_{n-1}-1}\,dt\sim a^{p_n/2}\),
since \(p_n/2<p_{n-1}\).
Using \(|g_{\tau,t}|\lesssim\min\{|g_\tau|,t\}\),
Minkowski first in the weighted average over \(U\) and then in the sum
over \(\tau\) gives
\begin{equation}\label{eq:coneweakstrongpacketaverage}
\begin{aligned}
 &\biggl[
 \int_0^\infty
 \biggl\|\Bigl(\sum_{\tau\subset\rho}|g_{\tau,t}|^2
 \Bigr)^{p_{n-1}/(2p_n)}
 \biggr\|_{L_\#^{p_n/2}(U)}^{p_n}
 t^{p_n/2-p_{n-1}-1}\,dt
 \biggr]^{1/2}\\
 &\le
 \fint_U
 \biggl[\int_0^\infty
 \Bigl(\sum_{\tau\subset\rho}|g_{\tau,t}(x)|^2
 \Bigr)^{p_{n-1}/2}
 t^{p_n/2-p_{n-1}-1}\,dt\biggr]^{1/2}\\
 &\lesssim
 \biggl\|\Bigl(\sum_{\tau\subset\rho}
 |g_\tau|^{p_n/p_{n-1}}\Bigr)^{p_{n-1}/(2p_n)}
 \biggr\|_{L_\#^{p_n/2}(U)}^{p_n/2}
 \le
 \biggl\|\Bigl(\sum_{\tau\subset\rho}|g_\tau|\Bigr)^{1/2}
 \biggr\|_{L_\#^{p_n/2}(U)}^{p_n/2}.
\end{aligned}
\end{equation}
The last inequality follows from
\((\sum_{\tau\subset\rho}|g_\tau|^{p_n/p_{n-1}})
^{p_{n-1}/p_n}\le\sum_{\tau\subset\rho}|g_\tau|\).
Squaring \eqref{eq:coneweakstrongpacketaverage} and substituting
into \eqref{eq:coneweakstrongintegrated}, we obtain
\[
 \int_{\alpha^2/2}^\infty t^{p_n/2-1}
 \biggl|\biggl\{\Bigl|\sum_\tau g_{\tau,t}\Bigr|
 >t/4\biggr\}\biggr|\,dt
 \lesssim_\epsilon \sigma^{-\epsilon}
 \sum_{\sigma\le s\le1}
 \sum_{\rho\in\Xi(s)}
 \sum_{\substack{
 U\parallel U_{\rho,\sigma^{-(n-1)}}\\
 U\in\mathcal G_\rho
 \bigl(\frac{\alpha}{C_\epsilon\sigma^{-\epsilon}}\bigr)}}
 |U|
 \biggl\|\Bigl(\sum_{\tau\subset\rho}|g_\tau|\Bigr)^{1/2}
 \biggr\|_{L_\#^{p_n/2}(U)}^{p_n}.
\]
Together with the estimate for
\(\sum_\tau(g_\tau-g_{\tau,t})\), this proves the lemma.
\end{proof}

\subsection{The estimate for Taylor cones}
\label{app:hcauxiliary}
\label{sec:hcanalytic}
\label{app:hcpackets}

\begin{proof}[Proof of Lemma~\ref{lem:hcpackettail}]
Use the wave packet decomposition of Appendix~\ref{subsec:momentconeinterpolation},
writing \(h_{\theta,t}\) for the sum of the wave packets of \(h_\theta\) with
\(\|\psi_{\theta,T}^{1/2}h_\theta\|_\infty\le t\). The wave packet cutoffs preserve the cap bounds
up to fixed constants. In particular, every \(\xi\in\supp\widehat{h_{\theta,t}}\) satisfies
\[
 \xi=\sum_{i=0}^d\lambda_i\phi_d^{(i)}(t_\theta),\qquad
 \max_{0\le i\le m+1}r^{i/d}|\lambda_i|\gtrsim1.
\]
Truncation can only decrease the number of active parents. As in
Appendix~\ref{subsec:momentconeinterpolation}, layer cake reduces the proof
to levels \(t\ge\alpha^2/2\), and Chebyshev's inequality together with
\eqref{eq:packetBernstein} bounds the contribution of the discarded wave packets
by the second sum in the lemma.

For the remaining wave packets, apply \eqref{eq:hcgeneralconvolution} at order
\(m+1\), scale \(r\), and amplitude \(t/4\). At \(m=d-3\), use
\eqref{eq:hcconstantestimate} and Proposition~\ref{prop:hcbase} instead.
Fix a box \(B\) selected by this estimate, with parent \(\rho\) of width \(u\).
Its weighted average is \(\gtrsim (t/(r^\epsilon N(u)))^{p_{d-m-1}/2}\).
To compare it with the average in the first sum of the lemma, observe that
\[
 |h_{\theta,t}|^2\lesssim \min\{|h_\theta|,t\}^2
 \le t^{2-p_{d-m}/p_{d-m-1}}|h_\theta|^{p_{d-m}/p_{d-m-1}}.
\]
All sums and convolutions are positive, and the powers in the average
telescope to \(p_{d-m-1}/4\). Thus the average with
\(|h_\theta|^{p_{d-m}/p_{d-m-1}}\) in place of \(|h_{\theta,t}|^2\)
is bounded below, up to a constant, by
\begin{equation}\label{eq:hcpacketselection}
 t^{p_{d-m}/4-p_{d-m-1}/2}
 \Bigl(\frac{t}{r^\epsilon N(u)}\Bigr)^{p_{d-m-1}/2}
 =\frac{t^{p_{d-m}/4}}{(r^\epsilon N(u))^{p_{d-m-1}/2}}
 \gtrsim \Bigl(\frac{\alpha}{r^\epsilon N(u)}\Bigr)^{p_{d-m}/2}.
\end{equation}
Here \(t\ge\alpha^2/2\) and \(p_{d-m-1}<p_{d-m}\). Taking a smaller
\(\epsilon\) in the input estimate absorbs the fixed constant in the
selection threshold. We may therefore sum over the boxes selected at
amplitude \(\alpha/r^\epsilon\) before integrating in \(t\).

Multiply the level set bound for
\(t^{p_{d-m-1}}|\{|\sum_\theta h_{\theta,t}|>t/4\}|\) by
\(d\mu(t)=t^{p_{d-m}/2-p_{d-m-1}-1}\,dt\). For each fixed box and chain
of scales, take the \(L_t^2(d\mu)\) norm of the average inside the outer
square. Minkowski moves this norm inside \(\fint_B\) and the sum over
subcaps. At each successive convolution, write
\(a=p_{d-m-\ell}/p_{d-m-\ell-1}\) for the power at that level. Pointwise
in the spatial variable, Minkowski gives
\[
 \big\|(F*w_{B_{\tau_\ell}})^a\big\|_{L_t^{p_{d-m-\ell-1}/2}(d\mu)}
 \le \Big(\|F\|_{L_t^{p_{d-m-\ell}/2}(d\mu)}*w_{B_{\tau_\ell}}\Big)^a
 \qquad(F\ge0).
\]
Starting at \(\ell=d-m-3\) and ending at \(\ell=1\), the time exponent
therefore increases from \(2=p_2/2\) to \(p_{d-m-1}/2\). At the innermost
sum, the scalar integral used in Appendix~\ref{subsec:momentconeinterpolation}
gives
\[
 \big\||h_{\theta,t}|^2\big\|_{L_t^{p_{d-m-1}/2}(d\mu)}
 \lesssim |h_\theta|^{p_{d-m}/p_{d-m-1}},
 \qquad p_{d-m}/2<p_{d-m-1}.
\]
When \(m=d-3\), the average has no intermediate convolutions and the time norm is
already \(L_t^2(t^{-3/2}\,dt)\), giving \(\fint_B\sum_{\theta\subset\rho}|h_\theta|^{7/4}\).
Squaring this bound for each average and summing over the boxes in
\eqref{eq:hcpacketselection} gives the first sum in the lemma.
\end{proof}

\section{From weak type to strong type}\label{subsec:curveweakstrong}

\begin{proof}[Proof of Lemma~\ref{lem:curveweakstrong}]
Fix \(r\), \(g\), and \(\alpha\) as in the lemma. By the layer cake formula,
\(\int_{|g|>\alpha}|g|^{p_n}\lesssim
\int_{\alpha/2}^\infty t^{p_n-1}|\{|g|>t\}|\,dt\).
For \(t>0\), set
\[m_{\tau,t}=\sum_{T\parallel\tau^*:\,
t/(CR)<\|\psi_{\tau,T}^{1/2}g_\tau\|_\infty\le t}\psi_{\tau,T}.\]
The wave packet partition of Section~\ref{sec:curvepruning} gives
\begin{equation}\label{eq:curveweakstrongdecomposition}
g=
\sum_{\substack{\tau,\ T\parallel\tau^*\\
\|\psi_{\tau,T}^{1/2}g_\tau\|_\infty\le t/(CR)}}
\psi_{\tau,T}g_\tau
+\sum_\tau g_\tau m_{\tau,t}
+\sum_{\substack{\tau,\ T\parallel\tau^*\\
\|\psi_{\tau,T}^{1/2}g_\tau\|_\infty>t}}
\psi_{\tau,T}g_\tau.
\end{equation}
The first sum on the right of \eqref{eq:curveweakstrongdecomposition} has absolute value at most
\(\frac{t}{CR}\sum_\tau\sum_T\psi_{\tau,T}^{1/2}\le t/4\),
since there are \(O(r^{1/n})\le O(R^{1/n})\) caps at scale \(r^{-1/n}\).
Thus the layer cake formula and \eqref{eq:curveweakstrongdecomposition} give
\begin{equation}\label{eq:curveweakstrongreduction}
\int_{|g|>\alpha}|g|^{p_n}
\lesssim
\int_{\alpha/2}^\infty t^{p_n-1}
\Biggl(
\biggl|\biggl\{\Bigl|\sum_\tau g_\tau m_{\tau,t}\Bigr|>t/4\biggr\}\biggr|
+
\biggl|\biggl\{\Bigl|
\sum_{\substack{\tau,\ T\parallel\tau^*\\
\|\psi_{\tau,T}^{1/2}g_\tau\|_\infty>t}}
\psi_{\tau,T}g_\tau
\Bigr|>t/4\biggr\}\biggr|
\Biggr)\,dt.
\end{equation}

We first estimate the first term in parentheses in \eqref{eq:curveweakstrongreduction} using the hypothesis of Lemma~\ref{lem:curveweakstrong} with \(g\) replaced by
\(\sum_\tau g_\tau m_{\tau,t}\) and \(\alpha\) replaced by \(t/4\).
Each \(g_\tau m_{\tau,t}\) has Fourier support in a fixed enlargement of \(\tau\).
Since \(g_\tau m_{\tau,t}=0\) whenever \(g_\tau=0\), the active counts in the weak estimate for \(\sum_\tau g_\tau m_{\tau,t}\) are at most \(N(s)\).
For each \(V\parallel V_{\rho,r}\) on its RHS, the selection condition, \(0\le m_{\tau,t}\le1\), and \(t\ge\alpha/2\) give
\begin{equation}\label{eq:curveweakstrongselection}
\begin{aligned}
 \fint_V
 \Bigl(\sum_{\tau\subset\rho}|g_\tau|^2*w_\tau\Bigr)^{p_n/4}
 &\ge
 \fint_V
 \Bigl(\sum_{\tau\subset\rho}
 (|g_\tau|^2m_{\tau,t}^2)*w_\tau\Bigr)^{p_n/4}\\
 &\gtrsim
 \Bigl(\frac{t}{\Lambda(R)N(s)}\Bigr)^{p_n/2}
 \gtrsim
 \Bigl(\frac{\alpha}{\Lambda(R)N(s)}\Bigr)^{p_n/2}.
\end{aligned}
\end{equation}
Hence \(V\in\mathcal G_\rho(\frac{\alpha}{R^\e\Lambda(R)})\).
Multiplying the weak estimate by \(t^{-1}\), integrating over \(t\ge\alpha/2\), and using \eqref{eq:curveweakstrongselection}, we obtain
\begin{equation}\label{eq:curveweakstrongintegrated}
\begin{aligned}
\int_{\alpha/2}^\infty t^{p_n-1}
\biggl|\biggl\{\Bigl|\sum_\tau g_\tau m_{\tau,t}\Bigr|>t/4\biggr\}\biggr|\,dt
&\lesssim
\Lambda(R)^{p_n}
\sum_{r^{-1/n}\le s\le1}\sum_{\rho\in\Theta(s)}
\sum_{\substack{V\parallel V_{\rho,r}\\
V\in\mathcal G_\rho(\frac{\alpha}{R^\e\Lambda(R)})}}|V|
\\
&\quad\times
\int_0^\infty
\biggl[
\fint_V
\bigl(\sum_{\tau\subset\rho}
(|g_\tau|^2m_{\tau,t}^2)*w_\tau\bigr)^{p_n/4}
\biggr]^2\frac{dt}{t}.
\end{aligned}
\end{equation}

We estimate the integral in \(t\) in each summand on the right of \eqref{eq:curveweakstrongintegrated} by the square of the corresponding weighted average over \(V\) with \(m_{\tau,t}\) removed.
By the definition of \(m_{\tau,t}\) and \(0\le m_{\tau,t}\le1\),
\begin{equation}\label{eq:curveweakstrongmask}
\int_0^\infty m_{\tau,t}(x)^{p_n}\frac{dt}{t}
\le
\sum_{T\parallel\tau^*}
\psi_{\tau,T}(x)
\int_{\|\psi_{\tau,T}^{1/2}g_\tau\|_\infty}
^{CR\|\psi_{\tau,T}^{1/2}g_\tau\|_\infty}
\frac{dt}{t}
\lessapprox1.
\end{equation}
Apply Minkowski first between \(L_t^2(dt/t)\) and
\(L_x^1(W_V\,dx/|V|)\), and then between
\(L_t^{p_n/2}(dt/t)\) and \(\ell_\tau^1L_y^1(w_\tau(y)\,dy)\).
By \eqref{eq:curveweakstrongmask}, for each \(V\) on the right of \eqref{eq:curveweakstrongintegrated},
\begin{equation}\label{eq:curveweakstrongpacketaverage}
\begin{aligned}
&\biggl\|
\fint_V
\bigl(\sum_{\tau\subset\rho}
(|g_\tau|^2m_{\tau,t}^2)*w_\tau\bigr)^{p_n/4}
\biggr\|_{L_t^2(dt/t)}\\
&\le\fint_V
\biggl\|
\sum_{\tau\subset\rho}\int
|g_\tau(x-y)|^2m_{\tau,t}(x-y)^2w_\tau(y)\,dy
\biggr\|_{L_t^{p_n/2}(dt/t)}^{p_n/4}\\
&\le\fint_V
\biggl[
\sum_{\tau\subset\rho}\int
|g_\tau(x-y)|^2
\|m_{\tau,t}(x-y)\|_{L_t^{p_n}(dt/t)}^2
w_\tau(y)\,dy
\biggr]^{p_n/4}\\
&\lessapprox
\fint_V
\bigl(\sum_{\tau\subset\rho}|g_\tau|^2*w_\tau\bigr)^{p_n/4}.
\end{aligned}
\end{equation}
Squaring \eqref{eq:curveweakstrongpacketaverage} and substituting into \eqref{eq:curveweakstrongintegrated} therefore gives
\begin{equation}\label{eq:curveweakstrongmiddle}
\begin{aligned}
\int_{\alpha/2}^\infty t^{p_n-1}
\biggl|\biggl\{\Bigl|\sum_\tau g_\tau m_{\tau,t}\Bigr|>t/4\biggr\}\biggr|\,dt
&\lesssim_\e R^\e\Lambda(R)^{p_n}
\sum_{r^{-1/n}\le s\le1}\sum_{\rho\in\Theta(s)}
\sum_{\substack{V\parallel V_{\rho,r}\\
V\in\mathcal G_\rho(\frac{\alpha}{R^\e\Lambda(R)})}}|V|
\biggl\|
\bigl(\sum_{\tau\subset\rho}|g_\tau|^2*w_\tau\bigr)^{1/2}
\biggr\|_{L^{p_n/2}_\#(V)}^{p_n}.
\end{aligned}
\end{equation}

We now estimate the second term in parentheses in \eqref{eq:curveweakstrongreduction} after integration in \(t\).
Chebyshev's inequality and
\(\|\psi_{\tau,T}g_\tau\|_1
\lesssim|T|\|\psi_{\tau,T}^{1/2}g_\tau\|_\infty\) give
\begin{equation}\label{eq:curveweakstronghigh}
\begin{aligned}
\int_{\alpha/2}^\infty t^{p_n-1}
\biggl|\biggl\{\Bigl|
\sum_{\substack{\tau,\ T\parallel\tau^*\\
\|\psi_{\tau,T}^{1/2}g_\tau\|_\infty>t}}
\psi_{\tau,T}g_\tau
\Bigr|>t/4\biggr\}\biggr|\,dt
&\lesssim
\sum_\tau\sum_{\substack{T\parallel\tau^*\\
\|\psi_{\tau,T}^{1/2}g_\tau\|_\infty>\alpha/2}}
|T|\|\psi_{\tau,T}^{1/2}g_\tau\|_\infty^{p_n}
\\
&\lesssim
\sum_\tau\sum_{\substack{T\parallel\tau^*\\
\fint_T|g_\tau|^{p_n/2}\gtrsim\alpha^{p_n/2}}}
|T|\biggl[\fint_T|g_\tau|^{p_n/2}\biggr]^2.
\end{aligned}
\end{equation}
The last inequality is \eqref{eq:packetBernstein} with exponent \(p_n/2\).
Since \(V_{\tau,r}=\tau^*\) and
\(|g_\tau|^2\lesssim|g_\tau|^2*w_\tau\), the last sum in \eqref{eq:curveweakstronghigh} is bounded by the \(s=r^{-1/n}\) terms on the right of \eqref{eq:curveweakstrongmiddle}.
Integrating the two terms in parentheses in \eqref{eq:curveweakstrongreduction}, we obtain in each case a bound by the right hand side of \eqref{eq:curveweakstrongmiddle}, proving the lemma.
\end{proof}

\bibliographystyle{alpha}
\bibliography{reference}

\end{document}